\documentclass{hha}
\usepackage[all]{xy}
\usepackage{amscd}

\newcommand{\si}{\psi}
\newcommand{\ot}{\otimes}

\begin{document}

\title{Hopf Cyclic Cohomology in Braided Monoidal Categories} 

\author{Masoud Khalkhali}             
\email{masoud@uwo.ca}
\address{Department of Mathematics,
         University of Western Ontario,
         London, Ontario, N6A 5B7,
         Canada}

\author{Arash Pourkia}
\email{apourkia@uwo.ca}
\address{Department of Mathematics,
         University of Western Ontario,
         London, Ontario, N6A 5B7,
         Canada}


\begin{abstract}
We  extend the formalism of Hopf cyclic cohomology
to  the context of braided   categories. For a Hopf algebra in a
braided monoidal abelian category we introduce the notion of stable
anti-Yetter-Drinfeld module. We associate a para-cocyclic and a
cocyclic object to a braided Hopf algebra endowed with a braided
modular pair in involution in the sense of Connes and Moscovici.
When the braiding is symmetric the full formalism of Hopf cyclic
cohomology with coefficients can be extended to our categorical setting.
\end{abstract}

\maketitle

%%*********************************************************************************************
\section{Introduction}
In ~\cite{cm2,cm3,cm4},  Connes and Moscovici, motivated by 
transverse index theory for foliations, defined a cohomology theory
of cyclic type for Hopf algebras by introducing the concept of a
modular pair in involution. This theory was later extended,  by the
first author and collaborators ~\cite{hkrs1,hkrs2,kr1,kr2},  by
introducing the notion of a stable anti-Yetter-Drinfeld module  as
coefficients for a cyclic cohomology theory of algebras or
coalgebras endowed with an action or coaction of a Hopf algebra.
Modular pairs in involution appeared naturally as one dimensional
stable anti-Yetter-Drinfeld modules in this latter theory. It is by
now clear that Hopf cyclic cohomology is the right noncommutative
analogue of group homology and Lie algebra homology. In particular
it allows an extension of Connes' theory of noncommutative
characteristic classes ~\cite{c2,cbook} to a setup involving an
action of a Hopf algebra or quantum group ~\cite{cm2,cm3}.

There are many examples of Hopf algebra-like objects that are close
to being a Hopf algebra but are not a Hopf algebra in the usual
sense. Examples include (differential graded) super Hopf algebras,
quasi-Hopf algebras, multiplier Hopf algebras, Hopf algebroids, and
locally compact quantum groups. In some cases, but certainly not
always,  these objects are Hopf algebras in an appropriate monoidal
category. Differential graded super Hopf algebras and quasi-Hopf
algebras are examples of this. Study of Hopf algebras in symmetric
monoidal categories goes back to ~\cite{saav}. Recent work on braided
Hopf algebras is mostly motivated by low dimensional topology
~\cite{maj2, maj, tura}.

In this paper we work in an arbitrary braided
monoidal abelian category and extend the notion of a stable
anti-Yetter-Drinfeld module over a Hopf algebra in such a category.
We show that to any braided Hopf algebra endowed with a braided
modular pair in involution one can associate a para-cocyclic object.
This para-cocyclic object is cyclic if the ambient category is
symmetric. In fact  Theorem ~\eqref{symm and ciclc cndtion} (cf.\ also
Remark ~\eqref{gnral taun}),  shows that this para-cocyclic object is
almost never cocyclic if the category is not symmetric. Of course,
as with any para-cocyclic object, by restricting to an appropriate
subspace we obtain a cocyclic object.

In the symmetric case one can go much further. In this case and for
an arbitrary braided stable anti-Yetter-Drinfeld module we obtain a
cocyclic object.  In fact in this case it  is no longer needed  to
restrict to Hopf algebras and one can work with module coalgebras.
As  a special case,  we define a Hopf cyclic cohomology for a
(differential graded) super Hopf algebra and  relate it to the
cohomology of super Lie algebras by considering the enveloping
algebra of a super Lie algebra.

The paper is organized as follows. In Section \ref{secbmc} we recall basic
notions of braided monoidal categories and \emph{braided Hopf algebras} which are, by definition,  Hopf algebra objects in such
 categories. In Section \ref{hopf cyclc for H C M} we define the notion of a  stable
 anti-Yetter-Drinfeld (SAYD) module for a braided Hopf algebra.  For symmetric braidings,  we define
 a Hopf cyclic cohomology theory for triples $(H, C, M)$ consisting
 of a coalgebra object endowed with an action of a braided Hopf
 algebra $H$, and a braided SAYD  module $M$. In Section \ref{brdd cm hpf thry} we show that when $C=H$ the Hopf cyclic complex
 of the triple $(C, C, {}^\sigma I_\delta)$ simplifies and we compute the resulting cocyclic object. 
In Sections \ref{super hpf alg} and \ref{super lie alg} we focus  on the special case of (differential graded) super Hopf algebras
 and show that, for the universal enveloping algebra of a (differential graded)  super Lie algebra,   
Hopf cyclic cohomology reduces to Lie   algebra homology.

The last section is in a sense the heart of this paper. Here we
 work in an arbitrary braided monoidal category and show that, for
 one dimensional SAYD's, one can always define a para-cocyclic object for a braided Hopf algebra.
 We compute  powers of the cyclic  operator and express it in terms of the braiding of the category
 (Theorem ~\eqref{symm and ciclc cndtion} and Remark ~\eqref{gnral taun}). 
Extending these results to more general braided  SAYD's  and  to braided triples $(H, C, M)$ is not straightforward
 and requires introducing extra structures. This more general case will be dealt with elsewhere.

We should mention that the cyclic cohomology of (\emph{ribbon})-algebras in braided monoidal 
abelian categories has been introduced and studied in ~\cite{akmaj}, motivated by non-associative geometry. 
\ack
We would like to thank the referee for very helpful comments on the paper.

%%*******************************************************************************************************
\section{Hopf algebras in braided monoidal categories}\label{secbmc}
Recall that a \emph{monoidal}, or \emph{tensor}, category $(\mathcal{C},\otimes, I, a, l, r)$ 
consists of a category $\mathcal{C}$, a functor $\otimes: \mathcal{C}\times \mathcal{C}\to \mathcal{C}$, 
an object $I\in \mathcal{C}$ (called  \emph{unit object}), and natural
isomorphisms, defined for all objects $A$, $B$, $C$, of $\mathcal{C}$,
\begin{gather*}
a=a_{A, B, C}: A\otimes (B\otimes C) \to (A\otimes B)\otimes C,\\
l=l_A: I\otimes A \to A, \quad \quad r=r_A: A\otimes I \to A,
\end{gather*}
(called the \emph{associativity} and \emph{unit constraints}, respectively),  such that the following
 \emph{pentagon} and \emph{triangle} diagrams commute ~\cite{ml,maj}:
\[
\xymatrix{ &((A\otimes B)\otimes C)\otimes D\ar[dl] \ar[dr]&\\
(A\otimes (B \otimes C))\otimes D\ar[d]& & (A\otimes B) \otimes (C \otimes D)\ar[d]\\
A\otimes ((B \otimes C)\otimes D )\ar[rr]& & A\otimes (B \otimes
(C\otimes D))}
\]
\[
\xymatrix{ (A\otimes I)\otimes B  \ar[rr] \ar[dr]& & A\otimes (I \otimes B)\ar[dl]\\
& A\otimes B& }
\]

The coherence theorem of Mac Lane ~\cite{ml}  asserts  that all
diagrams formed by $a, l, r$ by tensoring and composing, commute.
More precisely, it asserts that any two natural transformations
defined by $a, l, r$ between any two functors defined by  $\otimes$
and $I$ are equal.

A \emph{braided monoidal} category is a monoidal category
$\mathcal{C}$ endowed with a natural family of isomorphisms
\[\psi_{A, B}: A\otimes B \to B\otimes A,  \] called  \emph{braiding}
such that for all objects $A, B, C$ of $\mathcal{C}$ the following
diagrams commute (\emph{hexagon axioms}):
\[
\xymatrix{ & A\otimes (B\otimes C)\ar[r]^{\psi} & (B\otimes C)\otimes A\ar[dr]^{a^{-1}} &\\
(A\otimes B)\otimes C\ar[ur]^{a^{-1}} \ar[dr]^{\psi \otimes id}& & &B\otimes (C \otimes A) \\
&(B\otimes A)\otimes C\ar[r]^{a^{-1}}& B\otimes (A\otimes C) \ar[ru]^{id \otimes \psi} &}
\]
\[
\xymatrix{ & (A\otimes B) \otimes C\ar[r]^{\psi} & C \otimes (A \otimes B)\ar[dr]^{a} &\\
A\otimes (B\otimes C)\ar[ur]^{a} \ar[dr]^{id \otimes \psi }& & &(C \otimes A) \otimes B  \\
&A \otimes (C \otimes B)\ar[r]^{a}& (A\otimes C)\otimes B \ar[ru]^{\psi \otimes id}  &}
\]
If we show the braiding map $\psi_{A, B}$ by the following standard diagram as in ~\cite{akmaj, maj},
\begin{equation} \label{si} \notag
\xy /r2.5pc/:,  
    (0,0)="x1",       
  "x1"+(1,0)="x2",  , 
  "x1"+(0,-.5)="x3",
  "x2"+(0,-.5)="x4",
  "x1";"x2"**\dir{}?(.5)="m", 
  "m" + (0,-.25)="c", 
  "x1";"x4"**\crv{ "x1"+(0,-.25) & "c" & "x4"+(0,.25)},  
  "x2";"c"+<3pt,2pt>="c2"**\crv{"x2"+(0,-.25) & "c2"},
  "x3";"c"+<-3pt,-2pt>="c1"**\crv{"x3"+(0,.25) & "c1"}, 
  "c"+(0,.10)*!<0pt,-6pt>{\psi_{A,B}},  
  "x1"+(0,.5)*!<0pt,6pt>{A},  
  "x2"+(0,.5)*!<0pt,6pt>{B},    
  "x3"-(0,.5)*!<0pt,-6pt>{B},  
  "x4"-(0,.5)*!<0pt,-6pt>{A},  
\endxy  
\end{equation}
for any $A$ and $B$ in $\mathcal{C}$, then the naturality of $\psi$ can be visualized by the following identity:
\begin{equation} \label{naturalityofsi} 
\xy /r2.5pc/:, 
(3,0)="L11",
(4,0)="L12",
  "L11"="x1", 
  "x1"+(0,-.5)="x2",
  "x1";"x2"**\dir{}?(.5)="c", 
  "x1";"x2"**\dir{-},   
  "c"-(0,0)*!<-4pt,-4pt>{\bullet f },  
  "x1"+(0,.5)*!<0pt,6pt>{A},  
 "x2"="L21",
 "L12"="x1", 
  "x1"+(0,-.5)="x2",
  "x1";"x2"**\dir{}?(.5)="c", 
  "x1";"x2"**\dir{-},    
  "c"-(0,0)*!<-4pt,-2pt>{\bullet g },  
  "x1"+(0,.5)*!<0pt,6pt>{B},  
"x2"="L22",
 "L21"="x1",       
  "L22"="x2",  , 
  "x1"+(0,-.5)="x3",
  "x2"+(0,-.5)="x4",
  "x1";"x2"**\dir{}?(.5)="m", 
  "m" + (0,-.25)="c", 
  "x1";"x4"**\crv{ "x1"+(0,-.25) & "c" & "x4"+(0,.25)},  
  "x2";"c"+<3pt,2pt>="c2"**\crv{"x2"+(0,-.25) & "c2"},
  "x3";"c"+<-3pt,-2pt>="c1"**\crv{"x3"+(0,.25) & "c1"}, 
  "c"+(0,.10)*!<0pt,-6pt>{\psi_{A',B'}},  
  "x3"-(0,.5)*!<0pt,-6pt>{B'},  
 "x4"-(0,.5)*!<0pt,-6pt>{A'}, 
"L22"+(1,0)*!<0pt,0pt>{=}, 
    (6,0)="x1",       
  "x1"+(1,0)="x2",  , 
  "x1"+(0,-.5)="x3",
  "x2"+(0,-.5)="x4",
  "x1";"x2"**\dir{}?(.5)="m", 
  "m" + (0,-.25)="c", 
  "x1";"x4"**\crv{ "x1"+(0,-.25) & "c" & "x4"+(0,.25)},  
  "x2";"c"+<3pt,2pt>="c2"**\crv{"x2"+(0,-.25) & "c2"},
  "x3";"c"+<-3pt,-2pt>="c1"**\crv{"x3"+(0,.25) & "c1"}, 
  "c"+(0,.10)*!<0pt,-6pt>{\psi_{A,B}},  
  "x1"+(0,.5)*!<0pt,6pt>{A},  
  "x2"+(0,.5)*!<0pt,6pt>{B},    
   "x3"="x1", 
  "x1"+(0,-.5)="x2",
  "x1";"x2"**\dir{}?(.5)="c", 
  "x1";"x2"**\dir{-},    
  "c"-(0,0)*!<-4pt,0pt>{\bullet g },  
  "x2"-(0,.5)*!<0pt,-6pt>{B'},   
 "x4"="x1", 
  "x1"+(0,-.5)="x2",
  "x1";"x2"**\dir{}?(.5)="c", 
  "x1";"x2"**\dir{-},    
  "c"-(0,0)*!<-4pt,0pt>{\bullet f },  
  "x2"-(0,.5)*!<0pt,-6pt>{A'},   
 \endxy   
\end{equation} 
for any two morphisms $f: A \to A'$ and $f: B \to B'$ in $\mathcal{C}$. 
Notice that in special case when for example $B=I$, since  
\begin{equation}\label{when B is I}
\psi_{I,A} = \psi_{A,I} = id_A, \quad \quad \forall A \in \mathcal{C}  ,
\end{equation} 
we have
\begin{equation}\label{naturalityofsi when B is I}
\xy /r2.5pc/:, 
(1,0)="L11",
(2,0)="L12",
  "L11"="x1", 
  "x1"+(0,-.5)="x2",
  "x1";"x2"**\dir{}?(.5)="c", 
  "x1";"x2"**\dir{-},   
  "c"-(0,0)*!<-4pt,-2pt>{\bullet f },  
  "x1"+(0,.5)*!<0pt,6pt>{A},  
 "x2"="L21",
 "L12"="x1", 
  "x1"+(0,-.5)="x2",
  "x1";"x2"**\dir{}?(.5)="c", 
  "x1";"x2"**\dir{-},    
  "c"-(0,0)*!<-4pt,0pt>{\bullet g },  
  "x1"+(0,0)*!<0pt,-2pt>{\circ},  
 "x2"="L22",
 "L21"="x1",       
  "L22"="x2",  , 
  "x1"+(0,-.5)="x3",
  "x2"+(0,-.5)="x4",
  "x1";"x2"**\dir{}?(.5)="m", 
  "m" + (0,-.25)="c", 
  "x1";"x4"**\crv{ "x1"+(0,-.25) & "c" & "x4"+(0,.25)},  
  "x2";"c"+<3pt,2pt>="c2"**\crv{"x2"+(0,-.25) & "c2"},
  "x3";"c"+<-3pt,-2pt>="c1"**\crv{"x3"+(0,.25) & "c1"}, 
  "c"+(0,.10)*!<0pt,-6pt>{\psi_{A',B'}},  
  "x3"-(0,.5)*!<0pt,-6pt>{B'},  
 "x4"-(0,.5)*!<0pt,-6pt>{A'},  
  "L22"+(1,0)*!<0pt,0pt>{=}, 
(4,0)="L11",  
(5,0)="L12",   
  "L11"+(0,0)*!<0pt,-2pt>{\circ}, 
  "L12"+(0,.5)*!<0pt,6pt>{A}, 
 "L11"="x1", 
  "x1"+(0,-.5)="x2",
  "x1";"x2"**\dir{}?(.5)="c", 
  "x1";"x2"**\dir{-},   
  "c"-(0,0)*!<-4pt,0pt>{\bullet g },  
  "x2"-(0,.5)*!<0pt,-6pt>{B'},   
 "L12"="x1", 
  "x1"+(0,-.5)="x2",
  "x1";"x2"**\dir{}?(.5)="c", 
  "x1";"x2"**\dir{-},   
  "c"-(0,0)*!<-4pt,0pt>{\bullet f },  
  "x2"-(0,.5)*!<0pt,-6pt>{A'},   
 \endxy  
\end{equation}
In other cases when $A,A'$ or $B'$ is $I$ the naturality identities can be 
simplified in a similar way. 

A braiding is called a \emph{ symmetry} if we have
\[\psi_{B, A}\circ \psi_{A, B}=id_{A\otimes B}, \]
or in terms of braiding diagrams:

\begin{equation}\label{symmetry}
\xy /r2.5pc/:,
    (0,0)="x1",  "x1"+(1,0)="x2", "x1"+(0,-.5)="x3", "x2"+(0,-.5)="x4", 
  "x1";"x2"**\dir{}?(.5)="m", "m" + (0,-.25)="c", 
  "x1";"x4"**\crv{ "x1"+(0,-.25) & "c" & "x4"+(0,.25)},  
  "x2";"c"+<3pt,2pt>="c2"**\crv{"x2"+(0,-.25) & "c2"},
  "x3";"c"+<-3pt,-2pt>="c1"**\crv{"x3"+(0,.25) & "c1"}, 
  "c"+(0,.10)*!<0pt,-6pt>{\psi_{A,B}}, "x1"+(0,.5)*!<0pt,6pt>{A}, "x2"+(0,.5)*!<0pt,6pt>{B},    
  "x3" ="L11", "x4" ="L12", "L11"="x1", "L12"="x2", "x1"+(0,-.5)="x3", "x2"+(0,-.5)="x4",
  "x1";"x2"**\dir{}?(.5)="m", "m" + (0,-.25)="c", 
  "x1";"x4"**\crv{ "x1"+(0,-.25) & "c" & "x4"+(0,.25)},  
  "x2";"c"+<3pt,2pt>="c2"**\crv{"x2"+(0,-.25) & "c2"},
  "x3";"c"+<-3pt,-2pt>="c1"**\crv{"x3"+(0,.25) & "c1"}, 
  "c"+(0,.10)*!<0pt,-6pt>{\psi_{B,A}}, "x3"-(0,.5)*!<0pt,-6pt>{A}, "x4"-(0,.5)*!<0pt,-6pt>{B},  
"x3" ="L11", "x4" ="L12", 
  "L12"+(1,.5)*!<0pt,0pt>{=}, 
(3,0)="L11";(3,-1)="L21" **\dir{-},  (4,0)="L12";(4,-1)="L22" **\dir{-},  
  "L11"+(0,.5)*!<0pt,6pt>{A}, "L12"+(0,.5)*!<0pt,6pt>{B}, "L21"-(0,.5)*!<0pt,-6pt>{A}, "L22"-(0,.5)*!<0pt,-6pt>{B},  
\endxy   
\end{equation}
for all objects $A$ and $B$ of $\mathcal{C}$. Sometimes we just write $\psi^2
= id $ to signify the symmetry condition. A \emph{symmetric monoidal
category} is a monoidal category endowed with a symmetry.

A monoidal category is called \emph{strict} if its associativity and
unit isomorphisms  are in fact equalities. By  a theorem of  Mac
Lane ~\cite{ml} (cf.\ also ~\cite{maj}), any (braided) monoidal
category is \emph{ monoidal equivalent} to a  (braided) strict
monoidal category  in which $a,l$ and $r$ are just equalities and
the above commuting diagrams are reduced to the following
equalities:
\begin{gather*} 
(A \otimes B) \otimes C = A \otimes (B \otimes C)\\
I \otimes A = A \otimes I = A\\
\psi_{A,B \otimes C} =(id_B \otimes \psi_{A,C})(\psi_{A,B} \otimes id_C)\\
\psi_{A \otimes B, C} = (\psi_{A,C} \otimes id_B )(id_A \otimes \psi_{B,C}) 
\end{gather*}
for all objects $A,B,C$ of $\mathcal{C}$. The last identity in terms of braided diagrams is as follows:

\begin{equation}\label{si_AB C} \notag
\xy /r2.5pc/:,  
(3,0)="L11",
(4,0)="L12",
(5.5,0)="L13",
"L11"="a1",
"L12"="a2",
"L13"="a3",
"a1" + (0,-1)="b1",
"a2" + (0,-1)="b2",
"a3" + (0,-1)="b3",
"a1";"a2"**\dir{}?(.5)="m12", 
  "m12" + (0,-.25)="c12", 
"a2";"a3"**\dir{}?(.5)="m23", 
  "m23" + (0,-.25)="c23", 
  "a1";"b2"**\crv{ "a1"+(0,-.25) & "c12" & "b2"+(0,.25)},  
  "a2";"b3"**\crv{ "a2"+(0,-.25) & "c23" & "b3"+(0,.25)}, 
                "b1";"c12"+<-0pt,-11pt>="c21"**\crv{"b1"+(0,.25) & "c21"},         
               "a3";"c23"+<-0pt,-2pt>="c32"**\crv{"a3"+(0,-.25) & "c32"},
"c12" +(.35,-.20);"c23" +(-.25,-.10)**\dir{-}, 
"a1"+(0,.5)*!<0pt,6pt>{A}, 
"a2"+(0,.5)*!<0pt,6pt>{B}, 
"a3"+(0,.5)*!<0pt,6pt>{C}, 
"b1"-(0,.5)*!<0pt,-6pt>{C},
"b2"-(0,.5)*!<0pt,-6pt>{A}, 
"b3"-(0,.5)*!<0pt,-6pt>{B}, 
  "m23"-(0,0)*!<0pt,-0pt>{\si_{A\otimes B,C}},  
"a3"+(.5,-.25)*!<0pt,0pt>{=}, 
(6.5,0)="L11",
(7.5,0)="L12",
(8.5,0)="L13",
"L11"="a1",
"L12"="a2",
"L13"="a3",
"a1"+(0,.5)*!<0pt,6pt>{A}, 
"a2"+(0,.5)*!<0pt,6pt>{B}, 
"a3"+(0,.5)*!<0pt,6pt>{C}, 
"L11";"L11" +(0,-.5)= "L11"**\dir{-},
   "L12"="x1",       
  "L13"="x2",  , 
  "x1"+(0,-.5)="x3",
  "x2"+(0,-.5)="x4",
  "x1";"x2"**\dir{}?(.5)="m", 
  "m" + (0,-.25)="c", 
  "x1";"x4"**\crv{ "x1"+(0,-.25) & "c" & "x4"+(0,.25)},  
  "x2";"c"+<3pt,2pt>="c2"**\crv{"x2"+(0,-.25) & "c2"},
  "x3";"c"+<-3pt,-2pt>="c1"**\crv{"x3"+(0,.25) & "c1"}, 
  "c"+(0,.10)*!<0pt,-6pt>{\psi_{B,C}},  
"x3"="L12",
"x4"="L13",
   "L11"="x1",       
  "L12"="x2",  , 
  "x1"+(0,-.5)="x3",
  "x2"+(0,-.5)="x4",
  "x1";"x2"**\dir{}?(.5)="m", 
  "m" + (0,-.25)="c", 
  "x1";"x4"**\crv{ "x1"+(0,-.25) & "c" & "x4"+(0,.25)},  
  "x2";"c"+<3pt,2pt>="c2"**\crv{"x2"+(0,-.25) & "c2"},
  "x3";"c"+<-3pt,-2pt>="c1"**\crv{"x3"+(0,.25) & "c1"}, 
  "c"+(0,.10)*!<0pt,-6pt>{\psi_{A,C}},  
"x3"="L11",
"x4"="L12",
"L13";"L13" +(0,-.5)= "L13"**\dir{-},
"L11"-(0,.5)*!<0pt,-6pt>{C},
"L12"-(0,.5)*!<0pt,-6pt>{A}, 
"L13"-(0,.5)*!<0pt,-6pt>{B},   
\endxy  
\end{equation}
The other one is the same.

We note that if the original category is symmetric then its \emph{
strictification} is symmetric as well. This plays an important role
in our approach. In fact, using this result, we can safely assume
that our symmetric monoidal categories are strict and symmetric.
Working with strict categories drastically simplifies  the formalism
and that is what we shall do in this paper.

\begin{example} \label{ex cat of H,R mods} Let $(H,\,R =R_1 \otimes R_2)$ be a
quasitriangular Hopf algebra and $\mathcal{C}$ be the category of all left
$H$-modules. Then $\mathcal{C}$ is a braided monoidal abelian category. It
is symmetric if and only if  $R^{-1}=R_2 \otimes R_1$ ~\cite{maj}. Here
the monoidal  structure is defined by 
\[h \rhd (v \otimes w)=h^{(1)}\rhd v\otimes h^{(2)} \rhd w , \] 
and the braiding map $ \psi_{V \otimes W}$ acts by
\[ \psi_{V \otimes W}(v \otimes w):= (R_2 \rhd w \otimes R_1 \rhd  v) , \]  
for any $V$ and $W$ in $\mathcal{C}$, where $\rhd$ denotes the action of
$H$.  Throughout this paper we use  Sweedler's notation, with
summation understood, e.g. we write  $\Delta h =h^{(1)}\otimes h^{(2)}$ to
denote the comultiplication of Hopf algebras. 
\end{example}

\begin{example} In  a dual  manner if we  consider a co-quasitriangular Hopf
algebra $(H,\,R),$ then  the category of left  $H$-comodules is a
braided monoidal abelian category.  
\end{example}

\begin{example} \label{cat of CZ_2 mods} As a very special case of Example
~\eqref{ex cat of H,R mods}, let $H =\mathbb{C} \mathbb{Z}_2$ with
the non-trivial quasitriangular structure $R =R_1 \otimes R_2$ defined by

\begin{equation}\label{triang}\notag
 R := (\frac{1}{2})(1 \otimes 1 + 1 \otimes g + g \otimes 1 - g \otimes g) , 
\end{equation}
 where $g$ is the generator of the cyclic group $\mathbb{Z}_2$.
The category $\mathcal{C}=\mathbb{Z}_2$- Mod then is the category of super
vector spaces with even morphisms ~\cite{maj}. The braiding map $
\psi_{V \otimes W}$ for any $V=V_0 \oplus V_1$ and $W=W_0 \oplus W_1$ in
$\mathcal{C}$ acts as below:

\begin{equation}\label{brdsuper}
\psi_{V \otimes W}(v \otimes w) = (-1)^{|v||w|}~(w \otimes v).
\end{equation}

There is also a category of differential graded (DG) super vector
spaces whose objects are super complexes
$V_0\overset{d}{\underset{d}{\leftrightarrows}} V_1$, and its
morphisms are even chain maps. It is a braided monoidal category
with the same braiding map as ~\eqref{brdsuper}. 
\end{example}

\begin{remark} \label{cat of CZ_n mods} 
One can extend Example ~\eqref{cat of CZ_2 mods} to  $\mathbb{C} \mathbb{Z}_n$ for any $n > 2$ which
provides a good source of non-symmetric braided monoidal categories~\cite{maj}. 
\end{remark}

\begin{example} \label{YD cat} 
Let $H$ be a Hopf algebra over a field $k$ with
comultiplication $\Delta h =h^{(1)}\otimes h^{(2)}$ and a  bijective
antipode $S$. A left-left Yetter-Drinfeld (YD)
 H-module consist of a  vector space $V$,
a left  $H$-module structure on V ~\cite{rt,y}:
\[H \otimes V \to V, \quad \quad h\otimes v \mapsto hv,\]
and a  left  $H$-comodule structure on V:
\[V\to H\otimes V, \quad \quad v \mapsto  v_{(-1)}\otimes v_{(0)}.\]
The left action and coaction are supposed to satisfy the
Yetter-Drinfeld (YD)   compatibility condition:
\[(hv)_{(-1)} \otimes (hv)_{(0)} =h^{(1)}v_{(-1)}S(h^{(3)}) \otimes h^{(2)}v_{(0)},\]
for all $h\in H$ and $v\in V$.  The category of all YD H-modules is
called the Yetter-Drinfeld category of $H$, and is usually denoted
by ${}_{H}^{H}\mathcal{YD}$. It is a braided monoidal abelian
category with the braiding map:
\[\psi_{V \otimes W}(v \otimes w)=v_{(-1)}w \otimes v_{(0)}.\]
This category is in general not symmetric.  In fact the inverse of
the braiding is given by :
\[\psi^{-1}_{V \otimes W}(w \otimes v)=v_{(0)} \otimes S^{-1}(v_{(-1)})w.\]
\end{example}

\begin{definition}\label{Braided Hopf algebras} Let $\mathcal{C}$ be a strict braided
monoidal category. A Hopf algebra\\  $(H,m, \eta , \Delta , \varepsilon , S)$
in $\mathcal{C}$ consists of an object $H \in obj \mathcal{C}$,   morphisms $m: H
\otimes H \to H ~,~ \eta :I \to H ~,~ \Delta :H \to H \otimes H~,~ \varepsilon :H
\to I ~and~ S:H \to H$ called  multiplication, unit,
comultiplication, counit and antipode maps satisfying the relations:

\[ m(id\otimes m)=m(m\otimes id),\quad \text{associativity} \]
\begin{equation}\label{associativityofm} \notag
\xy /r2.5pc/:, 
(1,0)="L11",
(2,0)="L12",
(3,0)="L13",
"L11"+(0,.5)*!<0pt,6pt>{H}, 
"L12"+(0,.5)*!<0pt,6pt>{H},
"L13"+(0,.5)*!<0pt,6pt>{H},
  "L11";"L11" +(0,-.5)="L11"**\dir{-},
  "L12"="x1", 
  "L13"="x2",   
  "x1";"x2"**\dir{}?(.5)="m", 
  "m" + (0,-.25)="c", 
  "m" + (0,-.5)="x3", 
  "x1";"x3"**\crv{"x1"-(0,.25) & "c" & "x3"},
  "x2";"x3"**\crv{"x2"-(0,.25) & "c"  & "x3"}, 
"x3"="L12",
  "L11"="x1", 
  "L12"="x2",   
  "x1";"x2"**\dir{}?(.5)="m", 
  "m" + (0,-.25)="c", 
  "m" + (0,-.5)="x3", 
  "x1";"x3"**\crv{"x1"-(0,.25) & "c" & "x3"},
  "x2";"x3"**\crv{"x2"-(0,.25) & "c"  & "x3"}, 
"x3"="L11",
"L11"+(0,-.5)*!<0pt,-6pt>{H},  
  "x1"+(2.5,0)*!<0pt,0pt>{=}, 
(4,0)="L11",
(5,0)="L12",
(6,0)="L13",
"L11"+(0,.5)*!<0pt,6pt>{H}, 
"L12"+(0,.5)*!<0pt,6pt>{H},
"L13"+(0,.5)*!<0pt,6pt>{H},
  "L11"="x1", 
  "L12"="x2",   
  "x1";"x2"**\dir{}?(.5)="m", 
  "m" + (0,-.25)="c", 
  "m" + (0,-.5)="x3", 
  "x1";"x3"**\crv{"x1"-(0,.25) & "c" & "x3"},
  "x2";"x3"**\crv{"x2"-(0,.25) & "c"  & "x3"}, 
"x3"="L11",
  "L13";"L13" +(0,-.5)="L12"**\dir{-},
  "L11"="x1", 
  "L12"="x2",   
  "x1";"x2"**\dir{}?(.5)="m", 
  "m" + (0,-.25)="c", 
  "m" + (0,-.5)="x3", 
  "x1";"x3"**\crv{"x1"-(0,.25) & "c" & "x3"},
  "x2";"x3"**\crv{"x2"-(0,.25) & "c"  & "x3"}, 
"x3"="L11",
"L11"+(0,-.5)*!<0pt,-6pt>{H}, 
  "x1"+(2,0)*!<0pt,0pt>{=}, 
(7,0)="L11",
"L11"+(0,.5)*!<0pt,6pt>{H}, 
  "L11";"L11" +(0,-1)="L11"**\dir{-},
"L11"+(0,-.5)*!<0pt,-6pt>{H}, 
\endxy 
\end{equation}

\[m(\eta \otimes id)=m(id\otimes \eta)=id, \quad \text{unit}  \]
\begin{equation}\label{unitofm}\notag
\xy /r2.5pc/:, 
(1,0)="L11",
(2,0)="L12",
"L12"+(0,.5)*!<0pt,6pt>{H}, 
   "L11" + (0,-.25)="x1", 
  "x1"+(0,-.25)="x2",
  "x1";"x2"**\dir{}?(.5)="c", 
  "x1";"x2"**\dir{-},   
  "c"-(0,0)*!<-3pt,0pt>{\bullet \eta },  
  "x1"+(0,0)*!<0pt,-2pt>{\circ},
"x2"="L11",
  "L12";"L12" +(0,-.5)="L12"**\dir{-},
  "L11"="x1", 
  "L12"="x2",   
  "x1";"x2"**\dir{}?(.5)="m", 
  "m" + (0,-.25)="c", 
  "m" + (0,-.5)="x3", 
  "x1";"x3"**\crv{"x1"-(0,.25) & "c" & "x3"},
  "x2";"x3"**\crv{"x2"-(0,.25) & "c"  & "x3"}, 
"x3"="L11",
"L11"+(0,-.5)*!<0pt,-6pt>{H}, 
  "x1"+(1.5,0)*!<0pt,0pt>{=}, 
(3,0)="L11",
(4,0)="L12",
"L11"+(0,.5)*!<0pt,6pt>{H}, 
  "L11";"L11" +(0,-.5)="L11"**\dir{-},
   "L12" + (0,-.25)="x1", 
  "x1"+(0,-.25)="x2",
  "x1";"x2"**\dir{}?(.5)="c", 
  "x1";"x2"**\dir{-},   
   "c"-(0,0)*!<-3pt,0pt>{\bullet \eta },  
   "x1"+(0,0)*!<0pt,-2pt>{\circ},
"x2"="L12",
  "L11"="x1", 
  "L12"="x2",   
  "x1";"x2"**\dir{}?(.5)="m", 
  "m" + (0,-.25)="c", 
  "m" + (0,-.5)="x3", 
  "x1";"x3"**\crv{"x1"-(0,.25) & "c" & "x3"},
  "x2";"x3"**\crv{"x2"-(0,.25) & "c"  & "x3"}, 
"x3"="L11",
"L11"+(0,-.5)*!<0pt,-6pt>{H}, 
  "x1"+(1.5,0)*!<0pt,0pt>{=}, 
(5,0)="L11",
"L11"+(0,.5)*!<0pt,6pt>{H}, 
  "L11";"L11" +(0,-1)="L11"**\dir{-},
"L11"+(0,-.5)*!<0pt,-6pt>{H}, 
\endxy 
\end{equation}

\[(id \otimes \Delta) \Delta =(\Delta \otimes id)\Delta, \quad \text{coassociativity} \]
\begin{equation}\label{coassociativityofDel} \notag
\xy /r2.5pc/:, 
(1,0)="L11",
"L11"+(0,.5)*!<0pt,6pt>{H}, 
   "L11"="x1", 
  "x1"+(-.5,-.5)="x2",
  "x1"+(.5,-.5)="x3",
  "x1"+(0,-.25)="c",  
   "x1";"x2"**\crv{"x1" & "c" & "x2"+(0,.25)}, 
   "x1";"x3"**\crv{"x1" & "c" & "x3"+(0,.25)}, 
"x2"="L11",
"x3"="L12",
  "L11";"L11" +(0,-.5)="L11"**\dir{-},
   "L12"="x1", 
  "x1"+(-.5,-.5)="x2",
  "x1"+(.5,-.5)="x3",
  "x1"+(0,-.25)="c",  
   "x1";"x2"**\crv{"x1" & "c" & "x2"+(0,.25)}, 
   "x1";"x3"**\crv{"x1" & "c" & "x3"+(0,.25)}, 
"x2"="L12",
"x3"="L13",
"L11"+(0,-.5)*!<0pt,-6pt>{H}, 
"L12"+(0,-.5)*!<0pt,-6pt>{H}, 
"L13"+(0,-.5)*!<0pt,-6pt>{H}, 
  "x1"+(1,0)*!<0pt,0pt>{=}, 
(4,0)="L11",
"L11"+(0,.5)*!<0pt,6pt>{H}, 
   "L11"="x1", 
  "x1"+(-.5,-.5)="x2",
  "x1"+(.5,-.5)="x3",
  "x1"+(0,-.25)="c",  
   "x1";"x2"**\crv{"x1" & "c" & "x2"+(0,.25)}, 
   "x1";"x3"**\crv{"x1" & "c" & "x3"+(0,.25)}, 
"x2"="L11",
"x3"="L12",
   "L11"="x1", 
  "x1"+(-.5,-.5)="x2",
  "x1"+(.5,-.5)="x3",
  "x1"+(0,-.25)="c",  
   "x1";"x2"**\crv{"x1" & "c" & "x2"+(0,.25)}, 
   "x1";"x3"**\crv{"x1" & "c" & "x3"+(0,.25)}, 
"x2"="L21",
"x3"="L22",
  "L12";"L12" +(0,-.5)="L13"**\dir{-},
"L21"="L11",
"L22"="L12",
"L11"+(0,-.5)*!<0pt,-6pt>{H}, 
"L12"+(0,-.5)*!<0pt,-6pt>{H}, 
"L13"+(0,-.5)*!<0pt,-6pt>{H}, 	
\endxy 
\end{equation}

\[(\varepsilon \otimes id)\Delta =(id \otimes \varepsilon)\Delta =id, \quad \text{counit} \]
\begin{equation}\label{counityofDel} \notag 
\xy /r2.5pc/:, 
(1,0)="L11",
"L11"+(0,.5)*!<0pt,6pt>{H}, 
   "L11"="x1", 
  "x1"+(-.5,-.5)="x2",
  "x1"+(.5,-.5)="x3",
  "x1"+(0,-.25)="c",  
   "x1";"x2"**\crv{"x1" & "c" & "x2"+(0,.25)}, 
   "x1";"x3"**\crv{"x1" & "c" & "x3"+(0,.25)}, 
"x2"="L11",
"x3"="L12",
  "L11";"L11" +(0,-.5)="L11"**\dir{-},
   "L12"="x1", 
  "x1"+(0,-.25)="x2",
  "x1";"x2"**\dir{}?(.5)="c", 
  "x1";"x2"**\dir{-},   
  "c"-(0,0)*!<-3pt,0pt>{\bullet \varepsilon },  
  "x2"+(0,0)*!<0pt,2pt>{\circ},  
"L11"+(0,-.5)*!<0pt,-6pt>{H}, 
  "x1"+(.5,0)*!<0pt,0pt>{=}, 
(3,0)="L11",
"L11"+(0,.5)*!<0pt,6pt>{H}, 
   "L11"="x1", 
  "x1"+(-.5,-.5)="x2",
  "x1"+(.5,-.5)="x3",
  "x1"+(0,-.25)="c",  
   "x1";"x2"**\crv{"x1" & "c" & "x2"+(0,.25)}, 
   "x1";"x3"**\crv{"x1" & "c" & "x3"+(0,.25)}, 
"x2"="L11",
"x3"="L12",
   "L11"="x1", 
  "x1"+(0,-.25)="x2",
  "x1";"x2"**\dir{}?(.5)="c", 
  "x1";"x2"**\dir{-},   
  "c"-(0,0)*!<-3pt,0pt>{\bullet \varepsilon },  
  "x2"+(0,0)*!<0pt,2pt>{\circ},  
  "L12";"L12" +(0,-.5)="L11"**\dir{-},
"L11"+(0,-.5)*!<0pt,-6pt>{H}, 
  "x1"+(1.5,0)*!<0pt,0pt>{=}, 
(5,0)="L11",
"L11"+(0,.5)*!<0pt,6pt>{H}, 
  "L11";"L11" +(0,-1)="L11"**\dir{-},
"L11"+(0,-.5)*!<0pt,-6pt>{H}, 
\endxy 
\end{equation}

\[\Delta m = (m \otimes m)(id \otimes \psi \otimes id)(\Delta \otimes \Delta),\quad \text{compatibility}\]
\begin{equation}\label{compatibilityDelm}
\xy /r2.5pc/:, 
(1,0)="L11",
(2,0)="L12",
"L11"+(0,.5)*!<0pt,6pt>{H},
"L12"+(0,.5)*!<0pt,6pt>{H},
  "L11"="x1", 
  "L12"="x2",   
  "x1";"x2"**\dir{}?(.5)="m", 
  "m" + (0,-.25)="c", 
  "m" + (0,-.5)="x3", 
  "x1";"x3"**\crv{"x1"-(0,.25) & "c" & "x3"},
  "x2";"x3"**\crv{"x2"-(0,.25) & "c"  & "x3"}, 
"x3"="L11",
   "L11"="x1", 
  "x1"+(-.5,-.5)="x2",
  "x1"+(.5,-.5)="x3",
  "x1"+(0,-.25)="c",  
   "x1";"x2"**\crv{"x1" & "c" & "x2"+(0,.25)}, 
   "x1";"x3"**\crv{"x1" & "c" & "x3"+(0,.25)}, 
"x2"="L11",
"x3"="L12",
"L11"+(0,-.5)*!<0pt,-6pt>{H}, 
"L12"+(0,-.5)*!<0pt,-6pt>{H}, 
  "x1"+(1,0)*!<0pt,0pt>{=}, 
(4,0)="L11",
(6,0)="L12",
"L11"+(0,.5)*!<0pt,6pt>{H},
"L12"+(0,.5)*!<0pt,6pt>{H},
   "L11"="x1", 
  "x1"+(-.5,-.5)="x2",
  "x1"+(.5,-.5)="x3",
  "x1"+(0,-.25)="c",  
   "x1";"x2"**\crv{"x1" & "c" & "x2"+(0,.25)}, 
   "x1";"x3"**\crv{"x1" & "c" & "x3"+(0,.25)}, 
"x2"="L21",
"x3"="L22",
   "L12"="x1", 
  "x1"+(-.5,-.5)="x2",
  "x1"+(.5,-.5)="x3",
  "x1"+(0,-.25)="c",  
   "x1";"x2"**\crv{"x1" & "c" & "x2"+(0,.25)}, 
   "x1";"x3"**\crv{"x1" & "c" & "x3"+(0,.25)}, 
"x2"="L23",
"x3"="L24",
"L21"="L11",
"L22"="L12",
"L23"="L13",
"L24"="L14",
  "L11";"L11" +(0,-.5)="L11"**\dir{-},
  "L12"="x1",       
  "L13"="x2",  , 
  "x1"+(0,-.5)="x3",
  "x2"+(0,-.5)="x4",
  "x1";"x2"**\dir{}?(.5)="m", 
  "m" + (0,-.25)="c", 
  "x1";"x4"**\crv{ "x1"+(0,-.25) & "c" & "x4"+(0,.25)},  
  "x2";"c"+<3pt,2pt>="c2"**\crv{"x2"+(0,-.25) & "c2"},
  "x3";"c"+<-3pt,-2pt>="c1"**\crv{"x3"+(0,.25) & "c1"}, 
  "c"+(0,.10)*!<0pt,-6pt>{\psi},
"x3"="L12",
"x4"="L13",
  "L14";"L14" +(0,-.5)="L14"**\dir{-},
  "L11"="x1", 
  "L12"="x2",   
  "x1";"x2"**\dir{}?(.5)="m", 
  "m" + (0,-.25)="c", 
  "m" + (0,-.5)="x3", 
  "x1";"x3"**\crv{"x1"-(0,.25) & "c" & "x3"},
  "x2";"x3"**\crv{"x2"-(0,.25) & "c"  & "x3"}, 
"x3"="L11",
  "L13"="x1", 
  "L14"="x2",   
  "x1";"x2"**\dir{}?(.5)="m", 
  "m" + (0,-.25)="c", 
  "m" + (0,-.5)="x3", 
  "x1";"x3"**\crv{"x1"-(0,.25) & "c" & "x3"},
  "x2";"x3"**\crv{"x2"-(0,.25) & "c"  & "x3"}, 
"x3"="L12",
"L11"+(0,-.5)*!<0pt,-6pt>{H}, 
"L12"+(0,-.5)*!<0pt,-6pt>{H}, 
\endxy 
\end{equation}

\[\Delta \eta = \eta \otimes \eta ,~~\varepsilon m=\varepsilon \otimes \varepsilon ,~~  \varepsilon \eta =id_I  \]
\begin{equation}\label{compatibility epn eta} 
\xy /r2pc/:, 
(1,0)="L11",
   "L11" + (0,-.25)="x1", 
  "x1"+(0,-.25)="x2",
  "x1";"x2"**\dir{}?(.5)="c", 
  "x1";"x2"**\dir{-},   
  "c"-(0,0)*!<-3pt,0pt>{\bullet \eta },  
  "x1"+(0,0)*!<0pt,-2pt>{\circ},
"x2"="L11",
   "L11"="x1", 
  "x1"+(-.5,-.5)="x2",
  "x1"+(.5,-.5)="x3",
  "x1"+(0,-.25)="c",  
   "x1";"x2"**\crv{"x1" & "c" & "x2"+(0,.25)}, 
   "x1";"x3"**\crv{"x1" & "c" & "x3"+(0,.25)}, 
"x2"="L11",
"x3"="L12",
"L11"+(0,-.5)*!<0pt,-6pt>{H}, 
"L12"+(0,-.5)*!<0pt,-6pt>{H}, 
  "x1"+(1,0)*!<0pt,0pt>{=}, 
(3,0)="L11",
(4,0)="L12",
   "L11" + (0,-.25)="x1", 
  "x1"+(0,-.75)="x2",
  "x1";"x2"**\dir{}?(.5)="c", 
  "x1";"x2"**\dir{-},   
  "c"-(0,0)*!<-3pt,0pt>{\bullet \eta },  
  "x1"+(0,0)*!<0pt,-2pt>{\circ},
"x2"="L11",
   "L12" + (0,-.25)="x1", 
  "x1"+(0,-.75)="x2",
  "x1";"x2"**\dir{}?(.5)="c", 
  "x1";"x2"**\dir{-},   
  "c"-(0,0)*!<-3pt,0pt>{\bullet \eta },  
  "x1"+(0,0)*!<0pt,-2pt>{\circ},
"x2"="L12",
"L11"+(0,-.5)*!<0pt,-6pt>{H}, 
"L12"+(0,-.5)*!<0pt,-6pt>{H}, 
\endxy 
\quad \quad
\xy /r2pc/:, 
(1,0)="L11",
(2,0)="L12",
"L11"+(0,.5)*!<0pt,6pt>{H},
"L12"+(0,.5)*!<0pt,6pt>{H},
  "L11"="x1", 
  "L12"="x2",   
  "x1";"x2"**\dir{}?(.5)="m", 
  "m" + (0,-.25)="c", 
  "m" + (0,-.5)="x3", 
  "x1";"x3"**\crv{"x1"-(0,.25) & "c" & "x3"},
  "x2";"x3"**\crv{"x2"-(0,.25) & "c"  & "x3"}, 
"x3"="L11",
   "L11"="x1", 
  "x1"+(0,-.25)="x2",
  "x1";"x2"**\dir{}?(.5)="c", 
  "x1";"x2"**\dir{-},   
   "c"-(0,0)*!<-3pt,0pt>{\bullet \varepsilon },  
  "x2"+(0,0)*!<0pt,2pt>{\circ},  
  "x1"+(1,0)*!<0pt,0pt>{=}, 
(3,0)="L11",
(4,0)="L12",
"L11"+(0,.5)*!<0pt,6pt>{H},
"L12"+(0,.5)*!<0pt,6pt>{H},
   "L11"="x1", 
  "x1"+(0,-.75)="x2",
  "x1";"x2"**\dir{}?(.5)="c", 
  "x1";"x2"**\dir{-},   
  "c"-(0,0)*!<-3pt,0pt>{\bullet \varepsilon },  
  "x2"+(0,0)*!<0pt,2pt>{\circ},  
   "L12"="x1", 
  "x1"+(0,-.75)="x2",
  "x1";"x2"**\dir{}?(.5)="c", 
  "x1";"x2"**\dir{-},   
  "c"-(0,0)*!<-3pt,0pt>{\bullet \varepsilon },  
  "x2"+(0,0)*!<0pt,2pt>{\circ},  
\endxy 
\quad \quad
\xy /r2pc/:, 
(1,0)="L11",
   "L11" + (0,-.25)="x1", 
  "x1"+(0,-.25)="x2",
  "x1";"x2"**\dir{}?(.5)="c", 
  "x1";"x2"**\dir{-},   
  "c"-(0,0)*!<-3pt,0pt>{\bullet \eta },  
  "x1"+(0,0)*!<0pt,-2pt>{\circ},
"x2"="L11",
   "L11"="x1", 
  "x1"+(0,-.25)="x2",
  "x1";"x2"**\dir{}?(.5)="c", 
  "x1";"x2"**\dir{-},   
  "c"-(0,0)*!<-3pt,0pt>{\bullet \varepsilon },  
  "x2"+(0,0)*!<0pt,2pt>{\circ},  
  "x1"+(.5,0)*!<0pt,0pt>{=}, 
(2,-.25)="L11",
"L11"+(0,0)*!<0pt,-2pt>{\circ},
  "L11";"L11" +(0,-.5)="L11"**\dir{-},
"L11"+(0,0)*!<0pt,2pt>{\circ}, 												
\endxy
\end{equation}

\[m(S \otimes id)\Delta =m(id \otimes S)\Delta = \eta \varepsilon, \text{antipode} \]
\begin{equation}\label{antipode m(1,S)Del = eta epn}
\xy /r2.5pc/:, 
(1,0)="L11",
  "L11"+(0,.5)*!<0pt,6pt>{H},
   "L11"="x1", 
  "x1"+(-.5,-.5)="x2",
  "x1"+(.5,-.5)="x3",
  "x1"+(0,-.25)="c",  
   "x1";"x2"**\crv{"x1" & "c" & "x2"+(0,.25)}, 
   "x1";"x3"**\crv{"x1" & "c" & "x3"+(0,.25)}, 
"x2"= "L11",
"x3"= "L12",
 "L11"="x1", 
  "x1"+(0,-.5)="x2",
  "x1";"x2"**\dir{}?(.5)="c", 
  "x1";"x2"**\dir{-},  
  "c"-(0,0)*!<-4pt,0pt>{\bullet S }, 
"x2"= "L11",
"L12";"L12" + (0,-.5) = "L12" **\dir{-},  
   "L11" ="x1", 
    "L12"="x2",   
  "x1";"x2"**\dir{}?(.5)="m", 
  "m" + (0,-.25)="c", 
  "m" + (0,-.5)="x3", 
  "x1";"x3"**\crv{"x1"-(0,.25) & "c" & "x3"},
  "x2";"x3"**\crv{"x2"-(0,.25) & "c"  & "x3"}, 
"x3"="L11",
 "L11"-(0,.5)*!<0pt,-6pt>{H},  
  "L12"+(.5,.5)*!<0pt,0pt>{=}, 
 (3,0)="L11",
  "L11"+(0,.5)*!<0pt,6pt>{H},
   "L11"="x1", 
  "x1"+(-.5,-.5)="x2",
  "x1"+(.5,-.5)="x3",
  "x1"+(0,-.25)="c",  
   "x1";"x2"**\crv{"x1" & "c" & "x2"+(0,.25)}, 
   "x1";"x3"**\crv{"x1" & "c" & "x3"+(0,.25)}, 
"x2"= "L11",
"x3"= "L12",
"L11";"L11" + (0,-.5) = "L11" **\dir{-},  
"L12"="x1", 
  "x1"+(0,-.5)="x2",
  "x1";"x2"**\dir{}?(.5)="c", 
  "x1";"x2"**\dir{-},  
  "c"-(0,0)*!<-4pt,0pt>{\bullet S }, 
"x2"= "L12",
   "L11" ="x1", 
    "L12"="x2",   
  "x1";"x2"**\dir{}?(.5)="m", 
  "m" + (0,-.25)="c", 
  "m" + (0,-.5)="x3", 
  "x1";"x3"**\crv{"x1"-(0,.25) & "c" & "x3"},
  "x2";"x3"**\crv{"x2"-(0,.25) & "c"  & "x3"}, 
"x3"="L11",
 "L11"-(0,.5)*!<0pt,-6pt>{H},  
  "L12"+(.5,.5)*!<0pt,0pt>{=}, 
(5,0)="L11",
  "L11"+(0,.5)*!<0pt,6pt>{H},
  "L11"="x1", 
  "x1"+(0,-.25)="x2",
  "x1";"x2"**\dir{}?(.5)="c", 
  "x1";"x2"**\dir{-},   
  "c"-(0,0)*!<-3pt,0pt>{\bullet \varepsilon },  
   "x2"+(0,0)*!<0pt,2pt>{\circ},  
"x2" +(0,-.5)="x1", 
  "x1"+(0,-.25)="x2",
  "x1";"x2"**\dir{}?(.5)="c", 
  "x1";"x2"**\dir{-},   
  "c"-(0,0)*!<-3pt,0pt>{\bullet \eta },  
  "x1"+(0,0)*!<0pt,-2pt>{\circ}, 
"x2" ="L11",
 "L11"-(0,.5)*!<0pt,-6pt>{H},  
\endxy 
\end{equation}
By a braided Hopf algebra we mean a Hopf algebra in a braided
monoidal category.  \end{definition}

Notice that relation ~\eqref{compatibilityDelm}, which expresses the algebra 
property of the comultiplication, is the only relation that involves 
the braiding map $\psi$.

\begin{lemma} It is easy to prove that:
\begin{multline*}  
\Delta^2 m = (m \ot m\ot m)(id \ot \si \ot id \ot id \ot id)\\
(id \ot id \ot id \ot \si \ot id)(id \ot id \ot \si \ot id \ot id)(\Delta^2 \ot \Delta^2).
\end{multline*}
\begin{equation}\label{Del2 m}
\xy /r2.5pc/:, 
(0,0)="L11",
(1,0)="L12",
   "L11" ="x1", 
    "L12"="x2",   
  "x1";"x2"**\dir{}?(.5)="m", 
  "m" + (0,-.25)="c", 
  "m" + (0,-.5)="x3", 
  "x1";"x3"**\crv{"x1"-(0,.25) & "c" & "x3"},
  "x2";"x3"**\crv{"x2"-(0,.25) & "c"  & "x3"}, 
  "x1"+(0,.5)*!<0pt,6pt>{H},  
 "x2"+(0,.5)*!<0pt,6pt>{H},
"x3"="L31",
   "L31"="x1", 
   "x1"+(-,-.5)="x2",
  "x1"+(0,-.5)="x3",
  "x1"+(1,-.5)="x4",
  "x1"+(0,-.25)="c",  
   "x1";"x2"**\crv{"x1" & "c" & "x2"+(0,.25)}, 
   "x1";"x3"**\crv{"x1" & "c" & "x3"+(0,.25)}, 
   "x1";"x4"**\crv{"x1" & "c" & "x4"+(0,.25)}, 
"x2"="L31",
"x3"="L32",
"x4"="L33",
"L31"="L21",
"L32"="L22",
"L33"="L23",
  "L31"-(0,.5)*!<0pt,-6pt>{H},   
  "L32"-(0,.5)*!<0pt,-6pt>{H},  
  "L33"-(0,.5)*!<0pt,-6pt>{H}, 
  "L22"+(1.5,.5)*!<0pt,0pt>{=}, 
   (3,0)="x1", 
   "x1"+(-.25,-.5)="x2",
  "x1"+(0,-.5)="x3",
  "x1"+(.25,-.5)="x4",
  "x1"+(0,-.25)="c",  
   "x1";"x2"**\crv{"x1" & "c" & "x2"+(0,.25)}, 
   "x1";"x3"**\crv{"x1" & "c" & "x3"+(0,.25)}, 
   "x1";"x4"**\crv{"x1" & "c" & "x4"+(0,.25)}, 
  "x1"+(0,.5)*!<0pt,6pt>{H},  
"x2"="L31",
"x3"="L32",
"x4"="L33",
   (4,0)="x1", 
   "x1"+(-.25,-.5)="x2",
  "x1"+(0,-.5)="x3",
  "x1"+(.25,-.5)="x4",
  "x1"+(0,-.25)="c",  
   "x1";"x2"**\crv{"x1" & "c" & "x2"+(0,.25)}, 
   "x1";"x3"**\crv{"x1" & "c" & "x3"+(0,.25)}, 
   "x1";"x4"**\crv{"x1" & "c" & "x4"+(0,.25)}, 
  "x1"+(0,.5)*!<0pt,6pt>{H},  
"x2"="L34",
"x3"="L35",
"x4"="L36",
"L31"="L21",
"L32"="L22",
"L33"="L23",
"L34"="L24",
"L35"="L25",
"L36"="L26",
  "L21"="x1", 
  "x1"+(0,-.5)="x2",
  "x1";"x2"**\dir{-},  
"x2"= "L31",
  "L22"="x1", 
  "x1"+(0,-.5)="x2",
  "x1";"x2"**\dir{-},   
"x2"= "L32",
   "L23"="x1",       
   "L24"="x2",  
  "x1"+(0,-.5)="x3",
  "x2"+(0,-.5)="x4",
  "x1";"x2"**\dir{}?(.5)="m", 
  "m" + (0,-.25)="c", 
  "x1";"x4"**\crv{ "x1"+(0,-.25) & "c" & "x4"+(0,.25)},  
  "x2";"c"+<3pt,2pt>="c2"**\crv{"x2"+(0,-.25) & "c2"},
  "x3";"c"+<-3pt,-2pt>="c1"**\crv{"x3"+(0,.25) & "c1"}, 
 "c"+(0,.10)*!<0pt,-6pt>{\psi},  
"x3"= "L33",
"x4"= "L34",
  "L25"="x1", 
  "x1"+(0,-.5)="x2",
  "x1";"x2"**\dir{-},   
"x2"= "L35",
  "L26"="x1", 
  "x1"+(0,-.5)="x2",
  "x1";"x2"**\dir{-},   
"x2"= "L36",
"L31"="L21",
"L32"="L22",
"L33"="L23",
"L34"="L24",
"L35"="L25",
"L36"="L26",
  "L21"="x1", 
  "x1"+(0,-.5)="x2",
  "x1";"x2"**\dir{-},   
"x2"= "L31",
  "L22"="x1", 
  "x1"+(0,-.5)="x2",
  "x1";"x2"**\dir{-},   
"x2"= "L32",
  "L23"="x1", 
  "x1"+(0,-.5)="x2",
  "x1";"x2"**\dir{-},   
"x2"= "L33",
   "L24"="x1",      
   "L25"="x2",  
  "x1"+(0,-.5)="x3",
  "x2"+(0,-.5)="x4",
  "x1";"x2"**\dir{}?(.5)="m", 
  "m" + (0,-.25)="c", 
  "x1";"x4"**\crv{ "x1"+(0,-.25) & "c" & "x4"+(0,.25)},  
  "x2";"c"+<3pt,2pt>="c2"**\crv{"x2"+(0,-.25) & "c2"},
  "x3";"c"+<-3pt,-2pt>="c1"**\crv{"x3"+(0,.25) & "c1"}, 
  "c"+(0,.10)*!<0pt,-6pt>{\psi},  
"x3"= "L34",
"x4"= "L35",
  "L26"="x1", 
  "x1"+(0,-.5)="x2",
  "x1";"x2"**\dir{-},  
"x2"= "L36",
"L31"="L21",
"L32"="L22",
"L33"="L23",
"L34"="L24",
"L35"="L25",
"L36"="L26",
  "L21"="x1", 
  "x1"+(0,-.5)="x2",
  "x1";"x2"**\dir{-},   
"x2"= "L31",
   "L22"="x1",       
   "L23"="x2",  
  "x1"+(0,-.5)="x3",
  "x2"+(0,-.5)="x4",
  "x1";"x2"**\dir{}?(.5)="m", 
  "m" + (0,-.25)="c", 
  "x1";"x4"**\crv{ "x1"+(0,-.25) & "c" & "x4"+(0,.25)},  
  "x2";"c"+<3pt,2pt>="c2"**\crv{"x2"+(0,-.25) & "c2"},
  "x3";"c"+<-3pt,-2pt>="c1"**\crv{"x3"+(0,.25) & "c1"}, 
  "c"+(0,.10)*!<0pt,-6pt>{\psi},  
"x3"= "L32",
"x4"= "L33",
  "L24"="x1", 
  "x1"+(0,-.5)="x2",
  "x1";"x2"**\dir{-},  
"x2"= "L34",
  "L25"="x1", 
  "x1"+(0,-.5)="x2",
  "x1";"x2"**\dir{-},   
"x2"= "L35",
  "L26"="x1",
  "x1"+(0,-.5)="x2",
  "x1";"x2"**\dir{-},   
"x2"= "L36",
"L31"="L21",
"L32"="L22",
"L33"="L23",
"L34"="L24",
"L35"="L25",
"L36"="L26",
   "L21" ="x1", 
    "L22"="x2",   
  "x1";"x2"**\dir{}?(.5)="m", 
  "m" + (0,-.25)="c", 
  "m" + (0,-.5)="x3", 
  "x1";"x3"**\crv{"x1"-(0,.25) & "c" & "x3"},
  "x2";"x3"**\crv{"x2"-(0,.25) & "c"  & "x3"}, 
"x3"="L31",
   "L23" ="x1", 
    "L24"="x2",   
  "x1";"x2"**\dir{}?(.5)="m", 
  "m" + (0,-.25)="c", 
  "m" + (0,-.5)="x3", 
  "x1";"x3"**\crv{"x1"-(0,.25) & "c" & "x3"},
  "x2";"x3"**\crv{"x2"-(0,.25) & "c"  & "x3"}, 
"x3"="L32",
   "L25" ="x1", 
    "L26"="x2",   
  "x1";"x2"**\dir{}?(.5)="m", 
  "m" + (0,-.25)="c", 
  "m" + (0,-.5)="x3", 
  "x1";"x3"**\crv{"x1"-(0,.25) & "c" & "x3"},
  "x2";"x3"**\crv{"x2"-(0,.25) & "c"  & "x3"}, 
"x3"="L33",
  "L31"-(0,.5)*!<0pt,-6pt>{H},  
  "L32"-(0,.5)*!<0pt,-6pt>{H},  
  "L33"-(0,.5)*!<0pt,-6pt>{H}, 
\endxy 
\end{equation}
\end{lemma}

\begin{example} \label{def of superhopfalg} Any Hopf algebra in $\mathbb{Z}_2
$-Mod of Example ~\eqref{cat of CZ_2 mods} is called a super Hopf
algebra. See section ~\eqref{super hpf alg} for more details.
Similarly, a Hopf algebra in the category of differential graded
super vector spaces is a differential graded super Hopf algebra.
\end{example}

Notice that a super Hopf algebra is, in general, not a standard Hopf
algebra, since the comultiplication $\Delta$ is not an algebra map,
and the multiplication $m$ is not a coalgebra map in the standard
sense.

\begin{example} \label{hofalg in YD cat} For any $V$ in  ${}_{H}^{H}\mathcal{YD}$ of Example
~\eqref{YD cat},  the tensor algebra $T(V)$ is a braided Hopf algebra
in  ${}_{H}^{H}\mathcal{YD}$. Its   comultiplication,  counit, and antipode are
defined by $\Delta(v)=1\otimes v + v\otimes 1,$ $\varepsilon(v)=0,$ and $S(v)=-v$,
for all $v$ in $V$. 
\end{example}

The following proposition shows that standard properties of Hopf algebras hold for braided Hopf algebras.
\begin{proposition}  \label{prps of S}
If $(H,m, \eta , \Delta , \varepsilon , S)$ is a braided Hopf algebra (in $\mathcal{C}$), then:
\[Sm=m \psi (S \otimes S) =m(S \otimes S)\psi , \]
\begin{equation}\label{S m}
\xy /r2.5pc/:,   
(0,0)="L11",
(1,0)="L12",
  "L11"+(0,.5)*!<0pt,6pt>{H},
  "L12"+(0,.5)*!<0pt,6pt>{H}, 
  "L11"="x1", 
  "L12"="x2",  
  "x1";"x2"**\dir{}?(.5)="m", 
  "m" + (0,-.25)="c", 
  "m" + (0,-.5)="x3", 
  "x1";"x3"**\crv{"x1"-(0,.25) & "c" & "x3"},
  "x2";"x3"**\crv{"x2"-(0,.25) & "c"  & "x3"}, 
"x3"="L21",
  "L21"="x1",
  "x1"+(0,-.5)="x2",
  "x1";"x2"**\dir{}?(.5)="c", 
  "x1";"x2"**\dir{-},  
  "c"-(0,0)*!<-4pt,0pt>{\bullet S }, 
"x2"="L21",
  "x2"-(0,.5)*!<0pt,-6pt>{H}, 
 "x2" + (1.5,.25)*!<0pt,-6pt>{=}, 
(3,0)="L11",
(4,0)="L12",
  "L11"+(0,.5)*!<0pt,6pt>{H},
  "L12"+(0,.5)*!<0pt,6pt>{H}, 
  "L11"="x1", 
  "x1"+(0,-.5)="x2",
  "x1";"x2"**\dir{}?(.5)="c", 
  "x1";"x2"**\dir{-},   
  "c"-(0,0)*!<-4pt,0pt>{\bullet S }, 
"x2"="L21",
  "L12"="x1", 
  "x1"+(0,-.5)="x2",
  "x1";"x2"**\dir{}?(.5)="c", 
  "x1";"x2"**\dir{-},   
  "c"-(0,0)*!<-4pt,0pt>{\bullet S }, 
"x2"="L22",
  "L21"="x1",       
  "L22"="x2",  , 
  "x1"+(0,-.5)="x3",
  "x2"+(0,-.5)="x4",
  "x1";"x2"**\dir{}?(.5)="m", 
  "m" + (0,-.25)="c", 
  "x1";"x4"**\crv{ "x1"+(0,-.25) & "c" & "x4"+(0,.25)},  
  "x2";"c"+<3pt,2pt>="c2"**\crv{"x2"+(0,-.25) & "c2"},
  "x3";"c"+<-3pt,-2pt>="c1"**\crv{"x3"+(0,.25) & "c1"}, 
  "c"+(0,.10)*!<0pt,-6pt>{\psi}, 
"x3"="L21",
"x4"="L22",
  "L21"="x1", 
  "L22"="x2",  
  "x1";"x2"**\dir{}?(.5)="m", 
  "m" + (0,-.25)="c", 
  "m" + (0,-.5)="x3", 
  "x1";"x3"**\crv{"x1"-(0,.25) & "c" & "x3"},
  "x2";"x3"**\crv{"x2"-(0,.25) & "c"  & "x3"}, 
"x3"="L21",
  "x3"-(0,.5)*!<0pt,-6pt>{H}, 
 "x2" + (1.25,.25)*!<0pt,-6pt>{=}, 
(6,0)="L11",
(7,0)="L12",
  "L11"+(0,.5)*!<0pt,6pt>{H},
  "L12"+(0,.5)*!<0pt,6pt>{H}, 
  "L11"="x1",       
  "L12"="x2",  , 
  "x1"+(0,-.5)="x3",
  "x2"+(0,-.5)="x4",
  "x1";"x2"**\dir{}?(.5)="m", 
  "m" + (0,-.25)="c", 
  "x1";"x4"**\crv{ "x1"+(0,-.25) & "c" & "x4"+(0,.25)},  
  "x2";"c"+<3pt,2pt>="c2"**\crv{"x2"+(0,-.25) & "c2"},
  "x3";"c"+<-3pt,-2pt>="c1"**\crv{"x3"+(0,.25) & "c1"}, 
  "c"+(0,.10)*!<0pt,-6pt>{\psi}, 
"x3"="L21",
"x4"="L22",
  "L21"="x1", 
  "x1"+(0,-.5)="x2",
  "x1";"x2"**\dir{}?(.5)="c", 
  "x1";"x2"**\dir{-},   
  "c"-(0,0)*!<-4pt,0pt>{\bullet S }, 
"x2"="L21",
  "L22"="x1", 
  "x1"+(0,-.5)="x2",
  "x1";"x2"**\dir{}?(.5)="c", 
  "x1";"x2"**\dir{-},   
  "c"-(0,0)*!<-4pt,0pt>{\bullet S }, 
"x2"="L22",
  "L21"="x1", 
  "L22"="x2",   
  "x1";"x2"**\dir{}?(.5)="m", 
  "m" + (0,-.25)="c", 
  "m" + (0,-.5)="x3", 
  "x1";"x3"**\crv{"x1"-(0,.25) & "c" & "x3"},
  "x2";"x3"**\crv{"x2"-(0,.25) & "c"  & "x3"}, 
"x3"="L21",
  "x3"-(0,.5)*!<0pt,-6pt>{H}, 
\endxy 
\end{equation}

\[S \eta=\eta, \quad \quad \varepsilon S = \varepsilon   ,\]
\begin{equation}\label{epn S} 
\xy /r2.5pc/:,   
(1,0)="L11",
  "L11"+(0,.5)*!<0pt,6pt>{H},
   "L11"="x1", 
  "x1"+(0,-.25)="x2",
  "x1";"x2"**\dir{}?(.5)="c", 
  "x1";"x2"**\dir{-},   
  "c"-(0,0)*!<-3pt,0pt>{\bullet \eta },  
  "x1"+(0,0)*!<0pt,-2pt>{\circ}, 
"x2"="L11",
   "L11"="x1", 
  "x1"+(0,-.5)="x2",
  "x1";"x2"**\dir{}?(.5)="c", 
  "x1";"x2"**\dir{-},   
  "c"-(0,0)*!<-4pt,0pt>{\bullet S },  
"x2"="L11",
  "L11"-(0,.5)*!<0pt,-6pt>{H}, 
 "x2" + (.75,.25)*!<0pt,-6pt>{=}, 
(2.5,0)="L11",
  "L11"+(0,.5)*!<0pt,6pt>{H},
   "L11"="x1", 
  "x1"+(0,-.5)="x2",
  "x1";"x2"**\dir{}?(.5)="c", 
  "x1";"x2"**\dir{-},   
  "c"-(0,0)*!<-3pt,0pt>{\bullet \eta },  
  "x1"+(0,0)*!<0pt,-2pt>{\circ}, 
"x2"="L11",
  "L11"-(0,.5)*!<0pt,-6pt>{H}, 
\endxy 
\quad \quad
\xy /r2.5pc/:,   
(1,0)="L11",
  "L11"+(0,.5)*!<0pt,6pt>{H},
   "L11"="x1", 
  "x1"+(0,-.5)="x2",
  "x1";"x2"**\dir{}?(.5)="c", 
  "x1";"x2"**\dir{-},   
  "c"-(0,0)*!<-4pt,0pt>{\bullet S }, 
"x2" = "L11",
   "L11"="x1", 
  "x1"+(0,-.25)="x2",
  "x1";"x2"**\dir{}?(.5)="c", 
  "x1";"x2"**\dir{-},   
  "c"-(0,0)*!<-3pt,0pt>{\bullet \varepsilon },  
  "x2"+(0,0)*!<0pt,2pt>{\circ},  
  "x2"+(0,0)*!<0pt,2pt>{\circ},  
"x2" = "L11",
  "L11"+(.5,.5)*!<0pt,0pt>{=}, 
(2,0)="L11",
  "L11"+(0,.5)*!<0pt,6pt>{H},
   "L11"="x1", 
  "x1"+(0,-.55)="x2",
  "x1";"x2"**\dir{}?(.5)="c", 
  "x1";"x2"**\dir{-},   
  "c"-(0,0)*!<-3pt,0pt>{\bullet \varepsilon },  
  "x2"+(0,0)*!<0pt,2pt>{\circ},  
   "x2"+(0,0)*!<0pt,2pt>{\circ},  
"x2" = "L11",
\endxy  
\end{equation}

\[\Delta S = \psi (S \otimes S)\Delta = (S \otimes S) \psi \Delta   .\]
\begin{equation}\label{Del S}
\xy /r2.5pc/:,   
(1,0)="L11",
  "L11"+(0,.5)*!<0pt,6pt>{H},
   "L11"="x1", 
  "x1"+(0,-.5)="x2",
  "x1";"x2"**\dir{}?(.5)="c", 
  "x1";"x2"**\dir{-},   
  "c"-(0,0)*!<-4pt,0pt>{\bullet S }, 
"x2" = "L21",
   "L21"="x1", 
  "x1"+(-.5,-.5)="x2",
  "x1"+(.5,-.5)="x3",
  "x1"+(0,-.25)="c",  
   "x1";"x2"**\crv{"x1" & "c" & "x2"+(0,.25)}, 
   "x1";"x3"**\crv{"x1" & "c" & "x3"+(0,.25)}, 
"x2"= "L21",
"x3"= "L22",
  "L21"-(0,.5)*!<0pt,-6pt>{H},  
  "L22"-(0,.5)*!<0pt,-6pt>{H}, 
  "L22"+(.5,.5)*!<0pt,0pt>{=}, 
(3,0)="L11",
  "L11"+(0,.5)*!<0pt,6pt>{H},
   "L11"="x1", 
  "x1"+(-.5,-.5)="x2",
  "x1"+(.5,-.5)="x3",
  "x1"+(0,-.25)="c",  
   "x1";"x2"**\crv{"x1" & "c" & "x2"+(0,.25)}, 
   "x1";"x3"**\crv{"x1" & "c" & "x3"+(0,.25)}, 
"x2"= "L21",
"x3"= "L22",
   "L21"="x1", 
  "x1"+(0,-.5)="x2",
  "x1";"x2"**\dir{}?(.5)="c", 
  "x1";"x2"**\dir{-},  
  "c"-(0,0)*!<-4pt,0pt>{\bullet S }, 
"x2" = "L31",
   "L22"="x1", 
  "x1"+(0,-.5)="x2",
  "x1";"x2"**\dir{}?(.5)="c", 
  "x1";"x2"**\dir{-},   
  "c"-(0,0)*!<-4pt,0pt>{\bullet S }, 
"x2" = "L32",
    "L31"="x1",      
   "L32"="x2",  
  "x1"+(0,-.5)="x3",
  "x2"+(0,-.5)="x4",
  "x1";"x2"**\dir{}?(.5)="m", 
  "m" + (0,-.25)="c", 
  "x1";"x4"**\crv{ "x1"+(0,-.25) & "c" & "x4"+(0,.25)},  
  "x2";"c"+<3pt,2pt>="c2"**\crv{"x2"+(0,-.25) & "c2"},
  "x3";"c"+<-3pt,-2pt>="c1"**\crv{"x3"+(0,.25) & "c1"}, 
  "c"+(0,.10)*!<0pt,-6pt>{\psi}, 
"x3" = "L31",
"x4" = "L32",
  "L31"-(0,.5)*!<0pt,-6pt>{H},  
  "L32"-(0,.5)*!<0pt,-6pt>{H}, 
  "L22"+(.5,0)*!<0pt,0pt>{=}, 
(5,0)="L11",
  "L11"+(0,.5)*!<0pt,6pt>{H},
   "L11"="x1", 
  "x1"+(-.5,-.5)="x2",
  "x1"+(.5,-.5)="x3",
  "x1"+(0,-.25)="c",  
   "x1";"x2"**\crv{"x1" & "c" & "x2"+(0,.25)}, 
   "x1";"x3"**\crv{"x1" & "c" & "x3"+(0,.25)}, 
"x2"= "L21",
"x3"= "L22",
    "L21"="x1",       
   "L22"="x2",  
  "x1"+(0,-.5)="x3",
  "x2"+(0,-.5)="x4",
  "x1";"x2"**\dir{}?(.5)="m", 
  "m" + (0,-.25)="c", 
  "x1";"x4"**\crv{ "x1"+(0,-.25) & "c" & "x4"+(0,.25)},  
  "x2";"c"+<3pt,2pt>="c2"**\crv{"x2"+(0,-.25) & "c2"},
  "x3";"c"+<-3pt,-2pt>="c1"**\crv{"x3"+(0,.25) & "c1"}, 
  "c"+(0,.10)*!<0pt,-6pt>{\psi}, 
"x3" = "L31",
"x4" = "L32",
   "L31"="x1", 
  "x1"+(0,-.5)="x2",
  "x1";"x2"**\dir{}?(.5)="c", 
  "x1";"x2"**\dir{-},   
  "c"-(0,0)*!<-4pt,0pt>{\bullet S }, 
"x2" = "L31",
   "L32"="x1", 
  "x1"+(0,-.5)="x2",
  "x1";"x2"**\dir{}?(.5)="c", 
  "x1";"x2"**\dir{-},   
  "c"-(0,0)*!<-4pt,0pt>{\bullet S }, 
"x2" = "L32",
  "L31"-(0,.5)*!<0pt,-6pt>{H},  
  "L32"-(0,.5)*!<0pt,-6pt>{H}, 
\endxy  
\end{equation}
\end{proposition} 
 
\begin{proof}
~\cite{maj}
\end{proof}

\begin{definition} Let $H$ be a braided Hopf algebra in $\mathcal{C}$.  A right $H$-module
 is an object $M$ in $\mathcal{C}$ equipped with a morphism $\phi_M :M \otimes H
\to M$, called $H$ action, such that:
\[ (\phi)(id_M \otimes m_H) = (\phi)(\phi \otimes id_H), \quad \quad (\phi)(id_M \otimes \eta_H) = id_M. \]
\begin{equation} \label{H-Mod 2} \notag
\xy /r2.5pc/:,   
(1,0)="L11",
(2,0)="L12",
(3,0)="L13",
  "L11"+(0,.5)*!<0pt,6pt>{M},
  "L12"+(0,.5)*!<0pt,6pt>{H}, 
  "L13"+(0,.5)*!<0pt,6pt>{H},
  "L11";"L11" +(0,-.5)="L11"**\dir{-},
  "L12"="x1", 
  "L13"="x2",   
  "x1";"x2"**\dir{}?(.5)="m", 
  "m" + (0,-.25)="c", 
  "m" + (0,-.5)="x3", 
  "x1";"x3"**\crv{"x1"-(0,.25) & "c" & "x3"},
  "x2";"x3"**\crv{"x2"-(0,.25) & "c"  & "x3"}, 
  "c"+(0,.10)*!<0pt,-0pt>{m_H}, 
"x3"="L12",
  "L11"="x1", 
  "L12"="x2",   
  "x1";"x2"**\dir{}?(.5)="m", 
  "m" + (0,-.25)="c", 
  "m" + (0,-.5)="x3", 
  "x1";"x3"**\crv{"x1"-(0,.25) & "c" & "x3"},
  "x2";"x3"**\crv{"x2"-(0,.25) & "c"  & "x3"}, 
  "c"+(0,.10)*!<0pt,-0pt>{\phi_M}, 
"x3"="L11",
  "L11"-(0,.5)*!<0pt,-6pt>{M}, 
  "L11"+(1.5,.5)*!<0pt,0pt>{=}, 
(4,0)="L11",
(5,0)="L12",
(6,0)="L13",
  "L11"+(0,.5)*!<0pt,6pt>{M},
  "L12"+(0,.5)*!<0pt,6pt>{H}, 
  "L13"+(0,.5)*!<0pt,6pt>{H},
  "L11"="x1", 
  "L12"="x2",   
  "x1";"x2"**\dir{}?(.5)="m", 
  "m" + (0,-.25)="c", 
  "m" + (0,-.5)="x3", 
  "x1";"x3"**\crv{"x1"-(0,.25) & "c" & "x3"},
  "x2";"x3"**\crv{"x2"-(0,.25) & "c"  & "x3"}, 
  "c"+(0,.10)*!<0pt,-0pt>{\phi_M}, 
"x3"="L11",
  "L13";"L13" +(0,-.5)="L12"**\dir{-},
  "L11"="x1", 
  "L12"="x2",  
  "x1";"x2"**\dir{}?(.5)="m", 
  "m" + (0,-.25)="c", 
  "m" + (0,-.5)="x3", 
  "x1";"x3"**\crv{"x1"-(0,.25) & "c" & "x3"},
  "x2";"x3"**\crv{"x2"-(0,.25) & "c"  & "x3"}, 
  "c"+(0,.10)*!<0pt,-0pt>{\phi_M}, 
"x3"="L11",
  "L11"-(0,.5)*!<0pt,-6pt>{M}, 
\endxy  
\quad \quad
\xy /r2.5pc/:,  
(1,0)="L11",
(2,0)="L12",
  "L11"+(0,.5)*!<0pt,6pt>{M},
  "L11";"L11" +(0,-.5)="L11"**\dir{-},
   "L12" +(0, -.25)="x1", 
  "x1"+(0,-.25)="x2",
  "x1";"x2"**\dir{}?(.5)="c", 
  "x1";"x2"**\dir{-},   
  "c"-(0,0)*!<-3pt,0pt>{\bullet \eta },  
  "x1"+(0,0)*!<0pt,-2pt>{\circ},  
"x2"="L12",
  "L11"="x1", 
  "L12"="x2",   
  "x1";"x2"**\dir{}?(.5)="m", 
  "m" + (0,-.25)="c", 
  "m" + (0,-.5)="x3", 
  "x1";"x3"**\crv{"x1"-(0,.25) & "c" & "x3"},
  "x2";"x3"**\crv{"x2"-(0,.25) & "c"  & "x3"}, 
  "c"+(0,.10)*!<0pt,-0pt>{\phi_M}, 
"x3"="L11",
  "L11"-(0,.5)*!<0pt,-6pt>{M}, 
  "L11"+(1.5,.5)*!<0pt,0pt>{=}, 
(4,0)="L11",
  "L11"+(0,.5)*!<0pt,6pt>{M},
  "L11";"L11" +(0,-.5)="L11"**\dir{-},
  "L11"-(0,.5)*!<0pt,-6pt>{M}, 
\endxy 
\end{equation}

A  left $H$-comodule is an   object $M$ in $\mathcal{C}$ equipped with a
morphism $\rho_M :M \to H \otimes M$, called $H$ coaction,  such that:
\[(\Delta_H \otimes id_M)(\rho) = (id_M \otimes \rho)(\rho) , \quad \quad (\epsilon_H \otimes id_M)(\rho) = id_M.\]
\begin{equation}\label{H-COMod 2}\notag
\xy /r2.5pc/:,  
(1,0)="L11",
  "L11"+(0,.5)*!<0pt,6pt>{M},
    "L11"="x1", 
  "x1"+(-.5,-.5)="x2",
  "x1"+(.5,-.5)="x3",
  "x1"+(0,-.25)="c",  
   "x1";"x2"**\crv{"x1" & "c" & "x2"+(0,.25)}, 
   "x1";"x3"**\crv{"x1" & "c" & "x3"+(0,.25)}, 
  "c"+(0,0)*!<0pt,2pt>{\rho_M},  
"x2"="L11",
"x3"="L12",
    "L11"="x1", 
  "x1"+(-.5,-.5)="x2",
  "x1"+(.5,-.5)="x3",
  "x1"+(0,-.25)="c",  
   "x1";"x2"**\crv{"x1" & "c" & "x2"+(0,.25)}, 
   "x1";"x3"**\crv{"x1" & "c" & "x3"+(0,.25)}, 
  "c"+(0,0)*!<0pt,2pt>{\Delta_H},  
"x2"="L21",
"x3"="L22",
  "L12";"L12" +(0,-.5)="L13"**\dir{-},
"L21"="L11",
"L22"="L12",
  "L11"-(0,.5)*!<0pt,-6pt>{H}, 
  "L12"-(0,.5)*!<0pt,-6pt>{H}, 
  "L13"-(0,.5)*!<0pt,-6pt>{M}, 
  "L11"+(2,.5)*!<0pt,0pt>{=}, 
(3.5,0)="L11",
  "L11"+(0,.5)*!<0pt,6pt>{M},
    "L11"="x1", 
  "x1"+(-.5,-.5)="x2",
  "x1"+(.5,-.5)="x3",
  "x1"+(0,-.25)="c",  
   "x1";"x2"**\crv{"x1" & "c" & "x2"+(0,.25)}, 
   "x1";"x3"**\crv{"x1" & "c" & "x3"+(0,.25)}, 
  "c"+(0,0)*!<0pt,2pt>{\rho_M},  
"x2"="L11",
"x3"="L12",
  "L11";"L11" +(0,-.5)="L11"**\dir{-},
    "L12"="x1", 
  "x1"+(-.5,-.5)="x2",
  "x1"+(.5,-.5)="x3",
  "x1"+(0,-.25)="c",  
   "x1";"x2"**\crv{"x1" & "c" & "x2"+(0,.25)}, 
   "x1";"x3"**\crv{"x1" & "c" & "x3"+(0,.25)}, 
  "c"+(0,0)*!<0pt,2pt>{\rho_M},  
"x2"="L12",
"x3"="L13",
  "L11"-(0,.5)*!<0pt,-6pt>{H}, 
  "L12"-(0,.5)*!<0pt,-6pt>{H}, 
  "L13"-(0,.5)*!<0pt,-6pt>{M}, 
\endxy  
\quad \quad \quad
\xy /r2.5pc/:,   
(1,0)="L11",
  "L11"+(0,.5)*!<0pt,6pt>{M},
    "L11"="x1", 
  "x1"+(-.5,-.5)="x2",
  "x1"+(.5,-.5)="x3",
  "x1"+(0,-.25)="c",  
   "x1";"x2"**\crv{"x1" & "c" & "x2"+(0,.25)}, 
   "x1";"x3"**\crv{"x1" & "c" & "x3"+(0,.25)}, 
  "c"+(0,0)*!<0pt,2pt>{\rho_M},  
"x2"="L11",
"x3"="L12",
   "L11"="x1", 
  "x1"+(0,-.25)="x2",
  "x1";"x2"**\dir{}?(.5)="c", 
  "x1";"x2"**\dir{-},   
  "c"-(0,0)*!<-3pt,0pt>{\bullet \varepsilon },  
  "x2"+(0,0)*!<0pt,2pt>{\circ},  
  "x2"+(0,0)*!<0pt,2pt>{\circ},  
  "L12";"L12" +(0,-.5)="L11"**\dir{-},
  "L11"-(0,.5)*!<0pt,-6pt>{M}, 
  "L11"+(1,.5)*!<0pt,0pt>{=}, 
(3.5,0)="L11",
  "L11"+(0,.5)*!<0pt,6pt>{M},
    "L11"="x1", 
  "x1"+(-.5,-.5)="x2",
  "x1"+(.5,-.5)="x3",
  "x1"+(0,-.25)="c", 
   "x1";"x2"**\crv{"x1" & "c" & "x2"+(0,.25)}, 
   "x1";"x3"**\crv{"x1" & "c" & "x3"+(0,.25)}, 
  "c"+(0,0)*!<0pt,2pt>{\rho_M},  
"x2"="L11",
"x3"="L12",
  "L11"="x1",       
  "L12"="x2",  , 
  "x1"+(0,-.5)="x3",
  "x2"+(0,-.5)="x4",
  "x1";"x2"**\dir{}?(.5)="m", 
  "m" + (0,-.25)="c", 
  "x1";"x4"**\crv{ "x1"+(0,-.25) & "c" & "x4"+(0,.25)},  
  "x2";"c"+<3pt,2pt>="c2"**\crv{"x2"+(0,-.25) & "c2"},
  "x3";"c"+<-3pt,-2pt>="c1"**\crv{"x3"+(0,.25) & "c1"}, 
  "c"+(0,.10)*!<0pt,-3pt>{\psi_{H,M}},  
"x3"="L11",
"x4"="L12",
  "L11";"L11" +(0,-.5)="L11"**\dir{-},
   "L12"="x1", 
  "x1"+(0,-.25)="x2",
  "x1";"x2"**\dir{}?(.5)="c", 
  "x1";"x2"**\dir{-},   
  "c"-(0,0)*!<-3pt,0pt>{\bullet \varepsilon },  
  "x2"+(0,0)*!<0pt,2pt>{\circ},  
   "x2"+(0,0)*!<0pt,2pt>{\circ},  
  "L11"-(0,.5)*!<0pt,-6pt>{M}, 
  "L11"+(1.5,1)*!<0pt,0pt>{=}, 
(5.5,0)="L11",
  "L11"+(0,.5)*!<0pt,6pt>{M},
  "L11";"L11" +(0,-1)="L11"**\dir{-},
  "L11"-(0,.5)*!<0pt,-6pt>{M}, 
\endxy  
\end{equation}
\end{definition}

\begin{remark}
Let  $\mathcal{C}$  be  a  strict, braided, monoidal abelian category.
Throughout the paper we shall use the following  conventions to denote objects and morphisms of $\mathcal{C}$:
\item $\bullet$~$A^n$ for $A^ {\otimes n}$,
\item $\bullet$~$1$ for $id,$  e.g. we write $1_A$ or just $1$ for $id_A$, 
\item $\bullet$~$(f ,g)$ for $(f \otimes g)$,
\item $\bullet$~$id_n$ or just $1_n$ for $id_{A^n}$,
\item $\bullet$~$1_{A,B}$ for $1_A \otimes 1_B$,
\item $\bullet$~$\psi$ for $\psi_{A,A}$.

\noindent For example instead of writing 
$\Delta_H \, m_H = (m_H \otimes m_H)(id_H\otimes \psi_{H,H} \otimes id_H)(\Delta_H \otimes \Delta_H)$, which expresses
the fact that the comultiplication of a Hopf algebra  is an
anti-algebra map, we just write 
$\Delta m = (m,  m)(1 , \psi ,1)(\Delta , \Delta)$;
 or instead of $ \psi_{U,U \otimes U} =(id_U \otimes \psi_{U,U})(\psi_{U,U} \otimes id_U),$ 
 we simply write $\psi_{U,U^2} = (1,\psi)(\psi ,1)$ when there is no chance of confusion, and so on.
\end{remark}

%%***************************************************************************************************
\section{The cocyclic object of a braided triple ($H$,\,$C$,\,$M$)}\label{hopf cyclc for H C M}
In this section we extend the notion of a stable
anti-Yetter-Drinfeld (SAYD) module ~\cite{hkrs1,hkrs2} to braided
monoidal categories and define a cocyclic object for braided triples
($H$,\,$C$,\,$M$) in a \emph{symmetric} monoidal  abelian  category.
In the last section of this paper we treat the general non-symmetric
case which is much more subtle. Recall that, by definition,  in a
\emph{para-cocyclic object} all axioms of a cocyclic object are
satisfied except the relations $ \tau_n^{n+1}=id$. Given a
para-cocyclic object $X_n$, $n\geq 0$ in an abelian  category, we can always define a
cocyclic object by considering
\begin{equation}\label{pra to cyc}\notag
\overline{X}_n:= \text{ker}\,(id-\tau_n^{n+1}), 
\end{equation}
 and restricting the faces, degeneracies, and cyclic operators
to these subspaces. For general notion of cyclic and cocyclic objects 
 we refer to ~\cite{c1,cbook, l}

We fix a strict,  braided  monoidal  category $\mathcal{C}$, and a Hopf algebra $H$ in $\mathcal{C}$. For  the following definition $\mathcal{C}$ need not
be symmetric or additive.

\begin{definition} \label{r l SAYD mod}
A right-left braided stable anti-Yetter-Drinfeld (SAYD) $H$-module in $\mathcal{C}$ is an object $M$ in $\mathcal{C}$ such that:\\
(i) $M$ is a right $H$-module via an action $\phi_M :M \ot H \to M$, \\
(ii) $M$ is a left $H$-comodule via a coaction $\rho_M :M \to H \ot M$, \\
(iii) $M$ satisfies the braided anti-Yetter-Drinfeld condition, i.e.:
\begin{multline}\label{AYD}
(\rho)(\phi)= [(m)(S \ot m) \ot \phi][(\psi_{H^{\ot 2} ,H} \ot id_M \ot id_H)(id_{H^{\ot 2}} \ot \psi_{M ,H} \ot id_H)\\   (id_{H^{\ot 2}} \ot id_M \ot \psi_{H,H})(id_H \ot \psi_{M , H} \ot id_{H^{\ot 2}})][\rho \ot \Delta^2].
\end{multline}
(iv) $M$ is stable,  i.e. :
\[(\phi)(\psi_{H , M})(\rho) = id_M  .\]%\label{stable}
\end{definition}

\begin{remark} To deal with large expressions like ~\eqref{AYD} we break them
into two lines. 
\end{remark}

The above braided SAYD conditions $(iii)$ and $(iv)$ in terms of braided diagrams are as follows:
\begin{equation}\label{BSAYD} 
\xy /r2.5pc/:,   
  (0,0)="x1", 
  "x1"+(2,0)="x2",
  "x1";"x2"**\dir{}?(.5)="m", 
  "m"+(0,-.5)="c", 
  "m"+(0,-1)="x3", 
  "x1";"x3"**\crv{"x1"-(0,.5) & "c" & "x3"},
  "x2";"x3"**\crv{"x2"-(0,.5) & "c"  & "x3"}, 
  "c"+(0,.10)*!<0pt,-6pt>{\phi_M},   
  "x1"+(0,.5)*!<0pt,6pt>{M},  
  "x2"+(0,.5)*!<0pt,6pt>{H},   
  "x3"="x1", 
  "x1"+(-1,-1)="x2",
  "x1"+(1,-1)="x3",
  "x1"+(0,-.5)="c",  
   "x1";"x2"**\crv{"x1" & "c" & "x2"+(0,.5)}, 
   "x1";"x3"**\crv{"x1" & "c" & "x3"+(0,.5)}, 
  "x2"-(0,.5)*!<0pt,-6pt>{H}, 
  "x3"-(0,.5)*!<0pt,-6pt>{M},  
  "c"-(0,.10)*!<0pt,6pt>{\rho _M},   
  "x1"+(1,0)*!<0pt,0pt>{=}, 
(3,0)= "L11",
(6,0)= "L12",
  "L11"+(0,.5)*!<0pt,6pt>{M}, 
  "L12"+(0,.5)*!<0pt,6pt>{H}, 
  "L11"="x1", 
  "x1"+(-.5,-.5)="x2",
  "x1"+(.5,-.5)="x3",
  "x1"+(0,-.25)="c",  
   "x1";"x2"**\crv{"x1" & "c" & "x2"+(0,.25)}, 
   "x1";"x3"**\crv{"x1" & "c" & "x3"+(0,.25)}, 
  "c"+(0,0)*!<0pt,3pt>{\rho _M}, 
"x2" ="L21",
"x3" ="L22",
  "L12"="x1", 
  "x1"+(-1,-.5)="x2",
  "x1"+(0,-.5)="x3",
  "x1"+(1,-.5)="x4",
  "x1"+(0,-.25)="c",  
   "x1";"x2"**\crv{"x1" & "c" & "x2"+(0,.25)}, 
   "x1";"x3"**\crv{"x1" & "c" & "x3"+(0,.25)}, 
   "x1";"x4"**\crv{"x1" & "c" & "x4"+(0,.25)}, 
  "c"+(0,0)*!<8pt,5pt>{\Delta_H^2}, 
"L21" ="L11",
"L22" ="L12",
"x2"="L13",
"x3"="L14",
"x4"="L15",
  "L11";"L11" +(0,-.5)="L11"**\dir{-},
  "L12"="x1",       
  "L13"="x2",  , 
  "x1"+(0,-.5)="x3",
  "x2"+(0,-.5)="x4",
  "x1";"x2"**\dir{}?(.5)="m", 
  "m" + (0,-.25)="c", 
  "x1";"x4"**\crv{ "x1"+(0,-.25) & "c" & "x4"+(0,.25)},  
  "x2";"c"+<3pt,2pt>="c2"**\crv{"x2"+(0,-.25) & "c2"},
  "x3";"c"+<-3pt,-2pt>="c1"**\crv{"x3"+(0,.25) & "c1"}, 
  "c"+(0,.10)*!<0pt,-6pt>{\si_{M,H}},
"x3"="L12",
"x4"="L13",
  "L14";"L14" +(0,-.5)="L14"**\dir{-},
  "L15";"L15" +(0,-.5)="L15"**\dir{-},
  "L11";"L11" +(0,-.5)="L11"**\dir{-},
  "L12";"L12" +(0,-.5)="L12"**\dir{-},
  "L13";"L13" +(0,-.5)="L13"**\dir{-},
  "L14"="x1",       
  "L15"="x2",  , 
  "x1"+(0,-.5)="x3",
  "x2"+(0,-.5)="x4",
  "x1";"x2"**\dir{}?(.5)="m", 
  "m" + (0,-.25)="c", 
  "x1";"x4"**\crv{ "x1"+(0,-.25) & "c" & "x4"+(0,.25)},  
  "x2";"c"+<3pt,2pt>="c2"**\crv{"x2"+(0,-.25) & "c2"},
  "x3";"c"+<-3pt,-2pt>="c1"**\crv{"x3"+(0,.25) & "c1"}, 
  "c"+(0,.10)*!<0pt,-6pt>{\si_{H,H}},
"x3"="L14",
"x4"="L15",
"L11";"L11" +(0,-.5)="L11"**\dir{-},
"L12";"L12" +(0,-.5)="L12"**\dir{-},
  "L13"="x1",       
  "L14"="x2",  , 
  "x1"+(0,-.5)="x3",
  "x2"+(0,-.5)="x4",
  "x1";"x2"**\dir{}?(.5)="m", 
  "m" + (0,-.25)="c", 
  "x1";"x4"**\crv{ "x1"+(0,-.25) & "c" & "x4"+(0,.25)},  
  "x2";"c"+<3pt,2pt>="c2"**\crv{"x2"+(0,-.25) & "c2"},
  "x3";"c"+<-3pt,-2pt>="c1"**\crv{"x3"+(0,.25) & "c1"}, 
  "c"+(0,.10)*!<0pt,-6pt>{\si_{M,H}},
"x3"="L13",
"x4"="L14",
  "L15";"L15" +(0,-.5)="L15"**\dir{-},
"L11"="a1",
"L12"="a2",
"L13"="a3",
"a1" + (0,-1)="b1",
"a2" + (0,-1)="b2",
"a3" + (0,-1)="b3",
"a1";"a2"**\dir{}?(.5)="m12", 
  "m12" + (0,-.25)="c12", 
"a2";"a3"**\dir{}?(.5)="m23", 
  "m23" + (0,-.25)="c23", 
  "a1";"b2"**\crv{ "a1"+(0,-.25) & "c12" & "b2"+(0,.25)},  
  "a2";"b3"**\crv{ "a2"+(0,-.25) & "c23" & "b3"+(0,.25)}, 
                "b1";"c12"+<-0pt,-11pt>="c21"**\crv{"b1"+(0,.25) & "c21"},         
               "a3";"c23"+<-0pt,-2pt>="c32"**\crv{"a3"+(0,-.25) & "c32"},
"c12" +(.35,-.20);"c23" +(-.25,-.10)**\dir{-}, 
"c23"+(0,.10)*!<0pt,-6pt>{\si_{H^2,H}}, 
"b1"="L11",
"b2"="L12",
"b3"="L13",
  "L14";"L14" +(0,-1)="L14"**\dir{-},
  "L15";"L15" +(0,-1)="L15"**\dir{-},
  "L11"="x1", 
  "x1"+(0,-.5)="x2",
  "x1";"x2"**\dir{}?(.5)="c", 
  "x1";"x2"**\dir{-},   
  "c"-(0,0)*!<-4pt,0pt>{\bullet S_H },
"x2"="L11",
   "L12"="x1", 
  "L13"="x2",   
  "x1";"x2"**\dir{}?(.5)="m", 
  "m" + (0,-.25)="c", 
  "m" + (0,-.5)="x3", 
  "x1";"x3"**\crv{"x1"-(0,.25) & "c" & "x3"},
  "x2";"x3"**\crv{"x2"-(0,.25) & "c"  & "x3"}, 
  "c"+(0,.10)*!<0pt,-0pt>{m_H}, 
"x3"="L12",
   "L11"="x1", 
  "L12"="x2",   
  "x1";"x2"**\dir{}?(.5)="m", 
  "m" + (0,-.25)="c", 
  "m" + (0,-.5)="x3", 
  "x1";"x3"**\crv{"x1"-(0,.25) & "c" & "x3"},
  "x2";"x3"**\crv{"x2"-(0,.25) & "c"  & "x3"}, 
  "c"+(0,.10)*!<0pt,-0pt>{m_H}, 
"x3"="L11",
   "L14"="x1", 
  "L15"="x2",   
  "x1";"x2"**\dir{}?(.5)="m", 
  "m" + (0,-.25)="c", 
  "m" + (0,-.5)="x3", 
  "x1";"x3"**\crv{"x1"-(0,.25) & "c" & "x3"},
  "x2";"x3"**\crv{"x2"-(0,.25) & "c"  & "x3"}, 
  "c"+(0,.10)*!<0pt,-0pt>{\phi_M}, 
"x3"="L12",
  "L12";"L12" +(0,-.5)="L12"**\dir{-},
  "L11"-(0,.5)*!<0pt,-6pt>{H},
  "L12"-(0,.5)*!<0pt,-6pt>{M},
  \endxy
\quad \quad \quad \quad
\xy /r2.5pc/:,   
 (0,0)="x1", 
  "x1"+(-.5,-1)="x2",
  "x1"+(.5,-1)="x3",
  "x1"+(0,-.5)="c",  
   "x1";"x2"**\crv{"x1" & "c" & "x2"+(0,.5)}, 
   "x1";"x3"**\crv{"x1" & "c" & "x3"+(0,.5)}, 
  "x1"+(0,.5)*!<0pt,6pt>{M},  
   "c"-(0,.10)*!<0pt,0pt>{\rho _M},  
 "x2"="x1", 
  "x1"+(1,0)="x2",
  "x1"+(0,-1)="x3",
  "x1"+(1,-1)="x4",
  "x1"+(.5,-.5)="c",  
  "x1";"x4"**\crv{ "x1"+(0,-.5) & "c" & "x4"+(0,.5)},  
  "x2";"c"+<3pt,2pt>="c2"**\crv{"x2"+(0,-.5) & "c2"},
  "x3";"c"+<-3pt,-2pt>="c1"**\crv{"x3"+(0,.5) & "c1"}, 
  "c"+(0,.10)*!<0pt,-6pt>{\si_{H,M}},  
  "x3"="x1", 
  "x1"+(1,0)="x2",
  "x1"+(.5,-1)="x3",
  "x1"+(.5,-.5)="c", 
  "x1";"x3"**\crv{"x1"-(0,.5) & "c" & "x3"},
  "x2";"x3"**\crv{"x2"-(0,.5) & "c"  & "x3"}, 
  "c"+(0,.10)*!<0pt,-0pt>{\phi_M},   
  "x3"-(0,.5)*!<0pt,-6pt>{M},   
  "x2"+(.5,.5)*!<0pt,0pt>{=}, 
 (1.5,0)="x1";(1.5,-3)="x2" **\dir{-},  
  "x1"+(0,.5)*!<0pt,6pt>{M}, 
  "x2"-(0,.5)*!<0pt,-6pt>{M},  
  \endxy
\end{equation} 

\begin{definition} \label{H mod coalg}
A quadruple $(C,\, \Delta_C,\, \epsilon_C ,\, \phi_C)$ is called a left (braided) $H$-module-coalgebra in $\mathcal{C}$ if $(C,\Delta_C,\epsilon_C)$ is a coalgebra in $\mathcal{C}$, and  $C$ is a left $H$-module via an action $\phi_C :H \ot C \to C$ such that $\phi_C$ is a coalgebra map in $\mathcal{C}$ i.e. we have:
\[ \Delta_C \phi_C =(\phi_C \ot \phi_C)(id_H \ot \psi_{H,C} \ot id_C)(\Delta_H \ot \Delta_C), \quad \quad  \varepsilon_C \phi_C = \varepsilon_H \ot \varepsilon_C  \] 
\begin{equation}\label{H-module coalgebra 2}
\xy /r2.5pc/:, 
(1,0)="L11",
(2,0)="L12",
"L11"+(0,.5)*!<0pt,6pt>{H},
"L12"+(0,.5)*!<0pt,6pt>{C},
  "L11"="x1", 
  "L12"="x2",   
  "x1";"x2"**\dir{}?(.5)="m", 
  "m" + (0,-.25)="c", 
  "m" + (0,-.5)="x3", 
  "x1";"x3"**\crv{"x1"-(0,.25) & "c" & "x3"},
  "x2";"x3"**\crv{"x2"-(0,.25) & "c"  & "x3"}, 
  "c"+(0,.10)*!<0pt,-0pt>{\phi_C},  
"x3"="L11",
   "L11"="x1", 
  "x1"+(-.5,-.5)="x2",
  "x1"+(.5,-.5)="x3",
  "x1"+(0,-.25)="c",  
   "x1";"x2"**\crv{"x1" & "c" & "x2"+(0,.25)}, 
   "x1";"x3"**\crv{"x1" & "c" & "x3"+(0,.25)}, 
  "c"+(0,0)*!<0pt,4pt>{\Delta_C},  
"x2"="L11",
"x3"="L12",
"L11"+(0,-.5)*!<0pt,-6pt>{C}, 
"L12"+(0,-.5)*!<0pt,-6pt>{C}, 
  "x1"+(1,0)*!<0pt,0pt>{=}, 
(3.5,0)="L11",
(5.5,0)="L12",
"L11"+(0,.5)*!<0pt,6pt>{H},
"L12"+(0,.5)*!<0pt,6pt>{C},
   "L11"="x1", 
  "x1"+(-.5,-.5)="x2",
  "x1"+(.5,-.5)="x3",
  "x1"+(0,-.25)="c",  
   "x1";"x2"**\crv{"x1" & "c" & "x2"+(0,.25)}, 
   "x1";"x3"**\crv{"x1" & "c" & "x3"+(0,.25)}, 
  "c"+(0,0)*!<0pt,4pt>{\Delta_H},  
"x2"="L21",
"x3"="L22",
   "L12"="x1", 
  "x1"+(-.5,-.5)="x2",
  "x1"+(.5,-.5)="x3",
  "x1"+(0,-.25)="c",  
   "x1";"x2"**\crv{"x1" & "c" & "x2"+(0,.25)}, 
   "x1";"x3"**\crv{"x1" & "c" & "x3"+(0,.25)}, 
  "c"+(0,0)*!<0pt,4pt>{\Delta_C},  
"x2"="L23",
"x3"="L24",
"L21"="L11",
"L22"="L12",
"L23"="L13",
"L24"="L14",
  "L11";"L11" +(0,-.5)="L11"**\dir{-},
  "L12"="x1",       
  "L13"="x2",  , 
  "x1"+(0,-.5)="x3",
  "x2"+(0,-.5)="x4",
  "x1";"x2"**\dir{}?(.5)="m", 
  "m" + (0,-.25)="c", 
  "x1";"x4"**\crv{ "x1"+(0,-.25) & "c" & "x4"+(0,.25)},  
  "x2";"c"+<3pt,2pt>="c2"**\crv{"x2"+(0,-.25) & "c2"},
  "x3";"c"+<-3pt,-2pt>="c1"**\crv{"x3"+(0,.25) & "c1"}, 
  "c"+(0,.10)*!<0pt,-6pt>{\si_{H,C}},
"x3"="L12",
"x4"="L13",
  "L14";"L14" +(0,-.5)="L14"**\dir{-},
  "L11"="x1", 
  "L12"="x2",   
  "x1";"x2"**\dir{}?(.5)="m", 
  "m" + (0,-.25)="c", 
  "m" + (0,-.5)="x3", 
  "x1";"x3"**\crv{"x1"-(0,.25) & "c" & "x3"},
  "x2";"x3"**\crv{"x2"-(0,.25) & "c"  & "x3"}, 
  "c"+(0,.10)*!<0pt,-0pt>{\phi_C},  
"x3"="L11",
  "L13"="x1", 
  "L14"="x2",   
  "x1";"x2"**\dir{}?(.5)="m", 
  "m" + (0,-.25)="c", 
  "m" + (0,-.5)="x3", 
  "x1";"x3"**\crv{"x1"-(0,.25) & "c" & "x3"},
  "x2";"x3"**\crv{"x2"-(0,.25) & "c"  & "x3"}, 
  "c"+(0,.10)*!<0pt,-0pt>{\phi_C},  
"x3"="L12",
"L11"+(0,-.5)*!<0pt,-6pt>{C}, 
"L12"+(0,-.5)*!<0pt,-6pt>{C}, 
\endxy %\label{H-module coalgebra 1} 
\quad \quad
\xy /r2.5pc/:, 
(1,0)="L11",
(2,0)="L12",
"L11"+(0,.5)*!<0pt,6pt>{H},
"L12"+(0,.5)*!<0pt,6pt>{C},
  "L11"="x1", 
  "L12"="x2",   
  "x1";"x2"**\dir{}?(.5)="m", 
  "m" + (0,-.25)="c", 
  "m" + (0,-.5)="x3", 
  "x1";"x3"**\crv{"x1"-(0,.25) & "c" & "x3"},
  "x2";"x3"**\crv{"x2"-(0,.25) & "c"  & "x3"}, 
  "c"+(0,.10)*!<0pt,-0pt>{\phi_C},  
"x3"="L11",
   "L11"="x1", 
  "x1"+(0,-.25)="x2",
  "x1";"x2"**\dir{}?(.5)="c", 
  "x1";"x2"**\dir{-},   
  "c"-(0,0)*!<-6pt,0pt>{\bullet \varepsilon_C },  
  "x2"+(0,0)*!<0pt,2pt>{\circ},  
  "x1"+(1,0)*!<0pt,0pt>{=}, 
(3,0)="L11",
(4,0)="L12",
"L11"+(0,.5)*!<0pt,6pt>{H},
"L12"+(0,.5)*!<0pt,6pt>{C},
   "L11"="x1", 
  "x1"+(0,-.75)="x2",
  "x1";"x2"**\dir{}?(.5)="c", 
  "x1";"x2"**\dir{-},   
  "c"-(0,0)*!<-6pt,0pt>{\bullet \varepsilon_H },  
  "x2"+(0,0)*!<0pt,2pt>{\circ},  
   "L12"="x1", 
  "x1"+(0,-.75)="x2",
  "x1";"x2"**\dir{}?(.5)="c", 
  "x1";"x2"**\dir{-},   
  "c"-(0,0)*!<-6pt,0pt>{\bullet \varepsilon_C },  
  "x2"+(0,0)*!<0pt,2pt>{\circ},  
\endxy
\end{equation}
\end{definition}

\begin{definition} \label{diag act on Cn+1} Let $(C,\, \phi_C)$ be a left
$H$-module. The diagonal action of $H$ on $C^{n+1}:  =C^{\ot(n+1)}$
is defined by:
\[\phi_{C^{n+1}}:H \ot C^{n+1}  \to C^{n+1}  \]
\[\phi_{C^{n+1}} := \underbrace{ (\phi_C, \phi_C,...,\phi_C) }_{n+1~ times} (\mathcal{F}(\si_{H,C}))(\Delta_H^n \ot 1_{C^{n+1}}) , \]where,
\[\mathcal{F}(\psi_{H,C}): = \prod_{i=1}^n (id_{H^i}, \underbrace {\si_{H,C},\si_{H,C},...,\si_{H,C}}_{n+1-i~times} ,id_{C^i}) \label{F(si_H,C} .\]
\end{definition}
The following is the diagrammatic version of the diagonal action:
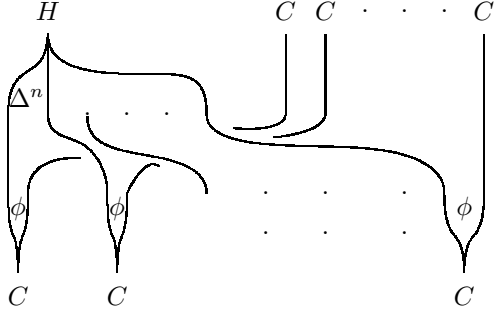
\begin{equation}\label{diagonalaction}\notag
\xy /r2.5pc/:,   
  (0,0)="x1", 
  "x1"+(-.5,-1)="x2",
  "x1"+(0,-1)="x3",
  "x1"+(2,-1)="x4",
  "x1"+(.5,-1)="a",
  "x1"+(1,-1)="b",
  "x1"+(1.5,-1)="d",
  "x1"+(0,-.5)="c",  
   "x1";"x2"**\crv{"x1" & "c" & "x2"+(0,.5)}, 
   "x1";"x3"**\crv{"x1" & "c" & "x3"+(0,.5)}, 
   "x1";"x4"**\crv{"x1"+(0,-.5) & "c"+(1.5,0) & "x4"+(0,.5)}, 
  "c"-(0,.10)*!<8pt,6pt>{\Delta^n},  
  "x1"+(0,.5)*!<0pt,6pt>{H},  
  "a"-(0,0)*!<0pt,0pt>{.},  
  "b"-(0,0)*!<0pt,0pt>{.},  
  "d"-(0,0)*!<0pt,0pt>{.},  
 "x1"="L11",
  (3,0)="L12"="y1",
  (3.5,0)="L13"="y2",
  (4,0)="L14"="y3",
  (4.5,0)="L15"="y4",
  (5,0)="L16"="y5",
  (5.5,0)="L17"="y6",
 "y1" ;"y1"+(0,-1) **\dir{-}, 
 "y2" ;"y2"+(0,-1) **\dir{-}, 
 "y6" ;"y6"+(0,-1) **\dir{-}, 
  "y1"+(0,.5)*!<0pt,6pt>{C},  
  "y2"+(0,.5)*!<0pt,6pt>{C},  
  "y3"+(0,.5)*!<0pt,6pt>{.},  
  "y4"+(0,.5)*!<0pt,6pt>{.},  
  "y5"+(0,.5)*!<0pt,6pt>{.},  
  "y6"+(0,.5)*!<0pt,6pt>{C}, 
 "x2" ="L21",
 "x3" ="L22",
 "b"="L23",
 "x4" ="L24",
 "y1"+(0,-1) ="L25",
 "y2"+(0,-1) ="L26",
 "y4"+(0,-1) ="L27",
 "y6"+(0,-1) ="L28",
  "x3";"a"+(.25,-1)**\crv{ "x3"+(0,-.3) & "a"+(0,-.25) & "a"+(.25,-1)+(0,.25)},   
  "x2"+(.25,-1);"a"+(0,-.5)+<-3pt,-2pt>="c1"**\crv{"x2"+(.255,-.5) & "c1"}, 
 "x2"+(.5,-1)+(1,0);"d"+(0,-.5)+<-3pt,-5pt>="c1"**\crv{"x2"+(.5,-.5)+(1.25,0) & "c1"}, 
  "x4";"y5"+(0,-2)**\crv{ "x4"+(0,-.5) & "y4"+(0,-1.25) & "y5"+(0,-2)+(0,.25)},  
  "a"+(0,-.05);"x4"+(0,-1)**\crv{ "a"+(0,-.45) & "d"+(0,-.5) & "x4"+(0,-1)+(0,.25)},  
  "y1"+(0,-1); "y1"+(-.75,-1.25)+<3pt,2pt>="c2"**\crv{ "y1"+(0,-1)+(0,-.25) & "c2"},
 "y2"+(0,-1); "y2"+(-.75,-1.5)+<3pt,6pt>="c2"**\crv{ "y2"+(0,-1)+(0,-.25) & "c2"},
 "y6" ;"y6"+(0,-2) **\dir{-}, 
"x2" ;"x2"+(0,-1) **\dir{-}, 
 "x4"+(.75,-1)*!<0pt,0pt>{.},  
 "x4"+(1.5,-1)*!<0pt,0pt>{.},  
 "x4"+(2.5,-1)*!<0pt,0pt>{.},  
 "x2" +(0,-1)="L31",
  "x2"+(.25,-1)="L32",
 "a"+(.25,-1)="L33",
  "x2"+(.5,-1)+(1,0)="L34",
   "x4"+(0,-1)="L35",
    "y5"+(0,-2)="L36",
    "y6"+(0,-2)="L37",
 "L31" ="x1",
   "L32"="x2",
  "x1";"x2"**\dir{}?(.5)="m", 
  "m" + (0,-.5)="c", 
  "m" + (0,-1)="x3", 
  "x1";"x3"**\crv{"x1"-(0,.5) & "c" & "x3"},
  "x2";"x3"**\crv{"x2"-(0,.5) & "c"  & "x3"}, 
  "c"+(0,.10)*!<0pt,-6pt>{\phi},   
  "x3"-(0,.5)*!<0pt,-6pt>{C},  
 "L31" ;"x1"**\dir{-}, 
 "L32" ;"x2"**\dir{-}, 
 "L33" ="x1",
"L34" ="x2",
  "x1";"x2"**\dir{}?(.5)="m", 
  "m" + (0,-.5)="c", 
  "m" + (0,-1)="x3", 
  "x1";"x3"**\crv{"x1"-(0,.5) & "c" & "x3"},
  "x2";"x3"**\crv{"x2"-(0,.5) & "c"  & "x3"}, 
  "c"+(0,.10)*!<0pt,-6pt>{\phi},   
  "x3"-(0,.5)*!<0pt,-6pt>{C},  
 "L36" ="x1",
"L37" ="x2",
  "x1";"x2"**\dir{}?(.5)="m", 
  "m" + (0,-.5)="c", 
  "m" + (0,-1)="x3", 
  "x1";"x3"**\crv{"x1"-(0,.5) & "c" & "x3"},
  "x2";"x3"**\crv{"x2"-(0,.5) & "c"  & "x3"}, 
  "c"+(0,.10)*!<0pt,-6pt>{\phi},   
  "x3"-(0,.5)*!<0pt,-6pt>{C},  
 "L36" ;"x1"**\dir{-}, 
 "L37" ;"x2"**\dir{-}, 
 "x4"+(.75,-1.5)*!<0pt,0pt>{.},  
 "x4"+(1.5,-1.5)*!<0pt,0pt>{.},  
 "x4"+(2.5,-1.5)*!<0pt,0pt>{.},  
 \endxy
\end{equation}

Now we are going to associate  a para-cocyclic object  to  any
triple $(H,C,M)$, where $H$ is a Hopf algebra, $C$ is a $H$-module
coalgebra and $M$ is a SAYD $H$-module, all in a symmetric monoidal
category $\mathcal{C}$. Notice that $\mathcal{C}$ need not be additive. For $n\geq
0$, let
\[ C^n=C^n(C,M) := M \ot C^{n+1}. \]
We
define faces $\delta_i:C^{n-1} \to C^n$, degeneracies $\sigma_i
:C^{n+1} \to C^n$ and cyclic maps $\tau_n :C^n \to C^n$ by:

\begin{eqnarray*}
\delta_i &=& \begin{cases}
 (1_M,1_{C^i}, \Delta_C, 1_{C^{n-i-1}}) \quad \quad \quad  \quad \quad \textrm{$0 \leq i < n$}\\
 (1_M , \si_{C , C^n})(1_M , \phi_C , 1_{C^n})(\si_{H,M} , 1_{C^{n+1}})(\rho _M , \Delta_C , 1_{C^{n-1}}) & \textrm{$i=n$}
  \end{cases}\\
\sigma_i &=& (1_M , 1_{C^{i+1}} , \varepsilon_C , 1_{C^{n-i}}), \quad \quad \quad \quad \quad  \textrm{$0 \leq i \leq n $}\\
 \tau_n &=& (1_M , \si_{C ,C^n})(1_M , \phi_C , 1_{C^n})(\si_{H,M} , 1_{C^{n+1}})(\rho _M ,1_{C^{n+1}})
\end{eqnarray*}

\begin{proposition} \label{main thm1} If $\mathcal{C}$ is a symmetric monoidal  category,
then $( C^\bullet ,\delta_i,\sigma_i, \tau )$ is a para-cocyclic
object in $\mathcal{C}$. 
\end{proposition}   

The idea of the proof of this proposition is very similar to and even easier than the proof of Theorem 
~\eqref{brdd ver of CM thry in nonsymmetric}, except that the symmetry condition is used in some steps.

Now let us assume in addition that $\mathcal{C}$ is an \emph{abelian}
category. Recall that given a right $H$- module $V$ via action $\phi_V$, 
and a left $H$-module $W$ via action $\phi_W$, the balanced tensor product 
$V\ot_H W$ is defined as the cokernel of the map:
\[\begin{CD}
 V \ot H \ot W @>(\phi_V \ot 1_W - 1_V \ot \phi_W) >> V \ot W 
\end{CD} \]
We form the balanced tensor products
\[ C^n_H=C^n_H(C,M) := M \ot_H C^{n+1}, \quad \quad n \geq 0,  \]
 with induced faces, degeneracies and cyclic maps denoted by $\widetilde{\delta_i}$, $\widetilde{\sigma_i}$ and $\widetilde{\tau_n}$.

\begin{theorem} \label{main thm} If $\mathcal{C}$ is a symmetric monoidal abelian
category,  then  $( C^\bullet_H ,\widetilde{\delta_i}, \widetilde{\sigma_i}, \widetilde{\tau_n})$ is a
cocyclic object  in $\mathcal{C}$. 
\end{theorem}

Some essential parts of the proof are visualized by braiding diagrams to help the reader for a better 
understanding of those parts involving long series of formulas. 

\begin{proof} The most difficult part is to show that the cyclic map
$\widetilde{\tau_n}$ is well defined on $C_H^n(C,M)$ for all $n$.  For this,
we use the following diagram:
\[\begin{CD}
 M \ot H \ot C^{n+1} @>f_{M,C^{n+1}} >> M \ot C^{n+1} @>\text{coker}(f_{M,C^{n+1}})>> M \ot_H C^{n+1}\\
   @.                 @VV\tau V                        @VV\widetilde{\tau} V                          \\
 M \ot H \ot C^{n+1} @>f_{M,C^{n+1}} >> M \ot C^{n+1} @>\text{coker}(f_{M,C^{n+1}})>> M \ot_H C^{n+1}
\end{CD} \]
where $f_{M,C^{n+1}}=(\phi_M \ot 1_{C^{n+1}} - 1_M \ot
\phi_{C^{n+1}})$. So we have to show that
\[\text{coker}(f_{M,C^{n+1}})~ (\tau)~ (f_{M,C^{n+1}}) = 0 ,\]
i.e.
\begin{multline}\label{cokerf tau f =0}
\text{coker}(f_{M,C^{n+1}})~[ (1_M , \si_{C , C^n})(1_M , \phi_C , 1_{C^n})(\si_{H,M} , 1_{C^{n+1}}) \\
(\rho _M , 1_{C^{n+1}})]~(\phi_M \ot 1_{C^{n+1}} - 1_M \ot \phi_{C^{n+1}}) =0  \quad \forall n .
\end{multline}

It is not hard to see that the equality ~\eqref{cokerf tau f =0} is
equivalent to: 
\begin{eqnarray} \label{cokerf (alf - bet)=0} \nonumber
\text{coker}(f_{M,H^2})~  [~[(1_{M,H}, m_H) (1_M, \si_{H^2, H})
(\si_{H,M}, 1_{H^2}) (\rho _M , \Delta_H)] \\ 
 - [(1_M, \eta_H, 1_H) (\si_{H,M}) (\rho _M)(\phi_M)]~]  =0 .
\end{eqnarray}
From ~\eqref{cokerf (alf - bet)=0} to ~\eqref{cokerf tau f =0} one uses the diagonal action 
of $H$ on $C^{n+1}$ and the coaction of $H$ on $M$. From  ~\eqref{cokerf tau f =0} 
to ~\eqref{cokerf (alf - bet)=0} one puts $n=1$ and then acts both sides of ~\eqref{cokerf tau f =0}  
on $(1_M, 1_H, \eta_H, \eta_H)$.

If we put: 
\[\alpha = [(1_{M,H}, m_H) (1_M, \si_{H^2, H}) (\si_{H,M}, 1_{H^2}) (\rho _M , \Delta_H)],\]
 and \[\beta = [(1_M, \eta_H, 1_H) (\si_{H,M}) (\rho _M)(\phi_M)],\] then ~\eqref{cokerf (alf - bet)=0} becomes:
\[\text{coker}(f_{M,H^2})~( \alpha - \beta ) = 0.\]

To prove this equality, we will define an isomorphism
\[\widetilde{\varphi} :M \ot_H H^2 \to M \ot_H H^2 \cong M \ot H \] and will show that,
\[(\widetilde{\varphi})~ \text{coker}(f_{M,H^2})~( \alpha - \beta ) = 0 .\]
More explicitly:

\noindent $\bullet$ \textbf{Step 1:} 

We claim that $\varphi$ defined as below is a $H$-linear isomorphism,
where the domain $H^2$ is considered as a $H$-module via diagonal
action ($\phi_{H^2}$), and the codomain $H^2$ is considered as a
$H$-module via multiplication in the first factor ($\phi'_{H^2}$).
\[ \varphi := (1_H, m_H)(1_H,S_H , 1_H))(\Delta_H , 1_H) :H^2 \to H^2  \label{def vfy } .\] It is easy to see that $\varphi$ is an isomorphism and in fact its inverse map is
\[ \varphi^{-1} := (1 ,m)(\Delta ,1)  \label{def vfy inv} .  \] To see that $\varphi$ is $H$-linear, we have to show that the following diagram commutes:
\[
\begin{CD}
H^2         @>\varphi>>     H^2\\
@A\phi_{H^2}AA              @AA\phi'_{H^2}A\\
H \ot H^2   @>>1 \ot \varphi>    H \ot H^2
\end{CD}
\]
i.e.,
\begin{equation}\label{vfy is H lin}
(\varphi)(\phi_{H^2}) = (\phi'_{H^2})(1 ,\varphi),
\end{equation}
 where  $\phi_{H^2} := (m ,m)(1, \si , 1)(\Delta, 1_{H^2})$ is the diagonal action of $H$ on the domain $H^2$, and $\phi'_{H^2} := (m ,1)$ is the other action of $H$ on the codomain $H^2$.
To prove that ~\eqref{vfy is H lin} is true we see that:
\begin{flalign*}
&RHS&\\
&= (m,1)~(1,1,m)(1,1, S,1)(1,\Delta ,1)&\\
&= (m,m)(1,1,S,1)(1, \Delta ,1),&
\end{flalign*}
and
\begin{flalign*}
&LHS&\\
&\overset{(1)}{=} (1,m)(1,S,1)(\Delta,1)~(m,m)(1,\si ,1)(\Delta,1,1)&\\
&\overset{(2)}{=} (1,m)(1,S,1)(1,1, m)(\Delta m,1,1)(1,\si ,1)(\Delta,1,1)&\\
&\overset{(3)}{=} (1,m)(1,1,m)(1,S,1,1)(m,m,1,1)(1,\si ,1,1,1)(\Delta,\Delta,1,1)(1,\si ,1)(\Delta,1,1)&\\
&\overset{(4)}{=} (1,m)(1,1,m)(m,1,1,1)(1,1,Sm,1,1)(1,\si ,1,1,1)(\Delta,\Delta,1,1)(1,\si ,1)(\Delta,1,1)&\\
&\overset{(5)}{=} (1,m)(1,1,m)(m,1,1,1)(1,1,m,1,1)(1,1,S,S,1,1)&\\
&~~~~~(1,1,\si ,1,1)(1,\si ,1,1,1)(\Delta,\Delta,1,1)(1,\si ,1)(\Delta,1,1)&\\
&\overset{(6)}{=} (1,m)(1,1,m)(m,1,1,1)(1,1,m,1,1)(1,1,S,S,1,1)&\\
&~~~~~(1,\si_{H,H^2} , 1,1)(1,1, \Delta ,1,1)(\Delta,1,1,1)(1,\si ,1)(\Delta,1,1)&\\
&\overset{(7)}{=} (1,m)(1,1,m)(m,1,1,1)(1,1,m,1,1)(1,1,S,S,1,1)&\\
&~~~~~(1,\Delta,1,1,1)(1,\si , 1,1)(1,1,\si,1)(\Delta,1,1,1)(\Delta,1,1)&\\
&\overset{(8)}{=} (1,m)(1,1,m)(m,1,1,1)(1,1,m,1,1)(1,1,S,S,1,1)&\\
&~~~~~(1,\Delta,1,1,1)(1, \si_{H^2,H})(1, \Delta,1,1)(\Delta,1,1)&\\
&\overset{(9)}{=} (1,m)(1,1,m)(m,1,1,1)(1,1,m,1,1)(1,1,S,S,1,1)&\\
&~~~~~(1,\Delta,1,1,1)(1,1,\Delta,1) (1, \si, 1)(\Delta,1,1)&\\
&\overset{(10)}{=} (1,m)(m,1,1)(1,1,1,m)(1,1,1,m,1)(1,1,S,S,1,1)&\\
&~~~~~(1,1,1,\Delta,1)(1,\Delta,1,1)(1, \si, 1)(\Delta,1,1)&\\
&\overset{(11)}{=} (1,m)(m,1,1)(1,1,1,m)(1,1,S,m(S,1)\Delta,1)(1,\Delta,1,1)(1, \si, 1)(\Delta,1,1)&\\
&\overset{(12)}{=} (1,m)(m,1,1)(1,1,1,m)(1,1,S,\eta ,1)(1,1,1,\varepsilon,1)(1,\Delta,1,1)(1, \si, 1)(\Delta,1,1)&\\
&\overset{(13)}{=} (1,m)(m,1,1)(1,1,S,m(\eta ,1))(1,\Delta,1)(1,1,\varepsilon,1)(1, \si, 1)(\Delta,1,1)&\\
&\overset{(14)}{=} (m,m)(1,1,S,1)(1,\Delta,1)(1,\varepsilon,1,1)(\Delta,1,1)&\\
&\overset{(14)}{=} (m,m)(1,1,S,1)(1,\Delta,1)((1,\varepsilon)\Delta,1,1)&\\
&\overset{(15)}{=} (m,m)(1,1,S,1)(1,\Delta,1)&
\end{flalign*}
Therefore $LHS= RHS$.

 Stages (1) to (15) are explained below:
\item (1) We use definitions $\varphi= (1, m)(1,S , 1))(\Delta , 1)$ and $\phi_{H^2}= (m ,m)(1, \si , 1)(\Delta, 1_{H^2})$.
\item (2) We use $(m,m) = (1,m)(m,1,1)$, commute $(1,m)$ and $(\Delta,1)$, i.e. $(\Delta,1)(1,m)=(1,1,m)(\Delta,1,1)$, and then use $(\Delta,1,1)(m,1,1)=(\Delta m,1,1)$.
\item (3) We commute $(1,S,1)$ and $(1,1, m)$, and use $\Delta m = (m \ot m)(id \ot \psi \ot id)(\Delta \ot \Delta)$.
\item (4) We use $(m,m,1,1) = (m,1,1,1)(1,1,m,1,1)$, commute $(m,1,1,1)$ and $(1,S,1,1)$, and then compose $(1,1,m,1,1)$ with $(1,1,S,1,1)$.
\item (5) We use $Sm=m(S,S)\si$.
\item (6) We use $\si_{H,H^2}=(1, \si)(\si ,1)$ and $(\Delta,\Delta,1,1)=(1,1, \Delta ,1,1)(\Delta,1,1,1)$.
\item (7) We use the naturality of $\si$ to commute $(1,\si_{H,H^2} ,1,1)$ and $(1,1, \Delta ,1,1)$, more explicitly we use the following
commuting diagram:
\[
\begin{CD}
H^2         @>(1,\Delta)>>     H \ot H^2\\
@V\si VV              @VV\si_{H,H^2}V\\
H^2   @>>(\Delta,1)>   H^2 \ot H 
\end{CD}
\]
also we commute $(\Delta,1,1,1)$ and $(1,\si ,1)$, i.e.  $(\Delta,1,1,1)(1, \si ,1) =(1,1, \si ,1)(\Delta ,1,1,1)$.
\item (8) We  use $\si_{H^2,H}=(\si ,1)(1, \si)$ and the coassociativity $(\Delta ,1)\Delta =(1, \Delta) \Delta$.
\item (9) We commute $(1, \si_{H^2,H})$ and $(1, \Delta,1,1)$ using the naturality of $\si$ as in step $(7)$.
\item (10) We use simple commutations and the associativity $m(1,m)=m(m,1)$, in four first parenthesis, also commute $(1,\Delta,1,1,1)$ and $(1,1,\Delta,1)$.
\item (11) We use  $(1,1,1,m,1)(1,1,S,S,1,1)(1,1,1,\Delta,1)=(1,1,S,m(S,1)\Delta,1)$.
\item (12) We use $m(S,1)\Delta =\eta \varepsilon$.
\item (13) We use $(1,1,1,\varepsilon,1)(1,\Delta,1,1) =(1,\Delta,1)(1,1,\varepsilon,1)$.
\item (14) We use $m(\eta ,1)=1$, also commute $(1,1,\varepsilon,1)$ and $(1,\si, 1)$ using naturality of $\si$. In fact we use this commuting
diagram:
\[
\begin{CD}
H^2         @>(\varepsilon,1)>>     I \ot H \\
@V\si VV              @VV\si_{I,H}=id_I V\\
H^2   @>>(1,\varepsilon)>   H \ot I 
\end{CD}
\]
\item (15) We use $(1,\varepsilon)\Delta=1$.

\noindent $\bullet$ \textbf{Step 2:}

Considering $1_M :M \to M$ as a $H$-linear isomorphism we have, $1_M \ot \varphi : M \ot H^2 \cong M \ot H^2 $ and so we have:
\begin{equation}\label{def wvfy}
\widetilde{\varphi} := 1_M \ot_H \varphi : M \ot_H H^2 \cong M \ot_H H^2 \cong (M \ot_H H) \ot H \cong M \ot H   
\end{equation}

\noindent $\bullet$ \textbf{Step 3:}

Now we  prove that $(\widetilde{\varphi})~ \text{coker}(f_{M,H^2})~( \alpha-\beta)= 0 $. For that, we look at this commuting diagram
\[
\begin{CD}
 M \ot H \ot H^2 @>f_{M,H^2} >> M \ot H^2 @>\text{coker}(f_{M,H^2})>> M \ot_H H^2\\
   @.                 @VV1_M \ot \varphi V                        @VV\widetilde{\varphi}V \\
 M \ot H \ot H^2 @>f'_{M,H^2} >> M \ot H^2 @>\text{coker}(f'_{M,H^2})>> M \ot_H H^2
\end{CD}
\]
where $f_{M,H^2}=(\phi_M \ot 1_{H^2} - 1_M \ot \phi_{H^2})$ and
$f'_{M,H^2}=(\phi_M \ot 1_{H^2} - 1_M \ot \phi'_{H^2})$, which shows,
$(\widetilde{\varphi})(\text{coker}(f)) = (\text{coker}(f'))(1_M \ot \varphi)$, 
 so we instead will prove $(\text{coker}(f'))(1_M \ot \varphi)~( \alpha -\beta ) =0$, i.e.,

\begin{equation}\label{cokrf'} 
(\text{coker}(f'))(1_M \ot \varphi)( \alpha ) =(\text{coker}(f'))(1_M \ot \varphi)~ \beta ).
\end{equation}

The LHS of ~\eqref{cokrf'} is equal to:
\begin{multline*}  
\overset{(1)}{=}(\text{coker}(f'))(1_M,1,m)(1_M,1,S,1)(1_M,\Delta ,1)(1_M,1,m)\\
(1_M,\si_{H^2,H})(\si_{H,M},1,1)(\rho,\Delta)=
\end{multline*}

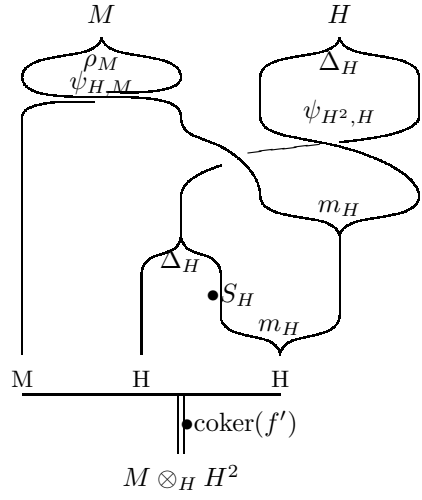
\[
\xy /r2.5pc/:,
(1,0)= "L11",
(4,0)= "L12",
  "L11"+(0,.5)*!<0pt,6pt>{M}, 
  "L12"+(0,.5)*!<0pt,6pt>{H}, 
  "L11"="x1", 
  "x1"+(-1,-.5)="x2",
  "x1"+(1,-.5)="x3",
  "x1"+(0,-.25)="c",  
   "x1";"x2"**\crv{"x1" & "c" & "x2"+(0,.25)}, 
   "x1";"x3"**\crv{"x1" & "c" & "x3"+(0,.25)}, 
  "c"+(0,0)*!<0pt,2pt>{\rho _M}, 
"x2" ="L21",
"x3" ="L22",
   "L12"="x1", 
  "x1"+(-1,-.5)="x2",
  "x1"+(1,-.5)="x3",
  "x1"+(0,-.25)="c",  
   "x1";"x2"**\crv{"x1" & "c" & "x2"+(0,.25)}, 
   "x1";"x3"**\crv{"x1" & "c" & "x3"+(0,.25)}, 
  "c"+(0,0)*!<0pt,2pt>{\Delta_H},  
"L21" ="L11",
"L22" ="L12",
"x2"="L13",
"x3"="L14",
  "L11"="x1",       
  "L12"="x2",  , 
  "x1"+(0,-.5)="x3",
  "x2"+(0,-.5)="x4",
  "x1";"x2"**\dir{}?(.5)="m", 
  "m" + (0,-.25)="c", 
  "x1";"x4"**\crv{ "x1"+(0,-.25) & "c" & "x4"+(0,.25)},  
  "x2";"c"+<3pt,2pt>="c2"**\crv{"x2"+(0,-.25) & "c2"},
  "x3";"c"+<-3pt,-2pt>="c1"**\crv{"x3"+(0,.25) & "c1"}, 
  "c"+(0,.10)*!<0pt,-2pt>{\si_{H,M}},
"x3"="L11",
"x4"="L12",
  "L13";"L13" +(0,-.5)="L13"**\dir{-},
  "L14";"L14" +(0,-.5)="L14"**\dir{-},
  "L11";"L11" +(0,-1)="L11"**\dir{-},
"L12"="a1",
"L13"="a2",
"L14"="a3",
"a1" + (0,-1)="b1",
"a2" + (0,-1)="b2",
"a3" + (0,-1)="b3",
"a1";"a2"**\dir{}?(.5)="m12", 
  "m12" + (0,-.25)="c12", 
"a2";"a3"**\dir{}?(.5)="m23", 
  "m23" + (0,-.25)="c23", 
  "a1";"b2"**\crv{ "a1"+(0,-.25) & "c12" & "b2"+(0,.25)},  
  "a2";"b3"**\crv{ "a2"+(0,-.25) & "c23" & "b3"+(0,.25)}, 
                "b1";"c12"+<-0pt,-11pt>="c21"**\crv{"b1"+(0,.25) & "c21"},         
               "a3";"c23"+<-0pt,-2pt>="c32"**\crv{"a3"+(0,-.25) & "c32"},
"c12" +(.35,-.20);"c23" +(-.25,-.10)**\dir{-}, 
"c23"+(0,.10)*!<0pt,-6pt>{\si_{H^2,H}}, 
"b1"="L12",
"b2"="L13",
"b3"="L14",
  "L11";"L11" +(0,-.5)="L11"**\dir{-},
  "L12";"L12" +(0,-.5)="L12"**\dir{-},
   "L13"="x1", 
  "L14"="x2",   
  "x1";"x2"**\dir{}?(.5)="m", 
  "m" + (0,-.25)="c", 
  "m" + (0,-.5)="x3", 
  "x1";"x3"**\crv{"x1"-(0,.25) & "c" & "x3"},
  "x2";"x3"**\crv{"x2"-(0,.25) & "c"  & "x3"}, 
  "c"+(0,.10)*!<0pt,-0pt>{m_H}, 
"x3"="L13",
  "L11";"L11" +(0,-.5)="L11"**\dir{-},
   "L12"="x1", 
  "x1"+(-.5,-.5)="x2",
  "x1"+(.5,-.5)="x3",
  "x1"+(0,-.25)="c",  
   "x1";"x2"**\crv{"x1" & "c" & "x2"+(0,.25)}, 
   "x1";"x3"**\crv{"x1" & "c" & "x3"+(0,.25)}, 
  "c"+(0,0)*!<0pt,2pt>{\Delta_H},  
"x2"="L21",
"x3"="L22",
  "L13";"L13" +(0,-.5)="L14"**\dir{-},
"L21" ="L12",
"L22" ="L13",
  "L11";"L11" +(0,-.5)="L11"**\dir{-},
  "L12";"L12" +(0,-.5)="L12"**\dir{-},
  "L13"="x1", 
  "x1"+(0,-.5)="x2",
  "x1";"x2"**\dir{}?(.5)="c", 
  "x1";"x2"**\dir{-},   
  "c"-(0,0)*!<-4pt,0pt>{\bullet S_H },
"x2"="L13",
  "L14";"L14" +(0,-.5)="L14"**\dir{-},
  "L11";"L11" +(0,-.5)="L11"**\dir{-},
  "L12";"L12" +(0,-.5)="L12"**\dir{-},
   "L13"="x1", 
  "L14"="x2",   
  "x1";"x2"**\dir{}?(.5)="m", 
  "m" + (0,-.25)="c", 
  "m" + (0,-.5)="x3", 
  "x1";"x3"**\crv{"x1"-(0,.25) & "c" & "x3"},
  "x2";"x3"**\crv{"x2"-(0,.25) & "c"  & "x3"}, 
  "c"+(0,.10)*!<0pt,-0pt>{m_H}, 
"x3"="L13",
  "L11"-(0,.5)*!<0pt,-6pt>{\text{\small{M}}},
  "L12"-(0,.5)*!<0pt,-6pt>{\text{\small{H}}},
  "L13"-(0,.5)*!<0pt,-6pt>{\text{\small{H}}},
"L11"-(0,.5) ;"L13"-(0,.5) **\dir{-},
  "L12" -(-.5,.5)="x1", 
  "x1"+(0,-.75)="x2",
  "x1";"x2"**\dir{}?(.5)="c", 
  "x1";"x2"**\dir{=},   
  "c"-(0,0)*!<-22pt,0pt>{\bullet \text{coker}(f') },
"x2"-(0,.5)*!<0pt,-6pt>{M \ot_H H^2 },
\endxy
\]

\begin{align*} 
&\overset{(1)}{=} (\text{coker}(f'))(1_M,1,m)(1_M,1,S,1)(1_M,\Delta,1)(1_M,1,m)\\
&~~~~~(1_M,\si_{H^2,H})(\si_{H,M},1,1)(\rho,\Delta)=\\
&\overset{(2)}{=} (\text{coker}(f'))~(1_M,1,m)(1_M,1,S,1)(1_M,\Delta ,1)(1_M, \si)(1_M,m ,1) (\si_{H,M}, 1,1) (\rho,\Delta)\\
&\overset{(3)}{=} (\text{coker}(f'))~(1_M,1,m)(1_M,1,S,1) (1_M, \si_{H, H^2})(1_M,1, \Delta )(1_M,m ,1)\\
&~~~~~(\si_{H,M}, 1,1) (\rho , \Delta)\\
&\overset{(4)}{=} (\text{coker}(f'))~(1_M,1,m)(1_M,1,S,1) (1_M, \si_{H, H^2})(1_M,m ,1,1) (\si_{H,M}, 1,1,1)\\ 
&~~~~~(1, 1_M,1, \Delta )(\rho , \Delta)\\
&\overset{(5)}{=}(\text{coker}(f'))(1_M,1,m)(1_M,1,S,1)(1_M,1,1,m)(1_M, \si_{H^2, H^2})(\si_{H,M},1,1,1)(\rho \Delta^2)\\
&\overset{(6)}{=} (\text{coker}(f'))~(1_M,1,m(S,m))(1_M, \si_{H^2, H^2})(\si_{H,M}, 1,1,1)(\rho , \Delta^2)\\
&\overset{(7)}{=} (\text{coker}(f'))~(1_M,m(1, \eta),m(S,m))(1_M, \si_{H^2, H^2})(\si_{H,M}, 1,1,1)(\rho , \Delta^2)\\
&\overset{(8)}{=} (\text{coker}(f'))~(1_M,m,1)(1_M,1, \eta,1) (1_M,1,m(S,m))(1_M, \si_{H^2, H^2})\\
&~~~~~(\si_{H,M}, 1,1,1)(\rho , \Delta^2)\\
\end{align*} 
\begin{align*}
&\overset{(9)}{=} (\text{coker}(f'))~(1_M,\phi'_{H^2})(1_M,1, \eta,1)(1_M,1,m(S,m))(1_M, \si_{H^2, H^2})\\
&~~~~~(\si_{H,M}, 1,1,1)(\rho , \Delta^2)\\
&\overset{(10)}{=} (\text{coker}(f'))~(\phi_M,1,1)(1_M,1, \eta,1) (1_M,1,m(S,m))\\
&~~~~~\underbrace{ (1_M, \si_{H^2, H^2})(\si_{H,M}, 1,1,1)} ~  (\rho , \Delta^2)=
\end{align*}

\begin{equation}\label{LHS2}  
\xy /r2.5pc/:,
(1,0)= "L11",
(4,0)= "L12",
  "L11"+(0,.5)*!<0pt,6pt>{M}, 
  "L12"+(0,.5)*!<0pt,6pt>{H}, 
  "L11"="x1", 
  "x1"+(-.75,-.5)="x2",
  "x1"+(.75,-.5)="x3",
  "x1"+(0,-.25)="c", 
   "x1";"x2"**\crv{"x1" & "c" & "x2"+(0,.25)}, 
   "x1";"x3"**\crv{"x1" & "c" & "x3"+(0,.25)}, 
  "c"+(0,0)*!<0pt,2pt>{\rho _M}, 
"x2" ="L21",
"x3" ="L22",
  "L12"="x1", 
  "x1"+(-1.25,-.5)="x2",
  "x1"+(0,-.5)="x3",
  "x1"+(1.25,-.5)="x4",
  "x1"+(0,-.25)="c",  
   "x1";"x2"**\crv{"x1" & "c" & "x2"+(0,.25)}, 
   "x1";"x3"**\crv{"x1" & "c" & "x3"+(0,.25)}, 
   "x1";"x4"**\crv{"x1" & "c" & "x4"+(0,.25)}, 
  "c"+(0,0)*!<8pt,5pt>{\Delta_H^2}, 
"L21" ="L11",
"L22" ="L12",
"x2"="L13",
"x3"="L14",
"x4"="L15",
					  "L11" + (-.5,0);"L15"+ (.5,0) **\dir{.},
  "L11"="x1",       
  "L12"="x2",  , 
  "x1"+(0,-.5)="x3",
  "x2"+(0,-.5)="x4",
  "x1";"x2"**\dir{}?(.5)="m", 
  "m" + (0,-.25)="c", 
  "x1";"x4"**\crv{ "x1"+(0,-.25) & "c" & "x4"+(0,.25)},  
  "x2";"c"+<3pt,2pt>="c2"**\crv{"x2"+(0,-.25) & "c2"},
  "x3";"c"+<-3pt,-2pt>="c1"**\crv{"x3"+(0,.25) & "c1"}, 
  "c"+(0,.10)*!<0pt,-2pt>{\si_{H,M}},
"x3"="L11",
"x4"="L12",
  "L13";"L13" +(0,-.5)="L13"**\dir{-},
  "L14";"L14" +(0,-.5)="L14"**\dir{-},
  "L15";"L15" +(0,-.5)="L15"**\dir{-},
  "L11";"L11" +(0,-1)="L11"**\dir{-},
"L12"="a1",
"L13"="a2",
"L14"="a3",
"L15"="a4",
"a1" + (0,-1)="b1",
"a2" + (0,-1)="b2",
"a3" + (0,-1)="b3",
"a4" + (0,-1)="b4",
"a1";"a3"**\dir{}?(.5)="m13", 
  "m13" + (0,-.25)="c13", 
"a1";"a2"**\dir{}?(.5)="m12", 
  "m12" + (0,-.25)="c12", 
"a2";"a3"**\dir{}?(.5)="m23", 
  "m23" + (0,-.25)="c23", 
"a2";"a4"**\dir{}?(.5)="m24", 
  "m24" + (0,-.25)="c24", 
  "a1";"b3"**\crv{ "a1"+(0,-.25) & "c13" & "b3"+(0,.25)},  
  "a2";"b4"**\crv{ "a2"+(0,-.25) & "c24" & "b4"+(0,.25)}, 
                "b1";"c12"+<5pt,-11pt>="c21"**\crv{"b1"+(0,.25) & "c21"},  
                "b2";"c23"+<0pt,-13pt>="c32"**\crv{"b2"+(0,.25) & "c32"},       
               "a3";"c23"+<3pt,0pt>="c32"**\crv{"a3"+(0,-.25) & "c32"},
               "a4";"c24"+<7pt,-4pt>="c42"**\crv{"a4"+(0,-.25) & "c42"},
"c23"+(0,.10)*!<0pt,-6pt>{\si_{H^2,H^2}}, 
"L11"="L11",
"b1"="L12",
"b2"="L13",
"b3"="L14",
"b4"="L15",
					  "L11" + (-.5,0);"L15"+ (.5,0) **\dir{.},
"L11";"L11" +(0,-1)="L11"**\dir{-},
"L12";"L12" +(0,-1)="L12"**\dir{-},
  "L13"="x1", 
  "x1"+(0,-.5)="x2",
  "x1";"x2"**\dir{}?(.5)="c", 
  "x1";"x2"**\dir{-},   
  "c"-(0,0)*!<-4pt,0pt>{\bullet S_H },
"x2"="L13",
   "L14"="x1", 
  "L15"="x2",   
  "x1";"x2"**\dir{}?(.5)="m", 
  "m" + (0,-.25)="c", 
  "m" + (0,-.5)="x3", 
  "x1";"x3"**\crv{"x1"-(0,.25) & "c" & "x3"},
  "x2";"x3"**\crv{"x2"-(0,.25) & "c"  & "x3"}, 
  "c"+(0,.10)*!<0pt,-0pt>{m_H}, 
"x3"="L14",
   "L13"="x1", 
  "L14"="x2",   
  "x1";"x2"**\dir{}?(.5)="m", 
  "m" + (0,-.25)="c", 
  "m" + (0,-.5)="x3", 
  "x1";"x3"**\crv{"x1"-(0,.25) & "c" & "x3"},
  "x2";"x3"**\crv{"x2"-(0,.25) & "c"  & "x3"}, 
  "c"+(0,.10)*!<0pt,-0pt>{m_H}, 
"x3"="L13",
   "L12" + (.5,0)="x1", 
  "x1"+(0,-.25)="x2",
  "x1";"x2"**\dir{}?(.5)="c", 
  "x1";"x2"**\dir{-},  
  "c"-(0,0)*!<-3pt,0pt>{\bullet \eta },  
  "x1"+(0,0)*!<0pt,-2pt>{\circ},
"L13"="L14",
"x2"="L13",
   "L11"="x1",
  "L12"="x2",   
  "x1";"x2"**\dir{}?(.5)="m", 
  "m" + (0,-.25)="c", 
  "m" + (0,-.5)="x3", 
  "x1";"x3"**\crv{"x1"-(0,.25) & "c" & "x3"},
  "x2";"x3"**\crv{"x2"-(0,.25) & "c"  & "x3"}, 
  "c"+(0,.10)*!<0pt,-0pt>{\phi_M}, 
"x3"="L11",
"L13";"L13" +(0,-.25)="L12"**\dir{-},
"L14";"L14" +(0,-.5)="L13"**\dir{-},
  "L11"-(0,.5)*!<0pt,-6pt>{\text{\small{M}}}, 
  "L12"-(0,.5)*!<0pt,-6pt>{\text{\small{H}}}, 
  "L13"-(0,.5)*!<0pt,-6pt>{\text{\small{H}}}, 
"L11"-(0,.5) ;"L13"-(0,.5) **\dir{-},
  "L12" -(0,.5)="x1", 
  "x1"+(0,-.75)="x2",
  "x1";"x2"**\dir{}?(.5)="c", 
  "x1";"x2"**\dir{=},   
  "c"-(0,0)*!<-22pt,0pt>{\bullet \text{coker}(f') },
"x2"-(0,.5)*!<0pt,-6pt>{M \ot_H H^2 },
\endxy 
\end{equation}
Stages (1) to (10) are explained below with (10) perhaps the most important one:
\item (1) We use $\varphi= (1, m)(1,S , 1))(\Delta , 1)$ and $\alpha = [(1_{M,H}, m_H) (1_M, \si_{H^2, H}) (\si_{H,M}, 1_{H^2}) (\rho , \Delta)]$.
\item (2) We use the naturality of $\si$ to commute $(1_M,1, m)$ and $(1_M, \si_{H^2, H})$,  i.e. we use the diagram:
\[
\begin{CD}
H^2 \ot H         @>(m,1)>>     H^2\\
@V\si_{H,H^2}VV              @VV\si V\\
H \ot H^2   @>>(1,m)>           H^2 
\end{CD}
\]
\item (3) We commute $(1_M,\Delta ,1)$ and $(1_M, \si)$, again using naturality of $\si$ :
\[
\begin{CD}
H^2          @>(1,\Delta)>>     H \ot H^2\\
@V\si VV              @VV\si_{H,H^2} V\\
H^2    @>>(\Delta,1)>      H^2 \ot H
\end{CD}
\]
\item (4) We commute $(1_M,1, \Delta )$ first with $(1_M,m ,1)$ and next with $ (\si_{H,M}, 1,1)$.
\item (5) We commute $ (1_M, \si_{H, H^2})$ and $(1_M,m ,1,1)$ using the diagram:
\[
\begin{CD}
H^2 \ot H^2         @>(m,1,1)>>     H \ot H^2\\
@V\si_{H^2,H^2} VV              @VV\si_{H,H^2} V\\
H^2 \ot  H^2   @>>(1,1,m)>      H^2 \ot H
\end{CD}
\]
Also we use $(1, 1_M,1, \Delta )(\rho , \Delta)=(\rho , \Delta^2)$ and use the definition $\Delta^2 := (1, \Delta )\Delta$.
\item (6) We use $(1_M,1,m)(1_M,1,S,1)(1_M,1,1,m)=(1_M,1,m(S,m))$.
\item (7) We use $1=m(1,\eta)$.
\item (8) We use $(1_M,m(1, \eta),m(S,m))=(1_M,m,1)(1_M,1, \eta,1) (1_M,1,m(S,m))$.
\item (9) We use the definition $\phi'_{H^2}=(m,1)$.
\item (10) We use $\text{coker}(f')~(1_M,\phi'_{H^2})=\text{coker}(f')~(\phi_M,1,1)$,
which comes from the definition of cokernel and this diagram:
\[
\begin{CD}
 M \ot H \ot H^2 @>{f'_{M,H^2}=(\phi_M \ot 1_{H^2} - 1_M \ot \phi'_{H^2})} >> M \ot H^2 @>\text{coker}(f')>> M \ot_H H^2
\end{CD}
\]

The RHS of ~\eqref{cokrf'} is equal to:
\begin{flalign*} 
&\overset{(1)}{=} (\text{coker}(f'))~(1_M,1,m)(1_M,1,S,1)(1_M,\Delta ,1)~(1_M, \eta, 1) (\si_{H,M}) (\rho )(\phi)&\\
&\overset{(2)}{=} (\text{coker}(f'))~(1_M,1,m)(1_M,1,S,1)(1_M, \eta,\eta,  1) (\si_{H,M}) (\rho )(\phi)&\\
&\overset{(3)}{=} (\text{coker}(f'))~(1_M,1,m)(1_M, \eta,\eta,  1) (\si_{H,M}) (\rho )(\phi)&\\
&\overset{(4)}{=} (\text{coker}(f'))~(1_M, \eta,1) (\si_{H,M}) (\rho )(\phi) =&
\end{flalign*}

\begin{equation}\label{pic1of the proofRHSpage13}
\xy /r2.5pc/:, 
(1,0)="L11",
(3,0)="L12",
  "L11"="x1",
  "L12"="x2",   
  "x1";"x2"**\dir{}?(.5)="m", 
  "m" + (0,-.25)="c", 
  "m" + (0,-.5)="x3", 
  "x1";"x3"**\crv{"x1"-(0,.25) & "c" & "x3"},
  "x2";"x3"**\crv{"x2"-(0,.25) & "c"  & "x3"}, 
  "c"+(0,.10)*!<0pt,-0pt>{\phi_{\text{\small{M}}}},  
  "x1"+(0,.5)*!<0pt,6pt>{\text{\small{M}}},  
  "x2"+(0,.5)*!<0pt,6pt>{\text{\small{H}}},  
"x3"="L11", 
   "L11"="x1", 
  "x1"+(-1,-.5)="x2",
  "x1"+(1,-.5)="x3",
  "x1"+(0,-.25)="c",  
   "x1";"x2"**\crv{"x1" & "c" & "x2"+(0,.25)}, 
   "x1";"x3"**\crv{"x1" & "c" & "x3"+(0,.25)}, 
  "c"+(0,0)*!<0pt,-0pt>{\rho _M},  
"x2"="L11", 
"x3"="L12", 
   "L11" ="x1",       
  "L12"  ="x2",  , 
  "x1"+(0,-.5)="x3",
  "x2"+(0,-.5)="x4",
  "x1";"x2"**\dir{}?(.5)="m", 
  "m" + (0,-.25)="c", 
  "x1";"x4"**\crv{ "x1"+(0,-.25) & "c" & "x4"+(0,.25)},  
  "x2";"c"+<3pt,2pt>="c2"**\crv{"x2"+(0,-.25) & "c2"},
  "x3";"c"+<-3pt,-2pt>="c1"**\crv{"x3"+(0,.25) & "c1"}, 
  "c"+(0,0)*!<0pt,-6pt>{\si_{H,M}}, 
"x3"="L11", 
"x4"="L13",
"L13" +(-1,0)="L12",  
  "L11" ;"L11" +(0,-.5)="L11" **\dir{-}, 
  "L13" ;"L13"+(0,-.5)="L13" **\dir{-}, 
   "L12" + (0,-.25)="x1", 
  "x1"+(0,-.25)="x2",
  "x1";"x2"**\dir{}?(.5)="c", 
  "x1";"x2"**\dir{-},   
  "c"-(0,0)*!<-3pt,0pt>{\bullet \eta },  
   "x1"+(0,0)*!<0pt,-2pt>{\circ},  
  "x2"-(0,.5)*!<0pt,-6pt>{\text{\small{H}}},  
  "L11"-(0,.5)*!<0pt,-6pt>{\text{\small{M}}},
  "L13"-(0,.5)*!<0pt,-6pt>{\text{\small{H}}},
 "x2"= "L12",
"L11"-(0,.5) ;"L13"-(0,.5) **\dir{-},
  "L12" -(0,.5)="x1", 
  "x1"+(0,-.75)="x2",
  "x1";"x2"**\dir{}?(.5)="c", 
  "x1";"x2"**\dir{=},   
  "c"-(0,0)*!<-22pt,0pt>{\bullet \text{coker}(f') },
"x2"-(0,.5)*!<0pt,-6pt>{M \ot_H H^2 },
  "L13"+(.75,.75)*!<0pt,0pt>{\overset{~\eqref{BSAYD}}{=}}, 
(5,0)= "L11",
(8,0)= "L12",
  "L11"+(0,.5)*!<0pt,6pt>{M}, 
  "L12"+(0,.5)*!<0pt,6pt>{H}, 
  "L11"="x1", 
  "x1"+(-.5,-.5)="x2",
  "x1"+(.5,-.5)="x3",
  "x1"+(0,-.25)="c",  
   "x1";"x2"**\crv{"x1" & "c" & "x2"+(0,.25)}, 
   "x1";"x3"**\crv{"x1" & "c" & "x3"+(0,.25)}, 
  "c"+(0,0)*!<0pt,3pt>{\rho _M}, 
"x2" ="L21",
"x3" ="L22",
  "L12"="x1", 
  "x1"+(-1,-.5)="x2",
  "x1"+(0,-.5)="x3",
  "x1"+(1,-.5)="x4",
  "x1"+(0,-.25)="c",  
   "x1";"x2"**\crv{"x1" & "c" & "x2"+(0,.25)}, 
   "x1";"x3"**\crv{"x1" & "c" & "x3"+(0,.25)}, 
   "x1";"x4"**\crv{"x1" & "c" & "x4"+(0,.25)}, 
  "c"+(0,0)*!<8pt,5pt>{\Delta_H^2}, 
"L21" ="L11",
"L22" ="L12",
"x2"="L13",
"x3"="L14",
"x4"="L15",
  "L11";"L11" +(0,-.5)="L11"**\dir{-},
  "L12"="x1",       
  "L13"="x2",  , 
  "x1"+(0,-.5)="x3",
  "x2"+(0,-.5)="x4",
  "x1";"x2"**\dir{}?(.5)="m", 
  "m" + (0,-.25)="c", 
  "x1";"x4"**\crv{ "x1"+(0,-.25) & "c" & "x4"+(0,.25)},  
  "x2";"c"+<3pt,2pt>="c2"**\crv{"x2"+(0,-.25) & "c2"},
  "x3";"c"+<-3pt,-2pt>="c1"**\crv{"x3"+(0,.25) & "c1"}, 
  "c"+(0,.10)*!<0pt,-6pt>{\si_{M,H}},
"x3"="L12",
"x4"="L13",
  "L14";"L14" +(0,-.5)="L14"**\dir{-},
  "L15";"L15" +(0,-.5)="L15"**\dir{-},
  "L11";"L11" +(0,-.5)="L11"**\dir{-},
  "L12";"L12" +(0,-.5)="L12"**\dir{-},
  "L13";"L13" +(0,-.5)="L13"**\dir{-},
  "L14"="x1",       
  "L15"="x2",  , 
  "x1"+(0,-.5)="x3",
  "x2"+(0,-.5)="x4",
  "x1";"x2"**\dir{}?(.5)="m", 
  "m" + (0,-.25)="c", 
  "x1";"x4"**\crv{ "x1"+(0,-.25) & "c" & "x4"+(0,.25)},  
  "x2";"c"+<3pt,2pt>="c2"**\crv{"x2"+(0,-.25) & "c2"},
  "x3";"c"+<-3pt,-2pt>="c1"**\crv{"x3"+(0,.25) & "c1"}, 
  "c"+(0,.10)*!<0pt,-6pt>{\si_{H,H}},
"x3"="L14",
"x4"="L15",
"L11";"L11" +(0,-.5)="L11"**\dir{-},
"L12";"L12" +(0,-.5)="L12"**\dir{-},
  "L13"="x1",       
  "L14"="x2",  , 
  "x1"+(0,-.5)="x3",
  "x2"+(0,-.5)="x4",
  "x1";"x2"**\dir{}?(.5)="m", 
  "m" + (0,-.25)="c", 
  "x1";"x4"**\crv{ "x1"+(0,-.25) & "c" & "x4"+(0,.25)},  
  "x2";"c"+<3pt,2pt>="c2"**\crv{"x2"+(0,-.25) & "c2"},
  "x3";"c"+<-3pt,-2pt>="c1"**\crv{"x3"+(0,.25) & "c1"}, 
  "c"+(0,.10)*!<0pt,-6pt>{\si_{M,H}},
"x3"="L13",
"x4"="L14",
  "L15";"L15" +(0,-.5)="L15"**\dir{-},
"L11"="a1",
"L12"="a2",
"L13"="a3",
"a1" + (0,-1)="b1",
"a2" + (0,-1)="b2",
"a3" + (0,-1)="b3",
"a1";"a2"**\dir{}?(.5)="m12", 
  "m12" + (0,-.25)="c12", 
"a2";"a3"**\dir{}?(.5)="m23", 
  "m23" + (0,-.25)="c23", 
  "a1";"b2"**\crv{ "a1"+(0,-.25) & "c12" & "b2"+(0,.25)},  
  "a2";"b3"**\crv{ "a2"+(0,-.25) & "c23" & "b3"+(0,.25)}, 
                "b1";"c12"+<-0pt,-11pt>="c21"**\crv{"b1"+(0,.25) & "c21"},         
               "a3";"c23"+<-0pt,-2pt>="c32"**\crv{"a3"+(0,-.25) & "c32"},
"c12" +(.35,-.20);"c23" +(-.25,-.10)**\dir{-}, 
"c23"+(0,.10)*!<0pt,-6pt>{\si_{H^2,H}}, 
"b1"="L11",
"b2"="L12",
"b3"="L13",
  "L14";"L14" +(0,-1)="L14"**\dir{-},
  "L15";"L15" +(0,-1)="L15"**\dir{-},
  "L11"="x1", 
  "x1"+(0,-.5)="x2",
  "x1";"x2"**\dir{}?(.5)="c", 
  "x1";"x2"**\dir{-},   
  "c"-(0,0)*!<-4pt,0pt>{\bullet S_H },
"x2"="L11",
   "L12"="x1", 
  "L13"="x2",  
  "x1";"x2"**\dir{}?(.5)="m", 
  "m" + (0,-.25)="c", 
  "m" + (0,-.5)="x3", 
  "x1";"x3"**\crv{"x1"-(0,.25) & "c" & "x3"},
  "x2";"x3"**\crv{"x2"-(0,.25) & "c"  & "x3"}, 
  "c"+(0,.10)*!<0pt,-0pt>{m_H}, 
"x3"="L12",
   "L11"="x1", 
  "L12"="x2",   
  "x1";"x2"**\dir{}?(.5)="m", 
  "m" + (0,-.25)="c", 
  "m" + (0,-.5)="x3", 
  "x1";"x3"**\crv{"x1"-(0,.25) & "c" & "x3"},
  "x2";"x3"**\crv{"x2"-(0,.25) & "c"  & "x3"}, 
  "c"+(0,.10)*!<0pt,-0pt>{m_H}, 
"x3"="L11",
   "L14"="x1", 
  "L15"="x2",   
  "x1";"x2"**\dir{}?(.5)="m", 
  "m" + (0,-.25)="c", 
  "m" + (0,-.5)="x3", 
  "x1";"x3"**\crv{"x1"-(0,.25) & "c" & "x3"},
  "x2";"x3"**\crv{"x2"-(0,.25) & "c"  & "x3"}, 
  "c"+(0,.10)*!<0pt,-0pt>{\phi_M}, 
"x3"="L12",
  "L12";"L12" +(0,-.5)="L12"**\dir{-},
   "L11" ="x1",       
  "L12"  ="x2",  , 
  "x1"+(0,-.5)="x3",
  "x2"+(0,-.5)="x4",
  "x1";"x2"**\dir{}?(.5)="m", 
  "m" + (0,-.25)="c", 
  "x1";"x4"**\crv{ "x1"+(0,-.25) & "c" & "x4"+(0,.25)},  
  "x2";"c"+<3pt,2pt>="c2"**\crv{"x2"+(0,-.25) & "c2"},
  "x3";"c"+<-3pt,-2pt>="c1"**\crv{"x3"+(0,.25) & "c1"}, 
  "c"+(0,0)*!<0pt,-6pt>{\si_{H,M}}, 
"x3"="L11", 
"x4"="L13",
"L13" +(-1,0)="L12",  
  "L11" ;"L11" +(0,-.5)="L11" **\dir{-}, 
  "L13" ;"L13"+(0,-.5)="L13" **\dir{-}, 
   "L12" + (0,-.25)="x1", 
  "x1"+(0,-.25)="x2",
  "x1";"x2"**\dir{}?(.5)="c", 
  "x1";"x2"**\dir{-},   
  "c"-(0,0)*!<-3pt,0pt>{\bullet \eta },  
  "x1"+(0,0)*!<0pt,-2pt>{\circ},  
  "x2"-(0,.5)*!<0pt,-6pt>{\text{\small{H}}},  
  "L11"-(0,.5)*!<0pt,-6pt>{\text{\small{M}}},
  "L13"-(0,.5)*!<0pt,-6pt>{\text{\small{H}}},
 "x2"= "L12",
"L11"-(0,.5) ;"L13"-(0,.5) **\dir{-},
  "L12" -(0,.5)="x1", 
  "x1"+(0,-.75)="x2",
  "x1";"x2"**\dir{}?(.5)="c", 
  "x1";"x2"**\dir{=},   
  "c"-(0,0)*!<-22pt,0pt>{\bullet \text{coker}(f') },
"x2"-(0,.5)*!<0pt,-6pt>{M \ot_H H^2 },
\endxy
\end{equation}

\begin{flalign*} 
&\overset{(5)}{=} (\text{coker}(f'))~(1_M, \eta,1) (\si_{H,M})~ (m(S,m) ,\phi)(\si_{H^2,H},1_M,1)(1,1, \si_{M,H},1)(1,1,1_M,\si)&\\
&~~~~~(1, \si_{M,H},1,1) (\rho , \Delta^2)&\\
&\overset{(6)}{=} (\text{coker}(f'))~(1_M, \eta,1) (\si_{H,M})~(1, \phi) (m(S,m) ,1_M ,1)(\si_{H^2,H},1_M,1)(1,1, \si_{M,H},1)&\\
&~~~~~(1,1,1_M,\si)(1, \si_{M,H},1,1) (\rho , \Delta^2)&\\
\end{flalign*}
\begin{flalign*}
&\overset{(7)}{=} (\text{coker}(f'))~(1_M, \eta,1) ( \phi,1) (\si_{H,M \ot H}) (m(S,m) ,1_M ,1)(\si_{H^2,H},1_M,1)(1,1, \si_{M,H},1)&\\
&~~~~~(1,1,1_M,\si)(1, \si_{M,H},1,1) (\rho , \Delta^2)&\\
&\overset{(8)}{=} (\text{coker}(f'))~(1_M, \eta,1) ( \phi,1)  (1_M ,1,m(S,m))(\si_{H^3,M \ot H})  (\si_{H^2,H},1_M,1)(1,1, \si_{M,H},1)&\\
&~~~~~(1,1,1_M,\si)(1, \si_{M,H},1,1) (\rho , \Delta^2)&\\
&\overset{(9)}{=} (\text{coker}(f'))~( \phi,1,1)  (1_M,1, \eta,1) (1_M ,1,m(S,m))&\\
&~~~~~\underbrace{ (\si_{H^3,M \ot H})  (\si_{H^2,H},1_M,1)(1,1, \si_{M,H},1)(1,1,1_M,\si)(1, \si_{M,H},1,1) } (\rho , \Delta^2)= &
\end{flalign*}

\begin{equation}\label{RHS2}  
\xy /r2.5pc/:,
(1,0)= "L11",
(4,0)= "L12",
  "L11"+(0,.5)*!<0pt,6pt>{M}, 
  "L12"+(0,.5)*!<0pt,6pt>{H}, 
  "L11"="x1", 
  "x1"+(-.75,-.5)="x2",
  "x1"+(.75,-.5)="x3",
  "x1"+(0,-.25)="c",  
   "x1";"x2"**\crv{"x1" & "c" & "x2"+(0,.25)}, 
   "x1";"x3"**\crv{"x1" & "c" & "x3"+(0,.25)}, 
  "c"+(0,0)*!<0pt,3pt>{\rho _M}, 
"x2" ="L21",
"x3" ="L22",
  "L12"="x1", 
  "x1"+(-1,-.5)="x2",
  "x1"+(0,-.5)="x3",
  "x1"+(1,-.5)="x4",
  "x1"+(0,-.25)="c",  
   "x1";"x2"**\crv{"x1" & "c" & "x2"+(0,.25)}, 
   "x1";"x3"**\crv{"x1" & "c" & "x3"+(0,.25)}, 
   "x1";"x4"**\crv{"x1" & "c" & "x4"+(0,.25)}, 
  "c"+(0,0)*!<8pt,5pt>{\Delta_H^2}, 
"L21" ="L11",
"L22" ="L12",
"x2"="L13",
"x3"="L14",
"x4"="L15",
					  "L11" + (-.5,0);"L15"+ (.5,0) **\dir{.},
  "L11";"L11" +(0,-.5)="L11"**\dir{-},
  "L12"="x1",       
  "L13"="x2",  , 
  "x1"+(0,-.5)="x3",
  "x2"+(0,-.5)="x4",
  "x1";"x2"**\dir{}?(.5)="m", 
  "m" + (0,-.25)="c", 
  "x1";"x4"**\crv{ "x1"+(0,-.25) & "c" & "x4"+(0,.25)},  
  "x2";"c"+<3pt,2pt>="c2"**\crv{"x2"+(0,-.25) & "c2"},
  "x3";"c"+<-3pt,-2pt>="c1"**\crv{"x3"+(0,.25) & "c1"}, 
  "c"+(0,.10)*!<0pt,-1pt>{\si_{M,H}},
"x3"="L12",
"x4"="L13",
  "L14";"L14" +(0,-.5)="L14"**\dir{-},
  "L15";"L15" +(0,-.5)="L15"**\dir{-},
  "L11";"L11" +(0,-.5)="L11"**\dir{-},
  "L12";"L12" +(0,-.5)="L12"**\dir{-},
  "L13";"L13" +(0,-.5)="L13"**\dir{-},
  "L14"="x1",       
  "L15"="x2",  , 
  "x1"+(0,-.5)="x3",
  "x2"+(0,-.5)="x4",
  "x1";"x2"**\dir{}?(.5)="m", 
  "m" + (0,-.25)="c", 
  "x1";"x4"**\crv{ "x1"+(0,-.25) & "c" & "x4"+(0,.25)},  
  "x2";"c"+<3pt,2pt>="c2"**\crv{"x2"+(0,-.25) & "c2"},
  "x3";"c"+<-3pt,-2pt>="c1"**\crv{"x3"+(0,.25) & "c1"}, 
  "c"+(0,.10)*!<0pt,-6pt>{\si_{H,H}},
"x3"="L14",
"x4"="L15",
"L11";"L11" +(0,-.5)="L11"**\dir{-},
"L12";"L12" +(0,-.5)="L12"**\dir{-},
  "L13"="x1",       
  "L14"="x2",  , 
  "x1"+(0,-.5)="x3",
  "x2"+(.5,-.5)="x4",
  "x1";"x2"**\dir{}?(.5)="m", 
  "m" + (0,-.25)="c", 
  "x1";"x4"**\crv{ "x1"+(0,-.25) & "c" & "x4"+(0,.25)},  
  "x2";"c"+<3pt,2pt>="c2"**\crv{"x2"+(0,-.25) & "c2"},
  "x3";"c"+<-3pt,-2pt>="c1"**\crv{"x3"+(0,.25) & "c1"}, 
  "c"+(0,.10)*!<0pt,-6pt>{\si_{M,H}},
"x3"="L13",
"x4"="L14",
  "L15";"L15" +(0,-.5)="L15"**\dir{-},
"L11"="a1",
"L12"="a2",
"L13"="a3",
"a1" + (0,-1)="b1",
"a2" + (0,-1)="b2",
"a3" + (0,-1)="b3",
"a1";"a2"**\dir{}?(.5)="m12", 
  "m12" + (0,-.25)="c12", 
"a2";"a3"**\dir{}?(.5)="m23", 
  "m23" + (0,-.25)="c23", 
  "a1";"b2"**\crv{ "a1"+(0,-.25) & "c12" & "b2"+(0,.25)},  
  "a2";"b3"**\crv{ "a2"+(0,-.25) & "c23" & "b3"+(0,.25)}, 
                "b1";"c12"+<-0pt,-11pt>="c21"**\crv{"b1"+(0,.25) & "c21"},         
               "a3";"c23"+<-0pt,-2pt>="c32"**\crv{"a3"+(0,-.25) & "c32"},
"c12" +(.35,-.20);"c23" +(-.25,-.10)**\dir{-}, 
"c23"+(0,.10)*!<0pt,-6pt>{\si_{H^2,H}}, 
"b1"="L11",
"b2"="L12",
"b3"="L13",
  "L14";"L14" +(0,-1)="L14"**\dir{-},
  "L15";"L15" +(0,-1)="L15"**\dir{-},
"L11"="a1",
"L12"="a2",
"L13"="a3",
"L14"="a4",
"L15"="a5",
"a1" + (0,-1)="b1",
"a2" + (0,-1)="b2",
"a3" + (0,-1)="b3",
"a4" + (0,-1)="b4",
"a5" + (0,-1)="b5",
"a1";"a3"**\dir{}?(.5)="m13", 
  "m13" + (0,-.25)="c13", 
"a1";"a2"**\dir{}?(.5)="m12", 
  "m12" + (0,-.25)="c12", 
"a2";"a3"**\dir{}?(.5)="m23", 
  "m23" + (0,-.25)="c23", 
"a2";"a4"**\dir{}?(.5)="m24", 
  "m24" + (0,-.25)="c24", 
"a3";"a4"**\dir{}?(.5)="m34", 
  "m34" + (0,-.25)="c34",
"a3";"a5"**\dir{}?(.5)="m35", 
  "m35" + (0,-.25)="c35",  
  "a1";"b3"**\crv{ "a1"+(0,-.25) & "c13" & "b3"+(0,.25)},  
  "a2";"b4"**\crv{ "a2"+(0,-.25) & "c24" & "b4"+(0,.25)}, 
  "a3";"b5"**\crv{ "a3"+(0,-.25) & "c35" & "b5"+(0,.25)}, 
               "b1";"c13"+<-0pt,-11pt>="c31"**\crv{"b1"+(0,.25) & "c31"}, 
               "b2";"c23"+<-0pt,-17pt>="c32"**\crv{"b2"+(0,.25) & "c32"}, 
               "a4";"c34"+<3pt,2pt>="c43"**\crv{"a4"+(0,-.25) & "c43"},
               "a5";"c35"+<3pt,-2pt>="c53"**\crv{"a5"+(0,-.25) & "c53"},
"c34"+(0,.10)*!<0pt,-6pt>{\si_{H^3,M\ot H}}, 
"b1" = "L11",
"b2" = "L12",
"b3" = "L13",
"b4" = "L14",
"b5" = "L15",
					  "L11" + (-.5,0);"L15"+ (.5,0) **\dir{.},
"L11";"L11" +(0,-1)="L11"**\dir{-},
"L12";"L12" +(0,-1)="L12"**\dir{-},
  "L13"="x1", 
  "x1"+(0,-.5)="x2",
  "x1";"x2"**\dir{}?(.5)="c", 
  "x1";"x2"**\dir{-},   
  "c"-(0,0)*!<-4pt,0pt>{\bullet S_H },
"x2"="L13",
   "L14"="x1", 
  "L15"="x2",   
  "x1";"x2"**\dir{}?(.5)="m", 
  "m" + (0,-.25)="c", 
  "m" + (0,-.5)="x3", 
  "x1";"x3"**\crv{"x1"-(0,.25) & "c" & "x3"},
  "x2";"x3"**\crv{"x2"-(0,.25) & "c"  & "x3"}, 
  "c"+(0,.10)*!<0pt,-0pt>{m_H}, 
"x3"="L14",
   "L13"="x1", 
  "L14"="x2",   
  "x1";"x2"**\dir{}?(.5)="m", 
  "m" + (0,-.25)="c", 
  "m" + (0,-.5)="x3", 
  "x1";"x3"**\crv{"x1"-(0,.25) & "c" & "x3"},
  "x2";"x3"**\crv{"x2"-(0,.25) & "c"  & "x3"}, 
  "c"+(0,.10)*!<0pt,-0pt>{m_H}, 
"x3"="L13",
   "L12" + (1,.25)="x1", 
  "x1"+(0,-.25)="x2",
  "x1";"x2"**\dir{}?(.5)="c", 
  "x1";"x2"**\dir{-},   
  "c"-(0,0)*!<-3pt,0pt>{\bullet \eta },  
  "x1"+(0,0)*!<0pt,-2pt>{\circ},
"L13"="L14",
"x2"="L13",
   "L11"="x1", 
  "L12"="x2",   
  "x1";"x2"**\dir{}?(.5)="m", 
  "m" + (0,-.25)="c", 
  "m" + (0,-.5)="x3", 
  "x1";"x3"**\crv{"x1"-(0,.25) & "c" & "x3"},
  "x2";"x3"**\crv{"x2"-(0,.25) & "c"  & "x3"}, 
  "c"+(0,.10)*!<0pt,-0pt>{\phi_M}, 
"x3"="L11",
"L13";"L13" +(0,-.5)="L12"**\dir{-},
"L14";"L14" +(0,-.5)="L13"**\dir{-},
  "L11"-(0,.5)*!<0pt,-6pt>{M}, 
  "L12"-(0,.5)*!<0pt,-6pt>{H}, 
  "L13"-(0,.5)*!<0pt,-6pt>{H}, 
"L11"-(0,.5) ;"L13"-(0,.5) **\dir{-},
  "L12" -(0,.5)="x1", 
  "x1"+(0,-.75)="x2",
  "x1";"x2"**\dir{}?(.5)="c", 
  "x1";"x2"**\dir{=},  
  "c"-(0,0)*!<-22pt,0pt>{\bullet \text{coker}(f') },
"x2"-(0,.5)*!<0pt,-6pt>{M \ot_H H^2 },
\endxy  
\end{equation}

Stages (1) to (9) are explained below.  Notice that as it is shown in picture ~\eqref{pic1of the proofRHSpage13}, in stage 
(5) the braided AYD condition ~\eqref{BSAYD} is used.
\item (1) We just use  the definitions of $\varphi= (1, m)(1,S , 1))(\Delta , 1)$ and, 
 $\beta = [(1_M, \eta, 1) (\si_{H,M}) (\rho )(\phi)]$.
\item (2) We use $(1_M,\Delta ,1)(1_M, \eta, 1)=(1_M,\Delta \eta, 1)$ and  $\Delta \eta = \eta \ot \eta$ .
\item (3) We use $(1_M,1,S,1)(1_M, \eta,\eta,  1)=(1_M, \eta,S \eta,  1)$ and  $S \eta= \eta$ .
\item (4) We use $(1_M,1,m)(1_M, \eta,\eta,  1)=(1_M, \eta,m(\eta,  1))$ and  $m(\eta,  1)=1$  .
\item (5) We use the AYD formula ~\eqref{BSAYD}
\item (6) We use $(m(S,m) ,\phi)=(1, \phi) (m(S,m) ,1_M ,1)$.
\item (7) We commute $(\si_{H,M})$ and $(1, \phi)$ using:
\[
\begin{CD}
H \ot (M \ot H)         @>(1,\phi)>>     H \ot M\\
@V\si_{H,M \ot H} VV              @VV\si_{H,M} V\\
(M \ot H) \ot H    @>>(\phi ,1)>  M \ot H  
\end{CD}
\]
\item (8) We commute $(\si_{H,M \ot H}) $ and $(m(S,m) ,1_M ,1)$ using:
\[\begin{CD}
H^3 \ot (M \ot H) @>(m(S,m),1_M,1)>>  H \ot (M \ot H)\\
@V\si_{H^3,M \ot H}VV @VV\si_{H,M \ot H}V \\
(M \ot H) \ot H^3 @>>(1_M,1,m(S,m))> (M \ot H) \ot H
\end{CD} 
\]
\item (9) We commute $(1_M, \eta,1)$ and $( \phi,1)$.

Thus, LHS $=$ RHS in ~\eqref{cokrf'} if:
\begin{multline*}  
(\si_{H^3 ,M \ot H}) (\si_{H^2 ,H}, 1_M,1) (1,1,\si_{M,H},1)(1,1,1_M ,\si) (1,\si_{M ,H}, 1,1)\\
=(1_M, \si_{H^2, H^2})(\si_{H,M}, 1,1,1),
\end{multline*}
i.e., if in Diagrams $~\eqref{LHS2}$ and $~\eqref{RHS2}$, the parts between doted lines be equal.
This is true since:

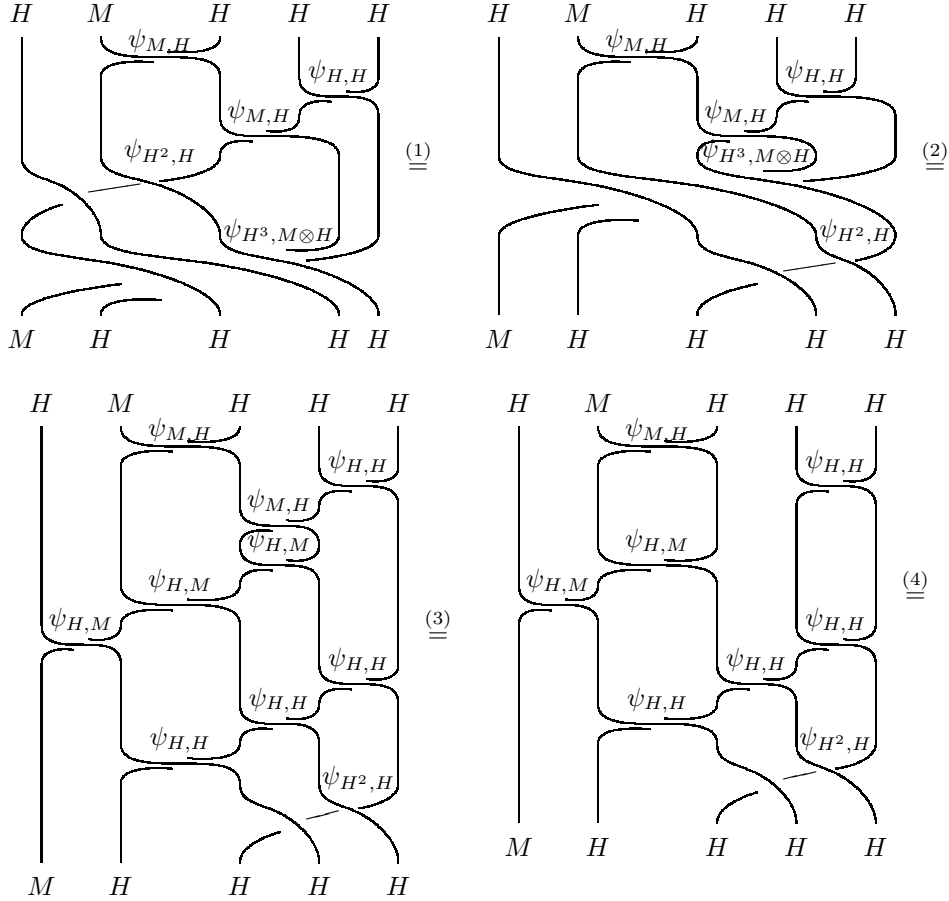
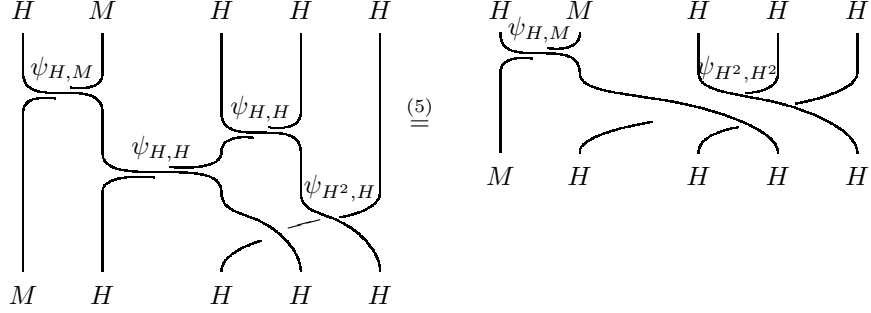
\[
\xy /r2.5pc/:,
(1,0)= "L11",
(2,0)= "L12",
(3.5,0)= "L13",
(4.5,0)= "L14",
(5.5,0)= "L15",
  "L11"+(0,.5)*!<0pt,6pt>{H}, 
  "L12"+(0,.5)*!<0pt,6pt>{M}, 
  "L13"+(0,.5)*!<0pt,6pt>{H}, 
  "L14"+(0,.5)*!<0pt,6pt>{H}, 
  "L15"+(0,.5)*!<0pt,6pt>{H},
  "L11";"L11" +(0,-.5)="L11"**\dir{-},
  "L12"="x1",       
  "L13"="x2",  , 
  "x1"+(0,-.5)="x3",
  "x2"+(0,-.5)="x4",
  "x1";"x2"**\dir{}?(.5)="m", 
  "m" + (0,-.25)="c", 
  "x1";"x4"**\crv{ "x1"+(0,-.25) & "c" & "x4"+(0,.25)},  
  "x2";"c"+<3pt,2pt>="c2"**\crv{"x2"+(0,-.25) & "c2"},
  "x3";"c"+<-3pt,-2pt>="c1"**\crv{"x3"+(0,.25) & "c1"}, 
  "c"+(0,.10)*!<0pt,-3pt>{\si_{M,H}},
"x3"="L12",
"x4"="L13",
  "L14";"L14" +(0,-.5)="L14"**\dir{-},
  "L15";"L15" +(0,-.5)="L15"**\dir{-},
  "L11";"L11" +(0,-.5)="L11"**\dir{-},
  "L12";"L12" +(0,-.5)="L12"**\dir{-},
  "L13";"L13" +(0,-.5)="L13"**\dir{-},
  "L14"="x1",       
  "L15"="x2",  , 
  "x1"+(0,-.5)="x3",
  "x2"+(0,-.5)="x4",
  "x1";"x2"**\dir{}?(.5)="m", 
  "m" + (0,-.25)="c", 
  "x1";"x4"**\crv{ "x1"+(0,-.25) & "c" & "x4"+(0,.25)},  
  "x2";"c"+<3pt,2pt>="c2"**\crv{"x2"+(0,-.25) & "c2"},
  "x3";"c"+<-3pt,-2pt>="c1"**\crv{"x3"+(0,.25) & "c1"}, 
  "c"+(0,.10)*!<0pt,-6pt>{\si_{H,H}},
"x3"="L14",
"x4"="L15",
"L11";"L11" +(0,-.5)="L11"**\dir{-},
"L12";"L12" +(0,-.5)="L12"**\dir{-},
  "L13"="x1",       
  "L14"="x2",  , 
  "x1"+(0,-.5)="x3",
  "x2"+(.5,-.5)="x4",
  "x1";"x2"**\dir{}?(.5)="m", 
  "m" + (0,-.25)="c", 
  "x1";"x4"**\crv{ "x1"+(0,-.25) & "c" & "x4"+(0,.25)},  
  "x2";"c"+<3pt,2pt>="c2"**\crv{"x2"+(0,-.25) & "c2"},
  "x3";"c"+<-3pt,-2pt>="c1"**\crv{"x3"+(0,.25) & "c1"}, 
  "c"+(0,.10)*!<0pt,-6pt>{\si_{M,H}},
"x3"="L13",
"x4"="L14",
  "L15";"L15" +(0,-.5)="L15"**\dir{-},
"L11"="a1",
"L12"="a2",
"L13"="a3",
"a1" + (0,-1)="b1",
"a2" + (0,-1)="b2",
"a3" + (0,-1)="b3",
"a1";"a2"**\dir{}?(.5)="m12", 
  "m12" + (0,-.25)="c12", 
"a2";"a3"**\dir{}?(.5)="m23", 
  "m23" + (0,-.25)="c23", 
  "a1";"b2"**\crv{ "a1"+(0,-.25) & "c12" & "b2"+(0,.25)},  
  "a2";"b3"**\crv{ "a2"+(0,-.25) & "c23" & "b3"+(0,.25)}, 
                "b1";"c12"+<-0pt,-11pt>="c21"**\crv{"b1"+(0,.25) & "c21"},         
               "a3";"c23"+<-0pt,-2pt>="c32"**\crv{"a3"+(0,-.25) & "c32"},
"c12" +(.35,-.20);"c23" +(-.25,-.10)**\dir{-}, %the between line !
"c23"+(0,.10)*!<0pt,-6pt>{\si_{H^2,H}}, 
"b1"="L11",
"b2"="L12",
"b3"="L13",
  "L14";"L14" +(0,-1)="L14"**\dir{-},
  "L15";"L15" +(0,-1)="L15"**\dir{-},
"L11"="a1",
"L12"="a2",
"L13"="a3",
"L14"="a4",
"L15"="a5",
"a1" + (0,-1)="b1",
"a2" + (0,-1)="b2",
"a3" + (0,-1)="b3",
"a4" + (0,-1)="b4",
"a5" + (0,-1)="b5",
"a1";"a3"**\dir{}?(.5)="m13", 
  "m13" + (0,-.25)="c13", 
"a1";"a2"**\dir{}?(.5)="m12", 
  "m12" + (0,-.25)="c12", 
"a2";"a3"**\dir{}?(.5)="m23", 
  "m23" + (0,-.25)="c23", 
"a2";"a4"**\dir{}?(.5)="m24", 
  "m24" + (0,-.25)="c24", 
"a3";"a4"**\dir{}?(.5)="m34", 
  "m34" + (0,-.25)="c34",
"a3";"a5"**\dir{}?(.5)="m35", 
  "m35" + (0,-.25)="c35",  
  "a1";"b3"**\crv{ "a1"+(0,-.25) & "c13" & "b3"+(0,.25)},  
  "a2";"b4"**\crv{ "a2"+(0,-.25) & "c24" & "b4"+(0,.25)}, 
  "a3";"b5"**\crv{ "a3"+(0,-.25) & "c35" & "b5"+(0,.25)}, 
               "b1";"c13"+<-0pt,-11pt>="c31"**\crv{"b1"+(0,.25) & "c31"}, 
               "b2";"c23"+<-0pt,-17pt>="c32"**\crv{"b2"+(0,.25) & "c32"}, 
               "a4";"c34"+<3pt,2pt>="c43"**\crv{"a4"+(0,-.25) & "c43"},
               "a5";"c35"+<3pt,-2pt>="c53"**\crv{"a5"+(0,-.25) & "c53"},
"c34"+(0,.10)*!<0pt,-6pt>{\si_{H^3,M\ot H}}, 
"b1" = "L11",
"b2" = "L12",
"b3" = "L13",
"b4" = "L14",
"b5" = "L15",
  "L11"-(0,.5)*!<0pt,-6pt>{M}, 
  "L12"-(0,.5)*!<0pt,-6pt>{H}, 
  "L13"-(0,.5)*!<0pt,-6pt>{H}, 
  "L14"-(0,.5)*!<0pt,-6pt>{H}, 
  "L15"-(0,.5)*!<0pt,-6pt>{H}, 
  "L15"+(.5,2)*!<0pt,0pt>{\overset{(1)}{=}}, 
\endxy 
%%%%%%%%%%%%%%%
\quad \quad
\xy /r2.5pc/:,
(1,0)= "L11",
(2,0)= "L12",
(3.5,0)= "L13",
(4.5,0)= "L14",
(5.5,0)= "L15",
  "L11"+(0,.5)*!<0pt,6pt>{H}, 
  "L12"+(0,.5)*!<0pt,6pt>{M}, 
  "L13"+(0,.5)*!<0pt,6pt>{H}, 
  "L14"+(0,.5)*!<0pt,6pt>{H}, 
  "L15"+(0,.5)*!<0pt,6pt>{H},
  "L11";"L11" +(0,-.5)="L11"**\dir{-},
  "L12"="x1",       
  "L13"="x2",  , 
  "x1"+(0,-.5)="x3",
  "x2"+(0,-.5)="x4",
  "x1";"x2"**\dir{}?(.5)="m", 
  "m" + (0,-.25)="c", 
  "x1";"x4"**\crv{ "x1"+(0,-.25) & "c" & "x4"+(0,.25)},  
  "x2";"c"+<3pt,2pt>="c2"**\crv{"x2"+(0,-.25) & "c2"},
  "x3";"c"+<-3pt,-2pt>="c1"**\crv{"x3"+(0,.25) & "c1"}, 
  "c"+(0,.10)*!<0pt,-3pt>{\si_{M,H}},
"x3"="L12",
"x4"="L13",
  "L14";"L14" +(0,-.5)="L14"**\dir{-},
  "L15";"L15" +(0,-.5)="L15"**\dir{-},
  "L11";"L11" +(0,-.5)="L11"**\dir{-},
  "L12";"L12" +(0,-.5)="L12"**\dir{-},
  "L13";"L13" +(0,-.5)="L13"**\dir{-},
  "L14"="x1",       
  "L15"="x2",  , 
  "x1"+(0,-.5)="x3",
  "x2"+(.5,-.5)="x4",
  "x1";"x2"**\dir{}?(.5)="m", 
  "m" + (0,-.25)="c", 
  "x1";"x4"**\crv{ "x1"+(0,-.25) & "c" & "x4"+(0,.25)},  
  "x2";"c"+<3pt,2pt>="c2"**\crv{"x2"+(0,-.25) & "c2"},
  "x3";"c"+<-3pt,-2pt>="c1"**\crv{"x3"+(0,.25) & "c1"}, 
  "c"+(0,.10)*!<0pt,-6pt>{\si_{H,H}},
"x3"="L14",
"x4"="L15",
"L11";"L11" +(0,-.5)="L11"**\dir{-},
"L12";"L12" +(0,-.5)="L12"**\dir{-},
  "L13"="x1",       
  "L14"="x2",  , 
  "x1"+(0,-.5)="x3",
  "x2"+(.5,-.5)="x4",
  "x1";"x2"**\dir{}?(.5)="m", 
  "m" + (0,-.25)="c", 
  "x1";"x4"**\crv{ "x1"+(0,-.25) & "c" & "x4"+(0,.25)},  
  "x2";"c"+<3pt,2pt>="c2"**\crv{"x2"+(0,-.25) & "c2"},
  "x3";"c"+<-3pt,-2pt>="c1"**\crv{"x3"+(0,.25) & "c1"}, 
  "c"+(0,.10)*!<0pt,-6pt>{\si_{M,H}},
"x3"="L13",
"x4"="L14",
  "L15";"L15" +(0,-.5)="L15"**\dir{-},
"L11"="a1",
"L12"="a2",
"L13"="a3",
"L14"="a4",
"L15"="a5",
"a1" + (0,-1)="b1",
"a2" + (0,-1)="b2",
"a3" + (0,-1)="b3",
"a4" + (0,-1)="b4",
"a5" + (0,-1)="b5",
"a1";"a3"**\dir{}?(.5)="m13", 
  "m13" + (0,-.25)="c13", 
"a1";"a2"**\dir{}?(.5)="m12", 
  "m12" + (0,-.25)="c12", 
"a2";"a3"**\dir{}?(.5)="m23", 
  "m23" + (0,-.25)="c23", 
"a2";"a4"**\dir{}?(.5)="m24", 
  "m24" + (0,-.25)="c24", 
"a3";"a4"**\dir{}?(.5)="m34", 
  "m34" + (0,-.25)="c34",
"a3";"a5"**\dir{}?(.5)="m35", 
  "m35" + (0,-.25)="c35",  
  "a1";"b3"**\crv{ "a1"+(0,-.25) & "c13" & "b3"+(0,.25)},  
  "a2";"b4"**\crv{ "a2"+(0,-.25) & "c24" & "b4"+(0,.25)}, 
  "a3";"b5"**\crv{ "a3"+(0,-.25) & "c35" & "b5"+(0,.25)}, 
               "b1";"c13"+<-0pt,-11pt>="c31"**\crv{"b1"+(0,.25) & "c31"}, 
               "b2";"c23"+<-0pt,-17pt>="c32"**\crv{"b2"+(0,.25) & "c32"}, 
               "a4";"c34"+<3pt,2pt>="c43"**\crv{"a4"+(0,-.25) & "c43"},
               "a5";"c35"+<3pt,-2pt>="c53"**\crv{"a5"+(0,-.25) & "c53"},
"c34"+(0,.10)*!<0pt,-6pt>{\si_{H^3,M\ot H}}, 
"b1" = "L11",
"b2" = "L12",
"b3" = "L13",
"b4" = "L14",
"b5" = "L15",
  "L11";"L11" +(0,-1)="L11"**\dir{-},
  "L12";"L12" +(0,-1)="L12"**\dir{-},
"L13"="a1",
"L14"="a2",
"L15"="a3",
"a1" + (0,-1)="b1",
"a2" + (0,-1)="b2",
"a3" + (0,-1)="b3",
"a1";"a2"**\dir{}?(.5)="m12", 
  "m12" + (0,-.25)="c12", 
"a2";"a3"**\dir{}?(.5)="m23", 
  "m23" + (0,-.25)="c23", 
  "a1";"b2"**\crv{ "a1"+(0,-.25) & "c12" & "b2"+(0,.25)},  
  "a2";"b3"**\crv{ "a2"+(0,-.25) & "c23" & "b3"+(0,.25)}, 
                "b1";"c12"+<-0pt,-11pt>="c21"**\crv{"b1"+(0,.25) & "c21"},         
               "a3";"c23"+<-0pt,-2pt>="c32"**\crv{"a3"+(0,-.25) & "c32"},
"c12" +(.35,-.20);"c23" +(-.25,-.10)**\dir{-}, %the between line !
"c23"+(0,.10)*!<0pt,-6pt>{\si_{H^2,H}}, 
"b1"="L13",
"b2"="L14",
"b3"="L15",
  "L11"-(0,.5)*!<0pt,-6pt>{M}, 
  "L12"-(0,.5)*!<0pt,-6pt>{H}, 
  "L13"-(0,.5)*!<0pt,-6pt>{H}, 
  "L14"-(0,.5)*!<0pt,-6pt>{H}, 
  "L15"-(0,.5)*!<0pt,-6pt>{H}, 
  "L15"+(.5,2)*!<0pt,0pt>{\overset{(2)}{=}}, 
\endxy
\]
%%%%%%%%%%%%%%%
\[
\xy /r2.5pc/:,
(1,0)= "L11",
(2,0)= "L12",
(3.5,0)= "L13",
(4.5,0)= "L14",
(5.5,0)= "L15",
  "L11"+(0,.5)*!<0pt,6pt>{H}, 
  "L12"+(0,.5)*!<0pt,6pt>{M}, 
  "L13"+(0,.5)*!<0pt,6pt>{H}, 
  "L14"+(0,.5)*!<0pt,6pt>{H}, 
  "L15"+(0,.5)*!<0pt,6pt>{H},
  "L11";"L11" +(0,-.5)="L11"**\dir{-},
  "L12"="x1",       
  "L13"="x2",  , 
  "x1"+(0,-.5)="x3",
  "x2"+(0,-.5)="x4",
  "x1";"x2"**\dir{}?(.5)="m", 
  "m" + (0,-.25)="c", 
  "x1";"x4"**\crv{ "x1"+(0,-.25) & "c" & "x4"+(0,.25)},  
  "x2";"c"+<3pt,2pt>="c2"**\crv{"x2"+(0,-.25) & "c2"},
  "x3";"c"+<-3pt,-2pt>="c1"**\crv{"x3"+(0,.25) & "c1"}, 
  "c"+(0,.10)*!<0pt,-3pt>{\si_{M,H}},
"x3"="L12",
"x4"="L13",
  "L14";"L14" +(0,-.5)="L14"**\dir{-},
  "L15";"L15" +(0,-.5)="L15"**\dir{-},
  "L11";"L11" +(0,-.5)="L11"**\dir{-},
  "L12";"L12" +(0,-.5)="L12"**\dir{-},
  "L13";"L13" +(0,-.5)="L13"**\dir{-},
  "L14"="x1",       
  "L15"="x2",  , 
  "x1"+(0,-.5)="x3",
  "x2"+(0,-.5)="x4",
  "x1";"x2"**\dir{}?(.5)="m", 
  "m" + (0,-.25)="c", 
  "x1";"x4"**\crv{ "x1"+(0,-.25) & "c" & "x4"+(0,.25)},  
  "x2";"c"+<3pt,2pt>="c2"**\crv{"x2"+(0,-.25) & "c2"},
  "x3";"c"+<-3pt,-2pt>="c1"**\crv{"x3"+(0,.25) & "c1"}, 
  "c"+(0,.10)*!<0pt,-6pt>{\si_{H,H}},
"x3"="L14",
"x4"="L15",
"L11";"L11" +(0,-.5)="L11"**\dir{-},
"L12";"L12" +(0,-.5)="L12"**\dir{-},
  "L13"="x1",       
  "L14"="x2",  , 
  "x1"+(0,-.5)="x3",
  "x2"+(0,-.5)="x4",
  "x1";"x2"**\dir{}?(.5)="m", 
  "m" + (0,-.25)="c", 
  "x1";"x4"**\crv{ "x1"+(0,-.25) & "c" & "x4"+(0,.25)},  
  "x2";"c"+<3pt,2pt>="c2"**\crv{"x2"+(0,-.25) & "c2"},
  "x3";"c"+<-3pt,-2pt>="c1"**\crv{"x3"+(0,.25) & "c1"}, 
  "c"+(0,.10)*!<0pt,-6pt>{\si_{M,H}},
"x3"="L13",
"x4"="L14",
  "L15";"L15" +(0,-.5)="L15"**\dir{-},
"L11";"L11" +(0,-.5)="L11"**\dir{-},
"L12";"L12" +(0,-.5)="L12"**\dir{-},
  "L13"="x1",       
  "L14"="x2",  , 
  "x1"+(0,-.5)="x3",
  "x2"+(0,-.5)="x4",
  "x1";"x2"**\dir{}?(.5)="m", 
  "m" + (0,-.25)="c", 
  "x1";"x4"**\crv{ "x1"+(0,-.25) & "c" & "x4"+(0,.25)},  
  "x2";"c"+<3pt,2pt>="c2"**\crv{"x2"+(0,-.25) & "c2"},
  "x3";"c"+<-3pt,-2pt>="c1"**\crv{"x3"+(0,.25) & "c1"}, 
  "c"+(0,.10)*!<0pt,-6pt>{\si_{H,M}},
"x3"="L13",
"x4"="L14",
  "L15";"L15" +(0,-.5)="L15"**\dir{-},
"L11";"L11" +(0,-.5)="L11"**\dir{-},
  "L12"="x1",       
  "L13"="x2",  , 
  "x1"+(0,-.5)="x3",
  "x2"+(0,-.5)="x4",
  "x1";"x2"**\dir{}?(.5)="m", 
  "m" + (0,-.25)="c", 
  "x1";"x4"**\crv{ "x1"+(0,-.25) & "c" & "x4"+(0,.25)},  
  "x2";"c"+<3pt,2pt>="c2"**\crv{"x2"+(0,-.25) & "c2"},
  "x3";"c"+<-3pt,-2pt>="c1"**\crv{"x3"+(0,.25) & "c1"}, 
  "c"+(0,.10)*!<0pt,-6pt>{\si_{H,M}},
"x3"="L12",
"x4"="L13",
  "L14";"L14" +(0,-.5)="L14"**\dir{-},
  "L15";"L15" +(0,-.5)="L15"**\dir{-},
  "L11"="x1",       
  "L12"="x2",  , 
  "x1"+(0,-.5)="x3",
  "x2"+(0,-.5)="x4",
  "x1";"x2"**\dir{}?(.5)="m", 
  "m" + (0,-.25)="c", 
  "x1";"x4"**\crv{ "x1"+(0,-.25) & "c" & "x4"+(0,.25)},  
  "x2";"c"+<3pt,2pt>="c2"**\crv{"x2"+(0,-.25) & "c2"},
  "x3";"c"+<-3pt,-2pt>="c1"**\crv{"x3"+(0,.25) & "c1"}, 
  "c"+(0,.10)*!<0pt,-6pt>{\si_{H,M}},
"x3"="L11",
"x4"="L12",
  "L13";"L13" +(0,-.5)="L13"**\dir{-},
  "L14";"L14" +(0,-.5)="L14"**\dir{-},
  "L15";"L15" +(0,-.5)="L15"**\dir{-},
  "L11";"L11" +(0,-.5)="L11"**\dir{-},
  "L12";"L12" +(0,-.5)="L12"**\dir{-},
  "L13";"L13" +(0,-.5)="L13"**\dir{-},
  "L14"="x1",       
  "L15"="x2",  , 
  "x1"+(0,-.5)="x3",
  "x2"+(0,-.5)="x4",
  "x1";"x2"**\dir{}?(.5)="m", 
  "m" + (0,-.25)="c", 
  "x1";"x4"**\crv{ "x1"+(0,-.25) & "c" & "x4"+(0,.25)},  
  "x2";"c"+<3pt,2pt>="c2"**\crv{"x2"+(0,-.25) & "c2"},
  "x3";"c"+<-3pt,-2pt>="c1"**\crv{"x3"+(0,.25) & "c1"}, 
  "c"+(0,.10)*!<0pt,-6pt>{\si_{H,H}},
"x3"="L14",
"x4"="L15",
  "L11";"L11" +(0,-.5)="L11"**\dir{-},
  "L12";"L12" +(0,-.5)="L12"**\dir{-},
  "L13"="x1",       
  "L14"="x2",  , 
  "x1"+(0,-.5)="x3",
  "x2"+(0,-.5)="x4",
  "x1";"x2"**\dir{}?(.5)="m", 
  "m" + (0,-.25)="c", 
  "x1";"x4"**\crv{ "x1"+(0,-.25) & "c" & "x4"+(0,.25)},  
  "x2";"c"+<3pt,2pt>="c2"**\crv{"x2"+(0,-.25) & "c2"},
  "x3";"c"+<-3pt,-2pt>="c1"**\crv{"x3"+(0,.25) & "c1"}, 
  "c"+(0,.10)*!<0pt,-6pt>{\si_{H,H}},
"x3"="L13",
"x4"="L14",
  "L15";"L15" +(0,-.5)="L15"**\dir{-},
  "L11";"L11" +(0,-.5)="L11"**\dir{-},
  "L12"="x1",       
  "L13"="x2",  , 
  "x1"+(0,-.5)="x3",
  "x2"+(0,-.5)="x4",
  "x1";"x2"**\dir{}?(.5)="m", 
  "m" + (0,-.25)="c", 
  "x1";"x4"**\crv{ "x1"+(0,-.25) & "c" & "x4"+(0,.25)},  
  "x2";"c"+<3pt,2pt>="c2"**\crv{"x2"+(0,-.25) & "c2"},
  "x3";"c"+<-3pt,-2pt>="c1"**\crv{"x3"+(0,.25) & "c1"}, 
  "c"+(0,.10)*!<0pt,-6pt>{\si_{H,H}},
"x3"="L12",
"x4"="L13",
  "L14";"L14" +(0,-.5)="L14"**\dir{-},
  "L15";"L15" +(0,-.5)="L15"**\dir{-},
  "L11";"L11" +(0,-1)="L11"**\dir{-},
  "L12";"L12" +(0,-1)="L12"**\dir{-},
"L13"="a1",
"L14"="a2",
"L15"="a3",
"a1" + (0,-1)="b1",
"a2" + (0,-1)="b2",
"a3" + (0,-1)="b3",
"a1";"a2"**\dir{}?(.5)="m12", 
  "m12" + (0,-.25)="c12", 
"a2";"a3"**\dir{}?(.5)="m23", 
  "m23" + (0,-.25)="c23", 
  "a1";"b2"**\crv{ "a1"+(0,-.25) & "c12" & "b2"+(0,.25)},  
  "a2";"b3"**\crv{ "a2"+(0,-.25) & "c23" & "b3"+(0,.25)}, 
                "b1";"c12"+<-0pt,-11pt>="c21"**\crv{"b1"+(0,.25) & "c21"},         
               "a3";"c23"+<-0pt,-2pt>="c32"**\crv{"a3"+(0,-.25) & "c32"},
"c12" +(.35,-.20);"c23" +(-.25,-.10)**\dir{-}, %the between line !
"c23"+(0,.10)*!<0pt,-6pt>{\si_{H^2,H}}, 
"b1"="L13",
"b2"="L14",
"b3"="L15",
  "L11"-(0,.5)*!<0pt,-6pt>{M}, 
  "L12"-(0,.5)*!<0pt,-6pt>{H}, 
  "L13"-(0,.5)*!<0pt,-6pt>{H}, 
  "L14"-(0,.5)*!<0pt,-6pt>{H}, 
  "L15"-(0,.5)*!<0pt,-6pt>{H}, 
  "L15"+(.5,3)*!<0pt,0pt>{\overset{(3)}{=}}, 
\endxy
%%%%%%%%%%%%%%%
\quad \quad
\xy /r2.5pc/:,
(1,0)= "L11",
(2,0)= "L12",
(3.5,0)= "L13",
(4.5,0)= "L14",
(5.5,0)= "L15",
  "L11"+(0,.5)*!<0pt,6pt>{H}, 
  "L12"+(0,.5)*!<0pt,6pt>{M}, 
  "L13"+(0,.5)*!<0pt,6pt>{H}, 
  "L14"+(0,.5)*!<0pt,6pt>{H}, 
  "L15"+(0,.5)*!<0pt,6pt>{H},
  "L11";"L11" +(0,-.5)="L11"**\dir{-},
  "L12"="x1",       
  "L13"="x2",  , 
  "x1"+(0,-.5)="x3",
  "x2"+(0,-.5)="x4",
  "x1";"x2"**\dir{}?(.5)="m", 
  "m" + (0,-.25)="c", 
  "x1";"x4"**\crv{ "x1"+(0,-.25) & "c" & "x4"+(0,.25)},  
  "x2";"c"+<3pt,2pt>="c2"**\crv{"x2"+(0,-.25) & "c2"},
  "x3";"c"+<-3pt,-2pt>="c1"**\crv{"x3"+(0,.25) & "c1"}, 
  "c"+(0,.10)*!<0pt,-3pt>{\si_{M,H}},
"x3"="L12",
"x4"="L13",
  "L14";"L14" +(0,-.5)="L14"**\dir{-},
  "L15";"L15" +(0,-.5)="L15"**\dir{-},
  "L11";"L11" +(0,-.5)="L11"**\dir{-},
  "L12";"L12" +(0,-.5)="L12"**\dir{-},
  "L13";"L13" +(0,-.5)="L13"**\dir{-},
  "L14"="x1",       
  "L15"="x2",  , 
  "x1"+(0,-.5)="x3",
  "x2"+(0,-.5)="x4",
  "x1";"x2"**\dir{}?(.5)="m", 
  "m" + (0,-.25)="c", 
  "x1";"x4"**\crv{ "x1"+(0,-.25) & "c" & "x4"+(0,.25)},  
  "x2";"c"+<3pt,2pt>="c2"**\crv{"x2"+(0,-.25) & "c2"},
  "x3";"c"+<-3pt,-2pt>="c1"**\crv{"x3"+(0,.25) & "c1"}, 
  "c"+(0,.10)*!<0pt,-6pt>{\si_{H,H}},
"x3"="L14",
"x4"="L15",
"L11";"L11" +(0,-.5)="L11"**\dir{-},
"L12";"L12" +(0,-.5)="L12"**\dir{-},
"L13";"L13" +(0,-.5)="L13"**\dir{-},
"L14";"L14" +(0,-.5)="L14"**\dir{-},
"L15";"L15" +(0,-.5)="L15"**\dir{-},
"L11";"L11" +(0,-.5)="L11"**\dir{-},
  "L12"="x1",       
  "L13"="x2",  , 
  "x1"+(0,-.5)="x3",
  "x2"+(0,-.5)="x4",
  "x1";"x2"**\dir{}?(.5)="m", 
  "m" + (0,-.25)="c", 
  "x1";"x4"**\crv{ "x1"+(0,-.25) & "c" & "x4"+(0,.25)},  
  "x2";"c"+<3pt,2pt>="c2"**\crv{"x2"+(0,-.25) & "c2"},
  "x3";"c"+<-3pt,-2pt>="c1"**\crv{"x3"+(0,.25) & "c1"}, 
  "c"+(0,.10)*!<0pt,-6pt>{\si_{H,M}},
"x3"="L12",
"x4"="L13",
  "L14";"L14" +(0,-.5)="L14"**\dir{-},
  "L15";"L15" +(0,-.5)="L15"**\dir{-},
  "L11"="x1",       
  "L12"="x2",  , 
  "x1"+(0,-.5)="x3",
  "x2"+(0,-.5)="x4",
  "x1";"x2"**\dir{}?(.5)="m", 
  "m" + (0,-.25)="c", 
  "x1";"x4"**\crv{ "x1"+(0,-.25) & "c" & "x4"+(0,.25)},  
  "x2";"c"+<3pt,2pt>="c2"**\crv{"x2"+(0,-.25) & "c2"},
  "x3";"c"+<-3pt,-2pt>="c1"**\crv{"x3"+(0,.25) & "c1"}, 
  "c"+(0,.10)*!<0pt,-6pt>{\si_{H,M}},
"x3"="L11",
"x4"="L12",
  "L13";"L13" +(0,-.5)="L13"**\dir{-},
  "L14";"L14" +(0,-.5)="L14"**\dir{-},
  "L15";"L15" +(0,-.5)="L15"**\dir{-},
  "L11";"L11" +(0,-.5)="L11"**\dir{-},
  "L12";"L12" +(0,-.5)="L12"**\dir{-},
  "L13";"L13" +(0,-.5)="L13"**\dir{-},
  "L14"="x1",       
  "L15"="x2",  , 
  "x1"+(0,-.5)="x3",
  "x2"+(0,-.5)="x4",
  "x1";"x2"**\dir{}?(.5)="m", 
  "m" + (0,-.25)="c", 
  "x1";"x4"**\crv{ "x1"+(0,-.25) & "c" & "x4"+(0,.25)},  
  "x2";"c"+<3pt,2pt>="c2"**\crv{"x2"+(0,-.25) & "c2"},
  "x3";"c"+<-3pt,-2pt>="c1"**\crv{"x3"+(0,.25) & "c1"}, 
  "c"+(0,.10)*!<0pt,-6pt>{\si_{H,H}},
"x3"="L14",
"x4"="L15",
  "L11";"L11" +(0,-.5)="L11"**\dir{-},
  "L12";"L12" +(0,-.5)="L12"**\dir{-},
  "L13"="x1",       
  "L14"="x2",  , 
  "x1"+(0,-.5)="x3",
  "x2"+(0,-.5)="x4",
  "x1";"x2"**\dir{}?(.5)="m", 
  "m" + (0,-.25)="c", 
  "x1";"x4"**\crv{ "x1"+(0,-.25) & "c" & "x4"+(0,.25)},  
  "x2";"c"+<3pt,2pt>="c2"**\crv{"x2"+(0,-.25) & "c2"},
  "x3";"c"+<-3pt,-2pt>="c1"**\crv{"x3"+(0,.25) & "c1"}, 
  "c"+(0,.10)*!<0pt,-6pt>{\si_{H,H}},
"x3"="L13",
"x4"="L14",
  "L15";"L15" +(0,-.5)="L15"**\dir{-},
  "L11";"L11" +(0,-.5)="L11"**\dir{-},
  "L12"="x1",       
  "L13"="x2",  , 
  "x1"+(0,-.5)="x3",
  "x2"+(0,-.5)="x4",
  "x1";"x2"**\dir{}?(.5)="m", 
  "m" + (0,-.25)="c", 
  "x1";"x4"**\crv{ "x1"+(0,-.25) & "c" & "x4"+(0,.25)},  
  "x2";"c"+<3pt,2pt>="c2"**\crv{"x2"+(0,-.25) & "c2"},
  "x3";"c"+<-3pt,-2pt>="c1"**\crv{"x3"+(0,.25) & "c1"}, 
  "c"+(0,.10)*!<0pt,-6pt>{\si_{H,H}},
"x3"="L12",
"x4"="L13",
  "L14";"L14" +(0,-.5)="L14"**\dir{-},
  "L15";"L15" +(0,-.5)="L15"**\dir{-},
  "L11";"L11" +(0,-1)="L11"**\dir{-},
  "L12";"L12" +(0,-1)="L12"**\dir{-},
"L13"="a1",
"L14"="a2",
"L15"="a3",
"a1" + (0,-1)="b1",
"a2" + (0,-1)="b2",
"a3" + (0,-1)="b3",
"a1";"a2"**\dir{}?(.5)="m12", 
  "m12" + (0,-.25)="c12", 
"a2";"a3"**\dir{}?(.5)="m23", 
  "m23" + (0,-.25)="c23", 
  "a1";"b2"**\crv{ "a1"+(0,-.25) & "c12" & "b2"+(0,.25)},  
  "a2";"b3"**\crv{ "a2"+(0,-.25) & "c23" & "b3"+(0,.25)}, 
                "b1";"c12"+<-0pt,-11pt>="c21"**\crv{"b1"+(0,.25) & "c21"},         
               "a3";"c23"+<-0pt,-2pt>="c32"**\crv{"a3"+(0,-.25) & "c32"},
"c12" +(.35,-.20);"c23" +(-.25,-.10)**\dir{-}, %the between line !
"c23"+(0,.10)*!<0pt,-6pt>{\si_{H^2,H}}, 
"b1"="L13",
"b2"="L14",
"b3"="L15",
  "L11"-(0,.5)*!<0pt,-6pt>{M}, 
  "L12"-(0,.5)*!<0pt,-6pt>{H}, 
  "L13"-(0,.5)*!<0pt,-6pt>{H}, 
  "L14"-(0,.5)*!<0pt,-6pt>{H}, 
  "L15"-(0,.5)*!<0pt,-6pt>{H}, 
  "L15"+(.5,3)*!<0pt,0pt>{\overset{(4)}{=}}, 
\endxy \]
%%%%%%%%%%%%%%%
\[
\xy /r2.5pc/:,
(1,0)= "L11",
(2,0)= "L12",
(3.5,0)= "L13",
(4.5,0)= "L14",
(5.5,0)= "L15",
  "L11"+(0,.5)*!<0pt,6pt>{H}, 
  "L12"+(0,.5)*!<0pt,6pt>{M}, 
  "L13"+(0,.5)*!<0pt,6pt>{H}, 
  "L14"+(0,.5)*!<0pt,6pt>{H}, 
  "L15"+(0,.5)*!<0pt,6pt>{H},
"L11";"L11" +(0,-.5)="L11"**\dir{-},
"L12";"L12" +(0,-.5)="L12"**\dir{-},
"L13";"L13" +(0,-.5)="L13"**\dir{-},
"L14";"L14" +(0,-.5)="L14"**\dir{-},
"L15";"L15" +(0,-.5)="L15"**\dir{-},
  "L11"="x1",       
  "L12"="x2",  , 
  "x1"+(0,-.5)="x3",
  "x2"+(0,-.5)="x4",
  "x1";"x2"**\dir{}?(.5)="m", 
  "m" + (0,-.25)="c", 
  "x1";"x4"**\crv{ "x1"+(0,-.25) & "c" & "x4"+(0,.25)},  
  "x2";"c"+<3pt,2pt>="c2"**\crv{"x2"+(0,-.25) & "c2"},
  "x3";"c"+<-3pt,-2pt>="c1"**\crv{"x3"+(0,.25) & "c1"}, 
  "c"+(0,.10)*!<0pt,-6pt>{\si_{H,M}},
"x3"="L11",
"x4"="L12",
  "L13";"L13" +(0,-.5)="L13"**\dir{-},
  "L14";"L14" +(0,-.5)="L14"**\dir{-},
  "L15";"L15" +(0,-.5)="L15"**\dir{-},
  "L11";"L11" +(0,-.5)="L11"**\dir{-},
  "L12";"L12" +(0,-.5)="L12"**\dir{-},
  "L13"="x1",       
  "L14"="x2",  , 
  "x1"+(0,-.5)="x3",
  "x2"+(0,-.5)="x4",
  "x1";"x2"**\dir{}?(.5)="m", 
  "m" + (0,-.25)="c", 
  "x1";"x4"**\crv{ "x1"+(0,-.25) & "c" & "x4"+(0,.25)},  
  "x2";"c"+<3pt,2pt>="c2"**\crv{"x2"+(0,-.25) & "c2"},
  "x3";"c"+<-3pt,-2pt>="c1"**\crv{"x3"+(0,.25) & "c1"}, 
  "c"+(0,.10)*!<0pt,-6pt>{\si_{H,H}},
"x3"="L13",
"x4"="L14",
  "L15";"L15" +(0,-.5)="L15"**\dir{-},
  "L11";"L11" +(0,-.5)="L11"**\dir{-},
  "L12"="x1",       
  "L13"="x2",  , 
  "x1"+(0,-.5)="x3",
  "x2"+(0,-.5)="x4",
  "x1";"x2"**\dir{}?(.5)="m", 
  "m" + (0,-.25)="c", 
  "x1";"x4"**\crv{ "x1"+(0,-.25) & "c" & "x4"+(0,.25)},  
  "x2";"c"+<3pt,2pt>="c2"**\crv{"x2"+(0,-.25) & "c2"},
  "x3";"c"+<-3pt,-2pt>="c1"**\crv{"x3"+(0,.25) & "c1"}, 
  "c"+(0,.10)*!<0pt,-6pt>{\si_{H,H}},
"x3"="L12",
"x4"="L13",
  "L14";"L14" +(0,-.5)="L14"**\dir{-},
  "L15";"L15" +(0,-.5)="L15"**\dir{-},
  "L11";"L11" +(0,-1)="L11"**\dir{-},
  "L12";"L12" +(0,-1)="L12"**\dir{-},
"L13"="a1",
"L14"="a2",
"L15"="a3",
"a1" + (0,-1)="b1",
"a2" + (0,-1)="b2",
"a3" + (0,-1)="b3",
"a1";"a2"**\dir{}?(.5)="m12", 
  "m12" + (0,-.25)="c12", 
"a2";"a3"**\dir{}?(.5)="m23", 
  "m23" + (0,-.25)="c23", 
  "a1";"b2"**\crv{ "a1"+(0,-.25) & "c12" & "b2"+(0,.25)},  
  "a2";"b3"**\crv{ "a2"+(0,-.25) & "c23" & "b3"+(0,.25)}, 
                "b1";"c12"+<-0pt,-11pt>="c21"**\crv{"b1"+(0,.25) & "c21"},         
               "a3";"c23"+<-0pt,-2pt>="c32"**\crv{"a3"+(0,-.25) & "c32"},
"c12" +(.35,-.20);"c23" +(-.25,-.10)**\dir{-}, 
"c23"+(0,.10)*!<0pt,-6pt>{\si_{H^2,H}}, 
"b1"="L13",
"b2"="L14",
"b3"="L15",
  "L11"-(0,.5)*!<0pt,-6pt>{M}, 
  "L12"-(0,.5)*!<0pt,-6pt>{H}, 
  "L13"-(0,.5)*!<0pt,-6pt>{H}, 
  "L14"-(0,.5)*!<0pt,-6pt>{H}, 
  "L15"-(0,.5)*!<0pt,-6pt>{H}, 
  "L15"+(.5,2)*!<0pt,0pt>{\overset{(5)}{=}}, 
\endxy
%%%%%%%%%%%%%%%
\quad \quad
\xy /r2.5pc/:,
(1,0)= "L11",
(2,0)= "L12",
(3.5,0)= "L13",
(4.5,0)= "L14",
(5.5,0)= "L15",
  "L11"+(0,.5)*!<0pt,6pt>{H}, 
  "L12"+(0,.5)*!<0pt,6pt>{M}, 
  "L13"+(0,.5)*!<0pt,6pt>{H}, 
  "L14"+(0,.5)*!<0pt,6pt>{H}, 
  "L15"+(0,.5)*!<0pt,6pt>{H},
  "L11"="x1",       
  "L12"="x2",  , 
  "x1"+(0,-.5)="x3",
  "x2"+(0,-.5)="x4",
  "x1";"x2"**\dir{}?(.5)="m", 
  "m" + (0,-.25)="c", 
  "x1";"x4"**\crv{ "x1"+(0,-.25) & "c" & "x4"+(0,.25)},  
  "x2";"c"+<3pt,2pt>="c2"**\crv{"x2"+(0,-.25) & "c2"},
  "x3";"c"+<-3pt,-2pt>="c1"**\crv{"x3"+(0,.25) & "c1"}, 
  "c"+(0,.10)*!<0pt,-6pt>{\si_{H,M}},
"x3"="L11",
"x4"="L12",
  "L13";"L13" +(0,-.5)="L13"**\dir{-},
  "L14";"L14" +(0,-.5)="L14"**\dir{-},
  "L15";"L15" +(0,-.5)="L15"**\dir{-},
  "L11";"L11" +(0,-1)="L11"**\dir{-},
"L12"="a1",
"L13"="a2",
"L14"="a3",
"L15"="a4",
"a1" + (0,-1)="b1",
"a2" + (0,-1)="b2",
"a3" + (0,-1)="b3",
"a4" + (0,-1)="b4",
"a1";"a3"**\dir{}?(.5)="m13", 
  "m13" + (0,-.25)="c13", 
"a1";"a2"**\dir{}?(.5)="m12", 
  "m12" + (0,-.25)="c12", 
"a2";"a3"**\dir{}?(.5)="m23", 
  "m23" + (0,-.25)="c23", 
"a2";"a4"**\dir{}?(.5)="m24", 
  "m24" + (0,-.25)="c24", 
  "a1";"b3"**\crv{ "a1"+(0,-.25) & "c13" & "b3"+(0,.25)},  
  "a2";"b4"**\crv{ "a2"+(0,-.25) & "c24" & "b4"+(0,.25)}, 
                "b1";"c12"+<5pt,-11pt>="c21"**\crv{"b1"+(0,.25) & "c21"},  
                "b2";"c23"+<0pt,-13pt>="c32"**\crv{"b2"+(0,.25) & "c32"},       
               "a3";"c23"+<3pt,0pt>="c32"**\crv{"a3"+(0,-.25) & "c32"},
               "a4";"c24"+<7pt,-4pt>="c42"**\crv{"a4"+(0,-.25) & "c42"},
"c23"+(0,.10)*!<0pt,-6pt>{\si_{H^2,H^2}}, 
"b1"="L12",
"b2"="L13",
"b3"="L14",
"b4"="L15",
  "L11"-(0,.5)*!<0pt,-6pt>{M}, 
  "L12"-(0,.5)*!<0pt,-6pt>{H}, 
  "L13"-(0,.5)*!<0pt,-6pt>{H}, 
  "L14"-(0,.5)*!<0pt,-6pt>{H}, 
  "L15"-(0,.5)*!<0pt,-6pt>{H}, 
\endxy
\]

Stages (1) to (5) are explained below. Notice that in stages $(3), (4)$ the symmetry condition $\si^2 = id$ 
is used.
\item (1) We commute $ (\si_{H^3 ,M \ot H})$ and $(\si_{H^2 ,H}, 1_M,1)$ using naturality of $\si$:

\[\begin{CD}
H^3 \ot (M \ot H) @>(\si_{H^2 ,H}, 1_M,1)>>  H^3 \ot (M \ot H)\\
@V\si_{H^3,M \ot H}VV @VV\si_{H^3,M \ot H}V \\
(M \ot H) \ot H^3 @>>(1_M,1, \si_{H^2 ,H})> (M \ot H) \ot H^3
\end{CD} \nonumber \]
\item (2) We  use 
\begin{multline*} 
(\si_{H^3 ,M \ot H}) =  (1_M, \si_{H,H},1,1)(1_M,1,\si_{H,H},1)(1_M,1,1,\si_{H,H})\\
(\si_{H,M},1,1,1)(1, \si_{H,M},1,1)(1,1, \si_{H,M},1)
\end{multline*}
\item (3) We use the symmetric property of $\si$ to put $ \si_{H,M} \si_{M,H} = 1$.
\item (4) We use the symmetric property of $\si$ to put $ \si_{H,M} \si_{M ,H} = 1 $ 
and $\si_{H,H} \si_{H,H}= 1 $.
\item (5) We  use $  \si_{H^2, H^2} =  (1,\si_{H^2 ,H})(\si_{H,H},1,1)(1,\si_{H,H},1)$.

As for the cyclic condition, considering the following diagram:
\[\begin{CD}
 M \ot H \ot C^{n+1} @>f_{M,C^{n+1}} >> M \ot C^{n+1} @>\text{coker}(f_{M,C^{n+1}})>> M \ot_H C^{n+1}\\
   @.                 @VV\tau V                        @VV\widetilde{\tau} V                                \\
 M \ot H \ot C^{n+1} @>f_{M,C^{n+1}} >> M \ot C^{n+1} @>\text{coker}(f_{M,C^{n+1}})>> M \ot_H C^{n+1}
\end{CD} \nonumber\]
where $f_{M,C^{n+1}}=(\phi_M \ot 1_{C^{n+1}} - 1_M \ot \phi_{C^{n+1}})$, and using the definition of cokernel, 
it can be verified that:
\begin{eqnarray*}
&(\widetilde{\tau}_n^{n+1})&\text{coker}(f_{M,C^{n+1}})=\text{coker}(f_{M,C^{n+1}}) (\tau_n^{n+1}) \\
&& =\text{coker}(f_{M,C^{n+1}}) (1_M \ot \phi_{C^{n+1}}) (\si_{H,M} \rho _M \ot 1_{C^{n+1}}) \\
&& =\text{coker}(f_{M,C^{n+1}}) (\phi_M \ot 1_{C^{n+1}}) (\si_{H,M} \rho _M \ot 1_{C^{n+1}}) \\
&& =\text{coker}(f_{M,C^{n+1}}) (\phi_M \si_{H,M} \rho _M \ot 1_{C^{n+1}}) .
\end{eqnarray*}
Now using the stability property of $M$, $(\phi_M)(\psi_{H , M})(\rho_M) = 1_M$, and
  the universal property of cokernel it is clear that:
\[ \widetilde{\tau}_n^{n+1} =id .  \]
 Also the other properties of a cocyclic object can be easily checked. This finishes the proof of Theorem ~\eqref{main thm}.
 \end{proof}

\begin{example} As a special case, if we put $C=H$ as a $H$-module coalgebra
over itself,  we obtain  a braided version of Connes-Moscovici's
Hopf cyclic theory ~\cite{cm2,cm3,cm4} in any symmetric monoidal
abelian category. We shall explain this example in more details in
the next section. \end{example}

%%**************************************************************************************************
\section{The braided version of Connes-Moscovici's Hopf cyclic cohomology}\label{brdd cm hpf thry}
Let $\mathcal{C}$ be a strict braided monoidal  category and
$(H,\,\Delta,\,\varepsilon,m,\,\eta,S)$ be a Hopf algebra in $\mathcal{C}$. Notice
that, except for Theorem ~\eqref{brdd ver of CM thry},   $\mathcal{C}$ is not
assumed to be symmetric or additive. 
\begin{definition} \label{char and cochar} A
character for $H$ is a morphism $\delta : H\to I$ in $\mathcal{C}$ which is
an algebra map, i.e, :
\[\delta m = \delta \ot \delta ~~~\text{and}~~~ \delta \eta = id_I . \]
A co-character for H is a morphism $\sigma : I \to H$ which is a
coalgebra map, i.e, :
\[\Delta \sigma = \sigma \ot \sigma ~~~\text{and}~~~ \varepsilon \sigma = id_I .\]
A pair $(\delta , \sigma)$ consisting of a character and a
co-character is called a braided modular pair if:
 \[\delta \sigma = id_I .\]
\end{definition}

\begin{definition}\label {(Braided Twised Antipode)} If $\delta$ is a character
for $H$, the corresponding $\delta$-twisted antipode $\widetilde{S}$
is defined by:
\[\widetilde{S} := (\delta \ot S)\Delta .\]

\begin{equation}\label{ws definition}
\xy /r2.5pc/:, 
(1,0)="L11",
"L11"+(0,.5)*!<0pt,6pt>{H},
   "L11"="x1", 
  "x1"+(-.5,-.5)="x2",
  "x1"+(.5,-.5)="x3",
  "x1"+(0,-.25)="c", 
   "x1";"x2"**\crv{"x1" & "c" & "x2"+(0,.25)}, 
   "x1";"x3"**\crv{"x1" & "c" & "x3"+(0,.25)}, 
  "c"+(0,0)*!<0pt,4pt>{\Delta},  
"x2"="L11",
"x3"="L12",
   "L11"="x1", 
  "x1"+(0,-.25)="x2",
  "x1";"x2"**\dir{}?(.5)="c", 
  "x1";"x2"**\dir{-},   
   "c"-(0,0)*!<-3pt,0pt>{\bullet \delta },  
  "x2"+(0,0)*!<0pt,2pt>{\circ},  
   "x2"+(0,0)*!<0pt,2pt>{\circ}, 
  "L12"="x1", 
  "x1"+(0,-.5)="x2",
  "x1";"x2"**\dir{}?(.5)="c", 
  "x1";"x2"**\dir{-},   
  "c"-(0,0)*!<-4pt,0pt>{\bullet S },
"x2"="L11",
"L11"+(0,-.5)*!<0pt,-6pt>{H}, 												
\endxy 
\end{equation}
\end{definition}

\begin{proposition} \label{propert of wS}
 If $\widetilde{S}$ is a $\delta$-twisted antipode for $H$ then we have:
\[\widetilde{S} m=m \psi (\widetilde{S} \ot \widetilde{S}) =m (\widetilde{S} \ot \widetilde{S}) \psi ,\]

\begin{equation}\label{ws m} 
\xy /r2.5pc/:,   
(0,0)="L11",
(1,0)="L12",
  "L11"+(0,.5)*!<0pt,6pt>{H},
  "L12"+(0,.5)*!<0pt,6pt>{H}, 
  "L11"="x1", 
  "L12"="x2",   
  "x1";"x2"**\dir{}?(.5)="m", 
  "m" + (0,-.25)="c", 
  "m" + (0,-.5)="x3", 
  "x1";"x3"**\crv{"x1"-(0,.25) & "c" & "x3"},
  "x2";"x3"**\crv{"x2"-(0,.25) & "c"  & "x3"}, 
"x3"="L21",
  "L21"="x1", 
  "x1"+(0,-.5)="x2",
  "x1";"x2"**\dir{}?(.5)="c", 
  "x1";"x2"**\dir{-},   
  "c"-(0,0)*!<-4pt,0pt>{\bullet \widetilde{\text\small{S}}  }, 
"x2"="L21",
  "x2"-(0,.5)*!<0pt,-6pt>{H}, 
 "x2" + (1.5,.25)*!<0pt,-6pt>{=}, 
(3,0)="L11",
(4,0)="L12",
  "L11"+(0,.5)*!<0pt,6pt>{H},
  "L12"+(0,.5)*!<0pt,6pt>{H}, 
  "L11"="x1", 
  "x1"+(0,-.5)="x2",
  "x1";"x2"**\dir{}?(.5)="c", 
  "x1";"x2"**\dir{-},   
  "c"-(0,0)*!<-4pt,0pt>{\bullet \widetilde{\text\small{S}}  }, 
"x2"="L21",
  "L12"="x1", 
  "x1"+(0,-.5)="x2",
  "x1";"x2"**\dir{}?(.5)="c", 
  "x1";"x2"**\dir{-},   
  "c"-(0,0)*!<-4pt,0pt>{\bullet \widetilde{\text\small{S}}  }, 
"x2"="L22",
  "L21"="x1",       
  "L22"="x2",  , 
  "x1"+(0,-.5)="x3",
  "x2"+(0,-.5)="x4",
  "x1";"x2"**\dir{}?(.5)="m", 
  "m" + (0,-.25)="c", 
  "x1";"x4"**\crv{ "x1"+(0,-.25) & "c" & "x4"+(0,.25)},  
  "x2";"c"+<3pt,2pt>="c2"**\crv{"x2"+(0,-.25) & "c2"},
  "x3";"c"+<-3pt,-2pt>="c1"**\crv{"x3"+(0,.25) & "c1"}, 
  "c"+(0,.10)*!<0pt,-6pt>{\si}, 
"x3"="L21",
"x4"="L22",
  "L21"="x1", 
  "L22"="x2",   
  "x1";"x2"**\dir{}?(.5)="m", 
  "m" + (0,-.25)="c", 
  "m" + (0,-.5)="x3", 
  "x1";"x3"**\crv{"x1"-(0,.25) & "c" & "x3"},
  "x2";"x3"**\crv{"x2"-(0,.25) & "c"  & "x3"}, 
"x3"="L21",
  "x3"-(0,.5)*!<0pt,-6pt>{H}, 
 "x2" + (1.25,.25)*!<0pt,-6pt>{=}, 
(6,0)="L11",
(7,0)="L12",
  "L11"+(0,.5)*!<0pt,6pt>{H},
  "L12"+(0,.5)*!<0pt,6pt>{H}, 
  "L11"="x1",       
  "L12"="x2",  , 
  "x1"+(0,-.5)="x3",
  "x2"+(0,-.5)="x4",
  "x1";"x2"**\dir{}?(.5)="m", 
  "m" + (0,-.25)="c", 
  "x1";"x4"**\crv{ "x1"+(0,-.25) & "c" & "x4"+(0,.25)},  
  "x2";"c"+<3pt,2pt>="c2"**\crv{"x2"+(0,-.25) & "c2"},
  "x3";"c"+<-3pt,-2pt>="c1"**\crv{"x3"+(0,.25) & "c1"}, 
  "c"+(0,.10)*!<0pt,-6pt>{\si}, 
"x3"="L21",
"x4"="L22",
  "L21"="x1", 
  "x1"+(0,-.5)="x2",
  "x1";"x2"**\dir{}?(.5)="c", 
  "x1";"x2"**\dir{-},   
  "c"-(0,0)*!<-4pt,0pt>{\bullet \widetilde{\text\small{S}}  }, 
"x2"="L21",
  "L22"="x1", 
  "x1"+(0,-.5)="x2",
  "x1";"x2"**\dir{}?(.5)="c", 
  "x1";"x2"**\dir{-},  
  "c"-(0,0)*!<-4pt,0pt>{\bullet \widetilde{\text\small{S}}  }, 
"x2"="L22",
  "L21"="x1", 
  "L22"="x2",   
  "x1";"x2"**\dir{}?(.5)="m", 
  "m" + (0,-.25)="c", 
  "m" + (0,-.5)="x3", 
  "x1";"x3"**\crv{"x1"-(0,.25) & "c" & "x3"},
  "x2";"x3"**\crv{"x2"-(0,.25) & "c"  & "x3"}, 
"x3"="L21",
  "x3"-(0,.5)*!<0pt,-6pt>{H}, 
\endxy 
\end{equation}

\[\widetilde{S} \eta=\eta \nonumber ,\]
\[\Delta \widetilde{S} = \psi (\widetilde{S} \ot S)\Delta = (S \ot \widetilde{S}) \psi \Delta ,\]

\begin{equation} \label{ Del ws}
\xy /r2.5pc/:,   
(0,0)="L11",
  "L11"+(0,.5)*!<0pt,6pt>{H},
   "L11"="x1", 
  "x1"+(-.5,-.5)="x2",
  "x1"+(.5,-.5)="x3",
  "x1"+(0,-.25)="c",  
   "x1";"x2"**\crv{"x1" & "c" & "x2"+(0,.25)}, 
   "x1";"x3"**\crv{"x1" & "c" & "x3"+(0,.25)}, 
"x2"= "L21",
"x3"= "L22",
   "L21"="x1", 
  "x1"+(0,-.5)="x2",
  "x1";"x2"**\dir{}?(.5)="c", 
  "x1";"x2"**\dir{-},  
  "c"-(0,0)*!<-4pt,0pt>{\bullet \widetilde{S} }, 
"x2" = "L31",
   "L22"="x1", 
  "x1"+(0,-.5)="x2",
  "x1";"x2"**\dir{}?(.5)="c", 
  "x1";"x2"**\dir{-},   
  "c"-(0,0)*!<-4pt,0pt>{\bullet S }, 
"x2" = "L32",
    "L31"="x1",       
   "L32"="x2",  
  "x1"+(0,-.5)="x3",
  "x2"+(0,-.5)="x4",
  "x1";"x2"**\dir{}?(.5)="m", 
  "m" + (0,-.25)="c", 
  "x1";"x4"**\crv{ "x1"+(0,-.25) & "c" & "x4"+(0,.25)},  
  "x2";"c"+<3pt,2pt>="c2"**\crv{"x2"+(0,-.25) & "c2"},
  "x3";"c"+<-3pt,-2pt>="c1"**\crv{"x3"+(0,.25) & "c1"}, 
  "c"+(0,.10)*!<0pt,-6pt>{\si}, 
"x3" = "L31",
"x4" = "L32",
  "L31"-(0,.5)*!<0pt,-6pt>{H},  
  "L32"-(0,.5)*!<0pt,-6pt>{H}, 
  "L22"+(.5,0)*!<0pt,0pt>{=}, 
"L11" + (2,0)="L11",
  "L11"+(0,.5)*!<0pt,6pt>{H},
   "L11"="x1", 
  "x1"+(-.5,-.5)="x2",
  "x1"+(.5,-.5)="x3",
  "x1"+(0,-.25)="c",  
   "x1";"x2"**\crv{"x1" & "c" & "x2"+(0,.25)}, 
   "x1";"x3"**\crv{"x1" & "c" & "x3"+(0,.25)}, 
"x2"= "L21",
"x3"= "L22",
    "L21"="x1",     
   "L22"="x2",  
  "x1"+(0,-.5)="x3",
  "x2"+(0,-.5)="x4",
  "x1";"x2"**\dir{}?(.5)="m", 
  "m" + (0,-.25)="c", 
  "x1";"x4"**\crv{ "x1"+(0,-.25) & "c" & "x4"+(0,.25)},  
  "x2";"c"+<3pt,2pt>="c2"**\crv{"x2"+(0,-.25) & "c2"},
  "x3";"c"+<-3pt,-2pt>="c1"**\crv{"x3"+(0,.25) & "c1"}, 
  "c"+(0,.10)*!<0pt,-6pt>{\si}, 
"x3" = "L31",
"x4" = "L32",
   "L31"="x1", 
  "x1"+(0,-.5)="x2",
  "x1";"x2"**\dir{}?(.5)="c", 
  "x1";"x2"**\dir{-},  
  "c"-(0,0)*!<-4pt,0pt>{\bullet S }, 
"x2" = "L31",
   "L32"="x1", 
  "x1"+(0,-.5)="x2",
  "x1";"x2"**\dir{}?(.5)="c", 
  "x1";"x2"**\dir{-},  
  "c"-(0,0)*!<-4pt,0pt>{\bullet \widetilde{S} }, 
"x2" = "L32",
  "L31"-(0,.5)*!<0pt,-6pt>{H},  
  "L32"-(0,.5)*!<0pt,-6pt>{H}, 
  "L22"+(.5,0)*!<0pt,0pt>{=}, 
"L11" + (2,0)="L11",
  "L11"+(0,.5)*!<0pt,6pt>{H},
   "L11"="x1", 
  "x1"+(0,-.5)="x2",
  "x1";"x2"**\dir{}?(.5)="c", 
  "x1";"x2"**\dir{-},  
  "c"-(0,0)*!<-4pt,0pt>{\bullet \widetilde{S} }, 
"x2" = "L21",
   "L21"="x1", 
  "x1"+(-.5,-.5)="x2",
  "x1"+(.5,-.5)="x3",
  "x1"+(0,-.25)="c",  
   "x1";"x2"**\crv{"x1" & "c" & "x2"+(0,.25)}, 
   "x1";"x3"**\crv{"x1" & "c" & "x3"+(0,.25)}, 
"x2"= "L21",
"x3"= "L22",
  "L21"-(0,.5)*!<0pt,-6pt>{H}, 
  "L22"-(0,.5)*!<0pt,-6pt>{H}, 
\endxy  
\end{equation}

\[\varepsilon \widetilde{S} = \delta ,~~~\delta \widetilde{S} =\varepsilon ,~~~\widetilde{S} \sigma = S \sigma \nonumber ,\]
\begin{equation}\label{s sgms = sgma -1}
m(\widetilde{S} \sigma \ot \sigma) = m(S \sigma \ot \sigma) = \eta .
\end{equation}
\end{proposition}

\begin{proof}  Here we give a diagrammatic proof of the relation ~\eqref{ws m}. Other relations 
can be easily proved with the same method. 

\[
\xy /r2.5pc/:, 
(1,0)="L11",
(2,0)="L12",
"L11"+(0,.5)*!<0pt,6pt>{H},
"L12"+(0,.5)*!<0pt,6pt>{H},
  "L11"="x1", 
  "L12"="x2",   
  "x1";"x2"**\dir{}?(.5)="m", 
  "m" + (0,-.25)="c", 
  "m" + (0,-.5)="x3", 
  "x1";"x3"**\crv{"x1"-(0,.25) & "c" & "x3"},
  "x2";"x3"**\crv{"x2"-(0,.25) & "c"  & "x3"}, 
"x3"="L11",
  "L11"="x1", 
  "x1"+(0,-.5)="x2",
  "x1";"x2"**\dir{}?(.5)="c", 
  "x1";"x2"**\dir{-},   
  "c"-(0,0)*!<-4pt,0pt>{\bullet \widetilde{S} },
"x2"="L11",
"L11"+(0,-.5)*!<0pt,-6pt>{H},		
  "L11"+(1,.5)*!<0pt,0pt>{=}, 
(3,0)="L11",
(4,0)="L12",
"L11"+(0,.5)*!<0pt,6pt>{H},
"L12"+(0,.5)*!<0pt,6pt>{H},
  "L11"="x1", 
  "L12"="x2",   
  "x1";"x2"**\dir{}?(.5)="m", 
  "m" + (0,-.25)="c", 
  "m" + (0,-.5)="x3", 
  "x1";"x3"**\crv{"x1"-(0,.25) & "c" & "x3"},
  "x2";"x3"**\crv{"x2"-(0,.25) & "c"  & "x3"}, 
"x3"="L11",
   "L11"="x1", 
  "x1"+(-.5,-.5)="x2",
  "x1"+(.5,-.5)="x3",
  "x1"+(0,-.25)="c", 
   "x1";"x2"**\crv{"x1" & "c" & "x2"+(0,.25)}, 
   "x1";"x3"**\crv{"x1" & "c" & "x3"+(0,.25)}, 
"x2"="L11",
"x3"="L12",
   "L11"="x1", 
  "x1"+(0,-.25)="x2",
  "x1";"x2"**\dir{}?(.5)="c", 
  "x1";"x2"**\dir{-},   
  "c"-(0,0)*!<-3pt,0pt>{\bullet \delta },  
  "x2"+(0,0)*!<0pt,2pt>{\circ},  
  "x2"+(0,0)*!<0pt,2pt>{\circ}, 
  "L12"="x1", 
  "x1"+(0,-.5)="x2",
  "x1";"x2"**\dir{}?(.5)="c", 
  "x1";"x2"**\dir{-},   
  "c"-(0,0)*!<-4pt,0pt>{\bullet S },
"x2"="L11",
"L11"+(0,-.5)*!<0pt,-6pt>{H},
  "L11"+(0,1)*!<0pt,0pt>{=}, 
(5,0)="L11",
(7,0)="L12",
"L11"+(0,.5)*!<0pt,6pt>{H},
"L12"+(0,.5)*!<0pt,6pt>{H},
   "L11"="x1", 
  "x1"+(-.5,-.5)="x2",
  "x1"+(.5,-.5)="x3",
  "x1"+(0,-.25)="c",  
   "x1";"x2"**\crv{"x1" & "c" & "x2"+(0,.25)}, 
   "x1";"x3"**\crv{"x1" & "c" & "x3"+(0,.25)}, 
"x2"="L21",
"x3"="L22",
   "L12"="x1", 
  "x1"+(-.5,-.5)="x2",
  "x1"+(.5,-.5)="x3",
  "x1"+(0,-.25)="c",  
   "x1";"x2"**\crv{"x1" & "c" & "x2"+(0,.25)}, 
   "x1";"x3"**\crv{"x1" & "c" & "x3"+(0,.25)}, 
"x2"="L23",
"x3"="L24",
"L21"="L11",
"L22"="L12",
"L23"="L13",
"L24"="L14",
  "L11";"L11" +(0,-.5)="L11"**\dir{-},
  "L12"="x1",       
  "L13"="x2",  , 
  "x1"+(0,-.5)="x3",
  "x2"+(0,-.5)="x4",
  "x1";"x2"**\dir{}?(.5)="m", 
  "m" + (0,-.25)="c", 
  "x1";"x4"**\crv{ "x1"+(0,-.25) & "c" & "x4"+(0,.25)},  
  "x2";"c"+<3pt,2pt>="c2"**\crv{"x2"+(0,-.25) & "c2"},
  "x3";"c"+<-3pt,-2pt>="c1"**\crv{"x3"+(0,.25) & "c1"}, 
  "c"+(0,.10)*!<0pt,-6pt>{\si},
"x3"="L12",
"x4"="L13",
  "L14";"L14" +(0,-.5)="L14"**\dir{-},
  "L11"="x1", 
  "L12"="x2",   
  "x1";"x2"**\dir{}?(.5)="m", 
  "m" + (0,-.25)="c", 
  "m" + (0,-.5)="x3", 
  "x1";"x3"**\crv{"x1"-(0,.25) & "c" & "x3"},
  "x2";"x3"**\crv{"x2"-(0,.25) & "c"  & "x3"}, 
"x3"="L11",
  "L13"="x1", 
  "L14"="x2",  
  "x1";"x2"**\dir{}?(.5)="m", 
  "m" + (0,-.25)="c", 
  "m" + (0,-.5)="x3", 
  "x1";"x3"**\crv{"x1"-(0,.25) & "c" & "x3"},
  "x2";"x3"**\crv{"x2"-(0,.25) & "c"  & "x3"}, 
"x3"="L12",
   "L11"="x1", 
  "x1"+(0,-.25)="x2",
  "x1";"x2"**\dir{}?(.5)="c", 
  "x1";"x2"**\dir{-},   
  "c"-(0,0)*!<-3pt,0pt>{\bullet \delta },  
  "x2"+(0,0)*!<0pt,2pt>{\circ},  
  "L12"="x1", 
  "x1"+(0,-.5)="x2",
  "x1";"x2"**\dir{}?(.5)="c", 
  "x1";"x2"**\dir{-},   
  "c"-(0,0)*!<-4pt,0pt>{\bullet S },
"x2"="L11",
"L11"+(0,-.5)*!<0pt,-6pt>{H},
  "L11" +(1,1.5)*!<0pt,0pt>{=}, 
(9,0)="L11",
(11,0)="L12",
"L11"+(0,.5)*!<0pt,6pt>{H},
"L12"+(0,.5)*!<0pt,6pt>{H},
   "L11"="x1", 
  "x1"+(-.5,-.5)="x2",
  "x1"+(.5,-.5)="x3",
  "x1"+(0,-.25)="c",  
   "x1";"x2"**\crv{"x1" & "c" & "x2"+(0,.25)}, 
   "x1";"x3"**\crv{"x1" & "c" & "x3"+(0,.25)}, 
"x2"="L21",
"x3"="L22",
   "L12"="x1", 
  "x1"+(-.5,-.5)="x2",
  "x1"+(.5,-.5)="x3",
  "x1"+(0,-.25)="c",  
   "x1";"x2"**\crv{"x1" & "c" & "x2"+(0,.25)}, 
   "x1";"x3"**\crv{"x1" & "c" & "x3"+(0,.25)}, 
"x2"="L23",
"x3"="L24",
"L21"="L11",
"L22"="L12",
"L23"="L13",
"L24"="L14",
  "L11";"L11" +(0,-.5)="L11"**\dir{-},
  "L12"="x1",       
  "L13"="x2",  , 
  "x1"+(0,-.5)="x3",
  "x2"+(0,-.5)="x4",
  "x1";"x2"**\dir{}?(.5)="m", 
  "m" + (0,-.25)="c", 
  "x1";"x4"**\crv{ "x1"+(0,-.25) & "c" & "x4"+(0,.25)},  
  "x2";"c"+<3pt,2pt>="c2"**\crv{"x2"+(0,-.25) & "c2"},
  "x3";"c"+<-3pt,-2pt>="c1"**\crv{"x3"+(0,.25) & "c1"}, 
"x3"="L12",
"x4"="L13",
  "L14";"L14" +(0,-.5)="L14"**\dir{-},
   "L11"="x1", 
  "x1"+(0,-.25)="x2",
  "x1";"x2"**\dir{}?(.5)="c", 
  "x1";"x2"**\dir{-},   
  "c"-(0,0)*!<-3pt,0pt>{\bullet \delta },  
  "x2"+(0,0)*!<0pt,2pt>{\circ},  
  "x2"+(0,0)*!<0pt,2pt>{\circ}, 
   "L12"="x1", 
  "x1"+(0,-.25)="x2",
  "x1";"x2"**\dir{}?(.5)="c", 
  "x1";"x2"**\dir{-},   
  "c"-(0,0)*!<-3pt,0pt>{\bullet \delta },  
  "x2"+(0,0)*!<0pt,2pt>{\circ},  
  "x2"+(0,0)*!<0pt,2pt>{\circ}, 
  "L13"="x1", 
  "x1"+(0,-.5)="x2",
  "x1";"x2"**\dir{}?(.5)="c", 
  "x1";"x2"**\dir{-},   
  "c"-(0,0)*!<-4pt,0pt>{\bullet S },
"x2"="L11",
  "L14"="x1", 
  "x1"+(0,-.5)="x2",
  "x1";"x2"**\dir{}?(.5)="c", 
  "x1";"x2"**\dir{-},  
  "c"-(0,0)*!<-4pt,0pt>{\bullet S },
"x2"="L12",
  "L11"="x1",       
  "L12"="x2",  , 
  "x1"+(0,-.5)="x3",
  "x2"+(0,-.5)="x4",
  "x1";"x2"**\dir{}?(.5)="m", 
  "m" + (0,-.25)="c", 
  "x1";"x4"**\crv{ "x1"+(0,-.25) & "c" & "x4"+(0,.25)},  
  "x2";"c"+<3pt,2pt>="c2"**\crv{"x2"+(0,-.25) & "c2"},
  "x3";"c"+<-3pt,-2pt>="c1"**\crv{"x3"+(0,.25) & "c1"}, 
  "c"+(0,.10)*!<0pt,-6pt>{\si},
"x3"="L11",
"x4"="L12",
  "L11"="x1", 
  "L12"="x2",   
  "x1";"x2"**\dir{}?(.5)="m", 
  "m" + (0,-.25)="c", 
  "m" + (0,-.5)="x3", 
  "x1";"x3"**\crv{"x1"-(0,.25) & "c" & "x3"},
  "x2";"x3"**\crv{"x2"-(0,.25) & "c"  & "x3"}, 
"x3"="L11",
"L11"+(0,-.5)*!<0pt,-6pt>{H},
\endxy
\]

\[
\xy /r2.5pc/:, 
  (0,-1)*!<0pt,0pt>{=}, 
(1,0)="L11",
(3,0)="L12",
"L11"+(0,.5)*!<0pt,6pt>{H},
"L12"+(0,.5)*!<0pt,6pt>{H},
   "L11"="x1", 
  "x1"+(-.5,-.5)="x2",
  "x1"+(.5,-.5)="x3",
  "x1"+(0,-.25)="c",  
   "x1";"x2"**\crv{"x1" & "c" & "x2"+(0,.25)}, 
   "x1";"x3"**\crv{"x1" & "c" & "x3"+(0,.25)}, 
"x2"="L21",
"x3"="L22",
   "L12"="x1", 
  "x1"+(-.5,-.5)="x2",
  "x1"+(.5,-.5)="x3",
  "x1"+(0,-.25)="c",  
   "x1";"x2"**\crv{"x1" & "c" & "x2"+(0,.25)}, 
   "x1";"x3"**\crv{"x1" & "c" & "x3"+(0,.25)}, 
"x2"="L23",
"x3"="L24",
"L21"="L11",
"L22"="L12",
"L23"="L13",
"L24"="L14",
   "L11"="x1", 
  "x1"+(0,-.25)="x2",
  "x1";"x2"**\dir{}?(.5)="c", 
  "x1";"x2"**\dir{-},   
  "c"-(0,0)*!<-3pt,0pt>{\bullet \delta },  
  "x2"+(0,0)*!<0pt,2pt>{\circ},  
  "x2"+(0,0)*!<0pt,2pt>{\circ}, 
  "L12"="x1", 
  "x1"+(0,-.5)="x2",
  "x1";"x2"**\dir{}?(.5)="c", 
  "x1";"x2"**\dir{-},   
  "c"-(0,0)*!<-4pt,0pt>{\bullet S },
"x2"="L11",
   "L13"="x1", 
  "x1"+(0,-.25)="x2",
  "x1";"x2"**\dir{}?(.5)="c", 
  "x1";"x2"**\dir{-},   
  "c"-(0,0)*!<-3pt,0pt>{\bullet \delta },  
  "x2"+(0,0)*!<0pt,2pt>{\circ},  
  "x2"+(0,0)*!<0pt,2pt>{\circ}, 
  "L14"="x1", 
  "x1"+(0,-.5)="x2",
  "x1";"x2"**\dir{}?(.5)="c", 
  "x1";"x2"**\dir{-},  
  "c"-(0,0)*!<-4pt,0pt>{\bullet S },
"x2"="L12",
  "L11"="x1",       
  "L12"="x2",  , 
  "x1"+(0,-.5)="x3",
  "x2"+(0,-.5)="x4",
  "x1";"x2"**\dir{}?(.5)="m", 
  "m" + (0,-.25)="c", 
  "x1";"x4"**\crv{ "x1"+(0,-.25) & "c" & "x4"+(0,.25)},  
  "x2";"c"+<3pt,2pt>="c2"**\crv{"x2"+(0,-.25) & "c2"},
  "x3";"c"+<-3pt,-2pt>="c1"**\crv{"x3"+(0,.25) & "c1"}, 
  "c"+(0,.10)*!<0pt,-2pt>{\si},
"x3"="L11",
"x4"="L12",
  "L11"="x1", 
  "L12"="x2",   
  "x1";"x2"**\dir{}?(.5)="m", 
  "m" + (0,-.25)="c", 
  "m" + (0,-.5)="x3", 
  "x1";"x3"**\crv{"x1"-(0,.25) & "c" & "x3"},
  "x2";"x3"**\crv{"x2"-(0,.25) & "c"  & "x3"}, 
"x3"="L11",
"L11"+(0,-.5)*!<0pt,-6pt>{H},
  "L11"+(1.5,1)*!<0pt,0pt>{=}, 
(4.5,0)="L11",
(5.5,0)="L12",
"L11"+(0,.5)*!<0pt,6pt>{H},
"L12"+(0,.5)*!<0pt,6pt>{H},
  "L11"="x1", 
  "x1"+(0,-.5)="x2",
  "x1";"x2"**\dir{}?(.5)="c", 
  "x1";"x2"**\dir{-},  
  "c"-(0,0)*!<-4pt,0pt>{\bullet \widetilde{S} },
"x2"="L11",
  "L12"="x1", 
  "x1"+(0,-.5)="x2",
  "x1";"x2"**\dir{}?(.5)="c", 
  "x1";"x2"**\dir{-},    
  "c"-(0,0)*!<-4pt,0pt>{\bullet \widetilde{S} },
"x2"="L12",
  "L11"="x1",       
  "L12"="x2",  , 
  "x1"+(0,-.5)="x3",
  "x2"+(0,-.5)="x4",
  "x1";"x2"**\dir{}?(.5)="m", 
  "m" + (0,-.25)="c", 
  "x1";"x4"**\crv{ "x1"+(0,-.25) & "c" & "x4"+(0,.25)},  
  "x2";"c"+<3pt,2pt>="c2"**\crv{"x2"+(0,-.25) & "c2"},
  "x3";"c"+<-3pt,-2pt>="c1"**\crv{"x3"+(0,.25) & "c1"}, 
  "c"+(0,.10)*!<0pt,-6pt>{\si},
"x3"="L11",
"x4"="L12",
  "L11"="x1", 
  "L12"="x2",  
  "x1";"x2"**\dir{}?(.5)="m", 
  "m" + (0,-.25)="c", 
  "m" + (0,-.5)="x3", 
  "x1";"x3"**\crv{"x1"-(0,.25) & "c" & "x3"},
  "x2";"x3"**\crv{"x2"-(0,.25) & "c"  & "x3"}, 
"x3"="L11",
"L11"+(0,-.5)*!<0pt,-6pt>{H},
 "x2" + (.5,-.25)*!<0pt,-6pt>{=}, 
(6.5,0)="L11",
(7.5,0)="L12",
  "L11"+(0,.5)*!<0pt,6pt>{H},
  "L12"+(0,.5)*!<0pt,6pt>{H}, 
  "L11"="x1",       
  "L12"="x2",  , 
  "x1"+(0,-.5)="x3",
  "x2"+(0,-.5)="x4",
  "x1";"x2"**\dir{}?(.5)="m", 
  "m" + (0,-.25)="c", 
  "x1";"x4"**\crv{ "x1"+(0,-.25) & "c" & "x4"+(0,.25)},  
  "x2";"c"+<3pt,2pt>="c2"**\crv{"x2"+(0,-.25) & "c2"},
  "x3";"c"+<-3pt,-2pt>="c1"**\crv{"x3"+(0,.25) & "c1"}, 
  "c"+(0,.10)*!<0pt,-6pt>{\si}, 
"x3"="L11",
"x4"="L12",
  "L11"="x1", 
  "x1"+(0,-.5)="x2",
  "x1";"x2"**\dir{}?(.5)="c", 
  "x1";"x2"**\dir{-},  
  "c"-(0,0)*!<-4pt,0pt>{\bullet \widetilde{\text\small{S}}  }, 
"x2"="L11",
  "L12"="x1", 
  "x1"+(0,-.5)="x2",
  "x1";"x2"**\dir{}?(.5)="c", 
  "x1";"x2"**\dir{-}, 
  "c"-(0,0)*!<-4pt,0pt>{\bullet \widetilde{\text\small{S}}  }, 
"x2"="L12",
  "L11"="x1", 
  "L12"="x2",   
  "x1";"x2"**\dir{}?(.5)="m", 
  "m" + (0,-.25)="c", 
  "m" + (0,-.5)="x3", 
  "x1";"x3"**\crv{"x1"-(0,.25) & "c" & "x3"},
  "x2";"x3"**\crv{"x2"-(0,.25) & "c"  & "x3"}, 
"x3"="L11",
  "L11"-(0,.5)*!<0pt,-6pt>{H}, 
\endxy 
\]

In the first identity we used the definition of $\widetilde{S}$, relation ~\eqref{ws definition}, in the 
second one we used relation ~\eqref{compatibilityDelm}, and in the third one relations 
~\eqref{compatibility epn eta} (replacing $\varepsilon$ by $\delta$) and ~\eqref{S m}. In the fourth identity
the naturality of $\si$, relation ~\eqref{naturalityofsi when B is I} is used. in the fifth identity again
 the definition of $\widetilde{S}$, relation ~\eqref{ws definition},  and in the last identity the naturality 
of $\si$,  relation ~\eqref{naturalityofsi}, is used.
\end{proof}

\begin{definition} \label{Braided MPI} A modular pair $(\delta, \sigma)$ for $H$
is called a braided modular pair in involution (BMPI) if:
\[m((m \ot id)(S \sigma \ot \widetilde{S}^2 \ot \sigma))=id .\]
Considering the fact that $ S \sigma = \widetilde{S} \sigma$, the BMPI condition can be shown 
by the following diagram.

\begin{equation}\label{BMPI}
\xy /r2.5pc/:, 
(0,0) = "L11",
 "L11"+(0,.5)*!<0pt,6pt>{H},
"L11" + (-1,0)= "L21",
"L11" = "L22",
"L11" + (1,0)= "L23",
   "L21"  + (0,-.5) ="x1",
  "x1"+(0,-1.5)="x2",
  "x1";"x2"**\dir{}?(.5)="c", 
  "x1";"x2"**\dir{-},  
  "c"-(0,0)*!<-3pt,-8pt>{\bullet \sigma }, 
  "c"+(0,0)*!<-3pt,10pt>{\bullet S }, 
  "x1"+(0,0)*!<0pt,-2pt>{\circ}, 
"x2" ="L21",
   "L22"="x1", 
  "x1"+(0,-2)="x2",
  "x1";"x2"**\dir{}?(.5)="c", 
  "x1";"x2"**\dir{-},  
  "c"-(0,0)*!<-4pt,0pt>{\bullet \widetilde{S}^2 }, 
"x2" ="L22",
   "L23"  + (0,-.5) ="x1", 
  "x1"+(0,-1.5)="x2",
  "x1";"x2"**\dir{}?(.5)="c", 
  "x1";"x2"**\dir{-},   
  "c"-(0,0)*!<-3pt,-8pt>{\bullet \sigma }, 
  "x1"+(0,0)*!<0pt,-2pt>{\circ}, 
"x2" ="L23",
  "L21"="x1", 
  "L22"="x2",  
  "x1";"x2"**\dir{}?(.5)="m", 
  "m" + (0,-.25)="c", 
  "m" + (0,-.5)="x3", 
  "x1";"x3"**\crv{"x1"-(0,.25) & "c" & "x3"},
  "x2";"x3"**\crv{"x2"-(0,.25) & "c"  & "x3"}, 
"x3"="L21",
"L23";"L23" + (0,-.5) = "L22" **\dir{-},  
  "L21"="x1",
  "L22"="x2",  
  "x1";"x2"**\dir{}?(.5)="m", 
  "m" + (0,-.25)="c", 
  "m" + (0,-.5)="x3", 
  "x1";"x3"**\crv{"x1"-(0,.25) & "c" & "x3"},
  "x2";"x3"**\crv{"x2"-(0,.25) & "c"  & "x3"}, 
"x3"="L21",
"L21"+(0,-.5)*!<0pt,-6pt>{H}, 
 "x3" + (1.5,1)*!<0pt,-6pt>{=},
(3.5,0) = "L11",
 "L11"+(0,.5)*!<0pt,6pt>{H},
"L11" + (-1,0)= "L21",
"L11" = "L22",
"L11" + (1,0)= "L23",
   "L21"  + (0,-.5) ="x1",
  "x1"+(0,-1.5)="x2",
  "x1";"x2"**\dir{}?(.5)="c", 
  "x1";"x2"**\dir{-},  
  "c"-(0,0)*!<-3pt,-8pt>{\bullet \sigma }, 
  "c"+(0,0)*!<-3pt,10pt>{\bullet \widetilde{S} }, 
  "x1"+(0,0)*!<0pt,-2pt>{\circ}, 
"x2" ="L21",
   "L22"="x1", 
  "x1"+(0,-2)="x2",
  "x1";"x2"**\dir{}?(.5)="c", 
  "x1";"x2"**\dir{-}, 
  "c"-(0,0)*!<-4pt,0pt>{\bullet \widetilde{S}^2 }, 
"x2" ="L22",
   "L23"  + (0,-.5) ="x1", 
  "x1"+(0,-1.5)="x2",
  "x1";"x2"**\dir{}?(.5)="c", 
  "x1";"x2"**\dir{-},  
  "c"-(0,0)*!<-3pt,-8pt>{\bullet \sigma }, 
  "x1"+(0,0)*!<0pt,-2pt>{\circ}, 
"x2" ="L23",
  "L21"="x1", 
  "L22"="x2",  
  "x1";"x2"**\dir{}?(.5)="m", 
  "m" + (0,-.25)="c", 
  "m" + (0,-.5)="x3", 
  "x1";"x3"**\crv{"x1"-(0,.25) & "c" & "x3"},
  "x2";"x3"**\crv{"x2"-(0,.25) & "c"  & "x3"}, 
"x3"="L21",
"L23";"L23" + (0,-.5) = "L22" **\dir{-},  
  "L21"="x1", 
  "L22"="x2",  
  "x1";"x2"**\dir{}?(.5)="m", 
  "m" + (0,-.25)="c", 
  "m" + (0,-.5)="x3", 
  "x1";"x3"**\crv{"x1"-(0,.25) & "c" & "x3"},
  "x2";"x3"**\crv{"x2"-(0,.25) & "c"  & "x3"}, 
"x3"="L21",
"L21"+(0,-.5)*!<0pt,-6pt>{H}, 
 "x3" + (1.5,1)*!<0pt,-6pt>{=}, 
(6,0)="L11",
(6,-3)="L12",
"L11"; "L12"**\dir{-},  
 "L11"+(0,.5)*!<0pt,6pt>{H},
"L12"+(0,-.5)*!<0pt,-6pt>{H}, 
\endxy   
\end{equation}
\end{definition}

\begin{example} \label{sig I del} 
One can easily check that, if I is considered as a right $H$-module via a character $\delta$:
\[ \phi_I = \delta : I \ot H =H \to I  ,\] and as a left $H$-comodule via a co-character $\sigma$:
\[ \rho_I = \sigma : I \to H \ot I =I ,\] then  $I$ is a braided SAYD module over $H$ if and only if $(\delta , \sigma)$ is a braided MPI. We denote this SAYD module by ${}^{\sigma} I_{\delta} $.
\end{example}

\begin{example} \label{ex ordnary hopfalg}
If $H$ is a Hopf algebra in $Vect_{\mathbb{C}}$, then the above definitions reduce to those of
Connes-Moscovici  ~\cite{cm2,cm3,cm4}. 
\end{example}

Now we are ready to give a braided version of Connes-Moscovici's Hopf cyclic theory.

\begin{theorem} \label{brdd ver of CM thry} 
Suppose $(H,\, (\delta,\sigma))$ is a braided Hopf algebra in a symmetric braided monoidal abelian
category $\mathcal{C}$, where $(\delta, \sigma)$ is a braided MPI. If we put $(C; ~ \phi_C, ~ \Delta_C) = (H; ~ m_H ,~ \Delta_H)$, and $M = {}^\sigma I_\delta$, then the cocyclic object of Theorem ~\eqref{main thm}
reduces to the following one which is a braided version of Connes-Moscovici's Hopf cyclic theory:
 \[C^0(H)=I~~~and~~~C^n(H) = H^n, \,  ~~  n \geq 1 , \] with faces, degeneracies, and cyclic maps given by:
\[
\delta_i = \left\{ \begin{array}{ll}
 (\eta, 1,1,...,1) & \textrm{ $i=0$}\\
 (1,1,...,1, \underset{ i-th}{\Delta}, 1,1,...1) & \textrm{$1 \leq i \leq n-1$}\\
 (1,1,...,1,\sigma) & \textrm{ $i=n$}
  \end{array} \right.
\]

\[\sigma_i =(1,1,..., \underset{(i+1)-th}{\varepsilon}, 1,1...,1), \,~~0 \leq i \leq n \]

\[
\tau_n = \left\{ \begin{array}{ll}
 id_I & \textrm{ $n=0$}\\
 (m_n)(\Delta^{n-1} \widetilde{S} , 1_{H^{n-1}} , \sigma) & \textrm{ $n \neq 0$}
  \end{array} \right.
\]
Here by $m_n$ we mean, $m_1=m$, and for $n \geq 2$:
\[m_n =m_{H^n} = (\underbrace{m,m,...,m}_{n~times})\mathcal{F}_n(\si)  \label{mn si} , \] where
\[\mathcal{F}_n(\si):=\prod_{j=1}^{n-1}~(1_{H^j}, \underbrace{\si, \si,...,\si}_{n-j~times},1_{H^j}).\]
\end{theorem}

for example for $n=2,3$, $\tau_2$ and $\tau_3$ are as follows:
\begin{flalign*}
\tau_2 &= (m_2)(\Delta \widetilde{S} ,1, \sigma) &\\
&=(m,m)(1,\si,1) ((S,\widetilde{S}) \si \Delta,1,\sigma) = (m,m)(1,\si,1) ( \si (\widetilde{S}, S) \Delta,1,\sigma)=&
\end{flalign*}

\begin{equation}\label{diagoftau2} 
\xy /r2.5pc/:,   
(0,0)="L11",
(2,0)="L12",
  "L11"+(0,.5)*!<0pt,6pt>{H},
  "L12"+(0,.5)*!<0pt,6pt>{H}, 
   "L11"="x1", 
  "x1"+(-.5,-.5)="x2",
  "x1"+(.5,-.5)="x3",
  "x1"+(0,-.25)="c",  
   "x1";"x2"**\crv{"x1" & "c" & "x2"+(0,.25)}, 
   "x1";"x3"**\crv{"x1" & "c" & "x3"+(0,.25)}, 
"x2"="L21",
"x3"="L22",
 "L12"="x1", 
  "x1"+(0,-.5)="x2",
  "x1";"x2"**\dir{-},   
"x2"="L23",
  "L12"+ (1,-.25)="x1", 
  "x1"+(0,-.25)="x2",
  "x1";"x2"**\dir{}?(.5)="c", 
  "x1";"x2"**\dir{-},   
  "c"-(0,0)*!<-3pt,0pt>{\bullet \sigma },  
  "x1"+(0,0)*!<0pt,-2pt>{\circ}, 
"x2"="L24",
   "L21"="x1",      
     "L22"="x2",  
  "x1"+(0,-.5)="x3",
  "x2"+(0,-.5)="x4",
  "x1";"x2"**\dir{}?(.5)="m", 
  "m" + (0,-.25)="c", 
  "x1";"x4"**\crv{ "x1"+(0,-.25) & "c" & "x4"+(0,.25)},  
  "x2";"c"+<3pt,2pt>="c2"**\crv{"x2"+(0,-.25) & "c2"},
  "x3";"c"+<-3pt,-2pt>="c1"**\crv{"x3"+(0,.25) & "c1"}, 
  "c"+(0,.10)*!<0pt,-6pt>{\si},
"x3"="L31",
"x4"="L32",
"L23"+(0,-0)="L33",
"L24"+(0,-0)="L34",
  "L31"="x1", 
  "x1"+(0,-.5)="x2",
  "x1";"x2"**\dir{}?(.5)="c", 
  "x1";"x2"**\dir{-},   
  "c"-(0,0)*!<-4pt,0pt>{\bullet \txt\small{S} }, 
"x2"="L41",
  "L32"="x1", 
  "x1"+(0,-.5)="x2",
  "x1";"x2"**\dir{}?(.5)="c", 
  "x1";"x2"**\dir{-},   
  "c"-(0,0)*!<-4pt,0pt>{\bullet \widetilde{\text\small{S}}},  
"x2"="L42",
  "L33"="x1", 
  "x1"+(0,-.5)="x2",
  "x1";"x2"**\dir{-},   
"x2"="L43",
  "L34"="x1", 
  "x1"+(0,-.5)="x2",
  "x1";"x2"**\dir{-},   
"x2"="L44",
"L43";"L43" + (0,-.5) = "L43" **\dir{-},  
"L44";"L44" + (0,-.5) = "L44" **\dir{-},  
"L41"="L21",
"L42"="L22",
"L43"="L23",
"L44"="L24",
"L21";"L21" + (0,-.5) = "L31" **\dir{-},  
  "L22"="x1",       
  "L23"="x2",  , 
  "x1"+(0,-.5)="x3",
  "x2"+(0,-.5)="x4",
  "x1";"x2"**\dir{}?(.5)="m", 
  "m" + (0,-.25)="c", 
  "x1";"x4"**\crv{ "x1"+(0,-.25) & "c" & "x4"+(0,.25)},  
  "x2";"c"+<3pt,2pt>="c2"**\crv{"x2"+(0,-.25) & "c2"},
  "x3";"c"+<-3pt,-2pt>="c1"**\crv{"x3"+(0,.25) & "c1"}, 
  "c"+(0,.10)*!<0pt,-6pt>{\si},
"x3"= "L32",
"x4"= "L33",
"L24";"L24" + (0,-.5) = "L34" **\dir{-},  
"L31"="L21",
"L32"="L22",
"L33"="L23",
"L34"="L24",
  "L21"="x1", 
  "L22"="x2",  
  "x1";"x2"**\dir{}?(.5)="m", 
  "m" + (0,-.25)="c", 
  "m" + (0,-.5)="x3", 
  "x1";"x3"**\crv{"x1"-(0,.25) & "c" & "x3"},
  "x2";"x3"**\crv{"x2"-(0,.25) & "c"  & "x3"}, 
  "x3"+(0,-.5)*!<0pt,-6pt>{H}, 
"x3"="L31",
  "L23"="x1", 
  "L24"="x2",   
  "x1";"x2"**\dir{}?(.5)="m", 
  "m" + (0,-.25)="c", 
  "m" + (0,-.5)="x3", 
  "x1";"x3"**\crv{"x1"-(0,.25) & "c" & "x3"},
  "x2";"x3"**\crv{"x2"-(0,.25) & "c"  & "x3"}, 
  "x3"+(0,-.5)*!<0pt,-6pt>{H},  
"x3"="L32",
"L31"="L21",
"L32"="L22",
 "L32" + (1,1)*!<0pt,0pt>{\overset{~\eqref{naturalityofsi}}{=}}, 
(4.5,0)="L11",
(6.5,0)="L12",
  "L11"+(0,.5)*!<0pt,6pt>{H},
  "L12"+(0,.5)*!<0pt,6pt>{H}, 
   "L11"="x1", 
  "x1"+(-.5,-.5)="x2",
  "x1"+(.5,-.5)="x3",
  "x1"+(0,-.25)="c",  
   "x1";"x2"**\crv{"x1" & "c" & "x2"+(0,.25)}, 
   "x1";"x3"**\crv{"x1" & "c" & "x3"+(0,.25)}, 
"x2"="L21",
"x3"="L22",
 "L12"="x1", 
  "x1"+(0,-.5)="x2",
  "x1";"x2"**\dir{-},  
"x2"="L23",
  "L12"+ (1,-.25)="x1", 
  "x1"+(0,-.25)="x2",
  "x1";"x2"**\dir{}?(.5)="c", 
  "x1";"x2"**\dir{-},  
  "c"-(0,0)*!<-3pt,0pt>{\bullet \sigma },  
  "x1"+(0,0)*!<0pt,-2pt>{\circ}, 
"x2"="L24",
  "L21"="x1", 
  "x1"+(0,-.5)="x2",
  "x1";"x2"**\dir{}?(.5)="c", 
  "x1";"x2"**\dir{-},  
  "c"-(0,0)*!<-4pt,0pt>{\bullet \widetilde{\text\small{S}} },
"x2"="L31",
  "L22"="x1",
  "x1"+(0,-.5)="x2",
  "x1";"x2"**\dir{}?(.5)="c", 
  "x1";"x2"**\dir{-},  
  "c"-(0,0)*!<-4pt,0pt>{\bullet \txt\small{S}}, 
"x2"="L32",
  "L23"="x1", 
  "x1"+(0,-.5)="x2",
  "x1";"x2"**\dir{-},   
"x2"="L33",
  "L24"="x1", 
  "x1"+(0,-.5)="x2",
  "x1";"x2"**\dir{-},   
"x2"="L34",
"L33";"L33" + (0,-.5) = "L33" **\dir{-},  
"L34";"L34" + (0,-.5) = "L34" **\dir{-},  
"L31"="L21",
"L32"="L22",
"L33"="L23",
"L34"="L24",
   "L21"="x1",       
     "L22"="x2",  
  "x1"+(0,-.5)="x3",
  "x2"+(0,-.5)="x4",
  "x1";"x2"**\dir{}?(.5)="m", 
  "m" + (0,-.25)="c", 
  "x1";"x4"**\crv{ "x1"+(0,-.25) & "c" & "x4"+(0,.25)},  
  "x2";"c"+<3pt,2pt>="c2"**\crv{"x2"+(0,-.25) & "c2"},
  "x3";"c"+<-3pt,-2pt>="c1"**\crv{"x3"+(0,.25) & "c1"}, 
  "c"+(0,.10)*!<0pt,-6pt>{\si},
"x3"="L31",
"x4"="L32",
"L23"+(0,-0)="L33",
"L24"+(0,-0)="L34",
"L31"="L21",
"L32"="L22",
"L33"="L23",
"L34"="L24",
"L21";"L21" + (0,-.5) = "L31" **\dir{-},  
  "L22"="x1",       
  "L23"="x2",  , 
  "x1"+(0,-.5)="x3",
  "x2"+(0,-.5)="x4",
  "x1";"x2"**\dir{}?(.5)="m", 
  "m" + (0,-.25)="c", 
  "x1";"x4"**\crv{ "x1"+(0,-.25) & "c" & "x4"+(0,.25)},  
  "x2";"c"+<3pt,2pt>="c2"**\crv{"x2"+(0,-.25) & "c2"},
  "x3";"c"+<-3pt,-2pt>="c1"**\crv{"x3"+(0,.25) & "c1"}, 
  "c"+(0,.10)*!<0pt,-6pt>{\si},
"x3"= "L32",
"x4"= "L33",
"L24";"L24" + (0,-.5) = "L34" **\dir{-},  
"L31"="L21",
"L32"="L22",
"L33"="L23",
"L34"="L24",
  "L21"="x1", 
  "L22"="x2",   
  "x1";"x2"**\dir{}?(.5)="m", 
  "m" + (0,-.25)="c", 
  "m" + (0,-.5)="x3", 
  "x1";"x3"**\crv{"x1"-(0,.25) & "c" & "x3"},
  "x2";"x3"**\crv{"x2"-(0,.25) & "c"  & "x3"}, 
  "x3"+(0,-.5)*!<0pt,-6pt>{H},
"x3"="L31",
  "L23"="x1", 
  "L24"="x2",  
  "x1";"x2"**\dir{}?(.5)="m", 
  "m" + (0,-.25)="c", 
  "m" + (0,-.5)="x3", 
  "x1";"x3"**\crv{"x1"-(0,.25) & "c" & "x3"},
  "x2";"x3"**\crv{"x2"-(0,.25) & "c"  & "x3"}, 
  "x3"+(0,-.5)*!<0pt,-6pt>{H}, 
"x3"="L32",
"L31"="L21",
"L32"="L22",
\endxy 
\end{equation}

\begin{flalign*} 
\tau_3 &= (m_3)(\Delta^2 \widetilde{S} , 1,1, \sigma) &\\
&=(m,m,m)(1,\si,3)(3,\si,1)(2,\si,2)((S,S,\widetilde{S})(1,\si)(\si,1)(1,\si)\Delta^2,1,1,\sigma) = &\\
&\overset{~\eqref{naturalityofsi}}{=}(m,m,m) (S,1,S,1,\widetilde{S},\sigma)(1,\si,2)(3,\si)(2,\si,1)(1,\si,2)(\si,3)(1,\si,2)(\Delta^2,2)&
\end{flalign*}

\begin{equation} \label{diagoftau3}
\xy /r2.5pc/:,  
(0,0)="L11",
(2,0)="L12",
(3,0)="L13",
  "L11"+(0,.5)*!<0pt,6pt>{H},
  "L12"+(0,.5)*!<0pt,6pt>{H}, 
  "L13"+(0,.5)*!<0pt,6pt>{H},
   "L11"="x1", 
  "x1"+(-1,-.5)="x2",
  "x1"+(0,-.5)="x3",
  "x1"+(1,-.5)="x4",
     "x1";"x3"**\dir{}?(.5)="c", 
   "x1";"x2"**\crv{"x1" & "c" & "x2"+(0,.5)}, 
   "x1";"x3"**\crv{"x1" & "c" & "x3"+(0,.5)}, 
   "x1";"x4"**\crv{"x1" & "c" & "x4"+(0,.5)}, 
"x2"="L21",
"x3"="L22",
"x4"="L23",
 "L12"="x1", 
  "x1"+(0,-.5)="x2",
  "x1";"x2"**\dir{-},   
"x2"="L24",
 "L13"="x1", 
  "x1"+(0,-.5)="x2",
  "x1";"x2"**\dir{-},   
"x2"="L25",
  "L13"+ (1,-.25)="x1", 
  "x1"+(0,-.25)="x2",
  "x1";"x2"**\dir{}?(.5)="c", 
  "x1";"x2"**\dir{-},   
  "c"-(0,0)*!<-3pt,0pt>{\bullet \sigma },  
  "x1"+(0,0)*!<0pt,-2pt>{\circ}, 
"x2"="L26",
  "L21"="x1", 
  "x1"+(0,-.5)="x2",
  "x1";"x2"**\dir{-},  
"x2"="L31",
   "L22"="x1",      
     "L23"="x2",  
  "x1"+(0,-.5)="x3",
  "x2"+(0,-.5)="x4",
  "x1";"x2"**\dir{}?(.5)="m", 
  "m" + (0,-.25)="c", 
  "x1";"x4"**\crv{ "x1"+(0,-.25) & "c" & "x4"+(0,.25)},  
  "x2";"c"+<3pt,2pt>="c2"**\crv{"x2"+(0,-.25) & "c2"},
  "x3";"c"+<-3pt,-2pt>="c1"**\crv{"x3"+(0,.25) & "c1"}, 
  "c"+(0,.10)*!<0pt,-6pt>{\si},
"x3"="L32",
"x4"="L33",
"L24"+(0,-0)="L34",
"L25"+(0,-0)="L35",
"L26"+(0,-0)="L36",
   "L31"="x1",      
     "L32"="x2",  
  "x1"+(0,-.5)="x3",
  "x2"+(0,-.5)="x4",
  "x1";"x2"**\dir{}?(.5)="m", 
  "m" + (0,-.25)="c", 
  "x1";"x4"**\crv{ "x1"+(0,-.25) & "c" & "x4"+(0,.25)},  
  "x2";"c"+<3pt,2pt>="c2"**\crv{"x2"+(0,-.25) & "c2"},
  "x3";"c"+<-3pt,-2pt>="c1"**\crv{"x3"+(0,.25) & "c1"}, 
  "c"+(0,.10)*!<0pt,-6pt>{\si},
"x3"="L41",
"x4"="L42",
  "L33"="x1", 
  "x1"+(0,-.5)="x2",
  "x1";"x2"**\dir{-},  
"x2"="L43",
  "L24"="x1", 
  "x1"+(0,-.5)="x2",
  "x1";"x2"**\dir{-},   
"x2"="L44",
  "L25"="x1", 
  "x1"+(0,-.5)="x2",
  "x1";"x2"**\dir{-},   
"x2"="L45",
  "L26"="x1", 
  "x1"+(0,-.5)="x2",
  "x1";"x2"**\dir{-},  
"x2"="L46",
  "L41"="x1", 
  "x1"+(0,-.5)="x2",
  "x1";"x2"**\dir{-}, 
"x2"="L51",
   "L42"="x1",       
     "L43"="x2",  
  "x1"+(0,-.5)="x3",
  "x2"+(0,-.5)="x4",
  "x1";"x2"**\dir{}?(.5)="m", 
  "m" + (0,-.25)="c", 
  "x1";"x4"**\crv{ "x1"+(0,-.25) & "c" & "x4"+(0,.25)},  
  "x2";"c"+<3pt,2pt>="c2"**\crv{"x2"+(0,-.25) & "c2"},
  "x3";"c"+<-3pt,-2pt>="c1"**\crv{"x3"+(0,.25) & "c1"}, 
  "c"+(0,.10)*!<0pt,-6pt>{\si},
"x3"="L52",
"x4"="L53",
  "L44"="x1", 
  "x1"+(0,-.5)="x2",
  "x1";"x2"**\dir{-},   
"x2"="L54",
  "L45"="x1", 
  "x1"+(0,-.5)="x2",
  "x1";"x2"**\dir{-},   
"x2"="L55",
  "L46"="x1", 
  "x1"+(0,-.5)="x2",
  "x1";"x2"**\dir{-},  
"x2"="L56",
  "L51"="x1", 
  "x1"+(0,-.5)="x2",
  "x1";"x2"**\dir{}?(.5)="c", 
  "x1";"x2"**\dir{-},   
  "c"-(0,0)*!<-4pt,0pt>{\bullet \text\small{S} },  
"x2"="L61",
  "L52"="x1",
  "x1"+(0,-.5)="x2",
  "x1";"x2"**\dir{}?(.5)="c", 
  "x1";"x2"**\dir{-},  
  "c"-(0,0)*!<-4pt,0pt>{\bullet \text\small{S}}, 
"x2"="L62",
  "L53"="x1", 
  "x1"+(0,-.5)="x2",
  "x1";"x2"**\dir{}?(.5)="c", 
  "x1";"x2"**\dir{-},  
  "c"-(0,0)*!<-4pt,0pt>{\bullet \widetilde{\text{\small{S}} } }, 
"x2"="L63",
  "L54"="x1", 
  "x1"+(0,-.5)="x2",
  "x1";"x2"**\dir{-},  
"x2"="L64",
  "L55"="x1", 
  "x1"+(0,-.5)="x2",
  "x1";"x2"**\dir{-},  
"x2"="L65",
  "L56"="x1", 
  "x1"+(0,-.5)="x2",
  "x1";"x2"**\dir{-}, 
"x2"="L66",
"L61"="L21",
"L62"="L22",
"L63"="L23",
  "L64";"L64" + (0,-.5)="L24"**\dir{-},  
  "L65";"L65" + (0,-.5)="L25"**\dir{-}, 
  "L66";"L66" + (0,-.5)="L26"**\dir{-}, 
  "L21"="x1", 
  "x1"+(0,-.5)="x2",
  "x1";"x2"**\dir{-},  
"x2"="L31",
  "L22"="x1", 
  "x1"+(0,-.5)="x2",
  "x1";"x2"**\dir{-},   
"x2"="L32",
   "L23"="x1",      
   "L24"="x2", 
  "x1"+(0,-.5)="x3",
  "x2"+(0,-.5)="x4",
  "x1";"x2"**\dir{}?(.5)="m", 
  "m" + (0,-.25)="c", 
  "x1";"x4"**\crv{ "x1"+(0,-.25) & "c" & "x4"+(0,.25)},  
  "x2";"c"+<3pt,2pt>="c2"**\crv{"x2"+(0,-.25) & "c2"},
  "x3";"c"+<-3pt,-2pt>="c1"**\crv{"x3"+(0,.25) & "c1"},   
  "c"+(0,.10)*!<0pt,-6pt>{\si},
"x3"="L33",
"x4"="L34",
  "L25"="x1", 
  "x1"+(0,-.5)="x2",
  "x1";"x2"**\dir{-},   
"x2"="L35",
  "L26"="x1", 
  "x1"+(0,-.5)="x2",
  "x1";"x2"**\dir{-},  
"x2"="L36",
"L31"="L21",
"L32"="L22",
"L33"="L23",
"L34"="L24",
"L35"="L25",
"L36"="L26",
  "L21"="x1", 
  "x1"+(0,-.5)="x2",
  "x1";"x2"**\dir{-},  
"x2"="L31",
  "L22"="x1", 
  "x1"+(0,-.5)="x2",
  "x1";"x2"**\dir{-},  
"x2"="L32",
  "L23"="x1", 
  "x1"+(0,-.5)="x2",
  "x1";"x2"**\dir{-},  
"x2"="L33",
   "L24"="x1",       
   "L25"="x2",  
  "x1"+(0,-.5)="x3",
  "x2"+(0,-.5)="x4",
  "x1";"x2"**\dir{}?(.5)="m", 
  "m" + (0,-.25)="c", 
  "x1";"x4"**\crv{ "x1"+(0,-.25) & "c" & "x4"+(0,.25)},  
  "x2";"c"+<3pt,2pt>="c2"**\crv{"x2"+(0,-.25) & "c2"},
  "x3";"c"+<-3pt,-2pt>="c1"**\crv{"x3"+(0,.25) & "c1"}, 
  "c"+(0,.10)*!<0pt,-6pt>{\si},
"x3"="L34",
"x4"="L35",
  "L26"="x1", 
  "x1"+(0,-.5)="x2",
  "x1";"x2"**\dir{-},  
"x2"="L36",
"L31"="L21",
"L32"="L22",
"L33"="L23",
"L34"="L24",
"L35"="L25",
"L36"="L26",
  "L21"="x1", 
  "x1"+(0,-.5)="x2",
  "x1";"x2"**\dir{-},  
"x2"="L31",
   "L22"="x1",      
   "L23"="x2",  
  "x1"+(0,-.5)="x3",
  "x2"+(0,-.5)="x4",
  "x1";"x2"**\dir{}?(.5)="m", 
  "m" + (0,-.25)="c", 
  "x1";"x4"**\crv{ "x1"+(0,-.25) & "c" & "x4"+(0,.25)},  
  "x2";"c"+<3pt,2pt>="c2"**\crv{"x2"+(0,-.25) & "c2"},
  "x3";"c"+<-3pt,-2pt>="c1"**\crv{"x3"+(0,.25) & "c1"}, 
  "c"+(0,.10)*!<0pt,-6pt>{\si},
"x3"="L32",
"x4"="L33",
  "L24"="x1", 
  "x1"+(0,-.5)="x2",
  "x1";"x2"**\dir{-},   
"x2"="L34",
  "L25"="x1", 
  "x1"+(0,-.5)="x2",
  "x1";"x2"**\dir{-}, 
"x2"="L35",
  "L26"="x1", 
  "x1"+(0,-.5)="x2",
  "x1";"x2"**\dir{-}, 
"x2"="L36",
"L31"="L21",
"L32"="L22",
"L33"="L23",
    "L34"="L24",
"L35"="L25",
"L36"="L26",
  "L21"="x1", 
  "L22"="x2",  
  "x1";"x2"**\dir{}?(.5)="m", 
  "m" + (0,-.25)="c", 
  "m" + (0,-.5)="x3", 
  "x1";"x3"**\crv{"x1"-(0,.25) & "c" & "x3"},
  "x2";"x3"**\crv{"x2"-(0,.25) & "c"  & "x3"}, 
"x3"= "L31",
  "L23"="x1", 
  "L24"="x2",  
  "x1";"x2"**\dir{}?(.5)="m", 
  "m" + (0,-.25)="c", 
  "m" + (0,-.5)="x3", 
  "x1";"x3"**\crv{"x1"-(0,.25) & "c" & "x3"},
  "x2";"x3"**\crv{"x2"-(0,.25) & "c"  & "x3"}, 
"x3"= "L32",
  "L25"="x1", 
  "L26"="x2",  
  "x1";"x2"**\dir{}?(.5)="m", 
  "m" + (0,-.25)="c", 
  "m" + (0,-.5)="x3", 
  "x1";"x3"**\crv{"x1"-(0,.25) & "c" & "x3"},
  "x2";"x3"**\crv{"x2"-(0,.25) & "c"  & "x3"}, 
"x3"= "L33",
"L31"="L21",
"L32"="L22",
"L33"="L23",
  "L21"-(0,.5)*!<0pt,-6pt>{H},  
  "L22"-(0,.5)*!<0pt,-6pt>{H}, 
  "L23"-(0,.5)*!<0pt,-6pt>{H}, 
  "L23"+(1,2)*!<0pt,-6pt>{\overset{~\eqref{naturalityofsi}}{=}}, 
(6,0)="L11",
(8,0)="L12",
(9,0)="L13",
  "L11"+(0,.5)*!<0pt,6pt>{H},
  "L12"+(0,.5)*!<0pt,6pt>{H}, 
  "L13"+(0,.5)*!<0pt,6pt>{H},
   "L11"="x1",
   "x1"+(-1,-.5)="x2",
  "x1"+(0,-.5)="x3",
  "x1"+(1,-.5)="x4",
  "x1"+(0,-.25)="c", 
   "x1";"x2"**\crv{"x1" & "c" & "x2"+(0,.25)}, 
   "x1";"x3"**\crv{"x1" & "c" & "x3"+(0,.25)}, 
   "x1";"x4"**\crv{"x1" & "c" & "x4"+(0,.25)}, 
"x2"="L21",
"x3"="L22",
"x4"="L23",
 "L12"="x1", 
  "x1"+(0,-.5)="x2",
  "x1";"x2"**\dir{-},  
"x2"="L24",
 "L13"="x1", 
  "x1"+(0,-.5)="x2",
  "x1";"x2"**\dir{-},  
"x2"="L25",
  "L21"="x1", 
  "x1"+(0,-.5)="x2",
  "x1";"x2"**\dir{-},  
"x2"= "L31",
   "L22"="x1",     
   "L23"="x2",  
  "x1"+(0,-.5)="x3",
  "x2"+(0,-.5)="x4",
  "x1";"x2"**\dir{}?(.5)="m", 
  "m" + (0,-.25)="c", 
  "x1";"x4"**\crv{ "x1"+(0,-.25) & "c" & "x4"+(0,.25)},  
  "x2";"c"+<3pt,2pt>="c2"**\crv{"x2"+(0,-.25) & "c2"},
  "x3";"c"+<-3pt,-2pt>="c1"**\crv{"x3"+(0,.25) & "c1"}, 
  "c"+(0,.10)*!<0pt,-6pt>{\si},  
"x3"= "L32",
"x4"= "L33",
  "L24"="x1", 
  "x1"+(0,-.5)="x2",
  "x1";"x2"**\dir{-},  
"x2"= "L34",
  "L25"="x1", 
  "x1"+(0,-.5)="x2",
  "x1";"x2"**\dir{-},   
"x2"= "L35",
"L31"="L21",
"L32"="L22",
"L33"="L23",
"L34"="L24",
"L35"="L25",
   "L21" ="x1",       
   "L22"="x2",  
  "x1"+(-0,-.5)="x3",
  "x2"+(0,-.5)="x4",
  "x1";"x2"**\dir{}?(.5)="m", 
  "m" + (0,-.25)="c", 
  "x1";"x4"**\crv{ "x1"+(0,-.25) & "c" & "x4"+(0,.25)},  
  "x2";"c"+<3pt,2pt>="c2"**\crv{"x2"+(0,-.25) & "c2"},
  "x3";"c"+<-3pt,-2pt>="c1"**\crv{"x3"+(0,.25) & "c1"}, 
  "c"+(0,.10)*!<0pt,-6pt>{\si}, 
"x3"= "L31",
"x4"= "L32",
  "L23"="x1", 
  "x1"+(0,-.5)="x2",
  "x1";"x2"**\dir{-},   
"x2"= "L33",
  "L24"="x1", 
  "x1"+(0,-.5)="x2",
  "x1";"x2"**\dir{-},   
"x2"= "L34",
  "L25"="x1", 
  "x1"+(0,-.5)="x2",
  "x1";"x2"**\dir{-},  
"x2"= "L35",
"L31"="L21",
"L32"="L22",
"L33"="L23",
"L34"="L24",
"L35"="L25",
  "L21"="x1", 
  "x1"+(0,-.5)="x2",
  "x1";"x2"**\dir{-},  
"x2"= "L31",
   "L22"="x1",      
   "L23"="x2",  
  "x1"+(0,-.5)="x3",
  "x2"+(0,-.5)="x4",
  "x1";"x2"**\dir{}?(.5)="m", 
  "m" + (0,-.25)="c", 
  "x1";"x4"**\crv{ "x1"+(0,-.25) & "c" & "x4"+(0,.25)},  
  "x2";"c"+<3pt,2pt>="c2"**\crv{"x2"+(0,-.25) & "c2"},
  "x3";"c"+<-3pt,-2pt>="c1"**\crv{"x3"+(0,.25) & "c1"}, 
  "c"+(0,.10)*!<0pt,-6pt>{\si}, 
"x3"= "L32",
"x4"= "L33",
  "L24"="x1", 
  "x1"+(0,-.5)="x2",
  "x1";"x2"**\dir{-},   
"x2"= "L34",
  "L25"="x1", 
  "x1"+(0,-.5)="x2",
  "x1";"x2"**\dir{-},  
"x2"= "L35",
"L31"="L21",
"L32"="L22",
"L33"="L23",
"L34"="L24",
"L35"="L25",
  "L21"="x1", 
  "x1"+(0,-.5)="x2",
  "x1";"x2"**\dir{-},  
"x2"= "L31",
  "L22"="x1", 
  "x1"+(0,-.5)="x2",
  "x1";"x2"**\dir{-},  
"x2"= "L32",
   "L23"="x1",       
   "L24"="x2",  
  "x1"+(0,-.5)="x3",
  "x2"+(0,-.5)="x4",
  "x1";"x2"**\dir{}?(.5)="m", 
  "m" + (0,-.25)="c", 
  "x1";"x4"**\crv{ "x1"+(0,-.25) & "c" & "x4"+(0,.25)},  
  "x2";"c"+<3pt,2pt>="c2"**\crv{"x2"+(0,-.25) & "c2"},
  "x3";"c"+<-3pt,-2pt>="c1"**\crv{"x3"+(0,.25) & "c1"}, 
  "c"+(0,.10)*!<0pt,-6pt>{\si},  
"x3"= "L33",
"x4"= "L34",
  "L25"="x1", 
  "x1"+(0,-.5)="x2",
  "x1";"x2"**\dir{-},  
"x2"= "L35",
"L31"="L21",
"L32"="L22",
"L33"="L23",
"L34"="L24",
"L35"="L25",
  "L21"="x1", 
  "x1"+(0,-.5)="x2",
  "x1";"x2"**\dir{-},   
"x2"= "L31",
  "L22"="x1", 
  "x1"+(0,-.5)="x2",
  "x1";"x2"**\dir{-},   
"x2"= "L32",
  "L23"="x1", 
  "x1"+(0,-.5)="x2",
  "x1";"x2"**\dir{-},   
"x2"= "L33",
   "L24"="x1",       
   "L25"="x2",  
  "x1"+(0,-.5)="x3",
  "x2"+(0,-.5)="x4",
  "x1";"x2"**\dir{}?(.5)="m", 
  "m" + (0,-.25)="c", 
  "x1";"x4"**\crv{ "x1"+(0,-.25) & "c" & "x4"+(0,.25)},  
  "x2";"c"+<3pt,2pt>="c2"**\crv{"x2"+(0,-.25) & "c2"},
  "x3";"c"+<-3pt,-2pt>="c1"**\crv{"x3"+(0,.25) & "c1"}, 
  "c"+(0,.10)*!<0pt,-6pt>{\si},  
"x3"= "L34",
"x4"= "L35",
"L31"="L21",
"L32"="L22",
"L33"="L23",
"L34"="L24",
"L35"="L25",
  "L21"="x1", 
  "x1"+(0,-.5)="x2",
  "x1";"x2"**\dir{-},   
"x2"= "L31",
   "L22"="x1",       
   "L23"="x2", 
  "x1"+(0,-.5)="x3",
  "x2"+(0,-.5)="x4",
  "x1";"x2"**\dir{}?(.5)="m", 
  "m" + (0,-.25)="c", 
  "x1";"x4"**\crv{ "x1"+(0,-.25) & "c" & "x4"+(0,.25)},  
  "x2";"c"+<3pt,2pt>="c2"**\crv{"x2"+(0,-.25) & "c2"},
  "x3";"c"+<-3pt,-2pt>="c1"**\crv{"x3"+(0,.25) & "c1"}, 
  "c"+(0,.10)*!<0pt,-6pt>{\si},  
"x3"= "L32",
"x4"= "L33",
  "L24"="x1", 
  "x1"+(0,-.5)="x2",
  "x1";"x2"**\dir{-},  
"x2"= "L34",
  "L25"="x1", 
  "x1"+(0,-.5)="x2",
  "x1";"x2"**\dir{-},   
"x2"= "L35",
"L31"="L21",
"L32"="L22",
"L33"="L23",
"L34"="L24",
"L35"="L25",
  "L21" ="x1", 
  "x1"+(0,-.5)="x2",
  "x1";"x2"**\dir{}?(.5)="c", 
  "x1";"x2"**\dir{-},  
  "c"-(0,0)*!<-4pt,0pt>{\bullet S }, 
"x2"= "L31",
  "L22"="x1", 
  "x1"+(0,-.5)="x2",
  "x1";"x2"**\dir{-},  
"x2"= "L32",
  "L23" ="x1", 
  "x1"+(0,-.5)="x2",
  "x1";"x2"**\dir{}?(.5)="c", 
  "x1";"x2"**\dir{-}, 
  "c"-(0,0)*!<-4pt,0pt>{\bullet S },  
"x2"= "L33",
  "L24"="x1", 
  "x1"+(0,-.5)="x2",
  "x1";"x2"**\dir{-},  
"x2"= "L34",
  "L25" ="x1", 
  "x1"+(0,-.5)="x2",
  "x1";"x2"**\dir{}?(.5)="c", 
  "x1";"x2"**\dir{-},  
  "c"-(0,0)*!<-4pt,0pt>{\bullet \widetilde{S} }, 
"x2"= "L35",
    "L25"+(1,-.25) ="x1", 
  "x1"+(0,-.25)="x2",
  "x1";"x2"**\dir{}?(.5)="c", 
  "x1";"x2"**\dir{-},  
  "c"-(0,0)*!<-3pt,0pt>{\bullet \sigma },  
  "x1"+(0,0)*!<0pt,-2pt>{\circ},  
"x2"= "L36",
"L31"="L21",
"L32"="L22",
"L33"="L23",
"L34"="L24",
"L35"="L25",
"L36"="L26",
   "L21" ="x1", 
    "L22"="x2",  
  "x1";"x2"**\dir{}?(.5)="m", 
  "m" + (0,-.25)="c", 
  "m" + (0,-.5)="x3", 
  "x1";"x3"**\crv{"x1"-(0,.25) & "c" & "x3"},
  "x2";"x3"**\crv{"x2"-(0,.25) & "c"  & "x3"}, 
"x3"="L31",
   "L23" ="x1", 
    "L24"="x2",  
  "x1";"x2"**\dir{}?(.5)="m", 
  "m" + (0,-.25)="c", 
  "m" + (0,-.5)="x3", 
  "x1";"x3"**\crv{"x1"-(0,.25) & "c" & "x3"},
  "x2";"x3"**\crv{"x2"-(0,.25) & "c"  & "x3"}, 
"x3"="L32",
   "L25" ="x1", 
    "L26"="x2",  
  "x1";"x2"**\dir{}?(.5)="m", 
  "m" + (0,-.25)="c", 
  "m" + (0,-.5)="x3", 
  "x1";"x3"**\crv{"x1"-(0,.25) & "c" & "x3"},
  "x2";"x3"**\crv{"x2"-(0,.25) & "c"  & "x3"}, 
"x3"="L33",
"L31"="L21",
"L32"="L22",
"L33"="L23",
  "L21"-(0,.5)*!<0pt,-6pt>{H},  
  "L22"-(0,.5)*!<0pt,-6pt>{H}, 
  "L23"-(0,.5)*!<0pt,-6pt>{H}, 
\endxy
\end{equation}

\begin{remark} In section ~\eqref{CM Hf cyclc in nonsymm} we discuss the
analogue of this theorem in a non-symmetric braided monoidal abelian
category. 
\end{remark}

\begin{remark} One can define the notions of $H$-module algebra and
$\delta$-invariant $\sigma$-trace and define a characteristic map as in
~\cite{cm2,cm3,cm4}. \end{remark}

\begin{example} {(Connes-Moscovici's theory)} If one puts $\mathcal{C}
=Vect_{\mathbb{C}}$ and $I = \mathbb{C}$, then the above formulas
reduce to those in ~\cite{cm2,cm3,cm4}. \end{example}
 As another example of
the above theory, we devote Section ~\eqref{super hpf alg} to provide
a Hopf cyclic theory for super Hopf algebras. But before that we
give some more results here:

\begin{lemma}  \label{H cocomm S2= id}
If $H$ is commutative in the sense that $m\si = m$ or cocommutative in the sense that $\Delta =  \si \Delta $,  then $ S^2 =id $ and thus $(\varepsilon ,\mathfrak{1})$ is an MPI for $H$.
\end{lemma}

\begin{proof} 
We have:
\begin{flalign*} 
1 = id_H &\overset{(1)}{=} m(\eta ,1)(\varepsilon,1)\Delta = m(\eta \varepsilon ,1)\Delta &\\
&\overset{(2)}{=} m(S \eta \varepsilon , 1)\Delta = m ( S(m(1,S)\Delta ,1 ) \Delta &\\
&\overset{(3)}{=} m ( m(S,S)\si(1,S) \Delta ,1 ) \Delta = m ( m(S^2,S)\si \Delta ,1 ) \Delta &\\
&\overset{(4)}{=} m ( m(S^2,S) \Delta ,1 ) \Delta = m(m,1)(S^2,S,1)(\Delta,1)\Delta &\\
&\overset{(5)}{=}  m(1,m)(S^2,S,1)(1, \Delta)\Delta =  m ( S^2,m(S,1)\Delta ) \Delta = m ( S^2, \eta \varepsilon ) \Delta&\\
&= m(1, \eta) (S^2, 1_I) (1,\varepsilon)\Delta = (S^2, 1_I) &\\
&\overset{(6)}{=} S^2&
\end{flalign*}

Stages (1) to (6) are explained below:
\item (1) We use $(\varepsilon,1)\Delta = m(\eta,1) = 1$.
\item (2) We use $S \eta = \eta$ and $m(1,S)\Delta = \eta \varepsilon$.
\item (3) We use $Sm = m (S,S) \si$ , commute $\si$ and $(1,S)$ using the naturality of $\si$, and use $(S,S)(S,1)=(S^2, S)$.
\item (4) In the commutative case we first use $(S,S)\si=\si (S,S)$ and then use $m\si=m$. In the cocommutative case we use $\Delta =  \si \Delta$, and then $( m(S^2,S) \Delta ,1 )=(m,1)(S^2,S,1)(\Delta,1)$.
\item (5) We use $m(m,1)=m(1,m)$, $(\Delta,1)\Delta =(1, \Delta)\Delta$, and $m(S,1)\Delta = \eta \varepsilon$.
\item (6) We use  $(\varepsilon,1)\Delta = m(\eta,1) = 1$ again.
\end{proof} 

\begin{remark} As one can see, in the proof of the above theorem we didn't
need the symmetric property $\si^2 = id$ for $\mathcal{C}$. \end{remark}

The following braided version of Theorem (4.2) in ~\cite{kr1}, can be
proved along the same lines.

\begin{theorem}  \label{ciclc hom=sum hoch hom} 
If $H$ is commutative, then
\begin{equation}\label{Hc =sum HH} 
HC_{(\varepsilon,\mathfrak{1})}^n (H)= \bigoplus_{i=n(mod~2)} HH^i (H,I), \forall n \geq 0,  
\end{equation}
where the left hand side is the cyclic cohomology of the braided Hopf algebra $H$ with MPI $(\varepsilon,\mathfrak{1})$ and the right hand side is the Hochschild cohomology of coalgebra $H$ with coefficients in $H$-bimodule $I$.
\end{theorem}   

%%************************************************************************************************************
\section{Hopf cyclic cohomology  for super Hopf algebras}\label{super hpf alg} 
In this section we give explicit formulas for the Hopf cyclic complex in the special case of super Hopf
algebras and compute it in the super commutative case. The degree of an element $a$  in a super vector 
space will be denoted by $|a|$.

\begin{definition} \label{def of superhopfalg 2} 
A Hopf algebra $H$ in $\mathbb{Z}_2 $-Mod is called a super Hopf algebra. Thus a super Hopf algebra
 $H = H_0 \oplus H_1$ is simultaneously a super algebra and a super  coalgebra and the two structures
 are compatible in the sense that  for all homogeneous elements $a, b~\in H$:
\begin{equation}\label{cmptblty for suprhopfalg}
\Delta(ab) = (-1)^{|a_{(2)}||b_{(1)}|}~ (a_{(1)} b_{(1)} \ot a_{(2)} b_{(2)}).
\end{equation}
 Furthermore, there is an even map  $S: H \to H$, the antipode,  such
 that:
 $ S(h^{(1)})h^{(2)} =h^{(1)}S(h^{(2)})=\varepsilon(h)1$, for all $h$ in $H$.
\end{definition}

One easily checks that formulas ~\eqref{S m} and ~\eqref{Del S} in proposition ~\eqref{prps of S} reduce to:
\[S(ab) = (-1)^{|a||b|}~S(b)S(a)  \label{frml dgr S ab} ,\]
\[\Delta S(a) = S(a)_{(1)} \ot S(a)_{(2)} = (-1)^{|a_{(1)}||a_{(2)}|}~S(a_{(2)}) \ot S(a_{(1)}) .\]
\begin{remark} We emphasis that condition ~\eqref{cmptblty for suprhopfalg} shows that a
super Hopf algebra is, in general, not a Hopf algebra in the category of
vector spaces. It is of course a Hopf algebra object in the category of super Hopf algebras 
(Example \ref{cat of CZ_2 mods})
\end{remark}

Let $H =H_0 \oplus H_1$  be a super Hopf algebra and $(\delta,
\sigma)$ be a super MPI for $H$. One checks that formulas ~\eqref{ws m}, 
~\eqref{ Del ws}, and ~\eqref{s sgms = sgma -1} in proposition ~\eqref{propert of wS} reduce to:
\[\widetilde{S}(ab) = (-1)^{|a||b|}~\widetilde{S}(b) \widetilde{S}(a) ,\]
\[\Delta \widetilde{S} (a)= \widetilde{S}(a)_{(1)} \ot \widetilde{S}(a)_{(2)} = (-1)^{|a_{(1)}||a_{(2)}|}~S(a_{(2)}) \ot \widetilde{S}(a_{(1)}),\]
\[S(\sigma) = \widetilde{S}(\sigma) = \sigma^{-1}.\]

\begin{proposition} \label{hpf cycl for sup hopfalg} Let $H =H_0 \oplus H_1$ be a
super Hopf algebra endowed with a braided modular pair in involution
$(\delta, \sigma)$. Then the complex, faces, degeneracies and cyclic
maps of Theorem ~\eqref{brdd ver of CM thry} can be written as :
 \[C^0(H)=\mathbb{C}~~~and~~~ C^n(H) = H^n, n \geq 1 , \]
\[
\delta_i(h_1,...,h_{n-1}) = \left\{ \begin{array}{ll}
 (1, h_1 ,h_2 ,...,h_{n-1}) & \textrm{ $i=0$}\\
 ( h_1 ,h_2 ,...,h_i^{(1)} ,h_i^{(2)},...,  h_{n-1}) & \textrm{$1 \leq i \leq n-1$}\\
 (h_1 ,h_2 ,...,h_{n-1}, \sigma) & \textrm{ $i=n$}
  \end{array} \right.
\]
\[ \sigma_i (h_1, h_2,..., h_{n+1})=\varepsilon(h_{i+1})~(h_1, h_2,...,h_i, h_{i+2},..., h_{n+1}), 0 \leq i \leq n, \]
\[ \tau_n (h_1,h_2,...,h_n) = \alpha \beta ~ ( S(h_1^{(n)}) h_2, S(h_1^{(n-1)}) h_3,...,S(h_1^{(2)}) h_n,\widetilde{S}(h_1^{(1)}) \sigma )  ,\] where $h_i$'s are homogeneous elements and:
\[ \alpha = \prod_{i=1}^{n-1}(-1)^{({|h_1^{(1)}|+...+|h_1^{(i)}|)(|h_1^{(i+1)}|)}} ,\]
\[   \beta = \prod_{j=1}^{n-1} (-1)^{|h_1^{(j)}| ( |h_2|+|h_3|+...+|h_{n-j+1}|) } . \]
\end{proposition}

The next lemma  is a corollary of Lemma ~\eqref{H cocomm S2= id}.
\begin{lemma} If a super Hopf algebra $H =H_0 \oplus H_1$ is super
commutative or super-cocommutative, 
 then $S^2 =1$ and thus $(\varepsilon,\mathfrak{1})$ is a super MPI for $H$.
\end{lemma}

The next proposition is a corollary of Theorem ~\eqref{ciclc hom=sum
hoch hom}

\begin{proposition} If a super Hopf algebra $H =H_0 \oplus H_1$ is super
commutative, then we have a decomposition:
\[ HC_{(\varepsilon,\mathfrak{1})}^n (H) = \bigoplus_{i=n(mod~2)} HH^i (H,\mathbb{C}), \forall n \geq 0. \]
\end{proposition}

\begin{remark} All the results of this section easily extend to the case of
differential graded super Hopf algebras. 
\end{remark}

%%*************************************************************************************************
\section{Hopf cyclic cohomology of the enveloping algebra of a super Lie algebra}\label{super lie alg}
A good reference for super Lie algebras, their enveloping
algebras, and a super analogue of the Poincar\'e-Birkhoff-Witt
theorem is ~\cite{ak6}. This latter result is specially important for
the proof of Lemma ~\eqref{antsummap} . Let $\mathfrak{g}=\mathfrak{g}_0 \oplus \mathfrak{g}_1$
be a super Lie algebra, let
\[\bigwedge \mathfrak{g}:=\frac{T(\mathfrak{g})}{(a\ot b +(-1)^{|a||b|} b \ot a)},\]
be the exterior algebra of $\mathfrak{g}$ and let
\[H=U(\mathfrak{g}) :=\frac{T(\mathfrak{g})}{([a,b] - a\ot b +(-1)^{|a||b|} b \ot a)},\]
be the enveloping algebra of $\mathfrak{g}$. Here $T(\mathfrak{g})$ is the tensor
algebra of $\mathfrak{g}$. All these constructions are done in the category
of super vector spaces ~\cite{ak6}. $U(\mathfrak{g})$ is a super cocommutative
super Hopf algebra. Our goal in this section is to show that,
analogous to the non-graded case ~\cite{cm2,cm3,cm4}, the relation
\[HP_{(\delta,1)}^*(U(\mathfrak{g}))=\bigoplus_{i=*(mod~2)} H_i (\mathfrak{g} ; \mathbb{C}_ \delta) ,\]
 holds, where $\delta$ is a character for $\mathfrak{g}$.  Here $HP_{(\delta,1)}^*(U(\mathfrak{g}))$ is
the periodic Hopf cyclic cohomology of the super Hopf algebra $H=U(\mathfrak{g})$, and $H_i (\mathfrak{g} ; \mathbb{C}_\delta) $ is the Lie algebra homology of $\mathfrak{g}$ with coefficient in the $\mathfrak{g}$-module $\mathbb{C}_ \delta$.

First, we notice that the Hochschild cohomology $HH^* (H,\mathbb{C}_\sigma)$ depends only on the coalgebra structure of $H$ and the grouplike elements $\sigma$ and $1$. In fact we have:

\begin{lemma}\label{hochofh}
$ HH^* (H,\mathbb{C}_ \sigma) = \text{Cotor}_H^*(\mathbb{C},\mathbb{C}_\sigma)$. 
\end{lemma} 

Let
\[S(\mathfrak{g}):=\frac{T(\mathfrak{g})}{(a\ot b -(-1)^{|a||b|} b \ot a)} ,\] denote
the symmetric algebra of the super vector space $\mathfrak{g}$. It is a
super cocommutative super Hopf algebra with the comultiplication
defined by $\Delta(x)= x \ot 1 + 1 \ot x$ for homogeneous elements $x$
of $\mathfrak{g}$. Since by the super Poincar\'e-Birkhoff-Witt theorem ~\cite{ak6},
$U(\mathfrak{g}) = S(\mathfrak{g})$ as super coalgebras one, using Lemma ~\eqref{hochofh},
 can prove the following lemma analogous to the non super case ~\cite{l}.

\begin{lemma}\label{antsummap}
 The antisymmetrization map
\[
\begin{CD}
A:\bigwedge^n \mathfrak{g}    @>>>  U(\mathfrak{g})^n ,
\end{CD}
\]
defined by
\[
A(x_1 \wedge...\wedge x_n)= (\sum_{\sigma \in S_n} (-1)^{\alpha_{\sigma}} sign(\sigma)(x_{\sigma(1)},...,x_{\sigma(n)}))/n!,\]%\label{antsymm map}
induces an isomorphism $ HH^* (U(\mathfrak{g}),\mathbb{C}) = \bigwedge^*(\mathfrak{g})$.  Here
\[\alpha_{\sigma}= \sum_{i=1}^n |x_{\sigma(i)}|(|x_1|+|x_2|+ \cdots \widehat{|x_{\sigma(j)}|}+\cdots+|x_{\sigma(i)-1}|),\]
 where $\widehat{|x_{\sigma(j)}|}$ means that if there are any of $|x_{\sigma(j)}|$'s for all $j<i$ they should be
 omitted. Simply $\alpha_{\sigma}$ contains $|a||b|$ for any two
 elements $a$ and $b$ of $x_i$'s, if they cross each other.
\end{lemma}  

The following complex is the super analogue of the
Chevalley-Eilenberg complex to compute the Lie algebra homology $H_
\bullet (\mathfrak{g} ,\, \mathbb{C}_\delta)$ of the super Lie algebra $\mathfrak{g}$:
\[
\begin{CD}
\bigwedge^0\mathfrak{g}  @<\delta<< \bigwedge^1\mathfrak{g} @<d<< \bigwedge^2\mathfrak{g} @<d<< \bigwedge^3\mathfrak{g} @<d<< \cdots
\end{CD}
\]
\begin{flalign*} 
d(x_1,...,x_n) &= (\sum_{i=1}^n (-1)^{i+1+\alpha_i} \delta(x_i)x_1 \wedge...\wedge \hat{x_i} \wedge...\wedge x_n)&\\ 
&+(\sum_{i<j}(-1)^{i+j+\alpha_i + \alpha_j - |x_i||x_j|} [x_i, x_j] \wedge x_1 \wedge...\wedge \hat{x_i} \wedge...\wedge \hat{x_j} \wedge...\wedge x_n) \nonumber,& 
\end{flalign*}
where $\alpha_1:=0$  and $\alpha_i = |x_i|(|x_1|+...+|x_{i-1}|)~;~i>1$.

Clearly we have a double complex
\[
\xymatrix{ \bigwedge^0\mathfrak{g}   \overset{d}{\underset{0}{\leftrightarrows}} \bigwedge^1\mathfrak{g}   \overset{d}{\underset{0}{\leftrightarrows}}\bigwedge^2\mathfrak{g}   \overset{d}{\underset{0}{\leftrightarrows}}... }
\]
whose total homology is $\bigoplus_{i=*(mod~2)} H_i (\mathfrak{g} ; \mathbb{C}_\delta)$.
On the other hand we have the $(b,B)$ double complex
\[
\xymatrix{ U(\mathfrak{g})^0   \overset{B}{\underset{b}{\leftrightarrows}} U(\mathfrak{g})^1   \overset{B}{\underset{b}{\leftrightarrows}} U(\mathfrak{g})^2  \overset{B}{\underset{b}{\leftrightarrows}}... }
\]
whose total homology is $HP_{(\delta,1)}^*(U(\mathfrak{g}))$.

One checks that the antisymmetrization map $A$ commutes with the
$B$-operator, i.e., $BA = Ad$. Now using Lemma ~\eqref{antsummap} we conclude that the
antisymmetrization map
\[
\begin{CD}
A:\bigwedge^n \mathfrak{g}    @>>>  U(\mathfrak{g})^n ,
\end{CD}
\]
defines a quasi-isomorphism between the above double
complexes. Summarizing everything we have: 

\begin{theorem}
\[HP_{(\delta,1)}^*(U(\mathfrak{g}))=\bigoplus_{i=*(mod~2)} H_i (\mathfrak{g} ; \mathbb{C}_ \delta) .\]
\end{theorem}

%%*************************************************************************************************
\section{Hopf cyclic cohomology  in non-symmetric monoidal categories}\label{CM Hf cyclc in nonsymm}

In Theorem ~\eqref{brdd ver of CM thry} we obtained a braided version
of Connes-Moscovici's Hopf cyclic theory in a symmetric  monoidal
category. This was obtained as a special case of a more general
result in Section ~\eqref{hopf cyclc for H C M}, for braided triples
($H$,\,$C$,\,$M$) in a symmetric monoidal category. In this section
we proceed to eliminate the restrictive symmetry condition $\si^2=id
$. For this, we will directly show that the complex of Theorem
~\eqref{brdd ver of CM thry} remains para-cocyclic in \emph{any}
braided abelian monoidal category $\mathcal{C}$, without any
 symmetry condition on the part of $\mathcal{C}$ (Theorem ~\eqref{brdd ver of CM thry in
nonsymmetric} below).  We shall also indicate  how the symmetry
condition on the braiding is related to  the cyclic condition
$\tau_n^{n+1} = id$ (Theorem ~\eqref{symm and ciclc cndtion} below).
The upshot is that  to obtain a cocyclic object and a Hopf cyclic
cohomology in non-symmetric monoidal categories one must inevitably
restrict to the subcomplex $\text{ker}(1-\tau_n^{n+1})$ (cf.\ formula
~\eqref{pra to cyc}). This procedure can, to some extent, be
generalized to braided triples as in Theorem ~\eqref{main thm}. This
more general case however  needs some additional structure and will
be dealt with elsewhere.

\begin{theorem} \label{brdd ver of CM thry in nonsymmetric} Let  $(H,\,
(\delta,\sigma))$  be  a braided Hopf algebra in a braided abelian
monoidal category $\mathcal{C}$, where $(\delta, \sigma)$ is a BMPI. The
following defines a  para-cocyclic object in $\mathcal{C}$:
 \[C^0(H)=I~~~and~~~C^n(H) = H^n, \,  ~~  n \geq 1,  \]
\[
\delta_i = \left\{ \begin{array}{ll}
 (\eta, 1,1,...,1) & \textrm{ $i=0$}\\
 (1,1,...,1, \underset{ i-th}{\Delta}, 1,1,...1) & \textrm{$1 \leq i \leq n-1$}\\
 (1,1,...,1,\sigma) & \textrm{ $i=n$}
  \end{array} \right.
\]
\[\sigma_i =(1,1,..., \underset{(i+1)-th}{\varepsilon}, 1,1...,1), \,~~0 \leq i \leq n~~\label{degen 0 2}   \nonumber\]
\[
\tau_n = \left\{ \begin{array}{ll}
 id_I & \textrm{ $n=0$}\\
 (m_n)(\Delta^{n-1} \widetilde{S} , 1_{H^{n-1}} , \sigma) & \textrm{ $n \neq 0$}
  \end{array} \right.
\]
Here by $m_n$ we mean, $m_1=m$, and for $n \geq 2$:
\[m_n =m_{H^n} = (\underbrace{m,m,...,m}_{n~times})\mathcal{F}_n(\si), \] where
\[\mathcal{F}_n(\si):=\prod_{j=1}^{n-1}~(1_{H^j}, \underbrace{\si, \si,...,\si}_{n-j~times},1_{H^j}).\]
\end{theorem}

Since the proof includes a very long series of long formulas, we try to visualize some essential steps by braiding diagrams. 
Of course it takes a very large amount of space to show all the steps by diagrams.
 
\begin{proof}  Among all  relations in a para-cocyclic object, only the
relations
 \[\tau_n \sigma_0 = \sigma_n \tau_{n+1}^2,  \quad  \tau_n \sigma_i = \sigma_{i-1} \tau_{n+1}, \]
 and others involving the cyclic operator
  $\tau_n$ are not obvious because the braiding map is involved. Here  we give a detailed proof of  the first  formula in  degree $n=2,$ i.e.
  we prove that $\tau_2 \sigma_0 =\sigma_2 \tau_3^2$.
 In the following, in addition to our previous conventions,  we shall  write $(3 ,\si, 2)$ for $(1,1,1, \si, 1,1)$, $\si_{23}$ for $\si_{H^2 \ot H^3},$ and so
 on. We have:

\begin{align*} 
\sigma_2 &\tau_3^2\\
\overset{~\eqref{diagoftau3}}{=}&(1,1,\varepsilon)(m_3)(\Delta^2 \widetilde{S} , 1,1, \sigma) (m_3)(\Delta^2 \widetilde{S} , 1,1, \sigma)\\
=&(1,1,\varepsilon)(m,m,m)(1,\si,3)(3,\si,1)(2,\si,2)((S,S,\widetilde{S})(1,\si)(\si,1)(1,\si)\Delta^2,1,1,\sigma)\\
~~~~~&(m,m,m)(1,\si,3)(3,\si,1)(2,\si,2)((S,S,\widetilde{S})(1,\si)(\si,1)(1,\si)\Delta^2,1,1,\sigma)\overset{~\eqref{diagoftau3}}{=}
\end{align*}

\[
\xy /r2pc/:,  
(0,0)="L11",
(2,0)="L12",
(3,0)="L13",
  "L11"+(0,.5)*!<0pt,6pt>{H},
  "L12"+(0,.5)*!<0pt,6pt>{H}, 
  "L13"+(0,.5)*!<0pt,6pt>{H},
   "L11"="x1", 
  "x1"+(-1,-.5)="x2",
  "x1"+(0,-.5)="x3",
  "x1"+(1,-.5)="x4",
     "x1";"x3"**\dir{}?(.5)="c", 
   "x1";"x2"**\crv{"x1" & "c" & "x2"+(0,.5)}, 
   "x1";"x3"**\crv{"x1" & "c" & "x3"+(0,.5)}, 
   "x1";"x4"**\crv{"x1" & "c" & "x4"+(0,.5)}, 
"x2"="L21",
"x3"="L22",
"x4"="L23",
 "L12"="x1", 
  "x1"+(0,-.5)="x2",
  "x1";"x2"**\dir{-},  
"x2"="L24",
 "L13"="x1", 
  "x1"+(0,-.5)="x2",
  "x1";"x2"**\dir{-},   
"x2"="L25",
  "L13"+ (1,-.25)="x1", 
  "x1"+(0,-.25)="x2",
  "x1";"x2"**\dir{}?(.5)="c", 
  "x1";"x2"**\dir{-},   
  "c"-(0,0)*!<-3pt,0pt>{\bullet \sigma },  
  "x1"+(0,0)*!<0pt,-2pt>{\circ},  
"x2"="L26",
  "L21"="x1", 
  "x1"+(0,-.5)="x2",
  "x1";"x2"**\dir{-},   
"x2"="L31",
   "L22"="x1",       
     "L23"="x2",  
  "x1"+(0,-.5)="x3",
  "x2"+(0,-.5)="x4",
  "x1";"x2"**\dir{}?(.5)="m", 
  "m" + (0,-.25)="c", 
  "x1";"x4"**\crv{ "x1"+(0,-.25) & "c" & "x4"+(0,.25)},  
  "x2";"c"+<3pt,2pt>="c2"**\crv{"x2"+(0,-.25) & "c2"},
  "x3";"c"+<-3pt,-2pt>="c1"**\crv{"x3"+(0,.25) & "c1"}, 
  "c"+(0,.10)*!<0pt,-4pt>{\si},
"x3"="L32",
"x4"="L33",
"L24"+(0,-0)="L34",
"L25"+(0,-0)="L35",
"L26"+(0,-0)="L36",
   "L31"="x1",       
     "L32"="x2", 
  "x1"+(0,-.5)="x3",
  "x2"+(0,-.5)="x4",
  "x1";"x2"**\dir{}?(.5)="m", 
  "m" + (0,-.25)="c", 
  "x1";"x4"**\crv{ "x1"+(0,-.25) & "c" & "x4"+(0,.25)},  
  "x2";"c"+<3pt,2pt>="c2"**\crv{"x2"+(0,-.25) & "c2"},
  "x3";"c"+<-3pt,-2pt>="c1"**\crv{"x3"+(0,.25) & "c1"}, 
  "c"+(0,.10)*!<0pt,-6pt>{\si},
"x3"="L41",
"x4"="L42",
  "L33"="x1", 
  "x1"+(0,-.5)="x2",
  "x1";"x2"**\dir{-}, 
"x2"="L43",
  "L24"="x1", 
  "x1"+(0,-.5)="x2",
  "x1";"x2"**\dir{-},  
"x2"="L44",
  "L25"="x1", 
  "x1"+(0,-.5)="x2",
  "x1";"x2"**\dir{-},   
"x2"="L45",
  "L26"="x1", 
  "x1"+(0,-.5)="x2",
  "x1";"x2"**\dir{-},  
"x2"="L46",
  "L41"="x1", 
  "x1"+(0,-.5)="x2",
  "x1";"x2"**\dir{-}, 
"x2"="L51",
   "L42"="x1",     
     "L43"="x2", 
  "x1"+(0,-.5)="x3",
  "x2"+(0,-.5)="x4",
  "x1";"x2"**\dir{}?(.5)="m", 
  "m" + (0,-.25)="c", 
  "x1";"x4"**\crv{ "x1"+(0,-.25) & "c" & "x4"+(0,.25)},  
  "x2";"c"+<3pt,2pt>="c2"**\crv{"x2"+(0,-.25) & "c2"},
  "x3";"c"+<-3pt,-2pt>="c1"**\crv{"x3"+(0,.25) & "c1"}, 
  "c"+(0,.10)*!<0pt,-6pt>{\si},
"x3"="L52",
"x4"="L53",
  "L44"="x1", 
  "x1"+(0,-.5)="x2",
  "x1";"x2"**\dir{-}, 
"x2"="L54",
  "L45"="x1", 
  "x1"+(0,-.5)="x2",
  "x1";"x2"**\dir{-},   
"x2"="L55",
  "L46"="x1", 
  "x1"+(0,-.5)="x2",
  "x1";"x2"**\dir{-},   
"x2"="L56",
  "L51"="x1", 
  "x1"+(0,-.5)="x2",
  "x1";"x2"**\dir{}?(.5)="c", 
  "x1";"x2"**\dir{-},  
  "c"-(0,0)*!<-4pt,0pt>{\bullet S },  
"x2"="L61",
  "L52"="x1", 
  "x1"+(0,-.5)="x2",
  "x1";"x2"**\dir{}?(.5)="c", 
  "x1";"x2"**\dir{-},   
  "c"-(0,0)*!<-4pt,0pt>{\bullet S},  
"x2"="L62",
  "L53"="x1", 
  "x1"+(0,-.5)="x2",
  "x1";"x2"**\dir{}?(.5)="c", 
  "x1";"x2"**\dir{-},  
  "c"-(0,0)*!<3pt,0pt>{\widetilde{S} \bullet}, 
"x2"="L63",
  "L54"="x1", 
  "x1"+(0,-.5)="x2",
  "x1";"x2"**\dir{-},  
"x2"="L64",
  "L55"="x1", 
  "x1"+(0,-.5)="x2",
  "x1";"x2"**\dir{-},  
"x2"="L65",
  "L56"="x1", 
  "x1"+(0,-.5)="x2",
  "x1";"x2"**\dir{-},  
"x2"="L66",
"L61"="L21",
"L62"="L22",
"L63"="L23",
  "L64";"L64" + (0,-.5)="L24"**\dir{-},  
  "L65";"L65" + (0,-.5)="L25"**\dir{-}, 
  "L66";"L66" + (0,-.5)="L26"**\dir{-},  
  "L21"="x1", 
  "x1"+(0,-.5)="x2",
  "x1";"x2"**\dir{-},   
"x2"="L31",
  "L22"="x1", 
  "x1"+(0,-.5)="x2",
  "x1";"x2"**\dir{-},   
"x2"="L32",
   "L23"="x1",      
   "L24"="x2",  
  "x1"+(0,-.5)="x3",
  "x2"+(0,-.5)="x4",
  "x1";"x2"**\dir{}?(.5)="m", 
  "m" + (0,-.25)="c", 
  "x1";"x4"**\crv{ "x1"+(0,-.25) & "c" & "x4"+(0,.25)},  
  "x2";"c"+<3pt,2pt>="c2"**\crv{"x2"+(0,-.25) & "c2"},
  "x3";"c"+<-3pt,-2pt>="c1"**\crv{"x3"+(0,.25) & "c1"},   
  "c"+(0,.10)*!<0pt,-6pt>{\si},
"x3"="L33",
"x4"="L34",
  "L25"="x1", 
  "x1"+(0,-.5)="x2",
  "x1";"x2"**\dir{-},  
"x2"="L35",
  "L26"="x1", 
  "x1"+(0,-.5)="x2",
  "x1";"x2"**\dir{-},  
"x2"="L36",
"L31"="L21",
"L32"="L22",
"L33"="L23",
"L34"="L24",
"L35"="L25",
"L36"="L26",
  "L21"="x1", 
  "x1"+(0,-.5)="x2",
  "x1";"x2"**\dir{-},   
"x2"="L31",
  "L22"="x1", 
  "x1"+(0,-.5)="x2",
  "x1";"x2"**\dir{-},   
"x2"="L32",
  "L23"="x1", 
  "x1"+(0,-.5)="x2",
  "x1";"x2"**\dir{-},   
"x2"="L33",
   "L24"="x1",   
   "L25"="x2",  
  "x1"+(0,-.5)="x3",
  "x2"+(0,-.5)="x4",
  "x1";"x2"**\dir{}?(.5)="m", 
  "m" + (0,-.25)="c", 
  "x1";"x4"**\crv{ "x1"+(0,-.25) & "c" & "x4"+(0,.25)},  
  "x2";"c"+<3pt,2pt>="c2"**\crv{"x2"+(0,-.25) & "c2"},
  "x3";"c"+<-3pt,-2pt>="c1"**\crv{"x3"+(0,.25) & "c1"}, 
  "c"+(0,.10)*!<0pt,-6pt>{\si},
"x3"="L34",
"x4"="L35",
  "L26"="x1", 
  "x1"+(0,-.5)="x2",
  "x1";"x2"**\dir{-},   
"x2"="L36",
"L31"="L21",
"L32"="L22",
"L33"="L23",
"L34"="L24",
"L35"="L25",
"L36"="L26",
  "L21"="x1", 
  "x1"+(0,-.5)="x2",
  "x1";"x2"**\dir{-},   
"x2"="L31",
   "L22"="x1",    
   "L23"="x2", 
  "x1"+(0,-.5)="x3",
  "x2"+(0,-.5)="x4",
  "x1";"x2"**\dir{}?(.5)="m", 
  "m" + (0,-.25)="c", 
  "x1";"x4"**\crv{ "x1"+(0,-.25) & "c" & "x4"+(0,.25)},  
  "x2";"c"+<3pt,2pt>="c2"**\crv{"x2"+(0,-.25) & "c2"},
  "x3";"c"+<-3pt,-2pt>="c1"**\crv{"x3"+(0,.25) & "c1"}, 
  "c"+(0,.10)*!<0pt,-4pt>{\si},
"x3"="L32",
"x4"="L33",
  "L24"="x1", 
  "x1"+(0,-.5)="x2",
  "x1";"x2"**\dir{-},   
"x2"="L34",
  "L25"="x1", 
  "x1"+(0,-.5)="x2",
  "x1";"x2"**\dir{-},   
"x2"="L35",
  "L26"="x1", 
  "x1"+(0,-.5)="x2",
  "x1";"x2"**\dir{-},   
"x2"="L36",
"L31"="L21",
"L32"="L22",
"L33"="L23",
    "L34"="L24",
"L35"="L25",
"L36"="L26",
  "L21"="x1", 
  "L22"="x2",  
  "x1";"x2"**\dir{}?(.5)="m", 
  "m" + (0,-.25)="c", 
  "m" + (0,-.5)="x3", 
  "x1";"x3"**\crv{"x1"-(0,.25) & "c" & "x3"},
  "x2";"x3"**\crv{"x2"-(0,.25) & "c"  & "x3"}, 
"x3"= "L31",
  "L23"="x1", 
  "L24"="x2",   
  "x1";"x2"**\dir{}?(.5)="m", 
  "m" + (0,-.25)="c", 
  "m" + (0,-.5)="x3", 
  "x1";"x3"**\crv{"x1"-(0,.25) & "c" & "x3"},
  "x2";"x3"**\crv{"x2"-(0,.25) & "c"  & "x3"}, 
"x3"= "L32",
  "L25"="x1", 
  "L26"="x2",  
  "x1";"x2"**\dir{}?(.5)="m", 
  "m" + (0,-.25)="c", 
  "m" + (0,-.5)="x3", 
  "x1";"x3"**\crv{"x1"-(0,.25) & "c" & "x3"},
  "x2";"x3"**\crv{"x2"-(0,.25) & "c"  & "x3"}, 
"x3"= "L33",
"L31"="L21",
"L32"="L22",
"L33"="L23",
"L21"="L11",
"L22"="L12",
"L23"="L13",
   "L11"="x1", 
  "x1"+(-1,-.5)="x2",
  "x1"+(0,-.5)="x3",
  "x1"+(1,-.5)="x4",
     "x1";"x3"**\dir{}?(.5)="c",
   "x1";"x2"**\crv{"x1" & "c" & "x2"+(0,.5)}, 
   "x1";"x3"**\crv{"x1" & "c" & "x3"+(0,.5)}, 
   "x1";"x4"**\crv{"x1" & "c" & "x4"+(0,.5)}, 
"x2"="L21",
"x3"="L22",
"x4"="L23",
 "L12"="x1", 
  "x1"+(0,-.5)="x2",
  "x1";"x2"**\dir{-},   
"x2"="L24",
 "L13"="x1", 
  "x1"+(0,-.5)="x2",
  "x1";"x2"**\dir{-},   
"x2"="L25",
  "L13"+ (1,-.25)="x1", 
  "x1"+(0,-.25)="x2",
  "x1";"x2"**\dir{}?(.5)="c", 
  "x1";"x2"**\dir{-},   
  "c"-(0,0)*!<-3pt,0pt>{\bullet \sigma },  
  "x1"+(0,0)*!<0pt,-2pt>{\circ},  
"x2"="L26",
  "L21"="x1", 
  "x1"+(0,-.5)="x2",
  "x1";"x2"**\dir{-},   
"x2"="L31",
   "L22"="x1",      
     "L23"="x2",  
  "x1"+(0,-.5)="x3",
  "x2"+(0,-.5)="x4",
  "x1";"x2"**\dir{}?(.5)="m", 
  "m" + (0,-.25)="c", 
  "x1";"x4"**\crv{ "x1"+(0,-.25) & "c" & "x4"+(0,.25)},  
  "x2";"c"+<3pt,2pt>="c2"**\crv{"x2"+(0,-.25) & "c2"},
  "x3";"c"+<-3pt,-2pt>="c1"**\crv{"x3"+(0,.25) & "c1"}, 
  "c"+(0,.10)*!<0pt,-6pt>{\si},
"x3"="L32",
"x4"="L33",
"L24"+(0,-0)="L34",
"L25"+(0,-0)="L35",
"L26"+(0,-0)="L36",
   "L31"="x1",      
     "L32"="x2",  
  "x1"+(0,-.5)="x3",
  "x2"+(0,-.5)="x4",
  "x1";"x2"**\dir{}?(.5)="m", 
  "m" + (0,-.25)="c", 
  "x1";"x4"**\crv{ "x1"+(0,-.25) & "c" & "x4"+(0,.25)},  
  "x2";"c"+<3pt,2pt>="c2"**\crv{"x2"+(0,-.25) & "c2"},
  "x3";"c"+<-3pt,-2pt>="c1"**\crv{"x3"+(0,.25) & "c1"}, 
  "c"+(0,.10)*!<0pt,-6pt>{\si},
"x3"="L41",
"x4"="L42",
  "L33"="x1", 
  "x1"+(0,-.5)="x2",
  "x1";"x2"**\dir{-},   
"x2"="L43",
  "L24"="x1", 
  "x1"+(0,-.5)="x2",
  "x1";"x2"**\dir{-},   
"x2"="L44",
  "L25"="x1", 
  "x1"+(0,-.5)="x2",
  "x1";"x2"**\dir{-},   
"x2"="L45",
  "L26"="x1", 
  "x1"+(0,-.5)="x2",
  "x1";"x2"**\dir{-},   
"x2"="L46",
  "L41"="x1", 
  "x1"+(0,-.5)="x2",
  "x1";"x2"**\dir{-},   
"x2"="L51",
   "L42"="x1",       
     "L43"="x2",  
  "x1"+(0,-.5)="x3",
  "x2"+(0,-.5)="x4",
  "x1";"x2"**\dir{}?(.5)="m", 
  "m" + (0,-.25)="c", 
  "x1";"x4"**\crv{ "x1"+(0,-.25) & "c" & "x4"+(0,.25)},  
  "x2";"c"+<3pt,2pt>="c2"**\crv{"x2"+(0,-.25) & "c2"},
  "x3";"c"+<-3pt,-2pt>="c1"**\crv{"x3"+(0,.25) & "c1"}, 
  "c"+(0,.10)*!<0pt,-6pt>{\si},
"x3"="L52",
"x4"="L53",
  "L44"="x1", 
  "x1"+(0,-.5)="x2",
  "x1";"x2"**\dir{-},   
"x2"="L54",
  "L45"="x1", 
  "x1"+(0,-.5)="x2",
  "x1";"x2"**\dir{-},   
"x2"="L55",
  "L46"="x1", 
  "x1"+(0,-.5)="x2",
  "x1";"x2"**\dir{-},   
"x2"="L56",
  "L51"="x1", 
  "x1"+(0,-.5)="x2",
  "x1";"x2"**\dir{}?(.5)="c", 
  "x1";"x2"**\dir{-},   
  "c"-(0,0)*!<-4pt,0pt>{\bullet S },  
"x2"="L61",
  "L52"="x1", 
  "x1"+(0,-.5)="x2",
  "x1";"x2"**\dir{}?(.5)="c", 
  "x1";"x2"**\dir{-},   
  "c"-(0,0)*!<-4pt,0pt>{\bullet S},  
"x2"="L62",
  "L53"="x1", 
  "x1"+(0,-.5)="x2",
  "x1";"x2"**\dir{}?(.5)="c", 
  "x1";"x2"**\dir{-},   
  "c"-(0,0)*!<3pt,0pt>{\widetilde{S} \bullet},  
"x2"="L63",
  "L54"="x1", 
  "x1"+(0,-.5)="x2",
  "x1";"x2"**\dir{-},  
"x2"="L64",
  "L55"="x1", 
  "x1"+(0,-.5)="x2",
  "x1";"x2"**\dir{-},   
"x2"="L65",
  "L56"="x1", 
  "x1"+(0,-.5)="x2",
  "x1";"x2"**\dir{-},   
"x2"="L66",
"L61"="L21",
"L62"="L22",
"L63"="L23",
  "L64";"L64" + (0,-.5)="L24"**\dir{-},  
  "L65";"L65" + (0,-.5)="L25"**\dir{-},  
  "L66";"L66" + (0,-.5)="L26"**\dir{-},  
  "L21"="x1", 
  "x1"+(0,-.5)="x2",
  "x1";"x2"**\dir{-},   
"x2"="L31",
  "L22"="x1", 
  "x1"+(0,-.5)="x2",
  "x1";"x2"**\dir{-},   
"x2"="L32",
   "L23"="x1",      
   "L24"="x2",  
  "x1"+(0,-.5)="x3",
  "x2"+(0,-.5)="x4",
  "x1";"x2"**\dir{}?(.5)="m", 
  "m" + (0,-.25)="c", 
  "x1";"x4"**\crv{ "x1"+(0,-.25) & "c" & "x4"+(0,.25)},  
  "x2";"c"+<3pt,2pt>="c2"**\crv{"x2"+(0,-.25) & "c2"},
  "x3";"c"+<-3pt,-2pt>="c1"**\crv{"x3"+(0,.25) & "c1"},   
  "c"+(0,.10)*!<0pt,-6pt>{\si},
"x3"="L33",
"x4"="L34",
  "L25"="x1", 
  "x1"+(0,-.5)="x2",
  "x1";"x2"**\dir{-},   
"x2"="L35",
  "L26"="x1", 
  "x1"+(0,-.5)="x2",
  "x1";"x2"**\dir{-},   
"x2"="L36",
"L31"="L21",
"L32"="L22",
"L33"="L23",
"L34"="L24",
"L35"="L25",
"L36"="L26",
  "L21"="x1", 
  "x1"+(0,-.5)="x2",
  "x1";"x2"**\dir{-},   
"x2"="L31",
  "L22"="x1", 
  "x1"+(0,-.5)="x2",
  "x1";"x2"**\dir{-},  
"x2"="L32",
  "L23"="x1", 
  "x1"+(0,-.5)="x2",
  "x1";"x2"**\dir{-},  
"x2"="L33",
   "L24"="x1",      
   "L25"="x2",  
  "x1"+(0,-.5)="x3",
  "x2"+(0,-.5)="x4",
  "x1";"x2"**\dir{}?(.5)="m", 
  "m" + (0,-.25)="c", 
  "x1";"x4"**\crv{ "x1"+(0,-.25) & "c" & "x4"+(0,.25)},  
  "x2";"c"+<3pt,2pt>="c2"**\crv{"x2"+(0,-.25) & "c2"},
  "x3";"c"+<-3pt,-2pt>="c1"**\crv{"x3"+(0,.25) & "c1"}, 
  "c"+(0,.10)*!<0pt,-6pt>{\si},
"x3"="L34",
"x4"="L35",
  "L26"="x1", 
  "x1"+(0,-.5)="x2",
  "x1";"x2"**\dir{-},   
"x2"="L36",
"L31"="L21",
"L32"="L22",
"L33"="L23",
"L34"="L24",
"L35"="L25",
"L36"="L26",
  "L21"="x1", 
  "x1"+(0,-.5)="x2",
  "x1";"x2"**\dir{-},   
"x2"="L31",
   "L22"="x1",      
   "L23"="x2",  
  "x1"+(0,-.5)="x3",
  "x2"+(0,-.5)="x4",
  "x1";"x2"**\dir{}?(.5)="m", 
  "m" + (0,-.25)="c", 
  "x1";"x4"**\crv{ "x1"+(0,-.25) & "c" & "x4"+(0,.25)},  
  "x2";"c"+<3pt,2pt>="c2"**\crv{"x2"+(0,-.25) & "c2"},
  "x3";"c"+<-3pt,-2pt>="c1"**\crv{"x3"+(0,.25) & "c1"}, 
  "c"+(0,.10)*!<0pt,-6pt>{\si},
"x3"="L32",
"x4"="L33",
  "L24"="x1", 
  "x1"+(0,-.5)="x2",
  "x1";"x2"**\dir{-},  
"x2"="L34",
  "L25"="x1", 
  "x1"+(0,-.5)="x2",
  "x1";"x2"**\dir{-},  
"x2"="L35",
  "L26"="x1", 
  "x1"+(0,-.5)="x2",
  "x1";"x2"**\dir{-},   
"x2"="L36",
"L31"="L21",
"L32"="L22",
"L33"="L23",
    "L34"="L24",
"L35"="L25",
"L36"="L26",
  "L21"="x1", 
  "L22"="x2",   
  "x1";"x2"**\dir{}?(.5)="m", 
  "m" + (0,-.25)="c", 
  "m" + (0,-.5)="x3", 
  "x1";"x3"**\crv{"x1"-(0,.25) & "c" & "x3"},
  "x2";"x3"**\crv{"x2"-(0,.25) & "c"  & "x3"}, 
"x3"= "L31",
  "L23"="x1", 
  "L24"="x2",   
  "x1";"x2"**\dir{}?(.5)="m", 
  "m" + (0,-.25)="c", 
  "m" + (0,-.5)="x3", 
  "x1";"x3"**\crv{"x1"-(0,.25) & "c" & "x3"},
  "x2";"x3"**\crv{"x2"-(0,.25) & "c"  & "x3"}, 
"x3"= "L32",
  "L25"="x1", 
  "L26"="x2",   
  "x1";"x2"**\dir{}?(.5)="m", 
  "m" + (0,-.25)="c", 
  "m" + (0,-.5)="x3", 
  "x1";"x3"**\crv{"x1"-(0,.25) & "c" & "x3"},
  "x2";"x3"**\crv{"x2"-(0,.25) & "c"  & "x3"}, 
"x3"= "L33",
"L31"="L21",
"L32"="L22",
"L33"="L23",
  "L21"-(0,.5)*!<0pt,-6pt>{H},  
  "L22"-(0,.5)*!<0pt,-6pt>{H}, 
   "L23"="x1", 
  "x1"+(0,-.5)="x2",
  "x1";"x2"**\dir{}?(.5)="c", 
  "x1";"x2"**\dir{-},  
  "c"-(0,0)*!<-3pt,0pt>{\bullet \varepsilon },  
  "x2"+(0,0)*!<0pt,2pt>{\circ},  
 \endxy 
\]

\begin{align*} 
\overset{~\eqref{naturalityofsi}}{=}&(m,m,\varepsilon m)(1,\si,3)(3,\si,1)(2,\si,2) (S,S,\widetilde{S},1,1,\sigma)\\
~~~~~&(1,\si,2)(\si,3)(1,\si,2)\\
~~~~~&(\Delta^2 m,m,m) (S,1,S,1,\widetilde{S},\sigma)(1,\si,2)(3,\si)(2,\si,1)\\
~~~~~&(1,\si,2)(\si,3)(1,\si,2)(\Delta^2,2)=
\end{align*}

\[
\xy /r2.5pc/:,  
(0,0)="L11",
(2,0)="L12",
(3,0)="L13",
  "L11"+(0,.5)*!<0pt,6pt>{H},
  "L12"+(0,.5)*!<0pt,6pt>{H}, 
  "L13"+(0,.5)*!<0pt,6pt>{H},
   "L11"="x1", 
   "x1"+(-1,-.5)="x2",
  "x1"+(0,-.5)="x3",
  "x1"+(1,-.5)="x4",
  "x1"+(0,-.25)="c", 
   "x1";"x2"**\crv{"x1" & "c" & "x2"+(0,.25)}, 
   "x1";"x3"**\crv{"x1" & "c" & "x3"+(0,.25)}, 
   "x1";"x4"**\crv{"x1" & "c" & "x4"+(0,.25)}, 
"x2"="L21",
"x3"="L22",
"x4"="L23",
 "L12"="x1",
  "x1"+(0,-.5)="x2",
  "x1";"x2"**\dir{-},   
"x2"="L24",
 "L13"="x1", 
  "x1"+(0,-.5)="x2",
  "x1";"x2"**\dir{-},   
"x2"="L25",
  "L21"="x1", 
  "x1"+(0,-.5)="x2",
  "x1";"x2"**\dir{-},   
"x2"= "L31",
   "L22"="x1",       
   "L23"="x2",  
  "x1"+(0,-.5)="x3",
  "x2"+(0,-.5)="x4",
  "x1";"x2"**\dir{}?(.5)="m", 
  "m" + (0,-.25)="c", 
  "x1";"x4"**\crv{ "x1"+(0,-.25) & "c" & "x4"+(0,.25)},  
  "x2";"c"+<3pt,2pt>="c2"**\crv{"x2"+(0,-.25) & "c2"},
  "x3";"c"+<-3pt,-2pt>="c1"**\crv{"x3"+(0,.25) & "c1"}, 
  "c"+(0,.10)*!<0pt,-4pt>{\si}, 
"x3"= "L32",
"x4"= "L33",
  "L24"="x1", 
  "x1"+(0,-.5)="x2",
  "x1";"x2"**\dir{-},   
"x2"= "L34",
  "L25"="x1", 
  "x1"+(0,-.5)="x2",
  "x1";"x2"**\dir{-},  
"x2"= "L35",
"L31"="L21",
"L32"="L22",
"L33"="L23",
"L34"="L24",
"L35"="L25",
   "L21" ="x1",     
   "L22"="x2",  
  "x1"+(-0,-.5)="x3",
  "x2"+(0,-.5)="x4",
  "x1";"x2"**\dir{}?(.5)="m", 
  "m" + (0,-.25)="c", 
  "x1";"x4"**\crv{ "x1"+(0,-.25) & "c" & "x4"+(0,.25)},  
  "x2";"c"+<3pt,2pt>="c2"**\crv{"x2"+(0,-.25) & "c2"},
  "x3";"c"+<-3pt,-2pt>="c1"**\crv{"x3"+(0,.25) & "c1"}, 
  "c"+(0,.10)*!<0pt,-6pt>{\si}, 
"x3"= "L31",
"x4"= "L32",
  "L23"="x1", 
  "x1"+(0,-.5)="x2",
  "x1";"x2"**\dir{-},  
"x2"= "L33",
  "L24"="x1", 
  "x1"+(0,-.5)="x2",
  "x1";"x2"**\dir{-},   
"x2"= "L34",
  "L25"="x1", 
  "x1"+(0,-.5)="x2",
  "x1";"x2"**\dir{-},   
"x2"= "L35",
"L31"="L21",
"L32"="L22",
"L33"="L23",
"L34"="L24",
"L35"="L25",
  "L21"="x1", 
  "x1"+(0,-.5)="x2",
  "x1";"x2"**\dir{-},  
"x2"= "L31",
   "L22"="x1",       
   "L23"="x2",  
  "x1"+(0,-.5)="x3",
  "x2"+(0,-.5)="x4",
  "x1";"x2"**\dir{}?(.5)="m", 
  "m" + (0,-.25)="c", 
  "x1";"x4"**\crv{ "x1"+(0,-.25) & "c" & "x4"+(0,.25)},  
  "x2";"c"+<3pt,2pt>="c2"**\crv{"x2"+(0,-.25) & "c2"},
  "x3";"c"+<-3pt,-2pt>="c1"**\crv{"x3"+(0,.25) & "c1"}, 
  "c"+(0,.10)*!<0pt,-6pt>{\si},  
"x3"= "L32",
"x4"= "L33",
  "L24"="x1", 
  "x1"+(0,-.5)="x2",
  "x1";"x2"**\dir{-},   
"x2"= "L34",
  "L25"="x1", 
  "x1"+(0,-.5)="x2",
  "x1";"x2"**\dir{-},  
"x2"= "L35",
"L31"="L21",
"L32"="L22",
"L33"="L23",
"L34"="L24",
"L35"="L25",
  "L21"="x1", 
  "x1"+(0,-.5)="x2",
  "x1";"x2"**\dir{-},   
"x2"= "L31",
  "L22"="x1", 
  "x1"+(0,-.5)="x2",
  "x1";"x2"**\dir{-},   
"x2"= "L32",
   "L23"="x1",      
   "L24"="x2",  
  "x1"+(0,-.5)="x3",
  "x2"+(0,-.5)="x4",
  "x1";"x2"**\dir{}?(.5)="m", 
  "m" + (0,-.25)="c", 
  "x1";"x4"**\crv{ "x1"+(0,-.25) & "c" & "x4"+(0,.25)},  
  "x2";"c"+<3pt,2pt>="c2"**\crv{"x2"+(0,-.25) & "c2"},
  "x3";"c"+<-3pt,-2pt>="c1"**\crv{"x3"+(0,.25) & "c1"}, 
  "c"+(0,.10)*!<0pt,-6pt>{\si},  
"x3"= "L33",
"x4"= "L34",
  "L25"="x1", 
  "x1"+(0,-.5)="x2",
  "x1";"x2"**\dir{-},   
"x2"= "L35",
"L31"="L21",
"L32"="L22",
"L33"="L23",
"L34"="L24",
"L35"="L25",
  "L21"="x1",
  "x1"+(0,-.5)="x2",
  "x1";"x2"**\dir{-},   
"x2"= "L31",
  "L22"="x1",
  "x1"+(0,-.5)="x2",
  "x1";"x2"**\dir{-},   
"x2"= "L32",
  "L23"="x1", 
  "x1"+(0,-.5)="x2",
  "x1";"x2"**\dir{-},  
"x2"= "L33",
   "L24"="x1",     
   "L25"="x2",  
  "x1"+(0,-.5)="x3",
  "x2"+(0,-.5)="x4",
  "x1";"x2"**\dir{}?(.5)="m", 
  "m" + (0,-.25)="c", 
  "x1";"x4"**\crv{ "x1"+(0,-.25) & "c" & "x4"+(0,.25)},  
  "x2";"c"+<3pt,2pt>="c2"**\crv{"x2"+(0,-.25) & "c2"},
  "x3";"c"+<-3pt,-2pt>="c1"**\crv{"x3"+(0,.25) & "c1"}, 
  "c"+(0,.10)*!<0pt,-6pt>{\si},  
"x3"= "L34",
"x4"= "L35",
"L31"="L21",
"L32"="L22",
"L33"="L23",
"L34"="L24",
"L35"="L25",
  "L21"="x1", 
  "x1"+(0,-.5)="x2",
  "x1";"x2"**\dir{-},   
"x2"= "L31",
   "L22"="x1",       
   "L23"="x2",  
  "x1"+(0,-.5)="x3",
  "x2"+(0,-.5)="x4",
  "x1";"x2"**\dir{}?(.5)="m", 
  "m" + (0,-.25)="c", 
  "x1";"x4"**\crv{ "x1"+(0,-.25) & "c" & "x4"+(0,.25)},  
  "x2";"c"+<3pt,2pt>="c2"**\crv{"x2"+(0,-.25) & "c2"},
  "x3";"c"+<-3pt,-2pt>="c1"**\crv{"x3"+(0,.25) & "c1"}, 
  "c"+(0,.10)*!<0pt,-6pt>{\si},  
"x3"= "L32",
"x4"= "L33",
  "L24"="x1", 
  "x1"+(0,-.5)="x2",
  "x1";"x2"**\dir{-},   
"x2"= "L34",
  "L25"="x1", 
  "x1"+(0,-.5)="x2",
  "x1";"x2"**\dir{-},   
"x2"= "L35",
"L31"="L21",
"L32"="L22",
"L33"="L23",
"L34"="L24",
"L35"="L25",
  "L21" ="x1", 
  "x1"+(0,-.5)="x2",
  "x1";"x2"**\dir{}?(.5)="c", 
  "x1";"x2"**\dir{-},   
  "c"-(0,0)*!<-4pt,0pt>{\bullet S},  
"x2"= "L31",
  "L22"="x1", 
  "x1"+(0,-.5)="x2",
  "x1";"x2"**\dir{-},   
"x2"= "L32",
  "L23" ="x1", 
  "x1"+(0,-.5)="x2",
  "x1";"x2"**\dir{}?(.5)="c", 
  "x1";"x2"**\dir{-},   
  "c"-(0,0)*!<-4pt,0pt>{\bullet S},  
"x2"= "L33",
  "L24"="x1", 
  "x1"+(0,-.5)="x2",
  "x1";"x2"**\dir{-},   
"x2"= "L34",
  "L25" ="x1", 
  "x1"+(0,-.5)="x2",
  "x1";"x2"**\dir{}?(.5)="c", 
  "x1";"x2"**\dir{-},   
  "c"-(0,0)*!<3pt,0pt>{\widetilde{S} \bullet}, 
"x2"= "L35",
    "L25"+(1,-.25) ="x1", 
  "x1"+(0,-.25)="x2",
  "x1";"x2"**\dir{}?(.5)="c", 
  "x1";"x2"**\dir{-},  
  "c"-(0,0)*!<-3pt,0pt>{\bullet \sigma },  
  "x1"+(0,0)*!<0pt,-2pt>{\circ},  
"x2"= "L36",
"L31"="L21",
"L32"="L22",
"L33"="L23",
"L34"="L24",
"L35"="L25",
"L36"="L26",
   "L21" ="x1", 
    "L22"="x2",   
  "x1";"x2"**\dir{}?(.5)="m", 
  "m" + (0,-.25)="c", 
  "m" + (0,-.5)="x3", 
  "x1";"x3"**\crv{"x1"-(0,.25) & "c" & "x3"},
  "x2";"x3"**\crv{"x2"-(0,.25) & "c"  & "x3"}, 
"x3"="L31",
   "L23" ="x1", 
    "L24"="x2",   
  "x1";"x2"**\dir{}?(.5)="m", 
  "m" + (0,-.25)="c", 
  "m" + (0,-.5)="x3", 
  "x1";"x3"**\crv{"x1"-(0,.25) & "c" & "x3"},
  "x2";"x3"**\crv{"x2"-(0,.25) & "c"  & "x3"}, 
"x3"="L34",
   "L25" ="x1", 
    "L26"="x2",   
  "x1";"x2"**\dir{}?(.5)="m", 
  "m" + (0,-.25)="c", 
  "m" + (0,-.5)="x3", 
  "x1";"x3"**\crv{"x1"-(0,.25) & "c" & "x3"},
  "x2";"x3"**\crv{"x2"-(0,.25) & "c"  & "x3"}, 
"x3"="L35",
   "L31"="x1", 
   "x1"+(-,-.5)="x2",
  "x1"+(0,-.5)="x3",
  "x1"+(1,-.5)="x4",
  "x1"+(0,-.25)="c",  
   "x1";"x2"**\crv{"x1" & "c" & "x2"+(0,.25)}, 
   "x1";"x3"**\crv{"x1" & "c" & "x3"+(0,.25)}, 
   "x1";"x4"**\crv{"x1" & "c" & "x4"+(0,.25)}, 
"x2"="L31",
"x3"="L32",
"x4"="L33",
  "L34";"L34" + (0,-.5) ="L34" **\dir{-}, 
  "L35";"L35" + (0,-.5) ="L35" **\dir{-}, 
"L31"="L21",
"L32"="L22",
"L33"="L23",
"L34"="L24",
"L35"="L25",
  "L21"="x1", 
  "x1"+(0,-.5)="x2",
  "x1";"x2"**\dir{-},   
"x2"= "L31",
   "L22"="x1",      
   "L23"="x2",  
  "x1"+(0,-.5)="x3",
  "x2"+(0,-.5)="x4",
  "x1";"x2"**\dir{}?(.5)="m", 
  "m" + (0,-.25)="c", 
  "x1";"x4"**\crv{ "x1"+(0,-.25) & "c" & "x4"+(0,.25)},  
  "x2";"c"+<3pt,2pt>="c2"**\crv{"x2"+(0,-.25) & "c2"},
  "x3";"c"+<-3pt,-2pt>="c1"**\crv{"x3"+(0,.25) & "c1"}, 
  "c"+(0,.10)*!<0pt,-6pt>{\si},  
"x3"= "L32",
"x4"= "L33",
  "L24"="x1", 
  "x1"+(0,-.5)="x2",
  "x1";"x2"**\dir{-},   
"x2"= "L34",
  "L25"="x1", 
  "x1"+(0,-.5)="x2",
  "x1";"x2"**\dir{-},   
"x2"= "L35",
"L31"="L21",
"L32"="L22",
"L33"="L23",
"L34"="L24",
"L35"="L25",
   "L21" ="x1",      
   "L22"="x2",  
  "x1"+(0,-.5)="x3",
  "x2"+(0,-.5)="x4",
  "x1";"x2"**\dir{}?(.5)="m", 
  "m" + (0,-.25)="c", 
  "x1";"x4"**\crv{ "x1"+(0,-.25) & "c" & "x4"+(0,.25)},  
  "x2";"c"+<3pt,2pt>="c2"**\crv{"x2"+(0,-.25) & "c2"},
  "x3";"c"+<-3pt,-2pt>="c1"**\crv{"x3"+(0,.25) & "c1"}, 
  "c"+(0,.10)*!<0pt,-6pt>{\si},  
"x3"= "L31",
"x4"= "L32",
  "L23"="x1", 
  "x1"+(0,-.5)="x2",
  "x1";"x2"**\dir{-},  
"x2"= "L33",
  "L24"="x1", 
  "x1"+(0,-.5)="x2",
  "x1";"x2"**\dir{-},   
"x2"= "L34",
  "L25"="x1", 
  "x1"+(0,-.5)="x2",
  "x1";"x2"**\dir{-},   
"x2"= "L35",
"L31"="L21",
"L32"="L22",
"L33"="L23",
"L34"="L24",
"L35"="L25",
  "L21"="x1", 
  "x1"+(0,-.5)="x2",
  "x1";"x2"**\dir{-},   
"x2"= "L31",
   "L22"="x1",       
   "L23"="x2",  
  "x1"+(0,-.5)="x3",
  "x2"+(0,-.5)="x4",
  "x1";"x2"**\dir{}?(.5)="m", 
  "m" + (0,-.25)="c", 
  "x1";"x4"**\crv{ "x1"+(0,-.25) & "c" & "x4"+(0,.25)},  
  "x2";"c"+<3pt,2pt>="c2"**\crv{"x2"+(0,-.25) & "c2"},
  "x3";"c"+<-3pt,-2pt>="c1"**\crv{"x3"+(0,.25) & "c1"}, 
  "c"+(0,.10)*!<0pt,-6pt>{\si},  
"x3"= "L32",
"x4"= "L33",
  "L24"="x1", 
  "x1"+(0,-.5)="x2",
  "x1";"x2"**\dir{-},   
"x2"= "L34",
  "L25"="x1", 
  "x1"+(0,-.5)="x2",
  "x1";"x2"**\dir{-},   
"x2"= "L35",
"L31"="L21",
"L32"="L22",
"L33"="L23",
"L34"="L24",
"L35"="L25",
  "L21" ="x1", 
  "x1"+(0,-.5)="x2",
  "x1";"x2"**\dir{}?(.5)="c", 
  "x1";"x2"**\dir{-},   
  "c"-(0,0)*!<-4pt,0pt>{\bullet S},  
"x2"= "L31",
  "L22" ="x1", 
  "x1"+(0,-.5)="x2",
  "x1";"x2"**\dir{}?(.5)="c", 
  "x1";"x2"**\dir{-},   
  "c"-(0,0)*!<-4pt,0pt>{\bullet S},  
"x2"= "L32",
  "L23" ="x1", 
  "x1"+(0,-.5)="x2",
  "x1";"x2"**\dir{}?(.5)="c", 
  "x1";"x2"**\dir{-},   
  "c"-(0,0)*!<3pt,0pt>{\widetilde{S} \bullet},  
"x2"= "L33",
  "L24"="x1", 
  "x1"+(0,-.5)="x2",
  "x1";"x2"**\dir{-},   
"x2"= "L34",
  "L25"="x1", 
  "x1"+(0,-.5)="x2",
  "x1";"x2"**\dir{-},   
"x2"= "L35",
  "L25"+ (1,-.25)="x1", 
  "x1"+(0,-.25)="x2",
  "x1";"x2"**\dir{}?(.5)="c", 
  "x1";"x2"**\dir{-},   
  "c"-(0,0)*!<-3pt,0pt>{\bullet \sigma }, 
  "x1"+(0,0)*!<0pt,-2pt>{\circ},  
"x2"= "L36",
"L31"="L21",
"L32"="L22",
"L33"="L23",
"L34"="L24",
"L35"="L25",
"L36"="L26",
  "L21"="x1", 
  "x1"+(0,-.5)="x2",
  "x1";"x2"**\dir{-},  
"x2"= "L31",
  "L22"="x1", 
  "x1"+(0,-.5)="x2",
  "x1";"x2"**\dir{-},   
"x2"= "L32",
   "L23"="x1",      
   "L24"="x2",  
  "x1"+(0,-.5)="x3",
  "x2"+(0,-.5)="x4",
  "x1";"x2"**\dir{}?(.5)="m", 
  "m" + (0,-.25)="c", 
  "x1";"x4"**\crv{ "x1"+(0,-.25) & "c" & "x4"+(0,.25)},  
  "x2";"c"+<3pt,2pt>="c2"**\crv{"x2"+(0,-.25) & "c2"},
  "x3";"c"+<-3pt,-2pt>="c1"**\crv{"x3"+(0,.25) & "c1"}, 
  "c"+(0,.10)*!<0pt,-6pt>{\si},  
"x3"= "L33",
"x4"= "L34",
  "L25"="x1", 
  "x1"+(0,-.5)="x2",
  "x1";"x2"**\dir{-}, 
"x2"= "L35",
  "L26"="x1", 
  "x1"+(0,-.5)="x2",
  "x1";"x2"**\dir{-},  
"x2"= "L36",
"L31"="L21",
"L32"="L22",
"L33"="L23",
"L34"="L24",
"L35"="L25",
"L36"="L26",
  "L21"="x1", 
  "x1"+(0,-.5)="x2",
  "x1";"x2"**\dir{-},  
"x2"= "L31",
  "L22"="x1", 
  "x1"+(0,-.5)="x2",
  "x1";"x2"**\dir{-},   
"x2"= "L32",
  "L23"="x1", 
  "x1"+(0,-.5)="x2",
  "x1";"x2"**\dir{-},   
"x2"= "L33",
   "L24"="x1",   
   "L25"="x2",  
  "x1"+(0,-.5)="x3",
  "x2"+(0,-.5)="x4",
  "x1";"x2"**\dir{}?(.5)="m", 
  "m" + (0,-.25)="c", 
  "x1";"x4"**\crv{ "x1"+(0,-.25) & "c" & "x4"+(0,.25)},  
  "x2";"c"+<3pt,2pt>="c2"**\crv{"x2"+(0,-.25) & "c2"},
  "x3";"c"+<-3pt,-2pt>="c1"**\crv{"x3"+(0,.25) & "c1"}, 
  "c"+(0,.10)*!<0pt,-6pt>{\si}, 
"x3"= "L34",
"x4"= "L35",
  "L26"="x1", 
  "x1"+(0,-.5)="x2",
  "x1";"x2"**\dir{-},   
"x2"= "L36",
"L31"="L21",
"L32"="L22",
"L33"="L23",
"L34"="L24",
"L35"="L25",
"L36"="L26",
  "L21"="x1", 
  "x1"+(0,-.5)="x2",
  "x1";"x2"**\dir{-},   
"x2"= "L31",
   "L22"="x1",      
   "L23"="x2",  
  "x1"+(0,-.5)="x3",
  "x2"+(0,-.5)="x4",
  "x1";"x2"**\dir{}?(.5)="m", 
  "m" + (0,-.25)="c", 
  "x1";"x4"**\crv{ "x1"+(0,-.25) & "c" & "x4"+(0,.25)},  
  "x2";"c"+<3pt,2pt>="c2"**\crv{"x2"+(0,-.25) & "c2"},
  "x3";"c"+<-3pt,-2pt>="c1"**\crv{"x3"+(0,.25) & "c1"}, 
  "c"+(0,.10)*!<0pt,-6pt>{\si}, 
"x3"= "L32",
"x4"= "L33",
  "L24"="x1", 
  "x1"+(0,-.5)="x2",
  "x1";"x2"**\dir{-},  
"x2"= "L34",
  "L25"="x1", 
  "x1"+(0,-.5)="x2",
  "x1";"x2"**\dir{-},   
"x2"= "L35",
  "L26"="x1", 
  "x1"+(0,-.5)="x2",
  "x1";"x2"**\dir{-},   
"x2"= "L36",
"L31"="L21",
"L32"="L22",
"L33"="L23",
"L34"="L24",
"L35"="L25",
"L36"="L26",
   "L21" ="x1", 
    "L22"="x2",  
  "x1";"x2"**\dir{}?(.5)="m", 
  "m" + (0,-.25)="c", 
  "m" + (0,-.5)="x3", 
  "x1";"x3"**\crv{"x1"-(0,.25) & "c" & "x3"},
  "x2";"x3"**\crv{"x2"-(0,.25) & "c"  & "x3"}, 
"x3"="L31",
   "L23" ="x1", 
    "L24"="x2",  
  "x1";"x2"**\dir{}?(.5)="m", 
  "m" + (0,-.25)="c", 
  "m" + (0,-.5)="x3", 
  "x1";"x3"**\crv{"x1"-(0,.25) & "c" & "x3"},
  "x2";"x3"**\crv{"x2"-(0,.25) & "c"  & "x3"}, 
"x3"="L32",
   "L25" ="x1", 
    "L26"="x2",  
  "x1";"x2"**\dir{}?(.5)="m", 
  "m" + (0,-.25)="c", 
  "m" + (0,-.5)="x3", 
  "x1";"x3"**\crv{"x1"-(0,.25) & "c" & "x3"},
  "x2";"x3"**\crv{"x2"-(0,.25) & "c"  & "x3"}, 
"x3"="L33",
   "L33"="x1", 
  "x1"+(0,-.25)="x2",
  "x1";"x2"**\dir{}?(.5)="c", 
  "x1";"x2"**\dir{-},   
  "c"-(0,0)*!<-3pt,0pt>{\bullet \varepsilon },  
  "x2"+(0,0)*!<0pt,2pt>{\circ},  
"x2"="L33",
  "L31";"L31" +(0,-.5)= "L31" **\dir{-}, 
  "L32";"L32" +(0,-.5)= "L32" **\dir{-},  
 "L31"-(0,.5)*!<0pt,-6pt>{H},   
 "L32"-(0,.5)*!<0pt,-6pt>{H}, 
\endxy
\]

\begin{align*} 
\overset{~\eqref{naturalityofsi}, ~\eqref{Del2 m}}{=}&(m,m,\varepsilon, \varepsilon)(S,1,S,1,\widetilde{S},\sigma)(1,\si,2)(3,\si)(2,\si,1)(1,\si,2)(\si,3)(1,\si,2)\\
~~~~~&((m,m,m)(1,\si,3)(3,\si,1)(2,\si,2)(\Delta^2 ,\Delta^2),m,m) (S,1,S,1,\widetilde{S},\sigma)\\
~~~~~&(1,\si,2)(3,\si)(2,\si,1)(1,\si,2)(\si,3)(1,\si,2)(\Delta^2,2) =
\end{align*}

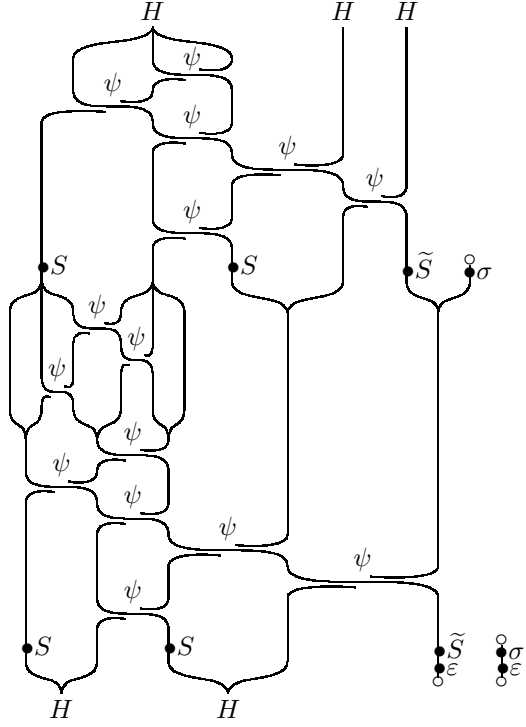
\[
\xy /r2pc/:,  
(0,0)="L11",
(3,0)="L12",
(4,0)="L13",
  "L11"+(0,.5)*!<0pt,6pt>{H},
  "L12"+(0,.5)*!<0pt,6pt>{H}, 
  "L13"+(0,.5)*!<0pt,6pt>{H},
   "L11"="x1", 
   "x1"+(-1.25,-.5)="x2",
  "x1"+(0,-.5)="x3",
  "x1"+(1.25,-.5)="x4",
  "x1"+(0,-.25)="c",  
   "x1";"x2"**\crv{"x1" & "c" & "x2"+(0,.25)}, 
   "x1";"x3"**\crv{"x1" & "c" & "x3"+(0,.25)}, 
   "x1";"x4"**\crv{"x1" & "c" & "x4"+(0,.25)}, 
"x2"="L21",
"x3"="L22",
"x4"="L23",
 "L12"="x1", 
  "x1"+(0,-.5)="x2",
  "x1";"x2"**\dir{-},   
"x2"="L24",
 "L13"="x1", 
  "x1"+(0,-.5)="x2",
  "x1";"x2"**\dir{-},  
"x2"="L25",
  "L21"="x1", 
  "x1"+(0,-.5)="x2",
  "x1";"x2"**\dir{-},   
"x2"= "L31",
   "L22"="x1",     
   "L23"="x2", 
  "x1"+(0,-.5)="x3",
  "x2"+(0,-.5)="x4",
  "x1";"x2"**\dir{}?(.5)="m", 
  "m" + (0,-.25)="c", 
  "x1";"x4"**\crv{ "x1"+(0,-.25) & "c" & "x4"+(0,.25)},  
  "x2";"c"+<3pt,2pt>="c2"**\crv{"x2"+(0,-.25) & "c2"},
  "x3";"c"+<-3pt,-2pt>="c1"**\crv{"x3"+(0,.25) & "c1"}, 
  "c"+(0,.10)*!<0pt,-4pt>{\si},  
"x3"= "L32",
"x4"= "L33",
  "L24"="x1", 
  "x1"+(0,-.5)="x2",
  "x1";"x2"**\dir{-},  
"x2"= "L34",
  "L25"="x1", 
  "x1"+(0,-.5)="x2",
  "x1";"x2"**\dir{-},  
"x2"= "L35",
"L31"="L21",
"L32"="L22",
"L33"="L23",
"L34"="L24",
"L35"="L25",
   "L21" ="x1",    
   "L22"="x2", 
  "x1"+(-.5,-.5)="x3",
  "x2"+(0,-.5)="x4",
  "x1";"x2"**\dir{}?(.5)="m", 
  "m" + (0,-.25)="c", 
  "x1";"x4"**\crv{ "x1"+(0,-.25) & "c" & "x4"+(0,.25)},  
  "x2";"c"+<3pt,2pt>="c2"**\crv{"x2"+(0,-.25) & "c2"},
  "x3";"c"+<-3pt,-2pt>="c1"**\crv{"x3"+(0,.25) & "c1"}, 
  "c"+(0,.10)*!<0pt,-6pt>{\si}, 
"x3"= "L31",
"x4"= "L32",
  "L23"="x1", 
  "x1"+(0,-.5)="x2",
  "x1";"x2"**\dir{-},   
"x2"= "L33",
  "L24"="x1", 
  "x1"+(0,-.5)="x2",
  "x1";"x2"**\dir{-},  
"x2"= "L34",
  "L25"="x1", 
  "x1"+(0,-.5)="x2",
  "x1";"x2"**\dir{-},   
"x2"= "L35",
"L31"="L21",
"L32"="L22",
"L33"="L23",
"L34"="L24",
"L35"="L25",
  "L21"="x1", 
  "x1"+(0,-.5)="x2",
  "x1";"x2"**\dir{-},   
"x2"= "L31",
   "L22"="x1",    
   "L23"="x2",  
  "x1"+(0,-.5)="x3",
  "x2"+(0,-.5)="x4",
  "x1";"x2"**\dir{}?(.5)="m", 
  "m" + (0,-.25)="c", 
  "x1";"x4"**\crv{ "x1"+(0,-.25) & "c" & "x4"+(0,.25)},  
  "x2";"c"+<3pt,2pt>="c2"**\crv{"x2"+(0,-.25) & "c2"},
  "x3";"c"+<-3pt,-2pt>="c1"**\crv{"x3"+(0,.25) & "c1"}, 
  "c"+(0,.10)*!<0pt,-6pt>{\si},  
"x3"= "L32",
"x4"= "L33",
  "L24"="x1", 
  "x1"+(0,-.5)="x2",
  "x1";"x2"**\dir{-},  
"x2"= "L34",
  "L25"="x1", 
  "x1"+(0,-.5)="x2",
  "x1";"x2"**\dir{-},  
"x2"= "L35",
"L31"="L21",
"L32"="L22",
"L33"="L23",
"L34"="L24",
"L35"="L25",
  "L21"="x1", 
  "x1"+(0,-.5)="x2",
  "x1";"x2"**\dir{-},   
"x2"= "L31",
  "L22"="x1", 
  "x1"+(0,-.5)="x2",
  "x1";"x2"**\dir{-},   
"x2"= "L32",
   "L23"="x1",      
   "L24"="x2",  
  "x1"+(0,-.5)="x3",
  "x2"+(0,-.5)="x4",
  "x1";"x2"**\dir{}?(.5)="m", 
  "m" + (0,-.25)="c", 
  "x1";"x4"**\crv{ "x1"+(0,-.25) & "c" & "x4"+(0,.25)},  
  "x2";"c"+<3pt,2pt>="c2"**\crv{"x2"+(0,-.25) & "c2"},
  "x3";"c"+<-3pt,-2pt>="c1"**\crv{"x3"+(0,.25) & "c1"}, 
  "c"+(0,.10)*!<0pt,-6pt>{\si}, 
"x3"= "L33",
"x4"= "L34",
  "L25"="x1", 
  "x1"+(0,-.5)="x2",
  "x1";"x2"**\dir{-},  
"x2"= "L35",
"L31"="L21",
"L32"="L22",
"L33"="L23",
"L34"="L24",
"L35"="L25",
  "L21"="x1", 
  "x1"+(0,-.5)="x2",
  "x1";"x2"**\dir{-},  
"x2"= "L31",
  "L22"="x1", 
  "x1"+(0,-.5)="x2",
  "x1";"x2"**\dir{-},  
"x2"= "L32",
  "L23"="x1", 
  "x1"+(0,-.5)="x2",
  "x1";"x2"**\dir{-},  
"x2"= "L33",
   "L24"="x1",   
   "L25"="x2",  
  "x1"+(0,-.5)="x3",
  "x2"+(0,-.5)="x4",
  "x1";"x2"**\dir{}?(.5)="m", 
  "m" + (0,-.25)="c", 
  "x1";"x4"**\crv{ "x1"+(0,-.25) & "c" & "x4"+(0,.25)},  
  "x2";"c"+<3pt,2pt>="c2"**\crv{"x2"+(0,-.25) & "c2"},
  "x3";"c"+<-3pt,-2pt>="c1"**\crv{"x3"+(0,.25) & "c1"}, 
  "c"+(0,.10)*!<0pt,-6pt>{\si},  
"x3"= "L34",
"x4"= "L35",
"L31"="L21",
"L32"="L22",
"L33"="L23",
"L34"="L24",
"L35"="L25",
  "L21"="x1", 
  "x1"+(0,-.5)="x2",
  "x1";"x2"**\dir{-},  
"x2"= "L31",
   "L22"="x1",    
   "L23"="x2",  
  "x1"+(0,-.5)="x3",
  "x2"+(0,-.5)="x4",
  "x1";"x2"**\dir{}?(.5)="m", 
  "m" + (0,-.25)="c", 
  "x1";"x4"**\crv{ "x1"+(0,-.25) & "c" & "x4"+(0,.25)},  
  "x2";"c"+<3pt,2pt>="c2"**\crv{"x2"+(0,-.25) & "c2"},
  "x3";"c"+<-3pt,-2pt>="c1"**\crv{"x3"+(0,.25) & "c1"}, 
  "c"+(0,.10)*!<0pt,-6pt>{\si}, 
"x3"= "L32",
"x4"= "L33",
  "L24"="x1", 
  "x1"+(0,-.5)="x2",
  "x1";"x2"**\dir{-},  
"x2"= "L34",
  "L25"="x1", 
  "x1"+(0,-.5)="x2",
  "x1";"x2"**\dir{-},  
"x2"= "L35",
"L31"="L21",
"L32"="L22",
"L33"="L23",
"L34"="L24",
"L35"="L25",
  "L21" ="x1", 
  "x1"+(0,-.5)="x2",
  "x1";"x2"**\dir{}?(.5)="c", 
  "x1";"x2"**\dir{-},  
  "c"-(0,0)*!<-4pt,0pt>{\bullet S }, 
"x2"= "L31",
  "L22"="x1", 
  "x1"+(0,-.5)="x2",
  "x1";"x2"**\dir{-},   
"x2"= "L32",
  "L23" ="x1", 
  "x1"+(0,-.5)="x2",
  "x1";"x2"**\dir{}?(.5)="c", 
  "x1";"x2"**\dir{-},   
  "c"-(0,0)*!<-4pt,0pt>{\bullet S },  
"x2"= "L33",
  "L24"="x1", 
  "x1"+(0,-.5)="x2",
  "x1";"x2"**\dir{-},  
"x2"= "L34",
  "L25" ="x1",
  "x1"+(0,-.5)="x2",
  "x1";"x2"**\dir{}?(.5)="c", 
  "x1";"x2"**\dir{-},   
  "c"-(0,0)*!<-4pt,0pt>{\bullet \widetilde{S} },  
"x2"= "L35",
    "L25"+(1,-.25) ="x1", 
  "x1"+(0,-.25)="x2",
  "x1";"x2"**\dir{}?(.5)="c", 
  "x1";"x2"**\dir{-},   
  "c"-(0,0)*!<-3pt,0pt>{\bullet \sigma }, 
  "x1"+(0,0)*!<0pt,-2pt>{\circ},  
"x2"= "L36",
"L31"="L21",
"L32"="L22",
"L33"="L23",
"L34"="L24",
"L35"="L25",
"L36"="L26",
   "L21"="x1", 
   "x1"+(-.5,-.5)="x2",
  "x1"+(0,-.5)="x3",
  "x1"+(.5,-.5)="x4",
  "x1"+(0,-.25)="c",  
   "x1";"x2"**\crv{"x1" & "c" & "x2"+(0,.25)}, 
   "x1";"x3"**\crv{"x1" & "c" & "x3"+(0,.25)}, 
   "x1";"x4"**\crv{"x1" & "c" & "x4"+(0,.25)}, 
"x2"="L31",
"x3"="L32",
"x4"="L33",
   "L22"="x1", 
   "x1"+(-.5,-.5)="x2",
  "x1"+(0,-.5)="x3",
  "x1"+(.5,-.5)="x4",
  "x1"+(0,-.25)="c",  
   "x1";"x2"**\crv{"x1" & "c" & "x2"+(0,.25)}, 
   "x1";"x3"**\crv{"x1" & "c" & "x3"+(0,.25)}, 
   "x1";"x4"**\crv{"x1" & "c" & "x4"+(0,.25)}, 
"x2"="L34",
"x3"="L35",
"x4"="L36",
   "L23" ="x1", 
    "L24"="x2",  
  "x1";"x2"**\dir{}?(.5)="m", 
  "m" + (0,-.25)="c", 
  "m" + (0,-.5)="x3", 
  "x1";"x3"**\crv{"x1"-(0,.25) & "c" & "x3"},
  "x2";"x3"**\crv{"x2"-(0,.25) & "c"  & "x3"}, 
"x3"="L37",
   "L25" ="x1", 
    "L26"="x2", 
  "x1";"x2"**\dir{}?(.5)="m", 
  "m" + (0,-.25)="c", 
  "m" + (0,-.5)="x3", 
  "x1";"x3"**\crv{"x1"-(0,.25) & "c" & "x3"},
  "x2";"x3"**\crv{"x2"-(0,.25) & "c"  & "x3"}, 
"x3"="L38",
"L31"="L21",
"L32"="L22",
"L33"="L23",
"L34"="L24",
"L35"="L25",
"L36"="L26",
"L37"="L27",
"L38"="L28",
  "L21"="x1", 
  "x1"+(0,-.5)="x2",
  "x1";"x2"**\dir{-},   
"x2"= "L31",
  "L22"="x1",
  "x1"+(0,-.5)="x2",
  "x1";"x2"**\dir{-},  
"x2"= "L32",
   "L23"="x1",   
   "L24"="x2",  
  "x1"+(0,-.5)="x3",
  "x2"+(0,-.5)="x4",
  "x1";"x2"**\dir{}?(.5)="m", 
  "m" + (0,-.25)="c", 
  "x1";"x4"**\crv{ "x1"+(0,-.25) & "c" & "x4"+(0,.25)},  
  "x2";"c"+<3pt,2pt>="c2"**\crv{"x2"+(0,-.25) & "c2"},
  "x3";"c"+<-3pt,-2pt>="c1"**\crv{"x3"+(0,.25) & "c1"}, 
 "c"+(0,.10)*!<0pt,-6pt>{\si}, 
"x3"= "L33",
"x4"= "L34",
  "L25"="x1", 
  "x1"+(0,-.5)="x2",
  "x1";"x2"**\dir{-},  
"x2"= "L35",
  "L26"="x1", 
  "x1"+(0,-.5)="x2",
  "x1";"x2"**\dir{-},  
"x2"= "L36",
  "L27"="x1", 
  "x1"+(0,-.5)="x2",
  "x1";"x2"**\dir{-},   
"x2"= "L37",
  "L28"="x1", 
  "x1"+(0,-.5)="x2",
  "x1";"x2"**\dir{-},  
"x2"= "L38",
"L31"="L21",
"L32"="L22",
"L33"="L23",
"L34"="L24",
"L35"="L25",
"L36"="L26",
"L37"="L27",
"L38"="L28",
  "L21"="x1", 
  "x1"+(0,-.5)="x2",
  "x1";"x2"**\dir{-},   
"x2"= "L31",
  "L22"="x1", 
  "x1"+(0,-.5)="x2",
  "x1";"x2"**\dir{-},  
"x2"= "L32",
  "L23"="x1", 
  "x1"+(0,-.5)="x2",
  "x1";"x2"**\dir{-},  
"x2"= "L33",
   "L24"="x1",  
   "L25"="x2", 
  "x1"+(0,-.5)="x3",
  "x2"+(0,-.5)="x4",
  "x1";"x2"**\dir{}?(.5)="m", 
  "m" + (0,-.25)="c", 
  "x1";"x4"**\crv{ "x1"+(0,-.25) & "c" & "x4"+(0,.25)},  
  "x2";"c"+<3pt,2pt>="c2"**\crv{"x2"+(0,-.25) & "c2"},
  "x3";"c"+<-3pt,-2pt>="c1"**\crv{"x3"+(0,.25) & "c1"}, 
  "c"+(0,.10)*!<0pt,-6pt>{\si},  
"x3"= "L34",
"x4"= "L35",
  "L26"="x1", 
  "x1"+(0,-.5)="x2",
  "x1";"x2"**\dir{-},   
"x2"= "L36",
  "L27"="x1", 
  "x1"+(0,-.5)="x2",
  "x1";"x2"**\dir{-},  
"x2"= "L37",
  "L28"="x1", 
  "x1"+(0,-.5)="x2",
  "x1";"x2"**\dir{-},  
"x2"= "L38",
"L31"="L21",
"L32"="L22",
"L33"="L23",
"L34"="L24",
"L35"="L25",
"L36"="L26",
"L37"="L27",
"L38"="L28",
  "L21"="x1", 
  "x1"+(0,-.5)="x2",
  "x1";"x2"**\dir{-}, 
"x2"= "L31",
   "L22"="x1",    
   "L23"="x2",  
  "x1"+(0,-.5)="x3",
  "x2"+(0,-.5)="x4",
  "x1";"x2"**\dir{}?(.5)="m", 
  "m" + (0,-.25)="c", 
  "x1";"x4"**\crv{ "x1"+(0,-.25) & "c" & "x4"+(0,.25)},  
  "x2";"c"+<3pt,2pt>="c2"**\crv{"x2"+(0,-.25) & "c2"},
  "x3";"c"+<-3pt,-2pt>="c1"**\crv{"x3"+(0,.25) & "c1"}, 
  "c"+(0,.10)*!<0pt,-6pt>{\si}, 
"x3"= "L32",
"x4"= "L33",
  "L24"="x1", 
  "x1"+(0,-.5)="x2",
  "x1";"x2"**\dir{-},  
"x2"= "L34",
  "L25"="x1", 
  "x1"+(0,-.5)="x2",
  "x1";"x2"**\dir{-},   
"x2"= "L35",
  "L26"="x1", 
  "x1"+(0,-.5)="x2",
  "x1";"x2"**\dir{-},  
"x2"= "L36",
  "L27"="x1", 
  "x1"+(0,-.5)="x2",
  "x1";"x2"**\dir{-},  
"x2"= "L37",
  "L28"="x1", 
  "x1"+(0,-.5)="x2",
  "x1";"x2"**\dir{-},   
"x2"= "L38",
"L31"="L21",
"L32"="L22",
"L33"="L23",
"L34"="L24",
"L35"="L25",
"L36"="L26",
"L37"="L27",
"L38"="L28",
   "L21" ="x1",
    "L22"="x2",  
  "x1";"x2"**\dir{}?(.5)="m", 
  "m" + (0,-.25)="c", 
  "m" + (0,-.5)="x3", 
  "x1";"x3"**\crv{"x1"-(0,.25) & "c" & "x3"},
  "x2";"x3"**\crv{"x2"-(0,.25) & "c"  & "x3"}, 
"x3"="L31",
   "L23" ="x1",
    "L24"="x2",  
  "x1";"x2"**\dir{}?(.5)="m", 
  "m" + (0,-.25)="c", 
  "m" + (0,-.5)="x3", 
  "x1";"x3"**\crv{"x1"-(0,.25) & "c" & "x3"},
  "x2";"x3"**\crv{"x2"-(0,.25) & "c"  & "x3"}, 
"x3"="L32",
   "L25" ="x1", 
    "L26"="x2",  
  "x1";"x2"**\dir{}?(.5)="m", 
  "m" + (0,-.25)="c", 
  "m" + (0,-.5)="x3", 
  "x1";"x3"**\crv{"x1"-(0,.25) & "c" & "x3"},
  "x2";"x3"**\crv{"x2"-(0,.25) & "c"  & "x3"}, 
"x3"="L33",
  "L27"="x1", 
  "x1"+(0,-.5)="x2",
  "x1";"x2"**\dir{-},   
"x2"= "L34",
  "L28"="x1", 
  "x1"+(0,-.5)="x2",
  "x1";"x2"**\dir{-},  
"x2"= "L35",
"L31"="L21",
"L32"="L22",
"L33"="L23",
"L34"="L24",
"L35"="L25",
  "L21"="x1",
  "x1"+(0,-.5)="x2",
  "x1";"x2"**\dir{-},  
"x2"= "L31",
   "L22"="x1",      
   "L23"="x2",  
  "x1"+(0,-.5)="x3",
  "x2"+(0,-.5)="x4",
  "x1";"x2"**\dir{}?(.5)="m", 
  "m" + (0,-.25)="c", 
  "x1";"x4"**\crv{ "x1"+(0,-.25) & "c" & "x4"+(0,.25)},  
  "x2";"c"+<3pt,2pt>="c2"**\crv{"x2"+(0,-.25) & "c2"},
  "x3";"c"+<-3pt,-2pt>="c1"**\crv{"x3"+(0,.25) & "c1"}, 
  "c"+(0,.10)*!<0pt,-6pt>{\si},  
"x3"= "L32",
"x4"= "L33",
  "L24"="x1", 
  "x1"+(0,-.5)="x2",
  "x1";"x2"**\dir{-},  
"x2"= "L34",
  "L25"="x1", 
  "x1"+(0,-.5)="x2",
  "x1";"x2"**\dir{-}, 
"x2"= "L35",
"L31"="L21",
"L32"="L22",
"L33"="L23",
"L34"="L24",
"L35"="L25",
   "L21" ="x1",   
   "L22"="x2",  
  "x1"+(0,-.5)="x3",
  "x2"+(0,-.5)="x4",
  "x1";"x2"**\dir{}?(.5)="m", 
  "m" + (0,-.25)="c", 
  "x1";"x4"**\crv{ "x1"+(0,-.25) & "c" & "x4"+(0,.25)},  
  "x2";"c"+<3pt,2pt>="c2"**\crv{"x2"+(0,-.25) & "c2"},
  "x3";"c"+<-3pt,-2pt>="c1"**\crv{"x3"+(0,.25) & "c1"}, 
  "c"+(0,.10)*!<0pt,-6pt>{\si},  
"x3"= "L31",
"x4"= "L32",
  "L23"="x1", 
  "x1"+(0,-.5)="x2",
  "x1";"x2"**\dir{-},  
"x2"= "L33",
  "L24"="x1", 
  "x1"+(0,-.5)="x2",
  "x1";"x2"**\dir{-},  
"x2"= "L34",
  "L25"="x1", 
  "x1"+(0,-.5)="x2",
  "x1";"x2"**\dir{-},   
"x2"= "L35",
"L31"="L21",
"L32"="L22",
"L33"="L23",
"L34"="L24",
"L35"="L25",
  "L21"="x1", 
  "x1"+(0,-.5)="x2",
  "x1";"x2"**\dir{-},  
"x2"= "L31",
   "L22"="x1",     
   "L23"="x2",  
  "x1"+(0,-.5)="x3",
  "x2"+(0,-.5)="x4",
  "x1";"x2"**\dir{}?(.5)="m", 
  "m" + (0,-.25)="c", 
  "x1";"x4"**\crv{ "x1"+(0,-.25) & "c" & "x4"+(0,.25)},  
  "x2";"c"+<3pt,2pt>="c2"**\crv{"x2"+(0,-.25) & "c2"},
  "x3";"c"+<-3pt,-2pt>="c1"**\crv{"x3"+(0,.25) & "c1"}, 
  "c"+(0,.10)*!<0pt,-6pt>{\si},  
"x3"= "L32",
"x4"= "L33",
  "L24"="x1", 
  "x1"+(0,-.5)="x2",
  "x1";"x2"**\dir{-},  
"x2"= "L34",
  "L25"="x1", 
  "x1"+(0,-.5)="x2",
  "x1";"x2"**\dir{-},   
"x2"= "L35",
"L31"="L21",
"L32"="L22",
"L33"="L23",
"L34"="L24",
"L35"="L25",
  "L21"="x1", 
  "x1"+(0,-.5)="x2",
  "x1";"x2"**\dir{-},   
"x2"= "L31",
  "L22"="x1", 
  "x1"+(0,-.5)="x2",
  "x1";"x2"**\dir{-},  
"x2"= "L32",
   "L23"="x1",    
   "L24"="x2",  
  "x1"+(0,-.5)="x3",
  "x2"+(0,-.5)="x4",
  "x1";"x2"**\dir{}?(.5)="m", 
  "m" + (0,-.25)="c", 
  "x1";"x4"**\crv{ "x1"+(0,-.25) & "c" & "x4"+(0,.25)},  
  "x2";"c"+<3pt,2pt>="c2"**\crv{"x2"+(0,-.25) & "c2"},
  "x3";"c"+<-3pt,-2pt>="c1"**\crv{"x3"+(0,.25) & "c1"}, 
  "c"+(0,.10)*!<0pt,-6pt>{\si}, 
"x3"= "L33",
"x4"= "L34",
  "L25"="x1", 
  "x1"+(0,-.5)="x2",
  "x1";"x2"**\dir{-},  
"x2"= "L35",
"L31"="L21",
"L32"="L22",
"L33"="L23",
"L34"="L24",
"L35"="L25",
  "L21"="x1", 
  "x1"+(0,-.5)="x2",
  "x1";"x2"**\dir{-},  
"x2"= "L31",
  "L22"="x1", 
  "x1"+(0,-.5)="x2",
  "x1";"x2"**\dir{-},  
"x2"= "L32",
  "L23"="x1", 
  "x1"+(0,-.5)="x2",
  "x1";"x2"**\dir{-}, 
"x2"= "L33",
   "L24"="x1",  
   "L25"="x2",  
  "x1"+(0,-.5)="x3",
  "x2"+(0,-.5)="x4",
  "x1";"x2"**\dir{}?(.5)="m", 
  "m" + (0,-.25)="c", 
  "x1";"x4"**\crv{ "x1"+(0,-.25) & "c" & "x4"+(0,.25)},  
  "x2";"c"+<3pt,2pt>="c2"**\crv{"x2"+(0,-.25) & "c2"},
  "x3";"c"+<-3pt,-2pt>="c1"**\crv{"x3"+(0,.25) & "c1"}, 
  "c"+(0,.10)*!<0pt,-6pt>{\si}, 
"x3"= "L34",
"x4"= "L35",
"L31"="L21",
"L32"="L22",
"L33"="L23",
"L34"="L24",
"L35"="L25",
  "L21"="x1", 
  "x1"+(0,-.5)="x2",
  "x1";"x2"**\dir{-},  
"x2"= "L31",
   "L22"="x1",     
   "L23"="x2",  
  "x1"+(0,-.5)="x3",
  "x2"+(0,-.5)="x4",
  "x1";"x2"**\dir{}?(.5)="m", 
  "m" + (0,-.25)="c", 
  "x1";"x4"**\crv{ "x1"+(0,-.25) & "c" & "x4"+(0,.25)},  
  "x2";"c"+<3pt,2pt>="c2"**\crv{"x2"+(0,-.25) & "c2"},
  "x3";"c"+<-3pt,-2pt>="c1"**\crv{"x3"+(0,.25) & "c1"}, 
  "c"+(0,.10)*!<0pt,-6pt>{\si},  
"x3"= "L32",
"x4"= "L33",
  "L24"="x1", 
  "x1"+(0,-.5)="x2",
  "x1";"x2"**\dir{-},  
"x2"= "L34",
  "L25"="x1", 
  "x1"+(0,-.5)="x2",
  "x1";"x2"**\dir{-},  
"x2"= "L35",
"L31"="L21",
"L32"="L22",
"L33"="L23",
"L34"="L24",
"L35"="L25",
  "L21" ="x1", 
  "x1"+(0,-.5)="x2",
  "x1";"x2"**\dir{}?(.5)="c", 
  "x1";"x2"**\dir{-},  
  "c"-(0,0)*!<-4pt,0pt>{\bullet S }, 
"x2"= "L31",
  "L22"="x1", 
  "x1"+(0,-.5)="x2",
  "x1";"x2"**\dir{-},   
"x2"= "L32",
  "L23" ="x1", 
  "x1"+(0,-.5)="x2",
  "x1";"x2"**\dir{}?(.5)="c", 
  "x1";"x2"**\dir{-},  
  "c"-(0,0)*!<-4pt,0pt>{\bullet S }, 
"x2"= "L33",
  "L24"="x1",
  "x1"+(0,-.5)="x2",
  "x1";"x2"**\dir{-},  
"x2"= "L34",
  "L25" ="x1", 
  "x1"+(0,-.5)="x2",
  "x1";"x2"**\dir{}?(.5)="c", 
  "x1";"x2"**\dir{-}, 
  "c"-(0,0)*!<-4pt,0pt>{\bullet \widetilde{S} },
"x2"= "L35",
    "L25"+(1,-.25) ="x1", 
  "x1"+(0,-.25)="x2",
  "x1";"x2"**\dir{}?(.5)="c", 
  "x1";"x2"**\dir{-},  
  "c"-(0,0)*!<-3pt,0pt>{\bullet \sigma },  
  "x1"+(0,0)*!<0pt,-2pt>{\circ},  
"x2"= "L36",
"L31"="L21",
"L32"="L22",
"L33"="L23",
"L34"="L24",
"L35"="L25",
"L36"="L26",
   "L21" ="x1", 
    "L22"="x2",  
  "x1";"x2"**\dir{}?(.5)="m", 
  "m" + (0,-.25)="c", 
  "m" + (0,-.5)="x3", 
  "x1";"x3"**\crv{"x1"-(0,.25) & "c" & "x3"},
  "x2";"x3"**\crv{"x2"-(0,.25) & "c"  & "x3"}, 
"x3"="L31",
   "L23" ="x1",
    "L24"="x2",  
  "x1";"x2"**\dir{}?(.5)="m", 
  "m" + (0,-.25)="c", 
  "m" + (0,-.5)="x3", 
  "x1";"x3"**\crv{"x1"-(0,.25) & "c" & "x3"},
  "x2";"x3"**\crv{"x2"-(0,.25) & "c"  & "x3"}, 
"x3"="L32",
    "L25" =  "x1", 
  "x1"+(0,-.25)="x2",
  "x1";"x2"**\dir{}?(.5)="c", 
  "x1";"x2"**\dir{-}, 
  "c"-(0,0)*!<-3pt,0pt>{\bullet \varepsilon }, 
  "x2"+(0,0)*!<0pt,2pt>{\circ}, 
    "L26" =  "x1",
  "x1"+(0,-.25)="x2",
  "x1";"x2"**\dir{}?(.5)="c", 
  "x1";"x2"**\dir{-},   
  "c"-(0,0)*!<-3pt,0pt>{\bullet \varepsilon },
  "x2"+(0,0)*!<0pt,2pt>{\circ}, 
 "L31"-(0,.5)*!<0pt,-6pt>{H},   
 "L32"-(0,.5)*!<0pt,-6pt>{H}, 
\endxy
\]

\begin{align*} 
=&(m,m)(S,1,S,1,\varepsilon\widetilde{S},\varepsilon\sigma)(1,\si,2)(3,\si)(2,\si,1)(1,\si,2)(\si,3)(1,\si,2)\\
~~~~~&(m,m,m,m,m)(1,\si,7)(3,\si,5)(2,\si,6)(\Delta^2,\Delta^2,4)(S,1,S,1,\widetilde{S},\sigma)\\
~~~~~&(1,\si,2)(3,\si)(2,\si,1)(1,\si,2)(\si,3)(1,\si,2)(\Delta^2,2)\\
=&(m,m)(S,1,S,1,\delta)(1,\si,2)(3,\si)(2,\si,1)(1,\si,2)(\si,3)(1,\si,2)\\
~~~~~&(m,m,m,m,m)(1,\si,7)(3,\si,5)(2,\si,6)(\Delta^2S,\Delta^2,S,1,\widetilde{S},\sigma)\\
~~~~~&(1,\si,2)(3,\si)(2,\si,1)(1,\si,2)(\si,3)(1,\si,2)(\Delta^2,2)\\
=&(m,m)(S,1,S,1,\delta)(m,m,m,m,m)(2,\si_{2,2},4)(6,\si_{2,2})(4,\si_{2,2},2)(2,\si_{2,2},4)\\
~~~~~&(\si_{2,2},6)(2,\si_{2,2},4)(1,\si,7)(3,\si,5)(2,\si,6)(S,S,S,3,S,1,\widetilde{S},\sigma)(1,\si,6)\\
~~~~~&(\si,7)(1,\si,6)(\Delta^2,\Delta^2,3)(1,\si,2)(3,\si)(2,\si,1)(1,\si,2)(\si,3)(1,\si,2)(\Delta^2,2)\\
=&(m,m)(Sm,m,Sm,m,\delta m)(2,\si_{2,2},4)(6,\si_{2,2})(4,\si_{2,2},2)(2,\si_{2,2},4)\\
~~~~~&(\si_{2,2},6)(2,\si_{2,2},4)(S,1,S,1,S,1,S,1,\widetilde{S},\sigma)(1,\si,6)(3,\si,4)(2,\si,5)(1,\si,6)\\
~~~~~&(\si,7)(1,\si,6)(\Delta^2,\Delta^2,3)(1,\si,2)(3,\si)(2,\si,1)(1,\si,2)(\si,3)(1,\si,2)(\Delta^2,2)\\
\end{align*}
\begin{align*}
=&(m,m)(m\si,m,m\si,m)(S,S,1,1,S,S,1,1,\delta,\delta)(2,\si_{2,2},4)(6,\si_{2,2})\\
~~~~~&(4,\si_{2,2},2)(2,\si_{2,2},4)(\si_{2,2},6)(2,\si_{2,2},4)(S,1,S,1,S,1,S,1,\widetilde{S},\sigma)\\
~~~~~&(1,\si,6)(3,\si,4)(2,\si,5) (1,\si,6)(\si,7)(1,\si,6)(\Delta^2,\Delta^2,3)\\
~~~~~&(1,\si,2)(3,\si)(2,\si,1)(1,\si,2)(\si,3)(1,\si,2)(\Delta^2,2)\\
=&(m,m)(m\si,m,m\si,m)(2,\si_{2,2},2)(\si_{2,2},4)(\delta,\delta,S,S,S,S,4)\\
~~~~~&(S,1,S,1,S,1,S,1,\widetilde{S},\sigma)(1,\si,6)(3,\si,4)(2,\si,5) (1,\si,6)(\si,7)(1,\si,6)(\Delta^2,\Delta^2,3)\\
~~~~~&(1,\si,2)(3,\si)(2,\si,1)(1,\si,2)(\si,3)(1,\si,2)(\Delta^2,2)\\
=&(m,m)(1,\si,1)(\si,2)(m\si,m\si,m,m)(\delta S,\delta,S^2,S,S^2,S,S,1,\widetilde{S},\sigma)\\
~~~~~&(1,\si,6)(3,\si,4)(2,\si,5)(1,\si,6)(\si,7)(1,\si,6)(\Delta^2,\Delta^2,3)  (1,\si,2)(3,\si)(2,\si,1)\\
~~~~~&(1,\si,2)(\si,3)(1,\si,2)(\Delta^2,2)\\
=&(m,m)(1,\si,1)(\si,2)(m\si,m\si,m,m)(1,\si,5)(\si,6)(S^2,S^2,\delta S,\delta,S,S,S,1,\widetilde{S},\sigma)\\
~~~~~&(\Delta^2,\Delta^2,3)  (1,\si,2)(3,\si)(2,\si,1)(1,\si,2)(\si,3)(1,\si,2)(\Delta^2,2)\\
=&(m,m)(1,\si,1)(\si,2)(m\si,m\si,m,m)(1,\si,5)(\si,6)(S^2,S^2,\delta S,\delta,S,S,S,1,\widetilde{S},\sigma)\\
~~~~~&(3,\Delta^2,3) (3,\si,2)(5,\si)(4,\si,1)(3,\si,2)(\Delta^2,4)(\si,3)(1,\si,2)(\Delta^2,2)=
\end{align*}

\[
\xy /r2pc/:,
(1,0)= "L11",
(4,0)= "L12",
(5,0)= "L13",
  "L11"+(0,.5)*!<0pt,6pt>{H}, 
  "L12"+(0,.5)*!<0pt,6pt>{H}, 
  "L13"+(0,.5)*!<0pt,6pt>{H}, 
  "L11"="x1", 
  "x1"+(-1.5,-.5)="x2",
  "x1"+(0,-.5)="x3",
  "x1"+(1.5,-.5)="x4",
  "x1"+(0,-.25)="c",  
   "x1";"x2"**\crv{"x1" & "c" & "x2"+(0,.25)}, 
   "x1";"x3"**\crv{"x1" & "c" & "x3"+(0,.25)}, 
   "x1";"x4"**\crv{"x1" & "c" & "x4"+(0,.25)}, 
"L12";"L12" +(0,-.5)="L14"**\dir{-},
"L13";"L13" +(0,-.5)="L15"**\dir{-},
"x2"="L11",
"x3"="L12",
"x4"="L13",
  "L11";"L11" +(0,-.5)="L11"**\dir{-},
  "L12"="x1",       
  "L13"="x2",  , 
  "x1"+(0,-.5)="x3",
  "x2"+(0,-.5)="x4",
  "x1";"x2"**\dir{}?(.5)="m", 
  "m" + (0,-.25)="c", 
  "x1";"x4"**\crv{ "x1"+(0,-.25) & "c" & "x4"+(0,.25)},  
  "x2";"c"+<3pt,2pt>="c2"**\crv{"x2"+(0,-.25) & "c2"},
  "x3";"c"+<-3pt,-2pt>="c1"**\crv{"x3"+(0,.25) & "c1"}, 
  "c"+(0,.10)*!<0pt,-4pt>{\si},
"x3"="L12",
"x4"="L13",
  "L14";"L14" +(0,-.5)="L14"**\dir{-},
  "L15";"L15" +(0,-.5)="L15"**\dir{-},
  "L11"="x1",       
  "L12"="x2",  , 
  "x1"+(-.25,-.5)="x3",
  "x2"+(.25,-.5)="x4",
  "x1";"x2"**\dir{}?(.5)="m", 
  "m" + (0,-.25)="c", 
  "x1";"x4"**\crv{ "x1"+(0,-.25) & "c" & "x4"+(0,.25)},  
  "x2";"c"+<3pt,2pt>="c2"**\crv{"x2"+(0,-.25) & "c2"},
  "x3";"c"+<-3pt,-2pt>="c1"**\crv{"x3"+(0,.25) & "c1"}, 
  "c"+(0,.10)*!<0pt,-6pt>{\si},
"x3"="L11",
"x4"="L12",
  "L13";"L13" +(0,-.5)="L13"**\dir{-},
  "L14";"L14" +(0,-.5)="L14"**\dir{-},
  "L15";"L15" +(0,-.5)="L15"**\dir{-},
  "L11"="x1", 
  "x1"+(-.75,-.5)="x2",
  "x1"+(0,-.5)="x3",
  "x1"+(.75,-.5)="x4",
  "x1"+(0,-.25)="c",  
   "x1";"x2"**\crv{"x1" & "c" & "x2"+(0,.25)}, 
   "x1";"x3"**\crv{"x1" & "c" & "x3"+(0,.25)}, 
   "x1";"x4"**\crv{"x1" & "c" & "x4"+(0,.25)}, 
"L12";"L12" +(0,-.5)="L24"**\dir{-},
"L13";"L13" +(0,-.5)="L25"**\dir{-},
"L14";"L14" +(0,-.5)="L26"**\dir{-},
"L15";"L15" +(0,-.5)="L27"**\dir{-},
"x2"="L11",
"x3"="L12",
"x4"="L13",
"L24"="L14",
"L25"="L15",
"L26"="L16",
"L27"="L17",
  "L11";"L11" +(0,-.5)="L11"**\dir{-},
  "L12";"L12" +(0,-.5)="L12"**\dir{-},
  "L13";"L13" +(0,-.5)="L13"**\dir{-},
  "L14"="x1",       
  "L15"="x2",  , 
  "x1"+(0,-.5)="x3",
  "x2"+(0,-.5)="x4",
  "x1";"x2"**\dir{}?(.5)="m", 
  "m" + (0,-.25)="c", 
  "x1";"x4"**\crv{ "x1"+(0,-.25) & "c" & "x4"+(0,.25)},  
  "x2";"c"+<3pt,2pt>="c2"**\crv{"x2"+(0,-.25) & "c2"},
  "x3";"c"+<-3pt,-2pt>="c1"**\crv{"x3"+(0,.25) & "c1"}, 
  "c"+(0,.10)*!<0pt,-6pt>{\si},
"x3"="L14",
"x4"="L15",
  "L16";"L16" +(0,-.5)="L16"**\dir{-},
  "L17";"L17" +(0,-.5)="L17"**\dir{-},
  "L11";"L11" +(0,-.5)="L11"**\dir{-},
  "L12";"L12" +(0,-.5)="L12"**\dir{-},
  "L13";"L13" +(0,-.5)="L13"**\dir{-},
  "L14";"L14" +(0,-.5)="L14"**\dir{-},
  "L15"="x1",       
  "L16"="x2",  , 
  "x1"+(0,-.5)="x3",
  "x2"+(0,-.5)="x4",
  "x1";"x2"**\dir{}?(.5)="m", 
  "m" + (0,-.25)="c", 
  "x1";"x4"**\crv{ "x1"+(0,-.25) & "c" & "x4"+(0,.25)},  
  "x2";"c"+<3pt,2pt>="c2"**\crv{"x2"+(0,-.25) & "c2"},
  "x3";"c"+<-3pt,-2pt>="c1"**\crv{"x3"+(0,.25) & "c1"}, 
  "c"+(0,.10)*!<0pt,-6pt>{\si},
"x3"="L15",
"x4"="L16",
  "L17";"L17" +(0,-.5)="L17"**\dir{-},
  "L11";"L11" +(0,-.5)="L11"**\dir{-},
  "L12";"L12" +(0,-.5)="L12"**\dir{-},
  "L13";"L13" +(0,-.5)="L13"**\dir{-},
  "L14";"L14" +(0,-.5)="L14"**\dir{-},
  "L15";"L15" +(0,-.5)="L15"**\dir{-},
  "L16"="x1",       
  "L17"="x2",  , 
  "x1"+(0,-.5)="x3",
  "x2"+(0,-.5)="x4",
  "x1";"x2"**\dir{}?(.5)="m", 
  "m" + (0,-.25)="c", 
  "x1";"x4"**\crv{ "x1"+(0,-.25) & "c" & "x4"+(0,.25)},  
  "x2";"c"+<3pt,2pt>="c2"**\crv{"x2"+(0,-.25) & "c2"},
  "x3";"c"+<-3pt,-2pt>="c1"**\crv{"x3"+(0,.25) & "c1"}, 
  "c"+(0,.10)*!<0pt,-6pt>{\si},
"x3"="L16",
"x4"="L17",
  "L11";"L11" +(0,-.5)="L11"**\dir{-},
  "L12";"L12" +(0,-.5)="L12"**\dir{-},
  "L13";"L13" +(0,-.5)="L13"**\dir{-},
  "L14"="x1",       
  "L15"="x2",  , 
  "x1"+(0,-.5)="x3",
  "x2"+(0,-.5)="x4",
  "x1";"x2"**\dir{}?(.5)="m", 
  "m" + (0,-.25)="c", 
  "x1";"x4"**\crv{ "x1"+(0,-.25) & "c" & "x4"+(0,.25)},  
  "x2";"c"+<3pt,2pt>="c2"**\crv{"x2"+(0,-.25) & "c2"},
  "x3";"c"+<-3pt,-2pt>="c1"**\crv{"x3"+(0,.25) & "c1"}, 
  "c"+(0,.10)*!<0pt,-6pt>{\si},
"x3"="L14",
"x4"="L15",
  "L16";"L16" +(0,-.5)="L16"**\dir{-},
  "L17";"L17" +(0,-.5)="L17"**\dir{-},
  "L11";"L11" +(0,-.5)="L11"**\dir{-},
  "L12";"L12" +(0,-.5)="L12"**\dir{-},
  "L13";"L13" +(0,-.5)="L13"**\dir{-},   
  "L14"="x1", 
  "x1"+(-.5,-.5)="x2",
  "x1"+(0,-.5)="x3",
  "x1"+(.5,-.5)="x4",
  "x1"+(0,-.25)="c",  
   "x1";"x2"**\crv{"x1" & "c" & "x2"+(0,.25)}, 
   "x1";"x3"**\crv{"x1" & "c" & "x3"+(0,.25)}, 
   "x1";"x4"**\crv{"x1" & "c" & "x4"+(0,.25)}, 
"L15";"L15" +(0,-.5)="L27"**\dir{-},
"L16";"L16" +(0,-.5)="L28"**\dir{-},
"L17";"L17" +(0,-.5)="L29"**\dir{-},
"x2"="L14",
"x3"="L15",
"x4"="L16",
"L27"="L17",
"L28"="L18",
"L29"="L19",
  "L11"-(0,.25)*!<3pt,-2pt>{S^2 \bullet }, 
  "L11";"L11" +(0,-.5)="L11"**\dir{-},
  "L12"-(0,.25)*!<3pt,-2pt>{S^2 \bullet }, 
  "L12";"L12" +(0,-.5)="L12"**\dir{-},
  "L13"-(0,.15)*!<3pt,0pt>{\delta S \bullet}, 
  "L13";"L13" +(0,-.3)="L13"**\dir{-},
  "L13"+(0,0)*!<0pt,2pt>{\circ}, 
  "L14"-(0,.15)*!<3pt,0pt>{\delta \bullet}, 
  "L14";"L14" +(0,-.3)="L14"**\dir{-},
  "L14"+(0,0)*!<0pt,2pt>{\circ},  
  "L15"-(0,.25)*!<-2pt,0pt>{\bullet S }, 
  "L15";"L15" +(0,-.5)="L15"**\dir{-}, 
  "L16"-(0,.25)*!<-2pt,0pt>{\bullet S}, 
  "L16";"L16" +(0,-.5)="L16"**\dir{-}, 
  "L17"-(0,.25)*!<-2pt,0pt>{\bullet S}, 
  "L17";"L17" +(0,-.5)="L17"**\dir{-}, 
  "L18";"L18" +(0,-.5)="L18"**\dir{-}, 
  "L19"-(0,.25)*!<-2pt,0pt>{\bullet \widetilde{S} }, 
  "L19";"L19" +(0,-.5)="L19"**\dir{-}, 
   "L19" + (1, .5)="x1", 
  "x1"+(0,-.5)="x2",
  "x1";"x2"**\dir{}?(.5)="c", 
  "x1";"x2"**\dir{-},   
  "c"-(0,0)*!<-3pt,0pt>{\bullet \sigma },  
  "x1"+(0,0)*!<0pt,-2pt>{\circ}, 
"x2"="L110",
  "L11";"L11" +(0,-.25)="L11"**\dir{-},
  "L12";"L12" +(0,-.25)="L12"**\dir{-},
  "L15";"L15" +(0,-.25)="L15"**\dir{-},
  "L16";"L16" +(0,-.25)="L16"**\dir{-},
  "L17";"L17" +(0,-.25)="L17"**\dir{-},
  "L18";"L18" +(0,-.25)="L18"**\dir{-},
  "L19";"L19" +(0,-.25)="L19"**\dir{-},
  "L110";"L110" +(0,-.25)="L110"**\dir{-},
  "L11"="x1",       
  "L12"="x2",  , 
  "x1"+(0,-.5)="x3",
  "x2"+(0,-.5)="x4",
  "x1";"x2"**\dir{}?(.5)="m", 
  "m" + (0,-.25)="c", 
  "x1";"x4"**\crv{ "x1"+(0,-.25) & "c" & "x4"+(0,.25)},  
  "x2";"c"+<3pt,2pt>="c2"**\crv{"x2"+(0,-.25) & "c2"},
  "x3";"c"+<-3pt,-2pt>="c1"**\crv{"x3"+(0,.25) & "c1"}, 
  "c"+(0,.10)*!<0pt,-4pt>{\si},
"x3"="L11",
"x4"="L12",
  "L15";"L15" +(0,-.5)="L15"**\dir{-},
  "L16";"L16" +(0,-.5)="L16"**\dir{-},
  "L17";"L17" +(0,-.5)="L17"**\dir{-},
  "L18";"L18" +(0,-.5)="L18"**\dir{-},
  "L19";"L19" +(0,-.5)="L19"**\dir{-},
  "L110";"L110" +(0,-.5)="L110"**\dir{-},
 "L11";"L11" +(0,-.5)="L11"**\dir{-},
  "L12"="x1",       
  "L15"="x2",  , 
  "x1"+(0,-.5)="x3",
  "x2"+(0,-.5)="x4",
  "x1";"x2"**\dir{}?(.5)="m", 
  "m" + (0,-.25)="c", 
  "x1";"x4"**\crv{ "x1"+(0,-.25) & "c" & "x4"+(0,.25)},  
  "x2";"c"+<3pt,2pt>="c2"**\crv{"x2"+(0,-.25) & "c2"},
  "x3";"c"+<-3pt,-2pt>="c1"**\crv{"x3"+(0,.25) & "c1"}, 
  "c"+(0,.10)*!<0pt,-6pt>{\si},
"x3"="L12",
"x4"="L15",
  "L16";"L16" +(0,-.5)="L16"**\dir{-},
  "L17";"L17" +(0,-.5)="L17"**\dir{-},
  "L18";"L18" +(0,-.5)="L18"**\dir{-},
  "L19";"L19" +(0,-.5)="L19"**\dir{-},
  "L110";"L110" +(0,-.5)="L110"**\dir{-},
  "L11"="x1",       
  "L12"="x2",  , 
  "x1"+(0,-.5)="x3",
  "x2"+(0,-.5)="x4",
  "x1";"x2"**\dir{}?(.5)="m", 
  "m" + (0,-.25)="c", 
  "x1";"x4"**\crv{ "x1"+(0,-.25) & "c" & "x4"+(0,.25)},  
  "x2";"c"+<3pt,2pt>="c2"**\crv{"x2"+(0,-.25) & "c2"},
  "x3";"c"+<-3pt,-2pt>="c1"**\crv{"x3"+(0,.25) & "c1"}, 
  "c"+(0,.10)*!<0pt,-6pt>{\si},
"x3"="L11",
"x4"="L12",
  "L15"="x1",       
  "L16"="x2",  , 
  "x1"+(0,-.5)="x3",
  "x2"+(0,-.5)="x4",
  "x1";"x2"**\dir{}?(.5)="m", 
  "m" + (0,-.25)="c", 
  "x1";"x4"**\crv{ "x1"+(0,-.25) & "c" & "x4"+(0,.25)},  
  "x2";"c"+<3pt,2pt>="c2"**\crv{"x2"+(0,-.25) & "c2"},
  "x3";"c"+<-3pt,-2pt>="c1"**\crv{"x3"+(0,.25) & "c1"}, 
  "c"+(0,.10)*!<0pt,-6pt>{\si},
"x3"="L15",
"x4"="L16",
  "L17";"L17" +(0,-.5)="L17"**\dir{-},
  "L18";"L18" +(0,-.5)="L18"**\dir{-},
  "L19";"L19" +(0,-.5)="L19"**\dir{-},
  "L110";"L110" +(0,-.5)="L110"**\dir{-},
   "L11"="x1", 
  "L12"="x2",   
  "x1";"x2"**\dir{}?(.5)="m", 
  "m" + (0,-.25)="c", 
  "m" + (0,-.5)="x3", 
  "x1";"x3"**\crv{"x1"-(0,.25) & "c" & "x3"},
  "x2";"x3"**\crv{"x2"-(0,.25) & "c"  & "x3"}, 
"x3"="L11",
   "L15"="x1", 
  "L16"="x2",   
  "x1";"x2"**\dir{}?(.5)="m", 
  "m" + (0,-.25)="c", 
  "m" + (0,-.5)="x3", 
  "x1";"x3"**\crv{"x1"-(0,.25) & "c" & "x3"},
  "x2";"x3"**\crv{"x2"-(0,.25) & "c"  & "x3"}, 
"x3"="L12",
   "L17"="x1", 
  "L18"="x2",   
  "x1";"x2"**\dir{}?(.5)="m", 
  "m" + (0,-.25)="c", 
  "m" + (0,-.5)="x3", 
  "x1";"x3"**\crv{"x1"-(0,.25) & "c" & "x3"},
  "x2";"x3"**\crv{"x2"-(0,.25) & "c"  & "x3"}, 
"x3"="L13",
   "L19"="x1", 
  "L110"="x2",   
  "x1";"x2"**\dir{}?(.5)="m", 
  "m" + (0,-.25)="c", 
  "m" + (0,-.5)="x3", 
  "x1";"x3"**\crv{"x1"-(0,.25) & "c" & "x3"},
  "x2";"x3"**\crv{"x2"-(0,.25) & "c"  & "x3"}, 
"x3"="L14",
  "L11"="x1",       
  "L12"="x2",  , 
  "x1"+(0,-.5)="x3",
  "x2"+(0,-.5)="x4",
  "x1";"x2"**\dir{}?(.5)="m", 
  "m" + (0,-.25)="c", 
  "x1";"x4"**\crv{ "x1"+(0,-.25) & "c" & "x4"+(0,.25)},  
  "x2";"c"+<3pt,2pt>="c2"**\crv{"x2"+(0,-.25) & "c2"},
  "x3";"c"+<-3pt,-2pt>="c1"**\crv{"x3"+(0,.25) & "c1"}, 
  "c"+(0,.10)*!<0pt,-6pt>{\si},
"x3"="L11",
"x4"="L12",
  "L13";"L13" +(0,-.5)="L13"**\dir{-},
  "L14";"L14" +(0,-.5)="L14"**\dir{-},
  "L11";"L11" +(0,-.5)="L11"**\dir{-},
  "L12"="x1",       
  "L13"="x2",  , 
  "x1"+(0,-.5)="x3",
  "x2"+(0,-.5)="x4",
  "x1";"x2"**\dir{}?(.5)="m", 
  "m" + (0,-.25)="c", 
  "x1";"x4"**\crv{ "x1"+(0,-.25) & "c" & "x4"+(0,.25)},  
  "x2";"c"+<3pt,2pt>="c2"**\crv{"x2"+(0,-.25) & "c2"},
  "x3";"c"+<-3pt,-2pt>="c1"**\crv{"x3"+(0,.25) & "c1"}, 
  "c"+(0,.10)*!<0pt,-6pt>{\si},
"x3"="L12",
"x4"="L13",
  "L14";"L14" +(0,-.5)="L14"**\dir{-},
   "L11"="x1", 
  "L12"="x2",   
  "x1";"x2"**\dir{}?(.5)="m", 
  "m" + (0,-.25)="c", 
  "m" + (0,-.5)="x3", 
  "x1";"x3"**\crv{"x1"-(0,.25) & "c" & "x3"},
  "x2";"x3"**\crv{"x2"-(0,.25) & "c"  & "x3"}, 
"x3"="L11",
   "L13"="x1", 
  "L14"="x2",  
  "x1";"x2"**\dir{}?(.5)="m", 
  "m" + (0,-.25)="c", 
  "m" + (0,-.5)="x3", 
  "x1";"x3"**\crv{"x1"-(0,.25) & "c" & "x3"},
  "x2";"x3"**\crv{"x2"-(0,.25) & "c"  & "x3"}, 
"x3"="L12",
  "L11"-(0,.5)*!<0pt,-6pt>{H}, 
  "L12"-(0,.5)*!<0pt,-6pt>{H}, 
\endxy 
\]

\begin{align*} 
=&(m,m)(1,\si,1)(\si,2)(m\si,m\si,m,m)(1,\si,5)(\si,6)(S^2,S^2,\delta S,\delta,S,S,S,1,\widetilde{S},\sigma)\\
~~~~~&(3,\Delta^2,3) (3,\si,2)(5,\si)(4,\si,1)(3,\si,2)(\si_{13},3)(1,\si_{13},2)(2,\Delta^2,2)(\Delta^2,2)\\
=&(m,m)(1,\si,1)(\si,2)(m\si,m\si,m,m)(1,\si,5)(\si,6)(S^2,S^2,\delta S,\delta,S,S,S,1,\widetilde{S},\sigma)\\
\end{align*}
\begin{align*}
~~~~~&(3,\si_{13},2)(7,\si)(4,\si_{13},1)(3,\si,4)(\si_{13},3)(1,\si_{13},2)(5,\Delta^2,1)(2,\Delta^2,2)(\Delta^2,2)\\
=&(m,m)(1,\si,1)(\si,2)(m\si,m\si,m,m)(1,\si,5)(\si,6)(2,\si_{12},3)(5,\si,3)(3,\si_{12},2)\\
~~~~~&(2,\si,4)(\si_{12},5)(1,\si_{12},4)(\widetilde{S},S,S^2,S^2,\delta S,\delta,S,S,1,\sigma)(5,\Delta^2,1)(2,\Delta^2,2)(\Delta^2,2)\\
=&(m,m)(1,\si,1)(\si,2)(m\si,m\si,m,m)(1,\si,5)(\si,6)(S^2,S^2,\delta,S,S,S,1, \widetilde{S},\delta S,\sigma)\\
~~~~~&(2,\si_{13},3)(6,\si_{21})(3,\si_{23},1)(2,\si,5)(\si_{12},6)(1,\si_{12},5)(5,\Delta^2,1)(2,\Delta^2,2)(\Delta^2,2)\\
=&(m,m)(1,\si,1)(\si,2)(m\si,m\si,m,m)(1,\si,5)(\si,6)(S^2,S^2,\delta,S,S,S,1, \widetilde{S},\sigma,\delta S)\\
~~~~~&(4,\si,3)(3,\si,4) (2,\si,5)(6,\si,1)(7,\si)(5,\si,2)(6,\si,1)(4,\si,3)(5,\si,2)\\
~~~~~&(3,\si,4)(4,\si,3)(2,\si,5)(1,\si,6)(\si,7)(2,\si,5)(1,\si,6)(5,\Delta^2,1)(2,\Delta^2,2)(\Delta^2,2)=\\
\end{align*}

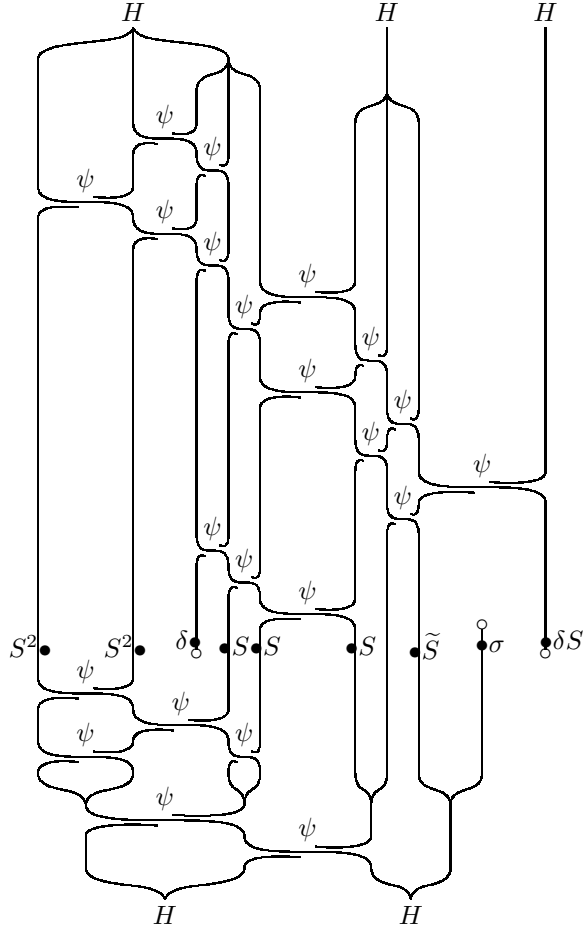
\[
\xy /r2pc/:,
(1,0)= "L11",
(5,0)= "L12",
(7.5,0)= "L13",
  "L11"+(0,.5)*!<0pt,6pt>{H}, 
  "L12"+(0,.5)*!<0pt,6pt>{H}, 
  "L13"+(0,.5)*!<0pt,6pt>{H}, 
  "L11"="x1", 
  "x1"+(-1.5,-.5)="x2",
  "x1"+(0,-.5)="x3",
  "x1"+(1.5,-.5)="x4",
  "x1"+(0,-.25)="c",  
   "x1";"x2"**\crv{"x1" & "c" & "x2"+(0,.25)}, 
   "x1";"x3"**\crv{"x1" & "c" & "x3"+(0,.25)}, 
   "x1";"x4"**\crv{"x1" & "c" & "x4"+(0,.25)}, 
"L12";"L12" +(0,-.5)="L14"**\dir{-},
"L13";"L13" +(0,-.5)="L15"**\dir{-},
"x2"="L11",
"x3"="L12",
"x4"="L13",
"L11";"L11" +(0,-.5)="L11"**\dir{-},
"L12";"L12" +(0,-.5)="L12"**\dir{-},
  "L13"="x1", 
  "x1"+(-.5,-.5)="x2",
  "x1"+(0,-.5)="x3",
  "x1"+(.5,-.5)="x4",
  "x1"+(0,-.25)="c",  
   "x1";"x2"**\crv{"x1" & "c" & "x2"+(0,.25)}, 
   "x1";"x3"**\crv{"x1" & "c" & "x3"+(0,.25)}, 
   "x1";"x4"**\crv{"x1" & "c" & "x4"+(0,.25)}, 
"x2"="L23",
"x3"="L24",
"x4"="L25",
"L14";"L14" +(0,-.5)="L16"**\dir{-},
"L15";"L15" +(0,-.5)="L17"**\dir{-},
"L23"="L13",
"L24"="L14",
"L25"="L15",
"L11";"L11" +(0,-.5)="L11"**\dir{-},
"L12";"L12" +(0,-.5)="L12"**\dir{-},
"L13";"L13" +(0,-.5)="L13"**\dir{-},
"L14";"L14" +(0,-.5)="L14"**\dir{-},
"L15";"L15" +(0,-.5)="L15"**\dir{-},
  "L16"="x1", 
  "x1"+(-.5,-.5)="x2",
  "x1"+(0,-.5)="x3",
  "x1"+(.5,-.5)="x4",
  "x1"+(0,-.25)="c",  
   "x1";"x2"**\crv{"x1" & "c" & "x2"+(0,.25)}, 
   "x1";"x3"**\crv{"x1" & "c" & "x3"+(0,.25)}, 
   "x1";"x4"**\crv{"x1" & "c" & "x4"+(0,.25)}, 
"x2"="L26",
"x3"="L27",
"x4"="L28",
"L17";"L17" +(0,-.5)="L19"**\dir{-},
"L26"="L16",
"L27"="L17",
"L28"="L18",
  "L11";"L11" +(0,-.5)="L11"**\dir{-},
  "L12"="x1",       
  "L13"="x2",  , 
  "x1"+(0,-.5)="x3",
  "x2"+(0,-.5)="x4",
  "x1";"x2"**\dir{}?(.5)="m", 
  "m" + (0,-.25)="c", 
  "x1";"x4"**\crv{ "x1"+(0,-.25) & "c" & "x4"+(0,.25)},  
  "x2";"c"+<3pt,2pt>="c2"**\crv{"x2"+(0,-.25) & "c2"},
  "x3";"c"+<-3pt,-2pt>="c1"**\crv{"x3"+(0,.25) & "c1"}, 
  "c"+(0,.10)*!<0pt,-6pt>{\si},
"x3"="L12",
"x4"="L13",
  "L14";"L14" +(0,-.5)="L14"**\dir{-},
  "L15";"L15" +(0,-.5)="L15"**\dir{-},
  "L16";"L16" +(0,-.5)="L16"**\dir{-},
  "L17";"L17" +(0,-.5)="L17"**\dir{-},
  "L18";"L18" +(0,-.5)="L18"**\dir{-},
  "L19";"L19" +(0,-.5)="L19"**\dir{-},
  "L11";"L11" +(0,-.5)="L11"**\dir{-},
  "L12";"L12" +(0,-.5)="L12"**\dir{-},
  "L13"="x1",       
  "L14"="x2",  , 
  "x1"+(0,-.5)="x3",
  "x2"+(0,-.5)="x4",
  "x1";"x2"**\dir{}?(.5)="m", 
  "m" + (0,-.25)="c", 
  "x1";"x4"**\crv{ "x1"+(0,-.25) & "c" & "x4"+(0,.25)},  
  "x2";"c"+<3pt,2pt>="c2"**\crv{"x2"+(0,-.25) & "c2"},
  "x3";"c"+<-3pt,-2pt>="c1"**\crv{"x3"+(0,.25) & "c1"}, 
  "c"+(0,.10)*!<0pt,-6pt>{\si},
"x3"="L13",
"x4"="L14",
  "L15";"L15" +(0,-.5)="L15"**\dir{-},
  "L16";"L16" +(0,-.5)="L16"**\dir{-},
  "L17";"L17" +(0,-.5)="L17"**\dir{-},
  "L18";"L18" +(0,-.5)="L18"**\dir{-},
  "L19";"L19" +(0,-.5)="L19"**\dir{-},
  "L11"="x1",       
  "L12"="x2",  , 
  "x1"+(0,-.5)="x3",
  "x2"+(0,-.5)="x4",
  "x1";"x2"**\dir{}?(.5)="m", 
  "m" + (0,-.25)="c", 
  "x1";"x4"**\crv{ "x1"+(0,-.25) & "c" & "x4"+(0,.25)},  
  "x2";"c"+<3pt,2pt>="c2"**\crv{"x2"+(0,-.25) & "c2"},
  "x3";"c"+<-3pt,-2pt>="c1"**\crv{"x3"+(0,.25) & "c1"}, 
  "c"+(0,.10)*!<0pt,-6pt>{\si},
"x3"="L11",
"x4"="L12",
  "L13";"L13" +(0,-.5)="L13"**\dir{-},
  "L14";"L14" +(0,-.5)="L14"**\dir{-},
  "L15";"L15" +(0,-.5)="L15"**\dir{-},
  "L16";"L16" +(0,-.5)="L16"**\dir{-},
  "L17";"L17" +(0,-.5)="L17"**\dir{-},
  "L18";"L18" +(0,-.5)="L18"**\dir{-},
  "L19";"L19" +(0,-.5)="L19"**\dir{-},
  "L11";"L11" +(0,-.5)="L11"**\dir{-},
  "L12"="x1",       
  "L13"="x2",  , 
  "x1"+(0,-.5)="x3",
  "x2"+(0,-.5)="x4",
  "x1";"x2"**\dir{}?(.5)="m", 
  "m" + (0,-.25)="c", 
  "x1";"x4"**\crv{ "x1"+(0,-.25) & "c" & "x4"+(0,.25)},  
  "x2";"c"+<3pt,2pt>="c2"**\crv{"x2"+(0,-.25) & "c2"},
  "x3";"c"+<-3pt,-2pt>="c1"**\crv{"x3"+(0,.25) & "c1"}, 
  "c"+(0,.10)*!<0pt,-6pt>{\si},
"x3"="L12",
"x4"="L13",
  "L14";"L14" +(0,-.5)="L14"**\dir{-},
  "L15";"L15" +(0,-.5)="L15"**\dir{-},
  "L16";"L16" +(0,-.5)="L16"**\dir{-},
  "L17";"L17" +(0,-.5)="L17"**\dir{-},
  "L18";"L18" +(0,-.5)="L18"**\dir{-},
  "L19";"L19" +(0,-.5)="L19"**\dir{-},
  "L11";"L11" +(0,-.5)="L11"**\dir{-},
  "L12";"L12" +(0,-.5)="L12"**\dir{-},
  "L13"="x1",       
  "L14"="x2",  , 
  "x1"+(0,-.5)="x3",
  "x2"+(0,-.5)="x4",
  "x1";"x2"**\dir{}?(.5)="m", 
  "m" + (0,-.25)="c", 
  "x1";"x4"**\crv{ "x1"+(0,-.25) & "c" & "x4"+(0,.25)},  
  "x2";"c"+<3pt,2pt>="c2"**\crv{"x2"+(0,-.25) & "c2"},
  "x3";"c"+<-3pt,-2pt>="c1"**\crv{"x3"+(0,.25) & "c1"}, 
  "c"+(0,.10)*!<0pt,-6pt>{\si},
"x3"="L13",
"x4"="L14",
  "L15";"L15" +(0,-.5)="L15"**\dir{-},
  "L16";"L16" +(0,-.5)="L16"**\dir{-},
  "L17";"L17" +(0,-.5)="L17"**\dir{-},
  "L18";"L18" +(0,-.5)="L18"**\dir{-},
  "L19";"L19" +(0,-.5)="L19"**\dir{-},
  "L11";"L11" +(0,-.5)="L11"**\dir{-},
  "L12";"L12" +(0,-.5)="L12"**\dir{-},
  "L13";"L13" +(0,-.5)="L13"**\dir{-},
  "L14";"L14" +(0,-.5)="L14"**\dir{-},
  "L15"="x1",       
  "L16"="x2",  , 
  "x1"+(0,-.5)="x3",
  "x2"+(0,-.5)="x4",
  "x1";"x2"**\dir{}?(.5)="m", 
  "m" + (0,-.25)="c", 
  "x1";"x4"**\crv{ "x1"+(0,-.25) & "c" & "x4"+(0,.25)},  
  "x2";"c"+<3pt,2pt>="c2"**\crv{"x2"+(0,-.25) & "c2"},
  "x3";"c"+<-3pt,-2pt>="c1"**\crv{"x3"+(0,.25) & "c1"}, 
  "c"+(0,.10)*!<0pt,-6pt>{\si},
"x3"="L15",
"x4"="L16",
  "L17";"L17" +(0,-.5)="L17"**\dir{-},
  "L18";"L18" +(0,-.5)="L18"**\dir{-},
  "L19";"L19" +(0,-.5)="L19"**\dir{-},
  "L11";"L11" +(0,-.5)="L11"**\dir{-},
  "L12";"L12" +(0,-.5)="L12"**\dir{-},
  "L13";"L13" +(0,-.5)="L13"**\dir{-},
  "L14"="x1",       
  "L15"="x2",  , 
  "x1"+(0,-.5)="x3",
  "x2"+(0,-.5)="x4",
  "x1";"x2"**\dir{}?(.5)="m", 
  "m" + (0,-.25)="c", 
  "x1";"x4"**\crv{ "x1"+(0,-.25) & "c" & "x4"+(0,.25)},  
  "x2";"c"+<3pt,2pt>="c2"**\crv{"x2"+(0,-.25) & "c2"},
  "x3";"c"+<-3pt,-2pt>="c1"**\crv{"x3"+(0,.25) & "c1"}, 
  "c"+(0,.10)*!<0pt,-6pt>{\si},
"x3"="L14",
"x4"="L15",
  "L16";"L16" +(0,-.5)="L16"**\dir{-},
  "L17";"L17" +(0,-.5)="L17"**\dir{-},
  "L18";"L18" +(0,-.5)="L18"**\dir{-},
  "L19";"L19" +(0,-.5)="L19"**\dir{-},
  "L11";"L11" +(0,-.5)="L11"**\dir{-},
  "L12";"L12" +(0,-.5)="L12"**\dir{-},
  "L13";"L13" +(0,-.5)="L13"**\dir{-},
  "L14";"L14" +(0,-.5)="L14"**\dir{-},
  "L15";"L15" +(0,-.5)="L15"**\dir{-},
  "L16"="x1",       
  "L17"="x2",  , 
  "x1"+(0,-.5)="x3",
  "x2"+(0,-.5)="x4",
  "x1";"x2"**\dir{}?(.5)="m", 
  "m" + (0,-.25)="c", 
  "x1";"x4"**\crv{ "x1"+(0,-.25) & "c" & "x4"+(0,.25)},  
  "x2";"c"+<3pt,2pt>="c2"**\crv{"x2"+(0,-.25) & "c2"},
  "x3";"c"+<-3pt,-2pt>="c1"**\crv{"x3"+(0,.25) & "c1"}, 
  "c"+(0,.10)*!<0pt,-6pt>{\si},
"x3"="L16",
"x4"="L17",
  "L18";"L18" +(0,-.5)="L18"**\dir{-},
  "L19";"L19" +(0,-.5)="L19"**\dir{-},
  "L11";"L11" +(0,-.5)="L11"**\dir{-},
  "L12";"L12" +(0,-.5)="L12"**\dir{-},
  "L13";"L13" +(0,-.5)="L13"**\dir{-},
  "L14";"L14" +(0,-.5)="L14"**\dir{-},
  "L15"="x1",       
  "L16"="x2",  , 
  "x1"+(0,-.5)="x3",
  "x2"+(0,-.5)="x4",
  "x1";"x2"**\dir{}?(.5)="m", 
  "m" + (0,-.25)="c", 
  "x1";"x4"**\crv{ "x1"+(0,-.25) & "c" & "x4"+(0,.25)},  
  "x2";"c"+<3pt,2pt>="c2"**\crv{"x2"+(0,-.25) & "c2"},
  "x3";"c"+<-3pt,-2pt>="c1"**\crv{"x3"+(0,.25) & "c1"}, 
  "c"+(0,.10)*!<0pt,-6pt>{\si},
"x3"="L15",
"x4"="L16",
  "L17";"L17" +(0,-.5)="L17"**\dir{-},
  "L18";"L18" +(0,-.5)="L18"**\dir{-},
  "L19";"L19" +(0,-.5)="L19"**\dir{-},
  "L11";"L11" +(0,-.5)="L11"**\dir{-},
  "L12";"L12" +(0,-.5)="L12"**\dir{-},
  "L13";"L13" +(0,-.5)="L13"**\dir{-},
  "L14";"L14" +(0,-.5)="L14"**\dir{-},
  "L15";"L15" +(0,-.5)="L15"**\dir{-},
  "L16";"L16" +(0,-.5)="L16"**\dir{-},
  "L17"="x1",       
  "L18"="x2",  , 
  "x1"+(0,-.5)="x3",
  "x2"+(0,-.5)="x4",
  "x1";"x2"**\dir{}?(.5)="m", 
  "m" + (0,-.25)="c", 
  "x1";"x4"**\crv{ "x1"+(0,-.25) & "c" & "x4"+(0,.25)},  
  "x2";"c"+<3pt,2pt>="c2"**\crv{"x2"+(0,-.25) & "c2"},
  "x3";"c"+<-3pt,-2pt>="c1"**\crv{"x3"+(0,.25) & "c1"}, 
  "c"+(0,.10)*!<0pt,-6pt>{\si},
"x3"="L17",
"x4"="L18",
  "L19";"L19" +(0,-.5)="L19"**\dir{-},
  "L11";"L11" +(0,-.5)="L11"**\dir{-},
  "L12";"L12" +(0,-.5)="L12"**\dir{-},
  "L13";"L13" +(0,-.5)="L13"**\dir{-},
  "L14";"L14" +(0,-.5)="L14"**\dir{-},
  "L15";"L15" +(0,-.5)="L15"**\dir{-},
  "L16"="x1",       
  "L17"="x2",  , 
  "x1"+(0,-.5)="x3",
  "x2"+(0,-.5)="x4",
  "x1";"x2"**\dir{}?(.5)="m", 
  "m" + (0,-.25)="c", 
  "x1";"x4"**\crv{ "x1"+(0,-.25) & "c" & "x4"+(0,.25)},  
  "x2";"c"+<3pt,2pt>="c2"**\crv{"x2"+(0,-.25) & "c2"},
  "x3";"c"+<-3pt,-2pt>="c1"**\crv{"x3"+(0,.25) & "c1"}, 
  "c"+(0,.10)*!<0pt,-6pt>{\si},
"x3"="L16",
"x4"="L17",
  "L18";"L18" +(0,-.5)="L18"**\dir{-},
  "L19";"L19" +(0,-.5)="L19"**\dir{-},
  "L11";"L11" +(0,-.5)="L11"**\dir{-},
  "L12";"L12" +(0,-.5)="L12"**\dir{-},
  "L13";"L13" +(0,-.5)="L13"**\dir{-},
  "L14";"L14" +(0,-.5)="L14"**\dir{-},
  "L15";"L15" +(0,-.5)="L15"**\dir{-},
  "L16";"L16" +(0,-.5)="L16"**\dir{-},
  "L17";"L17" +(0,-.5)="L17"**\dir{-},
  "L18"="x1",       
  "L19"="x2",  , 
  "x1"+(0,-.5)="x3",
  "x2"+(0,-.5)="x4",
  "x1";"x2"**\dir{}?(.5)="m", 
  "m" + (0,-.25)="c", 
  "x1";"x4"**\crv{ "x1"+(0,-.25) & "c" & "x4"+(0,.25)},  
  "x2";"c"+<3pt,2pt>="c2"**\crv{"x2"+(0,-.25) & "c2"},
  "x3";"c"+<-3pt,-2pt>="c1"**\crv{"x3"+(0,.25) & "c1"}, 
  "c"+(0,.10)*!<0pt,-6pt>{\si},
"x3"="L18",
"x4"="L19",
  "L11";"L11" +(0,-.5)="L11"**\dir{-},
  "L12";"L12" +(0,-.5)="L12"**\dir{-},
  "L13";"L13" +(0,-.5)="L13"**\dir{-},
  "L14";"L14" +(0,-.5)="L14"**\dir{-},
  "L15";"L15" +(0,-.5)="L15"**\dir{-},
  "L16";"L16" +(0,-.5)="L16"**\dir{-},
  "L17"="x1",       
  "L18"="x2",  , 
  "x1"+(0,-.5)="x3",
  "x2"+(0,-.5)="x4",
  "x1";"x2"**\dir{}?(.5)="m", 
  "m" + (0,-.25)="c", 
  "x1";"x4"**\crv{ "x1"+(0,-.25) & "c" & "x4"+(0,.25)},  
  "x2";"c"+<3pt,2pt>="c2"**\crv{"x2"+(0,-.25) & "c2"},
  "x3";"c"+<-3pt,-2pt>="c1"**\crv{"x3"+(0,.25) & "c1"}, 
  "c"+(0,.10)*!<0pt,-6pt>{\si},
"x3"="L17",
"x4"="L18",
  "L19";"L19" +(0,-.5)="L19"**\dir{-},
  "L11";"L11" +(0,-.5)="L11"**\dir{-},
  "L12";"L12" +(0,-.5)="L12"**\dir{-},
  "L13"="x1",       
  "L14"="x2",  , 
  "x1"+(0,-.5)="x3",
  "x2"+(0,-.5)="x4",
  "x1";"x2"**\dir{}?(.5)="m", 
  "m" + (0,-.25)="c", 
  "x1";"x4"**\crv{ "x1"+(0,-.25) & "c" & "x4"+(0,.25)},  
  "x2";"c"+<3pt,2pt>="c2"**\crv{"x2"+(0,-.25) & "c2"},
  "x3";"c"+<-3pt,-2pt>="c1"**\crv{"x3"+(0,.25) & "c1"}, 
  "c"+(0,.10)*!<0pt,-6pt>{\si},
"x3"="L13",
"x4"="L14",
  "L15";"L15" +(0,-.5)="L15"**\dir{-},
  "L16";"L16" +(0,-.5)="L16"**\dir{-},
  "L17";"L17" +(0,-.5)="L17"**\dir{-},
  "L18";"L18" +(0,-.5)="L18"**\dir{-},
  "L19";"L19" +(0,-.5)="L19"**\dir{-},
  "L11";"L11" +(0,-.5)="L11"**\dir{-},
  "L12";"L12" +(0,-.5)="L12"**\dir{-},
  "L13";"L13" +(0,-.5)="L13"**\dir{-},
  "L14"="x1",       
  "L15"="x2",  , 
  "x1"+(0,-.5)="x3",
  "x2"+(0,-.5)="x4",
  "x1";"x2"**\dir{}?(.5)="m", 
  "m" + (0,-.25)="c", 
  "x1";"x4"**\crv{ "x1"+(0,-.25) & "c" & "x4"+(0,.25)},  
  "x2";"c"+<3pt,2pt>="c2"**\crv{"x2"+(0,-.25) & "c2"},
  "x3";"c"+<-3pt,-2pt>="c1"**\crv{"x3"+(0,.25) & "c1"}, 
  "c"+(0,.10)*!<0pt,-6pt>{\si},
"x3"="L14",
"x4"="L15",
  "L16";"L16" +(0,-.5)="L16"**\dir{-},
  "L17";"L17" +(0,-.5)="L17"**\dir{-},
  "L18";"L18" +(0,-.5)="L18"**\dir{-},
  "L19";"L19" +(0,-.5)="L19"**\dir{-},
  "L11";"L11" +(0,-.5)="L11"**\dir{-},
  "L12";"L12" +(0,-.5)="L12"**\dir{-},
  "L13";"L13" +(0,-.5)="L13"**\dir{-},
  "L14";"L14" +(0,-.5)="L14"**\dir{-},
  "L15"="x1",       
  "L16"="x2",  , 
  "x1"+(0,-.5)="x3",
  "x2"+(0,-.5)="x4",
  "x1";"x2"**\dir{}?(.5)="m", 
  "m" + (0,-.25)="c", 
  "x1";"x4"**\crv{ "x1"+(0,-.25) & "c" & "x4"+(0,.25)},  
  "x2";"c"+<3pt,2pt>="c2"**\crv{"x2"+(0,-.25) & "c2"},
  "x3";"c"+<-3pt,-2pt>="c1"**\crv{"x3"+(0,.25) & "c1"}, 
  "c"+(0,.10)*!<0pt,-6pt>{\si},
"x3"="L15",
"x4"="L16",
  "L17";"L17" +(0,-.5)="L17"**\dir{-},
  "L18";"L18" +(0,-.5)="L18"**\dir{-},
  "L19";"L19" +(0,-.5)="L19"**\dir{-},
  "L11"-(0,.25)*!<3pt,0pt>{S^2 \bullet }, 
  "L11";"L11" +(0,-.5)="L11"**\dir{-},
  "L12"-(0,.25)*!<3pt,0pt>{S^2 \bullet }, 
  "L12";"L12" +(0,-.5)="L12"**\dir{-},
  "L13"-(0,.15)*!<3pt,0pt>{\delta \bullet}, 
  "L13";"L13" +(0,-.3)="L13"**\dir{-},
  "L13"+(0,0)*!<0pt,2pt>{\circ},  
  "L14"-(0,.25)*!<-2pt,0pt>{\bullet S }, 
  "L14";"L14" +(0,-.5)="L14"**\dir{-}, 
  "L15"-(0,.25)*!<-2pt,0pt>{\bullet S }, 
  "L15";"L15" +(0,-.5)="L15"**\dir{-}, 
  "L16"-(0,.25)*!<-2pt,0pt>{\bullet S }, 
  "L16";"L16" +(0,-.5)="L16"**\dir{-}, 
  "L17";"L17" +(0,-.5)="L17"**\dir{-}, 
  "L18"-(0,.25)*!<-2pt,0pt>{\bullet \widetilde{S} }, 
  "L18";"L18" +(0,-.5)="L18"**\dir{-}, 
   "L18" + (1, .5)="x1", 
  "x1"+(0,-.5)="x2",
  "x1";"x2"**\dir{}?(.5)="c", 
  "x1";"x2"**\dir{-},   
  "c"-(0,0)*!<-3pt,0pt>{\bullet \sigma },  
  "x1"+(0,0)*!<0pt,-2pt>{\circ}, 
"x2"="L29",
  "L19"-(0,.15)*!<-6pt,0pt>{\bullet \delta S}, 
  "L19";"L19" +(0,-.3)="L110"**\dir{-},
  "L110"+(0,0)*!<0pt,2pt>{\circ}, 
"L29"="L19",
  "L11";"L11" +(0,-.25)="L11"**\dir{-},
  "L12";"L12" +(0,-.25)="L12"**\dir{-},
  "L14";"L14" +(0,-.25)="L14"**\dir{-},
  "L15";"L15" +(0,-.25)="L15"**\dir{-},
  "L16";"L16" +(0,-.25)="L16"**\dir{-},
  "L17";"L17" +(0,-.25)="L17"**\dir{-},
  "L18";"L18" +(0,-.25)="L18"**\dir{-},
  "L19";"L19" +(0,-.25)="L19"**\dir{-},
  "L11"="x1",       
  "L12"="x2",  , 
  "x1"+(0,-.5)="x3",
  "x2"+(0,-.5)="x4",
  "x1";"x2"**\dir{}?(.5)="m", 
  "m" + (0,-.25)="c", 
  "x1";"x4"**\crv{ "x1"+(0,-.25) & "c" & "x4"+(0,.25)},  
  "x2";"c"+<3pt,2pt>="c2"**\crv{"x2"+(0,-.25) & "c2"},
  "x3";"c"+<-3pt,-2pt>="c1"**\crv{"x3"+(0,.25) & "c1"}, 
  "c"+(0,.10)*!<0pt,-6pt>{\si},
"x3"="L11",
"x4"="L12",
  "L14";"L14" +(0,-.5)="L13"**\dir{-},
  "L15";"L15" +(0,-.5)="L14"**\dir{-},
  "L16";"L16" +(0,-.5)="L15"**\dir{-},
  "L17";"L17" +(0,-.5)="L16"**\dir{-},
  "L18";"L18" +(0,-.5)="L17"**\dir{-},
  "L19";"L19" +(0,-.5)="L18"**\dir{-},
  "L11";"L11" +(0,-.5)="L11"**\dir{-},
  "L12"="x1",       
  "L13"="x2",  , 
  "x1"+(0,-.5)="x3",
  "x2"+(0,-.5)="x4",
  "x1";"x2"**\dir{}?(.5)="m", 
  "m" + (0,-.25)="c", 
  "x1";"x4"**\crv{ "x1"+(0,-.25) & "c" & "x4"+(0,.25)},  
  "x2";"c"+<3pt,2pt>="c2"**\crv{"x2"+(0,-.25) & "c2"},
  "x3";"c"+<-3pt,-2pt>="c1"**\crv{"x3"+(0,.25) & "c1"}, 
  "c"+(0,.10)*!<0pt,-6pt>{\si},
"x3"="L12",
"x4"="L13",
  "L14";"L14" +(0,-.5)="L14"**\dir{-},
  "L15";"L15" +(0,-.5)="L15"**\dir{-},
  "L16";"L16" +(0,-.5)="L16"**\dir{-},
  "L17";"L17" +(0,-.5)="L17"**\dir{-},
  "L18";"L18" +(0,-.5)="L18"**\dir{-},
  "L11"="x1",       
  "L12"="x2",  , 
  "x1"+(0,-.5)="x3",
  "x2"+(0,-.5)="x4",
  "x1";"x2"**\dir{}?(.5)="m", 
  "m" + (0,-.25)="c", 
  "x1";"x4"**\crv{ "x1"+(0,-.25) & "c" & "x4"+(0,.25)},  
  "x2";"c"+<3pt,2pt>="c2"**\crv{"x2"+(0,-.25) & "c2"},
  "x3";"c"+<-3pt,-2pt>="c1"**\crv{"x3"+(0,.25) & "c1"}, 
  "c"+(0,.10)*!<0pt,-6pt>{\si},
"x3"="L11",
"x4"="L12",
  "L13"="x1",       
  "L14"="x2",  , 
  "x1"+(0,-.5)="x3",
  "x2"+(0,-.5)="x4",
  "x1";"x2"**\dir{}?(.5)="m", 
  "m" + (0,-.25)="c", 
  "x1";"x4"**\crv{ "x1"+(0,-.25) & "c" & "x4"+(0,.25)},  
  "x2";"c"+<3pt,2pt>="c2"**\crv{"x2"+(0,-.25) & "c2"},
  "x3";"c"+<-3pt,-2pt>="c1"**\crv{"x3"+(0,.25) & "c1"}, 
  "c"+(0,.10)*!<0pt,-6pt>{\si},
"x3"="L13",
"x4"="L14",
  "L15";"L15" +(0,-.5)="L15"**\dir{-},
  "L16";"L16" +(0,-.5)="L16"**\dir{-},
  "L17";"L17" +(0,-.5)="L17"**\dir{-},
  "L18";"L18" +(0,-.5)="L18"**\dir{-},
   "L11"="x1", 
  "L12"="x2",  
  "x1";"x2"**\dir{}?(.5)="m", 
  "m" + (0,-.25)="c", 
  "m" + (0,-.5)="x3", 
  "x1";"x3"**\crv{"x1"-(0,.25) & "c" & "x3"},
  "x2";"x3"**\crv{"x2"-(0,.25) & "c"  & "x3"}, 
"x3"="L11",
   "L13"="x1", 
   "L14"="x2",   
  "x1";"x2"**\dir{}?(.5)="m", 
  "m" + (0,-.25)="c", 
  "m" + (0,-.5)="x3", 
  "x1";"x3"**\crv{"x1"-(0,.25) & "c" & "x3"},
  "x2";"x3"**\crv{"x2"-(0,.25) & "c"  & "x3"}, 
"x3"="L12",
   "L15"="x1", 
  "L16"="x2",   
  "x1";"x2"**\dir{}?(.5)="m", 
  "m" + (0,-.25)="c", 
  "m" + (0,-.5)="x3", 
  "x1";"x3"**\crv{"x1"-(0,.25) & "c" & "x3"},
  "x2";"x3"**\crv{"x2"-(0,.25) & "c"  & "x3"}, 
"x3"="L13",
   "L17"="x1", 
  "L18"="x2",  
  "x1";"x2"**\dir{}?(.5)="m", 
  "m" + (0,-.25)="c", 
  "m" + (0,-.5)="x3", 
  "x1";"x3"**\crv{"x1"-(0,.25) & "c" & "x3"},
  "x2";"x3"**\crv{"x2"-(0,.25) & "c"  & "x3"}, 
"x3"="L14",
  "L11"="x1",       
  "L12"="x2",  , 
  "x1"+(0,-.5)="x3",
  "x2"+(0,-.5)="x4",
  "x1";"x2"**\dir{}?(.5)="m", 
  "m" + (0,-.25)="c", 
  "x1";"x4"**\crv{ "x1"+(0,-.25) & "c" & "x4"+(0,.25)},  
  "x2";"c"+<3pt,2pt>="c2"**\crv{"x2"+(0,-.25) & "c2"},
  "x3";"c"+<-3pt,-2pt>="c1"**\crv{"x3"+(0,.25) & "c1"}, 
  "c"+(0,.10)*!<0pt,-6pt>{\si},
"x3"="L11",
"x4"="L12",
  "L13";"L13" +(0,-.5)="L13"**\dir{-},
  "L14";"L14" +(0,-.5)="L14"**\dir{-},
  "L11";"L11" +(0,-.5)="L11"**\dir{-},
  "L12"="x1",       
  "L13"="x2",  , 
  "x1"+(0,-.5)="x3",
  "x2"+(0,-.5)="x4",
  "x1";"x2"**\dir{}?(.5)="m", 
  "m" + (0,-.25)="c", 
  "x1";"x4"**\crv{ "x1"+(0,-.25) & "c" & "x4"+(0,.25)},  
  "x2";"c"+<3pt,2pt>="c2"**\crv{"x2"+(0,-.25) & "c2"},
  "x3";"c"+<-3pt,-2pt>="c1"**\crv{"x3"+(0,.25) & "c1"}, 
  "c"+(0,.10)*!<0pt,-6pt>{\si},
"x3"="L12",
"x4"="L13",
  "L14";"L14" +(0,-.5)="L14"**\dir{-},
   "L11"="x1", 
  "L12"="x2",  
  "x1";"x2"**\dir{}?(.5)="m", 
  "m" + (0,-.25)="c", 
  "m" + (0,-.5)="x3", 
  "x1";"x3"**\crv{"x1"-(0,.25) & "c" & "x3"},
  "x2";"x3"**\crv{"x2"-(0,.25) & "c"  & "x3"}, 
"x3"="L11",
   "L13"="x1", 
  "L14"="x2",   
  "x1";"x2"**\dir{}?(.5)="m", 
  "m" + (0,-.25)="c", 
  "m" + (0,-.5)="x3", 
  "x1";"x3"**\crv{"x1"-(0,.25) & "c" & "x3"},
  "x2";"x3"**\crv{"x2"-(0,.25) & "c"  & "x3"}, 
"x3"="L12",
  "L11"-(0,.5)*!<0pt,-6pt>{H}, 
  "L12"-(0,.5)*!<0pt,-6pt>{H}, 
\endxy 
\]

\begin{align*} 
=&(m,m)(1,\si,1)(\si,2)(m\si,m\si,m,m)(1,\si,5)(\si,6)(S^2,S^2,\delta,S,S,S,1, \widetilde{S},\delta S,\sigma) \\
~~~~~&(4,\si,3)(3,\si,4) (2,\si,5)(6,\si,1)(5,\si,2)(4,\si,3)(3,\si,4)(7,\si)(6,\si,1)\\
~~~~~&(5,\si,2)(4,\si,3)(2,\si,5)(1,\si,6)(\si,7)(2,\si,5)(1,\si,6)(5,\Delta^2,1)(2,\Delta^2,2)(\Delta^2,2)\\
=&(m,m)(1,\si,1)(\si,2)(m\si,m\si,m,m)(1,\si,5)(\si,6)\\
~~~~~&(S^2,S^2,\delta,S,S,S,1,\widetilde{S},\delta S,\sigma)(4,\si,3)(3,\si,4) (2,\si,5)(6,\si,1)\\
~~~~~&(5,\si,2)(4,\si,3)(3,\si,4)(2,\si,5)(1,\si,6)(\si,7)(2,\si,5)(1,\si,6)\\
~~~~~&(7,\si)(6,\si,1)(5,\si,2)(4,\si,3)(5,\Delta^2,1)(2,\Delta^2,2)(\Delta^2,2)\\
=&(m,m)(1,\si,1)(\si,2)(m,m,m,m)(\si,6)(2,\si,4)(1,\si,5)(\si,6)(3,\si,3)(2,\si,4)\\
~~~~~&(5,\si,1)(4,\si,2)(3,\si,3)(2,\si,4)(1,\si,5)(S^2,\widetilde{S},S^2,S,\delta,S,S,1,\sigma,\delta S)(\si,7)\\
~~~~~&(2,\si,5)(1,\si,6)(7,\si)(6,\si,1)(5,\si,2)(4,\si,3)(5,\Delta^2,1)(2,\Delta^2,2)(\Delta^2,2)\\
=&(m,m)(1,\si,1)(\si,2)(m,m,m,m)(\si,6)(2,\si,4)(1,\si,5)(3,\si,3)(2,\si,4)(5,\si,1)\\
~~~~~&(4,\si,2)(3,\si,3)(2,\si,4)(\si,6)(1,\si,5)(S^2,(\delta,S)\Delta,S^2,S,\delta,S,S,1,\sigma,\delta S)\\
~~~~~&(\si,7)(2,\si,5)(1,\si,6)(7,\si)(6,\si,1)(5,\si,2)(4,\si,3)(5,\Delta^2,1)(2,\Delta^2,2)(\Delta^2,2)\\
=&(m,m)(1,\si,1)(\si,2)(m,m,m,m)(\si,6)(2,\si,4)(1,\si,5)(3,\si,3)(2,\si,4)(5,\si,1)\\
~~~~~&(4,\si,2)(3,\si,3)(2,\si,4)(\si,6)(1,\si,5)(S^2,\delta,S,S^2,S,\delta,S,S,1,\sigma,\delta S)(1,\Delta,7)\\
~~~~~&(\si,7)(2,\si,5)(1,\si,6)(7,\si)(6,\si,1)(5,\si,2)(4,\si,3)(5,\Delta^2,1)(2,\Delta^2,2)(\Delta^2,2)\\
=&(1,m)(m,2)(\si_{12},1)(m,m,m,m)(\si,6)(2,\si,4)(1,\si,5)(3,\si,3)(2,\si,4)\\
~~~~~&(5,\si,1)(4,\si,2)(3,\si,3)(2,\si,4)(\si,6)(1,\si,5)(S^2,\delta,S,S^2,S,\delta,S,S,1,\sigma,\delta S)\\
~~~~~&(\si_{21},7)(3,\si,5)(2,\si,6)(8,\si)(7,\si,1)(6,\si,2)(5,\si,3)(\Delta,8)(5,\Delta,2)\\
~~~~~&(5,\Delta,1)(2,\Delta,3)(2,\Delta,2)(\Delta,3)(\Delta,2)\\
=&(1,m)(\si,1)(1,m,1)(m,2,m)(2,m,3)(4,m,2)(\si,6)(2,\si,4)(1,\si,5)(3,\si,3)\\
~~~~~&(2,\si,4)(5,\si,1)(4,\si,2)(3,\si,3)(2,\si,4)(\si,6)(1,\si,5)\\
~~~~~&(S^2,\delta,S,S^2,S,\delta,S,S,1,\sigma,\delta S)(\si_{21},7)(3,\si,5)(2,\si,6)(8,\si)(7,\si,1)\\
~~~~~&(5,\si_{12},2)(6,\Delta,2)(\Delta,7)(5,\Delta,1)(2,\Delta,3)(2,\Delta,2)(\Delta,3)(\Delta,2)\\
=&(1,m)(\si,1)(m,1,m)(2,m,2)(2,m,3)(4,m,2)(\si,6)(2,\si,4)(1,\si,5)(3,\si,3)\\
~~~~~&(2,\si,4)(5,\si,1)(4,\si,2)(3,\si,3)(2,\si,4)(\si,6)(1,\si,5)\\
~~~~~&(S^2,\delta,S,S^2,S,\delta,S,S,1,\sigma,\delta S)(\si_{21},7)(3,\si,5)(2,\si,6)(8,\si)(7,\si,1)(5,\si_{12},2)\\
~~~~~&(6,\Delta,2)(\Delta,7)(5,\Delta,1)(2,\Delta,3)(2,\Delta,2)(\Delta,3)(\Delta,2)\\
=&(1,m)(\si,1)(m,1,m)(2,m,2)(3,m,2)(4,m,2)(\si,6)(2,\si,4)(1,\si,5)(3,\si,3)\\
~~~~~&(2,\si,4)(5,\si,1)(4,\si,2)(3,\si,3)(2,\si,4)(\si,6)(1,\si,5)\\
~~~~~&(S^2,\delta,S,S^2,S,\delta,S,S,1,\sigma,\delta S)(\si,8)(5,\Delta,3)(1,\si,6)(3,\si,4)(2,\si,5)\\
~~~~~&(7,\si)(6,\si,1)(5,\si,2)(\Delta,7)(5,\Delta,1)(2,\Delta,3)(2,\Delta,2)(\Delta,3)(\Delta,2)\\
=&(1,m)(\si,1)(m,1,m)(2,m,2)(3,m,2)(3,m,3)(\si,6)(2,\si,4)(1,\si,5)(3,\si,3)\\
~~~~~&(2,\si,4)(5,\si,1)(4,\si,2)(3,\si,3)(2,\si,4)(\si,6)(1,\si,5)\\
~~~~~&(\delta,S^2,S,S^2,S,(\delta,S)\Delta,S,1,\sigma,\delta S)(1,\si,6)(3,\si,4)(2,\si,5)(7,\si)\\
~~~~~&(6,\si,1)(5,\si,2)(\Delta,7)(5,\Delta,1)(2,\Delta,3)(2,\Delta,2)(\Delta,3)(\Delta,2)=
\end{align*}

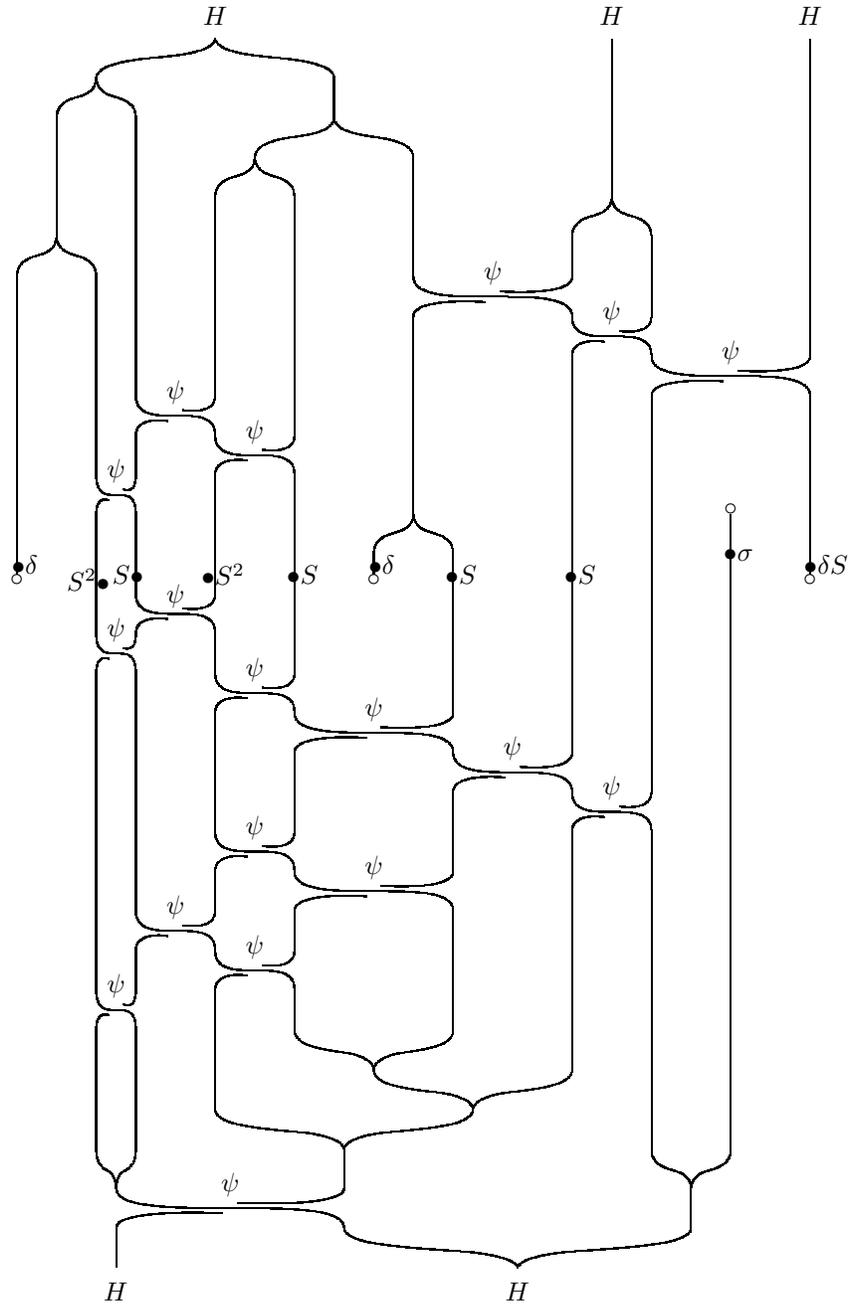
\[
\xy /r2.5pc/:,
(1,0)= "L11",
(6,0)= "L12",
(8.5,0)= "L13",
  "L11"+(0,.5)*!<0pt,6pt>{H}, 
  "L12"+(0,.5)*!<0pt,6pt>{H}, 
  "L13"+(0,.5)*!<0pt,6pt>{H}, 
   "L11"="x1", 
  "x1"+(-1.5,-.5)="x2",
  "x1"+(1.5,-.5)="x3",
  "x1"+(0,-.25)="c",  
   "x1";"x2"**\crv{"x1" & "c" & "x2"+(0,.25)}, 
   "x1";"x3"**\crv{"x1" & "c" & "x3"+(0,.25)}, 
"x2"= "L21",
"x3"= "L22",
  "L12";"L12" +(0,-.5)="L23"**\dir{-},
  "L13";"L13" +(0,-.5)="L24"**\dir{-},
"L21"= "L11",
"L22"= "L12",
"L23"= "L13",
"L24"= "L14",
   "L11"="x1", 
  "x1"+(-.5,-.5)="x2",
  "x1"+(.5,-.5)="x3",
  "x1"+(0,-.25)="c",  
   "x1";"x2"**\crv{"x1" & "c" & "x2"+(0,.25)}, 
   "x1";"x3"**\crv{"x1" & "c" & "x3"+(0,.25)}, 
"x2"= "L21",
"x3"= "L22",
  "L12";"L12" +(0,-.5)="L23"**\dir{-},
  "L13";"L13" +(0,-.5)="L24"**\dir{-},
  "L14";"L14" +(0,-.5)="L25"**\dir{-},
"L21"= "L11",
"L22"= "L12",
"L23"= "L13",
"L24"= "L14",
"L25"= "L15",
  "L11";"L11" +(0,-.5)="L11"**\dir{-},
  "L12";"L12" +(0,-.5)="L12"**\dir{-},
   "L13"="x1", 
  "x1"+(-1,-.5)="x2",
  "x1"+(1,-.5)="x3",
  "x1"+(0,-.25)="c",  
   "x1";"x2"**\crv{"x1" & "c" & "x2"+(0,.25)}, 
   "x1";"x3"**\crv{"x1" & "c" & "x3"+(0,.25)}, 
"x2"= "L23",
"x3"= "L24",
  "L14";"L14" +(0,-.5)="L25"**\dir{-},
  "L15";"L15" +(0,-.5)="L26"**\dir{-},
"L23"= "L13",
"L24"= "L14",
"L25"= "L15",
"L26"= "L16",
  "L11";"L11" +(0,-.5)="L11"**\dir{-},
  "L12";"L12" +(0,-.5)="L12"**\dir{-},
   "L13"="x1", 
  "x1"+(-.5,-.5)="x2",
  "x1"+(.5,-.5)="x3",
  "x1"+(0,-.25)="c",  
   "x1";"x2"**\crv{"x1" & "c" & "x2"+(0,.25)}, 
   "x1";"x3"**\crv{"x1" & "c" & "x3"+(0,.25)}, 
"x2"= "L23",
"x3"= "L24",
  "L14";"L14" +(0,-.5)="L25"**\dir{-},
  "L15";"L15" +(0,-.5)="L26"**\dir{-},
  "L16";"L16" +(0,-.5)="L27"**\dir{-},
"L23"= "L13",
"L24"= "L14",
"L25"= "L15",
"L26"= "L16",
"L27"= "L17",
  "L11";"L11" +(0,-.5)="L11"**\dir{-},
  "L12";"L12" +(0,-.5)="L12"**\dir{-},
  "L13";"L13" +(0,-.5)="L13"**\dir{-},
  "L14";"L14" +(0,-.5)="L14"**\dir{-},
  "L15";"L15" +(0,-.5)="L15"**\dir{-},
   "L16"="x1", 
  "x1"+(-.5,-.5)="x2",
  "x1"+(.5,-.5)="x3",
  "x1"+(0,-.25)="c",  
   "x1";"x2"**\crv{"x1" & "c" & "x2"+(0,.25)}, 
   "x1";"x3"**\crv{"x1" & "c" & "x3"+(0,.25)}, 
"x2"= "L26",
"x3"= "L27",
  "L17";"L17" +(0,-.5)="L28"**\dir{-},
"L26"= "L16",
"L27"= "L17",
"L28"= "L18",
   "L11"="x1", 
  "x1"+(-.5,-.5)="x2",
  "x1"+(.5,-.5)="x3",
  "x1"+(0,-.25)="c",  
   "x1";"x2"**\crv{"x1" & "c" & "x2"+(0,.25)}, 
   "x1";"x3"**\crv{"x1" & "c" & "x3"+(0,.25)}, 
"x2"= "L21",
"x3"= "L22",
  "L12";"L12" +(0,-.5)="L23"**\dir{-},
  "L13";"L13" +(0,-.5)="L24"**\dir{-},
  "L14";"L14" +(0,-.5)="L25"**\dir{-},
  "L15";"L15" +(0,-.5)="L26"**\dir{-},
  "L16";"L16" +(0,-.5)="L27"**\dir{-},
  "L17";"L17" +(0,-.5)="L28"**\dir{-},
  "L18";"L18" +(0,-.5)="L29"**\dir{-},
"L21"= "L11",
"L22"= "L12",
"L23"= "L13",
"L24"= "L14",
"L25"= "L15",
"L26"= "L16",
"L27"= "L17",
"L28"= "L18",
"L29"= "L19",
  "L11";"L11" +(0,-.5)="L11"**\dir{-},
  "L12";"L12" +(0,-.5)="L12"**\dir{-},
  "L13";"L13" +(0,-.5)="L13"**\dir{-},
  "L14";"L14" +(0,-.5)="L14"**\dir{-},
  "L15";"L15" +(0,-.5)="L15"**\dir{-},
  "L16"="x1",       
  "L17"="x2",  , 
  "x1"+(0,-.5)="x3",
  "x2"+(0,-.5)="x4",
  "x1";"x2"**\dir{}?(.5)="m", 
  "m" + (0,-.25)="c", 
  "x1";"x4"**\crv{ "x1"+(0,-.25) & "c" & "x4"+(0,.25)},  
  "x2";"c"+<3pt,2pt>="c2"**\crv{"x2"+(0,-.25) & "c2"},
  "x3";"c"+<-3pt,-2pt>="c1"**\crv{"x3"+(0,.25) & "c1"}, 
  "c"+(0,.10)*!<0pt,-6pt>{\si},
"x3"="L16",
"x4"="L17",
  "L18";"L18" +(0,-.5)="L18"**\dir{-},
  "L19";"L19" +(0,-.5)="L19"**\dir{-},
  "L11";"L11" +(0,-.5)="L11"**\dir{-},
  "L12";"L12" +(0,-.5)="L12"**\dir{-},
  "L13";"L13" +(0,-.5)="L13"**\dir{-},
  "L14";"L14" +(0,-.5)="L14"**\dir{-},
  "L15";"L15" +(0,-.5)="L15"**\dir{-},
  "L16";"L16" +(0,-.5)="L16"**\dir{-},
  "L17"="x1",       
  "L18"="x2",  , 
  "x1"+(0,-.5)="x3",
  "x2"+(0,-.5)="x4",
  "x1";"x2"**\dir{}?(.5)="m", 
  "m" + (0,-.25)="c", 
  "x1";"x4"**\crv{ "x1"+(0,-.25) & "c" & "x4"+(0,.25)},  
  "x2";"c"+<3pt,2pt>="c2"**\crv{"x2"+(0,-.25) & "c2"},
  "x3";"c"+<-3pt,-2pt>="c1"**\crv{"x3"+(0,.25) & "c1"}, 
  "c"+(0,.10)*!<0pt,-6pt>{\si},
"x3"="L17",
"x4"="L18",
  "L19";"L19" +(0,-.5)="L19"**\dir{-},
  "L11";"L11" +(0,-.5)="L11"**\dir{-},
  "L12";"L12" +(0,-.5)="L12"**\dir{-},
  "L13";"L13" +(0,-.5)="L13"**\dir{-},
  "L14";"L14" +(0,-.5)="L14"**\dir{-},
  "L15";"L15" +(0,-.5)="L15"**\dir{-},
  "L16";"L16" +(0,-.5)="L16"**\dir{-},
  "L17";"L17" +(0,-.5)="L17"**\dir{-},
  "L18"="x1",       
  "L19"="x2",  , 
  "x1"+(0,-.5)="x3",
  "x2"+(0,-.5)="x4",
  "x1";"x2"**\dir{}?(.5)="m", 
  "m" + (0,-.25)="c", 
  "x1";"x4"**\crv{ "x1"+(0,-.25) & "c" & "x4"+(0,.25)},  
  "x2";"c"+<3pt,2pt>="c2"**\crv{"x2"+(0,-.25) & "c2"},
  "x3";"c"+<-3pt,-2pt>="c1"**\crv{"x3"+(0,.25) & "c1"}, 
  "c"+(0,.10)*!<0pt,-6pt>{\si},
"x3"="L18",
"x4"="L19",
  "L11";"L11" +(0,-.5)="L11"**\dir{-},
  "L12";"L12" +(0,-.5)="L12"**\dir{-},
  "L13"="x1",       
  "L14"="x2",  , 
  "x1"+(0,-.5)="x3",
  "x2"+(0,-.5)="x4",
  "x1";"x2"**\dir{}?(.5)="m", 
  "m" + (0,-.25)="c", 
  "x1";"x4"**\crv{ "x1"+(0,-.25) & "c" & "x4"+(0,.25)},  
  "x2";"c"+<3pt,2pt>="c2"**\crv{"x2"+(0,-.25) & "c2"},
  "x3";"c"+<-3pt,-2pt>="c1"**\crv{"x3"+(0,.25) & "c1"}, 
  "c"+(0,.10)*!<0pt,-6pt>{\si},
"x3"="L13",
"x4"="L14",
  "L15";"L15" +(0,-.5)="L15"**\dir{-},
  "L16";"L16" +(0,-.5)="L16"**\dir{-},
  "L17";"L17" +(0,-.5)="L17"**\dir{-},
  "L18";"L18" +(0,-.5)="L18"**\dir{-},
  "L19";"L19" +(0,-.5)="L19"**\dir{-},
  "L11";"L11" +(0,-.5)="L11"**\dir{-},
  "L12";"L12" +(0,-.5)="L12"**\dir{-},
  "L13";"L13" +(0,-.5)="L13"**\dir{-},
  "L14"="x1",       
  "L15"="x2",  , 
  "x1"+(0,-.5)="x3",
  "x2"+(0,-.5)="x4",
  "x1";"x2"**\dir{}?(.5)="m", 
  "m" + (0,-.25)="c", 
  "x1";"x4"**\crv{ "x1"+(0,-.25) & "c" & "x4"+(0,.25)},  
  "x2";"c"+<3pt,2pt>="c2"**\crv{"x2"+(0,-.25) & "c2"},
  "x3";"c"+<-3pt,-2pt>="c1"**\crv{"x3"+(0,.25) & "c1"}, 
  "c"+(0,.10)*!<0pt,-6pt>{\si},
"x3"="L14",
"x4"="L15",
  "L16";"L16" +(0,-.5)="L16"**\dir{-},
  "L17";"L17" +(0,-.5)="L17"**\dir{-},
  "L18";"L18" +(0,-.5)="L18"**\dir{-},
  "L19";"L19" +(0,-.5)="L19"**\dir{-},
  "L11";"L11" +(0,-.5)="L11"**\dir{-},
  "L12"="x1",       
  "L13"="x2",  , 
  "x1"+(0,-.5)="x3",
  "x2"+(0,-.5)="x4",
  "x1";"x2"**\dir{}?(.5)="m", 
  "m" + (0,-.25)="c", 
  "x1";"x4"**\crv{ "x1"+(0,-.25) & "c" & "x4"+(0,.25)},  
  "x2";"c"+<3pt,2pt>="c2"**\crv{"x2"+(0,-.25) & "c2"},
  "x3";"c"+<-3pt,-2pt>="c1"**\crv{"x3"+(0,.25) & "c1"}, 
  "c"+(0,.10)*!<0pt,-6pt>{\si},
"x3"="L12",
"x4"="L13",
  "L14";"L14" +(0,-.5)="L14"**\dir{-},
  "L15";"L15" +(0,-.5)="L15"**\dir{-},
  "L16";"L16" +(0,-.5)="L16"**\dir{-},
  "L17";"L17" +(0,-.5)="L17"**\dir{-},
  "L18";"L18" +(0,-.5)="L18"**\dir{-},
  "L19";"L19" +(0,-.5)="L19"**\dir{-},
  "L11";"L11" +(0,-.5)="L11"**\dir{-},
  "L12";"L12" +(0,-.5)="L12"**\dir{-},
  "L13";"L13" +(0,-.5)="L13"**\dir{-},
  "L14";"L14" +(0,-.5)="L14"**\dir{-},
  "L15";"L15" +(0,-.5)="L15"**\dir{-},
   "L16"="x1", 
  "x1"+(-.5,-.5)="x2",
  "x1"+(.5,-.5)="x3",
  "x1"+(0,-.25)="c",  
   "x1";"x2"**\crv{"x1" & "c" & "x2"+(0,.25)}, 
   "x1";"x3"**\crv{"x1" & "c" & "x3"+(0,.25)}, 
"x2"= "L26",
"x3"= "L27",
  "L17";"L17" +(0,-.5)="L28"**\dir{-},
  "L18";"L18" +(0,-.5)="L29"**\dir{-},
  "L19";"L19" +(0,-.5)="L210"**\dir{-},
"L26"= "L16",
"L27"= "L17",
"L28"= "L18",
"L29"= "L19",
"L210"= "L110",
  "L11"="x1", 
  "x1"+(0,-.25)="x2",
  "x1";"x2"**\dir{}?(.5)="c", 
  "x1";"x2"**\dir{-},   
"c"-(0,0)*!<-3pt,0pt>{\bullet \delta  },  
"x2"+(0,0)*!<0pt,2pt>{\circ},  
  "L12"-(0,.25)*!<3pt,2pt>{S^2  \bullet }, 
  "L12";"L12" +(0,-.5)="L11"**\dir{-}, 
  "L13"-(0,.25)*!<3pt,0pt>{S\bullet}, 
  "L13";"L13" +(0,-.5)="L12"**\dir{-},
  "L14"-(0,.25)*!<-3pt,0pt>{\bullet S^2  }, 
  "L14";"L14" +(0,-.5)="L13"**\dir{-}, 
  "L15"-(0,.25)*!<-3pt,0pt>{\bullet S }, 
  "L15";"L15" +(0,-.5)="L14"**\dir{-},
   "L16"="x1", 
  "x1"+(0,-.25)="x2",
  "x1";"x2"**\dir{}?(.5)="c", 
  "x1";"x2"**\dir{-},   
"c"-(0,0)*!<-3pt,0pt>{\bullet \delta },  
"x2"+(0,0)*!<0pt,2pt>{\circ},  
  "L17"-(0,.25)*!<-3pt,0pt>{\bullet S }, 
  "L17";"L17" +(0,-.5)="L15"**\dir{-}, 
  "L18"-(0,.25)*!<-3pt,0pt>{\bullet S }, 
  "L18";"L18" +(0,-.5)="L16"**\dir{-},
  "L19";"L19" +(0,-.5)="L17"**\dir{-},
   "L19" + (1, .5)="x1", 
  "x1"+(0,-1)="x2",
  "x1";"x2"**\dir{}?(.5)="c", 
  "x1";"x2"**\dir{-},   
  "c"-(0,0)*!<-3pt,0pt>{\bullet \sigma },  
  "x1"+(0,0)*!<0pt,-2pt>{\circ}, 
"x2"="L18",
  "L110"="x1", 
  "x1"+(0,-.25)="x2",
  "x1";"x2"**\dir{}?(.5)="c", 
  "x1";"x2"**\dir{-},   
"c"-(0,0)*!<-6pt,0pt>{\bullet \delta S },  
"x2"+(0,0)*!<0pt,2pt>{\circ},  
  "L11";"L11" +(0,-.5)="L11"**\dir{-},
  "L12"="x1",       
  "L13"="x2",  , 
  "x1"+(0,-.5)="x3",
  "x2"+(0,-.5)="x4",
  "x1";"x2"**\dir{}?(.5)="m", 
  "m" + (0,-.25)="c", 
  "x1";"x4"**\crv{ "x1"+(0,-.25) & "c" & "x4"+(0,.25)},  
  "x2";"c"+<3pt,2pt>="c2"**\crv{"x2"+(0,-.25) & "c2"},
  "x3";"c"+<-3pt,-2pt>="c1"**\crv{"x3"+(0,.25) & "c1"}, 
  "c"+(0,.10)*!<0pt,-4pt>{\si},
"x3"="L12",
"x4"="L13",
  "L14";"L14" +(0,-.5)="L14"**\dir{-},
  "L15";"L15" +(0,-.5)="L15"**\dir{-},
  "L16";"L16" +(0,-.5)="L16"**\dir{-},
  "L17";"L17" +(0,-.5)="L17"**\dir{-},
  "L18";"L18" +(0,-.5)="L18"**\dir{-},
  "L11"="x1",       
  "L12"="x2",  , 
  "x1"+(0,-.5)="x3",
  "x2"+(0,-.5)="x4",
  "x1";"x2"**\dir{}?(.5)="m", 
  "m" + (0,-.25)="c", 
  "x1";"x4"**\crv{ "x1"+(0,-.25) & "c" & "x4"+(0,.25)},  
  "x2";"c"+<3pt,2pt>="c2"**\crv{"x2"+(0,-.25) & "c2"},
  "x3";"c"+<-3pt,-2pt>="c1"**\crv{"x3"+(0,.25) & "c1"}, 
  "c"+(0,.10)*!<0pt,-6pt>{\si},
"x3"="L11",
"x4"="L12",
  "L13";"L13" +(0,-.5)="L13"**\dir{-},
  "L14";"L14" +(0,-.5)="L14"**\dir{-},
  "L15";"L15" +(0,-.5)="L15"**\dir{-},
  "L16";"L16" +(0,-.5)="L16"**\dir{-},
  "L17";"L17" +(0,-.5)="L17"**\dir{-},
  "L18";"L18" +(0,-.5)="L18"**\dir{-},
  "L11";"L11" +(0,-.5)="L11"**\dir{-},
  "L12";"L12" +(0,-.5)="L12"**\dir{-},
  "L13"="x1",       
  "L14"="x2",  , 
  "x1"+(0,-.5)="x3",
  "x2"+(0,-.5)="x4",
  "x1";"x2"**\dir{}?(.5)="m", 
  "m" + (0,-.25)="c", 
  "x1";"x4"**\crv{ "x1"+(0,-.25) & "c" & "x4"+(0,.25)},  
  "x2";"c"+<3pt,2pt>="c2"**\crv{"x2"+(0,-.25) & "c2"},
  "x3";"c"+<-3pt,-2pt>="c1"**\crv{"x3"+(0,.25) & "c1"}, 
  "c"+(0,.10)*!<0pt,-6pt>{\si},
"x3"="L13",
"x4"="L14",
  "L15";"L15" +(0,-.5)="L15"**\dir{-},
  "L16";"L16" +(0,-.5)="L16"**\dir{-},
  "L17";"L17" +(0,-.5)="L17"**\dir{-},
  "L18";"L18" +(0,-.5)="L18"**\dir{-},
  "L11";"L11" +(0,-.5)="L11"**\dir{-},
  "L12";"L12" +(0,-.5)="L12"**\dir{-},
  "L13";"L13" +(0,-.5)="L13"**\dir{-},
  "L14"="x1",       
  "L15"="x2",  , 
  "x1"+(0,-.5)="x3",
  "x2"+(0,-.5)="x4",
  "x1";"x2"**\dir{}?(.5)="m", 
  "m" + (0,-.25)="c", 
  "x1";"x4"**\crv{ "x1"+(0,-.25) & "c" & "x4"+(0,.25)},  
  "x2";"c"+<3pt,2pt>="c2"**\crv{"x2"+(0,-.25) & "c2"},
  "x3";"c"+<-3pt,-2pt>="c1"**\crv{"x3"+(0,.25) & "c1"}, 
  "c"+(0,.10)*!<0pt,-6pt>{\si},
"x3"="L14",
"x4"="L15",
  "L16";"L16" +(0,-.5)="L16"**\dir{-},
  "L17";"L17" +(0,-.5)="L17"**\dir{-},
  "L18";"L18" +(0,-.5)="L18"**\dir{-},
  "L11";"L11" +(0,-.5)="L11"**\dir{-},
  "L12";"L12" +(0,-.5)="L12"**\dir{-},
  "L13";"L13" +(0,-.5)="L13"**\dir{-},
  "L14";"L14" +(0,-.5)="L14"**\dir{-},
  "L15"="x1",       
  "L16"="x2",  , 
  "x1"+(0,-.5)="x3",
  "x2"+(0,-.5)="x4",
  "x1";"x2"**\dir{}?(.5)="m", 
  "m" + (0,-.25)="c", 
  "x1";"x4"**\crv{ "x1"+(0,-.25) & "c" & "x4"+(0,.25)},  
  "x2";"c"+<3pt,2pt>="c2"**\crv{"x2"+(0,-.25) & "c2"},
  "x3";"c"+<-3pt,-2pt>="c1"**\crv{"x3"+(0,.25) & "c1"}, 
  "c"+(0,.10)*!<0pt,-6pt>{\si},
"x3"="L15",
"x4"="L16",
  "L17";"L17" +(0,-.5)="L17"**\dir{-},
  "L18";"L18" +(0,-.5)="L18"**\dir{-},
  "L11";"L11" +(0,-.5)="L11"**\dir{-},
  "L12";"L12" +(0,-.5)="L12"**\dir{-},
  "L13";"L13" +(0,-.5)="L13"**\dir{-},
  "L14";"L14" +(0,-.5)="L14"**\dir{-},
  "L15";"L15" +(0,-.5)="L15"**\dir{-},
  "L16"="x1",       
  "L17"="x2",  , 
  "x1"+(0,-.5)="x3",
  "x2"+(0,-.5)="x4",
  "x1";"x2"**\dir{}?(.5)="m", 
  "m" + (0,-.25)="c", 
  "x1";"x4"**\crv{ "x1"+(0,-.25) & "c" & "x4"+(0,.25)},  
  "x2";"c"+<3pt,2pt>="c2"**\crv{"x2"+(0,-.25) & "c2"},
  "x3";"c"+<-3pt,-2pt>="c1"**\crv{"x3"+(0,.25) & "c1"}, 
  "c"+(0,.10)*!<0pt,-6pt>{\si},
"x3"="L16",
"x4"="L17",
  "L18";"L18" +(0,-.5)="L18"**\dir{-},
  "L11";"L11" +(0,-.5)="L11"**\dir{-},
  "L12";"L12" +(0,-.5)="L12"**\dir{-},
  "L13"="x1",       
  "L14"="x2",  , 
  "x1"+(0,-.5)="x3",
  "x2"+(0,-.5)="x4",
  "x1";"x2"**\dir{}?(.5)="m", 
  "m" + (0,-.25)="c", 
  "x1";"x4"**\crv{ "x1"+(0,-.25) & "c" & "x4"+(0,.25)},  
  "x2";"c"+<3pt,2pt>="c2"**\crv{"x2"+(0,-.25) & "c2"},
  "x3";"c"+<-3pt,-2pt>="c1"**\crv{"x3"+(0,.25) & "c1"}, 
  "c"+(0,.10)*!<0pt,-6pt>{\si},
"x3"="L13",
"x4"="L14",
  "L15";"L15" +(0,-.5)="L15"**\dir{-},
  "L16";"L16" +(0,-.5)="L16"**\dir{-},
  "L17";"L17" +(0,-.5)="L17"**\dir{-},
  "L18";"L18" +(0,-.5)="L18"**\dir{-},
  "L11";"L11" +(0,-.5)="L11"**\dir{-},
  "L12";"L12" +(0,-.5)="L12"**\dir{-},
  "L13";"L13" +(0,-.5)="L13"**\dir{-},
  "L14"="x1",       
  "L15"="x2",  , 
  "x1"+(0,-.5)="x3",
  "x2"+(0,-.5)="x4",
  "x1";"x2"**\dir{}?(.5)="m", 
  "m" + (0,-.25)="c", 
  "x1";"x4"**\crv{ "x1"+(0,-.25) & "c" & "x4"+(0,.25)},  
  "x2";"c"+<3pt,2pt>="c2"**\crv{"x2"+(0,-.25) & "c2"},
  "x3";"c"+<-3pt,-2pt>="c1"**\crv{"x3"+(0,.25) & "c1"}, 
  "c"+(0,.10)*!<0pt,-6pt>{\si},
"x3"="L14",
"x4"="L15",
  "L16";"L16" +(0,-.5)="L16"**\dir{-},
  "L17";"L17" +(0,-.5)="L17"**\dir{-},
  "L18";"L18" +(0,-.5)="L18"**\dir{-},
 "L11";"L11" +(0,-.5)="L11"**\dir{-},
  "L12"="x1",       
  "L13"="x2",  , 
  "x1"+(0,-.5)="x3",
  "x2"+(0,-.5)="x4",
  "x1";"x2"**\dir{}?(.5)="m", 
  "m" + (0,-.25)="c", 
  "x1";"x4"**\crv{ "x1"+(0,-.25) & "c" & "x4"+(0,.25)},  
  "x2";"c"+<3pt,2pt>="c2"**\crv{"x2"+(0,-.25) & "c2"},
  "x3";"c"+<-3pt,-2pt>="c1"**\crv{"x3"+(0,.25) & "c1"}, 
  "c"+(0,.10)*!<0pt,-6pt>{\si},
"x3"="L12",
"x4"="L13",
  "L14";"L14" +(0,-.5)="L14"**\dir{-},
  "L15";"L15" +(0,-.5)="L15"**\dir{-},
  "L16";"L16" +(0,-.5)="L16"**\dir{-},
  "L17";"L17" +(0,-.5)="L17"**\dir{-},
  "L18";"L18" +(0,-.5)="L18"**\dir{-},
  "L11";"L11" +(0,-.5)="L11"**\dir{-},
  "L12";"L12" +(0,-.5)="L12"**\dir{-},
  "L13"="x1",       
  "L14"="x2",  , 
  "x1"+(0,-.5)="x3",
  "x2"+(0,-.5)="x4",
  "x1";"x2"**\dir{}?(.5)="m", 
  "m" + (0,-.25)="c", 
  "x1";"x4"**\crv{ "x1"+(0,-.25) & "c" & "x4"+(0,.25)},  
  "x2";"c"+<3pt,2pt>="c2"**\crv{"x2"+(0,-.25) & "c2"},
  "x3";"c"+<-3pt,-2pt>="c1"**\crv{"x3"+(0,.25) & "c1"}, 
  "c"+(0,.10)*!<0pt,-6pt>{\si},
"x3"="L13",
"x4"="L14",
  "L15";"L15" +(0,-.5)="L15"**\dir{-},
  "L16";"L16" +(0,-.5)="L16"**\dir{-},
  "L17";"L17" +(0,-.5)="L17"**\dir{-},
  "L18";"L18" +(0,-.5)="L18"**\dir{-},
  "L11"="x1",       
  "L12"="x2",  , 
  "x1"+(0,-.5)="x3",
  "x2"+(0,-.5)="x4",
  "x1";"x2"**\dir{}?(.5)="m", 
  "m" + (0,-.25)="c", 
  "x1";"x4"**\crv{ "x1"+(0,-.25) & "c" & "x4"+(0,.25)},  
  "x2";"c"+<3pt,2pt>="c2"**\crv{"x2"+(0,-.25) & "c2"},
  "x3";"c"+<-3pt,-2pt>="c1"**\crv{"x3"+(0,.25) & "c1"}, 
  "c"+(0,.10)*!<0pt,-6pt>{\si},
"x3"="L11",
"x4"="L12",
  "L13";"L13" +(0,-.5)="L13"**\dir{-},
  "L14";"L14" +(0,-.5)="L14"**\dir{-},
  "L15";"L15" +(0,-.5)="L15"**\dir{-},
  "L16";"L16" +(0,-.5)="L16"**\dir{-},
  "L17";"L17" +(0,-.5)="L17"**\dir{-},
  "L18";"L18" +(0,-.5)="L18"**\dir{-},
  "L11";"L11" +(0,-.5)="L11"**\dir{-},
  "L12";"L12" +(0,-.5)="L12"**\dir{-},
  "L13";"L13" +(0,-.5)="L13"**\dir{-},
  "L14"="x1", 
  "L15"="x2",   
  "x1";"x2"**\dir{}?(.5)="m", 
  "m" + (0,-.25)="c", 
  "m" + (0,-.5)="x3", 
  "x1";"x3"**\crv{"x1"-(0,.25) & "c" & "x3"},
  "x2";"x3"**\crv{"x2"-(0,.25) & "c"  & "x3"}, 
"x3"="L14",
  "L16";"L16" +(0,-.5)="L15"**\dir{-},
  "L17";"L17" +(0,-.5)="L16"**\dir{-},
  "L18";"L18" +(0,-.5)="L17"**\dir{-},
  "L11";"L11" +(0,-.5)="L11"**\dir{-},
  "L12";"L12" +(0,-.5)="L12"**\dir{-},
  "L13";"L13" +(0,-.5)="L13"**\dir{-},
  "L14"="x1", 
  "L15"="x2",  
  "x1";"x2"**\dir{}?(.5)="m", 
  "m" + (0,-.25)="c", 
  "m" + (0,-.5)="x3", 
  "x1";"x3"**\crv{"x1"-(0,.25) & "c" & "x3"},
  "x2";"x3"**\crv{"x2"-(0,.25) & "c"  & "x3"}, 
"x3"="L14",
  "L16";"L16" +(0,-.5)="L15"**\dir{-},
  "L17";"L17" +(0,-.5)="L16"**\dir{-},
  "L11";"L11" +(0,-.5)="L11"**\dir{-},
  "L12";"L12" +(0,-.5)="L12"**\dir{-},
  "L13"="x1", 
  "L14"="x2",  
  "x1";"x2"**\dir{}?(.5)="m", 
  "m" + (0,-.25)="c", 
  "m" + (0,-.5)="x3", 
  "x1";"x3"**\crv{"x1"-(0,.25) & "c" & "x3"},
  "x2";"x3"**\crv{"x2"-(0,.25) & "c"  & "x3"}, 
"x3"="L13",
  "L15";"L15" +(0,-.5)="L14"**\dir{-},
  "L16";"L16" +(0,-.5)="L15"**\dir{-},
  "L11"="x1", 
  "L12"="x2",  
  "x1";"x2"**\dir{}?(.5)="m", 
  "m" + (0,-.25)="c", 
  "m" + (0,-.5)="x3", 
  "x1";"x3"**\crv{"x1"-(0,.25) & "c" & "x3"},
  "x2";"x3"**\crv{"x2"-(0,.25) & "c"  & "x3"}, 
"x3"="L11",
  "L13";"L13" +(0,-.5)="L12"**\dir{-},
  "L14"="x1", 
  "L15"="x2",  
  "x1";"x2"**\dir{}?(.5)="m", 
  "m" + (0,-.25)="c", 
  "m" + (0,-.5)="x3", 
  "x1";"x3"**\crv{"x1"-(0,.25) & "c" & "x3"},
  "x2";"x3"**\crv{"x2"-(0,.25) & "c"  & "x3"}, 
"x3"="L13",
  "L11"="x1",       
  "L12"="x2",  , 
  "x1"+(0,-.5)="x3",
  "x2"+(0,-.5)="x4",
  "x1";"x2"**\dir{}?(.5)="m", 
  "m" + (0,-.25)="c", 
  "x1";"x4"**\crv{ "x1"+(0,-.25) & "c" & "x4"+(0,.25)},  
  "x2";"c"+<3pt,2pt>="c2"**\crv{"x2"+(0,-.25) & "c2"},
  "x3";"c"+<-3pt,-2pt>="c1"**\crv{"x3"+(0,.25) & "c1"}, 
  "c"+(0,.10)*!<0pt,-6pt>{\si},
"x3"="L11",
"x4"="L12",
  "L13";"L13" +(0,-.5)="L13"**\dir{-},
  "L11";"L11" +(0,-.5)="L11"**\dir{-},
  "L12"="x1", 
  "L13"="x2",   
  "x1";"x2"**\dir{}?(.5)="m", 
  "m" + (0,-.25)="c", 
  "m" + (0,-.5)="x3", 
  "x1";"x3"**\crv{"x1"-(0,.25) & "c" & "x3"},
  "x2";"x3"**\crv{"x2"-(0,.25) & "c"  & "x3"}, 
"x3"="L12",
  "L11"-(0,.5)*!<0pt,-6pt>{H}, 
  "L12"-(0,.5)*!<0pt,-6pt>{H}, 
\endxy
\]

\begin{align*} 
=&(1,m)(\si,1)(m,1,m)(2,m,2)(3,m,2)(3,m,3)(\si,6)(2,\si,4)(1,\si,5)(3,\si,3)\\
~~~~~&(2,\si,4)(5,\si,1)(4,\si,2)(3,\si,3)(2,\si,4)(\si,6)(1,\si,5)(\si,6)\\
~~~~~&(\delta,S,S^2,S^2,S,\widetilde{S},S,1,\sigma,\delta S)(3,\si,4)(2,\si,5)(7,\si)(6,\si,1)(5,\si,2)\\
~~~~~&(\Delta,7)(5,\Delta,1)(2,\Delta,3)(2,\Delta,2)(\Delta,3)(\Delta,2)\\
=&(1,m)(\si,1)(m,1,m)(2,m,2)(3,m,2)(\si,5)(3,m,3)(2,\si,)(3,\si,3)(1,\si,5)\\
~~~~~&(2,\si,4)(5,\si,1)(4,\si,2)(3,\si,3)(2,\si,4)(1,\si,5)(\si,6)(1,\si,5)\\
~~~~~&(\delta,S,S^2,S^2,S,\widetilde{S},S,1,\sigma,\delta S)(3,\si,4)(2,\si,5)(7,\si)(6,\si,1)(5,\si,2)\\
~~~~~&(\Delta,7)(5,\Delta,1)(2,\Delta,3) (2,\Delta,2)(\Delta,3)(\Delta,2)\\
=&(1,m)(\si,1)(m,1,m)(2,m,2)(3,m,2)(\si,5)(3,m,3)(2,\si_{21},3)(1,\si_{21},4)(5,\si,1)\\
~~~~~&(4,\si,2)(3,\si,3)(1,\si_{12},4)(\si,6)(\delta,S,S^2,S^2,S,\widetilde{S},S,1,\sigma,\delta S)(2,\si,5)(3,\si,4)\\
~~~~~&(2,\si,5)(7,\si)(6,\si,1)(5,\si,2)(\Delta,7)(5,\Delta,1)(2,\Delta,3)(2,\Delta,2)(\Delta,3)(\Delta,2)\\
=&(1,m)(\si,1)(m,1,m)(2,m,2)(3,m,2)(\si,5)(2,\si,3)(1,\si,4)(1,m,5)(5,\si,1)\\
~~~~~&(4,\si,2)(3,\si,3)(1,\si_{12},4)(\si,6)(\delta,S,S^2,S^2,S,\widetilde{S},S,1,\sigma,\delta S)(3,\si,4)(2,\si,5)\\
~~~~~&(3,\si,4)(7,\si)(6,\si,1)(5,\si,2)(\Delta,7)(5,\Delta,1)(2,\Delta,3)(2,\Delta,2)(\Delta,3)(\Delta,2)\\
=&(1,m)(\si,1)(m,1,m)(2,m,2)(3,m,2)(\si,5)(2,\si,3)(1,\si,4)(4,\si,1)(3,\si,2)\\
~~~~~&(2,\si,3)(1,m,5)(1,\si_{12},4)(\si,6)(\delta,S,S^2,S^2,S,\widetilde{S},S,1,\sigma,\delta S)(3,\si,4)(2,\si,5)\\
~~~~~&(3,\si,4)(7,\si)(6,\si,1)(5,\si,2)(\Delta,7)(5,\Delta,1)(2,\Delta,3)(2,\Delta,2)(\Delta,3)(\Delta,2)\\
=&(1,m)(\si,1)(m,1,m)(2,m,2)(3,m,2)(\si,5)(2,\si,3)(1,\si,4)(4,\si,1)(3,\si,2)\\
~~~~~&(2,\si,3)(1,\si,4)(2,m,4)(\si,6)(\delta,S,S^2,S^2,S,\widetilde{S},S,1,\sigma,\delta S)(3,\si,4)\\
~~~~~&(2,\si_{21},4)(7,\si)(6,\si,1)(5,\si,2)(\Delta,7)(5,\Delta,1)(2,\Delta,3)(2,\Delta,2)(1,\Delta,2)(\Delta,2)\\
=&(1,m)(\si,1)(m,1,m)(2,m,2)(3,m,2)(\si,5)(2,\si,3)(1,\si,4)(4,\si,1)(3,\si,2)\\
~~~~~&(2,\si,3)(1,\si,4)(\si,5)(2,m,4)(\delta,S,S^2,S^2,S,\widetilde{S},S,1,\sigma,\delta S)(3,\si,4)\\
~~~~~&(2,\si_{21},4)(7,\si)(6,\si,1)(5,\si,2)(\Delta,7)(5,\Delta,1)(1,\Delta,4)(1,\Delta,3)(1,\Delta,2)(\Delta,2)\\
=&(1,m)(\si,1)(m,1,m)(2,m,2)(3,m,2)(\si,5)(2,\si,3)(1,\si,4)(4,\si,1)(3,\si,2)\\
~~~~~&(2,\si,3)(1,\si,4)(\si,5)(2,m,4)(\delta,S,S^2,S^2,S,\widetilde{S},S,1,\sigma,\delta S)(3,\si,4)\\
~~~~~&(2,\si_{21},4)(2,\Delta,5)(6,\si)(5,\si,1)(4,\si,2)(\Delta,6)(4,\Delta,1)(1,\Delta,3)(1,\Delta,2)(\Delta,2)\\
=&(1,m)(\si,1)(m,1,m)(2,m,2)(3,m,2)(\si,5)(2,\si,3)(1,\si,4)(4,\si,1)(3,\si,2)\\
~~~~~&(2,\si,3)(1,\si,4)(\si,5)(2,m,4)(\delta,S,S^2,S^2,S,\widetilde{S},S,1,\sigma,\delta S)(3,\si,4)\\
~~~~~&(3,\Delta,4)(2,\si,4)(6,\si)(5,\si,1)(4,\si,2)(\Delta,6)(4,\Delta,1)(1,\Delta,3)(1,\Delta,2)(\Delta,2)\\
=&(1,m)(\si,1)(m,1,m)(2,m,2)(3,m,2)(\si,5)(2,\si,3)(1,\si,4)(4,\si,1)(3,\si,2)\\
~~~~~&(2,\si,3)(1,\si,4)(\si,5)(\delta,S,S^2,m(S^2,S)\si\Delta,\widetilde{S},S,1,\sigma,\delta S)(2,\si,4)\\
~~~~~&(6,\si)(5,\si,1)(4,\si,2)(\Delta,6)(4,\Delta,1)(1,\Delta,3)(1,\Delta,2)(\Delta,2)=
\end{align*}

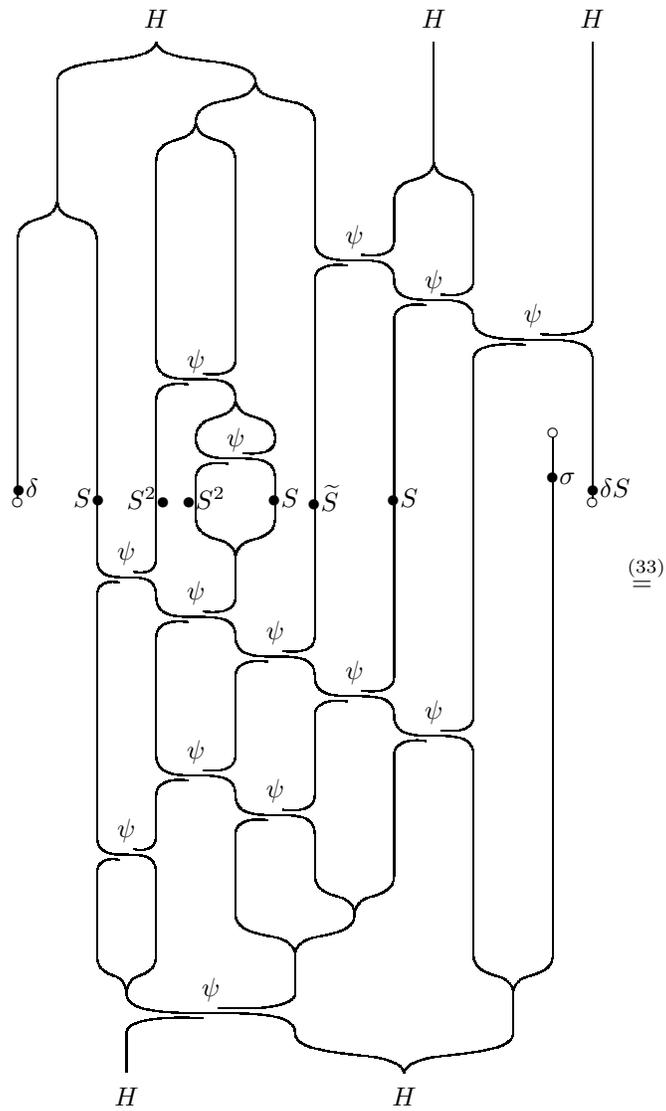
\[
\xy /r2.5pc/:,
(1,0)= "L11",
(4.5,0)= "L12",
(6.5,0)= "L13",
  "L11"+(0,.5)*!<0pt,6pt>{H}, 
  "L12"+(0,.5)*!<0pt,6pt>{H}, 
  "L13"+(0,.5)*!<0pt,6pt>{H}, 
   "L11"="x1", 
  "x1"+(-1.25,-.5)="x2",
  "x1"+(1.25,-.5)="x3",
  "x1"+(0,-.25)="c",  
   "x1";"x2"**\crv{"x1" & "c" & "x2"+(0,.25)}, 
   "x1";"x3"**\crv{"x1" & "c" & "x3"+(0,.25)}, 
"x2"= "L21",
"x3"= "L22",
  "L12";"L12" +(0,-.5)="L23"**\dir{-},
  "L13";"L13" +(0,-.5)="L24"**\dir{-},
"L21"= "L11",
"L22"= "L12",
"L23"= "L13",
"L24"= "L14",
  "L11";"L11" +(0,-.5)="L11"**\dir{-},
   "L12"="x1", 
  "x1"+(-.75,-.5)="x2",
  "x1"+(.75,-.5)="x3",
  "x1"+(0,-.25)="c",  
   "x1";"x2"**\crv{"x1" & "c" & "x2"+(0,.25)}, 
   "x1";"x3"**\crv{"x1" & "c" & "x3"+(0,.25)}, 
"x2"= "L22",
"x3"= "L23",
  "L13";"L13" +(0,-.5)="L24"**\dir{-},
  "L14";"L14" +(0,-.5)="L25"**\dir{-},
"L22"= "L12",
"L23"= "L13",
"L24"= "L14",
"L25"= "L15",
  "L11";"L11" +(0,-.5)="L11"**\dir{-},
   "L12"="x1", 
  "x1"+(-.5,-.5)="x2",
  "x1"+(.5,-.5)="x3",
  "x1"+(0,-.25)="c",  
   "x1";"x2"**\crv{"x1" & "c" & "x2"+(0,.25)}, 
   "x1";"x3"**\crv{"x1" & "c" & "x3"+(0,.25)}, 
"x2"= "L22",
"x3"= "L23",
  "L13";"L13" +(0,-.5)="L24"**\dir{-},
  "L14";"L14" +(0,-.5)="L25"**\dir{-},
  "L15";"L15" +(0,-.5)="L26"**\dir{-},
"L22"= "L12",
"L23"= "L13",
"L24"= "L14",
"L25"= "L15",
"L26"= "L16",
  "L11";"L11" +(0,-.5)="L11"**\dir{-},
  "L12";"L12" +(0,-.5)="L12"**\dir{-},
  "L13";"L13" +(0,-.5)="L13"**\dir{-},
  "L14";"L14" +(0,-.5)="L14"**\dir{-},
  "L15"="x1", 
  "x1"+(-.5,-.5)="x2",
  "x1"+(.5,-.5)="x3",
  "x1"+(0,-.25)="c",  
   "x1";"x2"**\crv{"x1" & "c" & "x2"+(0,.25)}, 
   "x1";"x3"**\crv{"x1" & "c" & "x3"+(0,.25)}, 
"x2"= "L25",
"x3"= "L26",
   "L16";"L16" +(0,-.5)="L27"**\dir{-},
"L25"= "L15",
"L26"= "L16",
"L27"= "L17",
   "L11"="x1", 
  "x1"+(-.5,-.5)="x2",
  "x1"+(.5,-.5)="x3",
  "x1"+(0,-.25)="c",  
   "x1";"x2"**\crv{"x1" & "c" & "x2"+(0,.25)}, 
   "x1";"x3"**\crv{"x1" & "c" & "x3"+(0,.25)}, 
"x2"= "L21",
"x3"= "L22",
  "L12";"L12" +(0,-.5)="L23"**\dir{-},
  "L13";"L13" +(0,-.5)="L24"**\dir{-},
  "L14";"L14" +(0,-.5)="L25"**\dir{-},
  "L15";"L15" +(0,-.5)="L26"**\dir{-},
  "L16";"L16" +(0,-.5)="L27"**\dir{-},
  "L17";"L17" +(0,-.5)="L28"**\dir{-},
"L21"= "L11",
"L22"= "L12",
"L23"= "L13",
"L24"= "L14",
"L25"= "L15",
"L26"= "L16",
"L27"= "L17",
"L28"= "L18",
  "L11";"L11" +(0,-.5)="L11"**\dir{-},
  "L12";"L12" +(0,-.5)="L12"**\dir{-},
  "L13";"L13" +(0,-.5)="L13"**\dir{-},
  "L14";"L14" +(0,-.5)="L14"**\dir{-},
  "L15"="x1",       
  "L16"="x2",  , 
  "x1"+(0,-.5)="x3",
  "x2"+(0,-.5)="x4",
  "x1";"x2"**\dir{}?(.5)="m", 
  "m" + (0,-.25)="c", 
  "x1";"x4"**\crv{ "x1"+(0,-.25) & "c" & "x4"+(0,.25)},  
  "x2";"c"+<3pt,2pt>="c2"**\crv{"x2"+(0,-.25) & "c2"},
  "x3";"c"+<-3pt,-2pt>="c1"**\crv{"x3"+(0,.25) & "c1"}, 
  "c"+(0,.10)*!<0pt,-6pt>{\si},
"x3"="L15",
"x4"="L16",
  "L17";"L17" +(0,-.5)="L17"**\dir{-},
  "L18";"L18" +(0,-.5)="L18"**\dir{-},
  "L11";"L11" +(0,-.5)="L11"**\dir{-},
  "L12";"L12" +(0,-.5)="L12"**\dir{-},
  "L13";"L13" +(0,-.5)="L13"**\dir{-},
  "L14";"L14" +(0,-.5)="L14"**\dir{-},
  "L15";"L15" +(0,-.5)="L15"**\dir{-},
  "L16"="x1",       
  "L17"="x2",  , 
  "x1"+(0,-.5)="x3",
  "x2"+(0,-.5)="x4",
  "x1";"x2"**\dir{}?(.5)="m", 
  "m" + (0,-.25)="c", 
  "x1";"x4"**\crv{ "x1"+(0,-.25) & "c" & "x4"+(0,.25)},  
  "x2";"c"+<3pt,2pt>="c2"**\crv{"x2"+(0,-.25) & "c2"},
  "x3";"c"+<-3pt,-2pt>="c1"**\crv{"x3"+(0,.25) & "c1"}, 
  "c"+(0,.10)*!<0pt,-4pt>{\si},
"x3"="L16",
"x4"="L17",
  "L18";"L18" +(0,-.5)="L18"**\dir{-},
  "L11";"L11" +(0,-.5)="L11"**\dir{-},
  "L12";"L12" +(0,-.5)="L12"**\dir{-},
  "L13";"L13" +(0,-.5)="L13"**\dir{-},
  "L14";"L14" +(0,-.5)="L14"**\dir{-},
  "L15";"L15" +(0,-.5)="L15"**\dir{-},
  "L16";"L16" +(0,-.5)="L16"**\dir{-},
  "L17"="x1",       
  "L18"="x2",  , 
  "x1"+(0,-.5)="x3",
  "x2"+(0,-.5)="x4",
  "x1";"x2"**\dir{}?(.5)="m", 
  "m" + (0,-.25)="c", 
  "x1";"x4"**\crv{ "x1"+(0,-.25) & "c" & "x4"+(0,.25)},  
  "x2";"c"+<3pt,2pt>="c2"**\crv{"x2"+(0,-.25) & "c2"},
  "x3";"c"+<-3pt,-2pt>="c1"**\crv{"x3"+(0,.25) & "c1"}, 
  "c"+(0,.10)*!<0pt,-6pt>{\si},
"x3"="L17",
"x4"="L18",
  "L11";"L11" +(0,-.5)="L11"**\dir{-},
  "L12";"L12" +(0,-.5)="L12"**\dir{-},
  "L13"="x1",       
  "L14"="x2",  , 
  "x1"+(0,-.5)="x3",
  "x2"+(0,-.5)="x4",
  "x1";"x2"**\dir{}?(.5)="m", 
  "m" + (0,-.25)="c", 
  "x1";"x4"**\crv{ "x1"+(0,-.25) & "c" & "x4"+(0,.25)},  
  "x2";"c"+<3pt,2pt>="c2"**\crv{"x2"+(0,-.25) & "c2"},
  "x3";"c"+<-3pt,-2pt>="c1"**\crv{"x3"+(0,.25) & "c1"}, 
  "c"+(0,.10)*!<0pt,-4pt>{\si},
"x3"="L13",
"x4"="L14",
  "L15";"L15" +(0,-.5)="L15"**\dir{-},
  "L16";"L16" +(0,-.5)="L16"**\dir{-},
  "L17";"L17" +(0,-.5)="L17"**\dir{-},
  "L18";"L18" +(0,-.5)="L18"**\dir{-},
  "L11";"L11" +(0,-.5)="L11"**\dir{-},
  "L12";"L12" +(0,-.5)="L12"**\dir{-},
  "L13";"L13" +(0,-.5)="L13"**\dir{-},
   "L14"="x1", 
  "x1"+(-.5,-.5)="x2",
  "x1"+(.5,-.5)="x3",
  "x1"+(0,-.25)="c",  
   "x1";"x2"**\crv{"x1" & "c" & "x2"+(0,.25)}, 
   "x1";"x3"**\crv{"x1" & "c" & "x3"+(0,.25)}, 
"x2"= "L24",
"x3"= "L25",
   "L15";"L15" +(0,-.5)="L26"**\dir{-},
   "L16";"L16" +(0,-.5)="L27"**\dir{-},
   "L17";"L17" +(0,-.5)="L28"**\dir{-},
   "L18";"L18" +(0,-.5)="L29"**\dir{-},
"L24"= "L14",
"L25"= "L15",
"L26"= "L16",
"L27"= "L17",
"L28"= "L18",
"L29"= "L19",
  "L11";"L11" +(0,-.5)="L11"**\dir{-},
  "L12";"L12" +(0,-.5)="L12"**\dir{-},
  "L13";"L13" +(0,-.5)="L13"**\dir{-},
  "L14"="x1",       
  "L15"="x2",  , 
  "x1"+(0,-.5)="x3",
  "x2"+(0,-.5)="x4",
  "x1";"x2"**\dir{}?(.5)="m", 
  "m" + (0,-.25)="c", 
  "x1";"x4"**\crv{ "x1"+(0,-.25) & "c" & "x4"+(0,.25)},  
  "x2";"c"+<3pt,2pt>="c2"**\crv{"x2"+(0,-.25) & "c2"},
  "x3";"c"+<-3pt,-2pt>="c1"**\crv{"x3"+(0,.25) & "c1"}, 
  "c"+(0,.10)*!<0pt,-4pt>{\si},
"x3"="L14",
"x4"="L15",
  "L16";"L16" +(0,-.5)="L16"**\dir{-},
  "L17";"L17" +(0,-.5)="L17"**\dir{-},
  "L18";"L18" +(0,-.5)="L18"**\dir{-},
  "L19";"L19" +(0,-.5)="L19"**\dir{-},
  "L11"="x1", 
  "x1"+(0,-.25)="x2",
  "x1";"x2"**\dir{}?(.5)="c", 
  "x1";"x2"**\dir{-},   
"c"-(0,0)*!<-3pt,0pt>{\bullet \delta  },  
"x2"+(0,0)*!<0pt,2pt>{\circ},  
  "L12"-(0,.25)*!<3pt,0pt>{S \bullet }, 
  "L12";"L12" +(0,-.5)="L11"**\dir{-}, 
  "L13"-(0,.25)*!<3pt,0pt>{S^2 \bullet}, 
  "L13";"L13" +(0,-.5)="L12"**\dir{-},
  "L14"-(0,.25)*!<-3pt,0pt>{\bullet S^2  }, 
  "L14";"L14" +(0,-.5)="L13"**\dir{-}, 
  "L15"-(0,.25)*!<-3pt,0pt>{\bullet S }, 
  "L15";"L15" +(0,-.5)="L14"**\dir{-},
  "L16"-(0,.25)*!<-3pt,0pt>{\bullet \widetilde{S}}, 
  "L16";"L16" +(0,-.5)="L15"**\dir{-},
  "L17"-(0,.25)*!<-3pt,0pt>{\bullet S }, 
  "L17";"L17" +(0,-.5)="L16"**\dir{-},
  "L18";"L18" +(0,-.5)="L17"**\dir{-},
   "L18" + (1, .5)="x1", 
  "x1"+(0,-1)="x2",
  "x1";"x2"**\dir{}?(.5)="c", 
  "x1";"x2"**\dir{-},   
  "c"-(0,0)*!<-3pt,0pt>{\bullet \sigma },  
  "x1"+(0,0)*!<0pt,-2pt>{\circ}, 
"x2"="L18",
  "L19"="x1", 
  "x1"+(0,-.25)="x2",
  "x1";"x2"**\dir{}?(.5)="c", 
  "x1";"x2"**\dir{-},   
"c"-(0,0)*!<-6pt,0pt>{\bullet \delta S},  
"x2"+(0,0)*!<0pt,2pt>{\circ},  
  "L11";"L11" +(0,-.5)="L11"**\dir{-},
  "L12";"L12" +(0,-.5)="L12"**\dir{-},
  "L13"="x1", 
  "L14"="x2",  
  "x1";"x2"**\dir{}?(.5)="m", 
  "m" + (0,-.25)="c", 
  "m" + (0,-.5)="x3", 
  "x1";"x3"**\crv{"x1"-(0,.25) & "c" & "x3"},
  "x2";"x3"**\crv{"x2"-(0,.25) & "c"  & "x3"}, 
"x3"="L13",
  "L15";"L15" +(0,-.5)="L14"**\dir{-},
  "L16";"L16" +(0,-.5)="L15"**\dir{-},
  "L17";"L17" +(0,-.5)="L16"**\dir{-},
  "L18";"L18" +(0,-.5)="L17"**\dir{-},
  "L11"="x1",       
  "L12"="x2",  , 
  "x1"+(0,-.5)="x3",
  "x2"+(0,-.5)="x4",
  "x1";"x2"**\dir{}?(.5)="m", 
  "m" + (0,-.25)="c", 
  "x1";"x4"**\crv{ "x1"+(0,-.25) & "c" & "x4"+(0,.25)},  
  "x2";"c"+<3pt,2pt>="c2"**\crv{"x2"+(0,-.25) & "c2"},
  "x3";"c"+<-3pt,-2pt>="c1"**\crv{"x3"+(0,.25) & "c1"}, 
  "c"+(0,.10)*!<0pt,-6pt>{\si},
"x3"="L11",
"x4"="L12",
  "L13";"L13" +(0,-.5)="L13"**\dir{-},
  "L14";"L14" +(0,-.5)="L14"**\dir{-},
  "L15";"L15" +(0,-.5)="L15"**\dir{-},
  "L16";"L16" +(0,-.5)="L16"**\dir{-},
  "L17";"L17" +(0,-.5)="L17"**\dir{-},
  "L11";"L11" +(0,-.5)="L11"**\dir{-},
  "L12"="x1",       
  "L13"="x2",  , 
  "x1"+(0,-.5)="x3",
  "x2"+(0,-.5)="x4",
  "x1";"x2"**\dir{}?(.5)="m", 
  "m" + (0,-.25)="c", 
  "x1";"x4"**\crv{ "x1"+(0,-.25) & "c" & "x4"+(0,.25)},  
  "x2";"c"+<3pt,2pt>="c2"**\crv{"x2"+(0,-.25) & "c2"},
  "x3";"c"+<-3pt,-2pt>="c1"**\crv{"x3"+(0,.25) & "c1"}, 
  "c"+(0,.10)*!<0pt,-6pt>{\si},
"x3"="L12",
"x4"="L13",
  "L14";"L14" +(0,-.5)="L14"**\dir{-},
  "L15";"L15" +(0,-.5)="L15"**\dir{-},
  "L16";"L16" +(0,-.5)="L16"**\dir{-},
  "L17";"L17" +(0,-.5)="L17"**\dir{-},
  "L11";"L11" +(0,-.5)="L11"**\dir{-},
  "L12";"L12" +(0,-.5)="L12"**\dir{-},
  "L13"="x1",       
  "L14"="x2",  , 
  "x1"+(0,-.5)="x3",
  "x2"+(0,-.5)="x4",
  "x1";"x2"**\dir{}?(.5)="m", 
  "m" + (0,-.25)="c", 
  "x1";"x4"**\crv{ "x1"+(0,-.25) & "c" & "x4"+(0,.25)},  
  "x2";"c"+<3pt,2pt>="c2"**\crv{"x2"+(0,-.25) & "c2"},
  "x3";"c"+<-3pt,-2pt>="c1"**\crv{"x3"+(0,.25) & "c1"}, 
  "c"+(0,.10)*!<0pt,-6pt>{\si},
"x3"="L13",
"x4"="L14",
  "L15";"L15" +(0,-.5)="L15"**\dir{-},
  "L16";"L16" +(0,-.5)="L16"**\dir{-},
  "L17";"L17" +(0,-.5)="L17"**\dir{-},
  "L11";"L11" +(0,-.5)="L11"**\dir{-},
  "L12";"L12" +(0,-.5)="L12"**\dir{-},
  "L13";"L13" +(0,-.5)="L13"**\dir{-},
  "L14"="x1",       
  "L15"="x2",  , 
  "x1"+(0,-.5)="x3",
  "x2"+(0,-.5)="x4",
  "x1";"x2"**\dir{}?(.5)="m", 
  "m" + (0,-.25)="c", 
  "x1";"x4"**\crv{ "x1"+(0,-.25) & "c" & "x4"+(0,.25)},  
  "x2";"c"+<3pt,2pt>="c2"**\crv{"x2"+(0,-.25) & "c2"},
  "x3";"c"+<-3pt,-2pt>="c1"**\crv{"x3"+(0,.25) & "c1"}, 
  "c"+(0,.10)*!<0pt,-6pt>{\si},
"x3"="L14",
"x4"="L15",
  "L16";"L16" +(0,-.5)="L16"**\dir{-},
  "L17";"L17" +(0,-.5)="L17"**\dir{-},
  "L11";"L11" +(0,-.5)="L11"**\dir{-},
  "L12";"L12" +(0,-.5)="L12"**\dir{-},
  "L13";"L13" +(0,-.5)="L13"**\dir{-},
  "L14";"L14" +(0,-.5)="L14"**\dir{-},
  "L15"="x1",       
  "L16"="x2",  , 
  "x1"+(0,-.5)="x3",
  "x2"+(0,-.5)="x4",
  "x1";"x2"**\dir{}?(.5)="m", 
  "m" + (0,-.25)="c", 
  "x1";"x4"**\crv{ "x1"+(0,-.25) & "c" & "x4"+(0,.25)},  
  "x2";"c"+<3pt,2pt>="c2"**\crv{"x2"+(0,-.25) & "c2"},
  "x3";"c"+<-3pt,-2pt>="c1"**\crv{"x3"+(0,.25) & "c1"}, 
  "c"+(0,.10)*!<0pt,-6pt>{\si},
"x3"="L15",
"x4"="L16",
  "L17";"L17" +(0,-.5)="L17"**\dir{-},
  "L11";"L11" +(0,-.5)="L11"**\dir{-},
  "L12"="x1",       
  "L13"="x2",  , 
  "x1"+(0,-.5)="x3",
  "x2"+(0,-.5)="x4",
  "x1";"x2"**\dir{}?(.5)="m", 
  "m" + (0,-.25)="c", 
  "x1";"x4"**\crv{ "x1"+(0,-.25) & "c" & "x4"+(0,.25)},  
  "x2";"c"+<3pt,2pt>="c2"**\crv{"x2"+(0,-.25) & "c2"},
  "x3";"c"+<-3pt,-2pt>="c1"**\crv{"x3"+(0,.25) & "c1"}, 
  "c"+(0,.10)*!<0pt,-6pt>{\si},
"x3"="L12",
"x4"="L13",
  "L14";"L14" +(0,-.5)="L14"**\dir{-},
  "L15";"L15" +(0,-.5)="L15"**\dir{-},
  "L16";"L16" +(0,-.5)="L16"**\dir{-},
  "L17";"L17" +(0,-.5)="L17"**\dir{-},
  "L11";"L11" +(0,-.5)="L11"**\dir{-},
  "L12";"L12" +(0,-.5)="L12"**\dir{-},
  "L13"="x1",       
  "L14"="x2",  , 
  "x1"+(0,-.5)="x3",
  "x2"+(0,-.5)="x4",
  "x1";"x2"**\dir{}?(.5)="m", 
  "m" + (0,-.25)="c", 
  "x1";"x4"**\crv{ "x1"+(0,-.25) & "c" & "x4"+(0,.25)},  
  "x2";"c"+<3pt,2pt>="c2"**\crv{"x2"+(0,-.25) & "c2"},
  "x3";"c"+<-3pt,-2pt>="c1"**\crv{"x3"+(0,.25) & "c1"}, 
  "c"+(0,.10)*!<0pt,-6pt>{\si},
"x3"="L13",
"x4"="L14",
  "L15";"L15" +(0,-.5)="L15"**\dir{-},
  "L16";"L16" +(0,-.5)="L16"**\dir{-},
  "L17";"L17" +(0,-.5)="L17"**\dir{-},
  "L11"="x1",       
  "L12"="x2",  , 
  "x1"+(0,-.5)="x3",
  "x2"+(0,-.5)="x4",
  "x1";"x2"**\dir{}?(.5)="m", 
  "m" + (0,-.25)="c", 
  "x1";"x4"**\crv{ "x1"+(0,-.25) & "c" & "x4"+(0,.25)},  
  "x2";"c"+<3pt,2pt>="c2"**\crv{"x2"+(0,-.25) & "c2"},
  "x3";"c"+<-3pt,-2pt>="c1"**\crv{"x3"+(0,.25) & "c1"}, 
  "c"+(0,.10)*!<0pt,-6pt>{\si},
"x3"="L11",
"x4"="L12",
  "L13";"L13" +(0,-.5)="L13"**\dir{-},
  "L14";"L14" +(0,-.5)="L14"**\dir{-},
  "L15";"L15" +(0,-.5)="L15"**\dir{-},
  "L16";"L16" +(0,-.5)="L16"**\dir{-},
  "L17";"L17" +(0,-.5)="L17"**\dir{-},
  "L11";"L11" +(0,-.5)="L11"**\dir{-},
  "L12";"L12" +(0,-.5)="L12"**\dir{-},
  "L13";"L13" +(0,-.5)="L13"**\dir{-},
  "L14"="x1",
  "L15"="x2",  
  "x1";"x2"**\dir{}?(.5)="m", 
  "m" + (0,-.25)="c", 
  "m" + (0,-.5)="x3", 
  "x1";"x3"**\crv{"x1"-(0,.25) & "c" & "x3"},
  "x2";"x3"**\crv{"x2"-(0,.25) & "c"  & "x3"}, 
"x3"="L14",
  "L16";"L16" +(0,-.5)="L15"**\dir{-},
  "L17";"L17" +(0,-.5)="L16"**\dir{-},
  "L11";"L11" +(0,-.5)="L11"**\dir{-},
  "L12";"L12" +(0,-.5)="L12"**\dir{-},
  "L13"="x1", 
  "L14"="x2",   
  "x1";"x2"**\dir{}?(.5)="m", 
  "m" + (0,-.25)="c", 
  "m" + (0,-.5)="x3", 
  "x1";"x3"**\crv{"x1"-(0,.25) & "c" & "x3"},
  "x2";"x3"**\crv{"x2"-(0,.25) & "c"  & "x3"}, 
"x3"="L13",
  "L15";"L15" +(0,-.5)="L14"**\dir{-},
  "L16";"L16" +(0,-.5)="L15"**\dir{-},
  "L11"="x1", 
  "L12"="x2",  
  "x1";"x2"**\dir{}?(.5)="m", 
  "m" + (0,-.25)="c", 
  "m" + (0,-.5)="x3", 
  "x1";"x3"**\crv{"x1"-(0,.25) & "c" & "x3"},
  "x2";"x3"**\crv{"x2"-(0,.25) & "c"  & "x3"}, 
"x3"="L11",
  "L13";"L13" +(0,-.5)="L12"**\dir{-},
  "L14"="x1", 
  "L15"="x2",   
  "x1";"x2"**\dir{}?(.5)="m", 
  "m" + (0,-.25)="c", 
  "m" + (0,-.5)="x3", 
  "x1";"x3"**\crv{"x1"-(0,.25) & "c" & "x3"},
  "x2";"x3"**\crv{"x2"-(0,.25) & "c"  & "x3"}, 
"x3"="L13",
  "L11"="x1",       
  "L12"="x2",  , 
  "x1"+(0,-.5)="x3",
  "x2"+(0,-.5)="x4",
  "x1";"x2"**\dir{}?(.5)="m", 
  "m" + (0,-.25)="c", 
  "x1";"x4"**\crv{ "x1"+(0,-.25) & "c" & "x4"+(0,.25)},  
  "x2";"c"+<3pt,2pt>="c2"**\crv{"x2"+(0,-.25) & "c2"},
  "x3";"c"+<-3pt,-2pt>="c1"**\crv{"x3"+(0,.25) & "c1"}, 
  "c"+(0,.10)*!<0pt,-6pt>{\si},
"x3"="L11",
"x4"="L12",
  "L13";"L13" +(0,-.5)="L13"**\dir{-},
  "L11";"L11" +(0,-.5)="L11"**\dir{-},
  "L12"="x1", 
  "L13"="x2",   
  "x1";"x2"**\dir{}?(.5)="m", 
  "m" + (0,-.25)="c", 
  "m" + (0,-.5)="x3", 
  "x1";"x3"**\crv{"x1"-(0,.25) & "c" & "x3"},
  "x2";"x3"**\crv{"x2"-(0,.25) & "c"  & "x3"}, 
"x3"="L12",
  "L11"-(0,.5)*!<0pt,-6pt>{H}, 
  "L12"-(0,.5)*!<0pt,-6pt>{H}, 
 "L12"+(3,6.5)*!<0pt,6pt>{\overset{~\eqref{m(S2,ws)si Del}}{=}},
\endxy 
\]

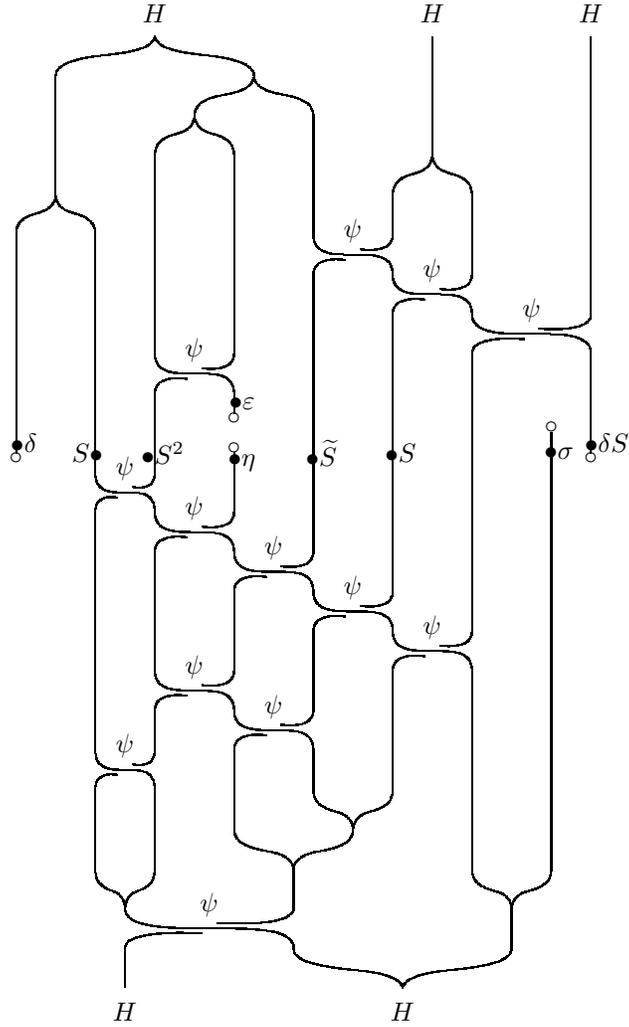
\[
\xy /r2.5pc/:,
(0,0)= "L11",
(3.5,0)= "L12",
(5.5,0)= "L13",
  "L11"+(0,.5)*!<0pt,6pt>{H}, 
  "L12"+(0,.5)*!<0pt,6pt>{H}, 
  "L13"+(0,.5)*!<0pt,6pt>{H}, 
   "L11"="x1", 
  "x1"+(-1.25,-.5)="x2",
  "x1"+(1.25,-.5)="x3",
  "x1"+(0,-.25)="c",  
   "x1";"x2"**\crv{"x1" & "c" & "x2"+(0,.25)}, 
   "x1";"x3"**\crv{"x1" & "c" & "x3"+(0,.25)}, 
"x2"= "L21",
"x3"= "L22",
  "L12";"L12" +(0,-.5)="L23"**\dir{-},
  "L13";"L13" +(0,-.5)="L24"**\dir{-},
"L21"= "L11",
"L22"= "L12",
"L23"= "L13",
"L24"= "L14",
  "L11";"L11" +(0,-.5)="L11"**\dir{-},
   "L12"="x1", 
  "x1"+(-.75,-.5)="x2",
  "x1"+(.75,-.5)="x3",
  "x1"+(0,-.25)="c",  
   "x1";"x2"**\crv{"x1" & "c" & "x2"+(0,.25)}, 
   "x1";"x3"**\crv{"x1" & "c" & "x3"+(0,.25)}, 
"x2"= "L22",
"x3"= "L23",
  "L13";"L13" +(0,-.5)="L24"**\dir{-},
  "L14";"L14" +(0,-.5)="L25"**\dir{-},
"L22"= "L12",
"L23"= "L13",
"L24"= "L14",
"L25"= "L15",
  "L11";"L11" +(0,-.5)="L11"**\dir{-},
   "L12"="x1", 
  "x1"+(-.5,-.5)="x2",
  "x1"+(.5,-.5)="x3",
  "x1"+(0,-.25)="c",  
   "x1";"x2"**\crv{"x1" & "c" & "x2"+(0,.25)}, 
   "x1";"x3"**\crv{"x1" & "c" & "x3"+(0,.25)}, 
"x2"= "L22",
"x3"= "L23",
  "L13";"L13" +(0,-.5)="L24"**\dir{-},
  "L14";"L14" +(0,-.5)="L25"**\dir{-},
  "L15";"L15" +(0,-.5)="L26"**\dir{-},
"L22"= "L12",
"L23"= "L13",
"L24"= "L14",
"L25"= "L15",
"L26"= "L16",
  "L11";"L11" +(0,-.5)="L11"**\dir{-},
  "L12";"L12" +(0,-.5)="L12"**\dir{-},
  "L13";"L13" +(0,-.5)="L13"**\dir{-},
  "L14";"L14" +(0,-.5)="L14"**\dir{-},
   "L15"="x1", 
  "x1"+(-.5,-.5)="x2",
  "x1"+(.5,-.5)="x3",
  "x1"+(0,-.25)="c",  
   "x1";"x2"**\crv{"x1" & "c" & "x2"+(0,.25)}, 
   "x1";"x3"**\crv{"x1" & "c" & "x3"+(0,.25)}, 
"x2"= "L25",
"x3"= "L26",
   "L16";"L16" +(0,-.5)="L27"**\dir{-},
"L25"= "L15",
"L26"= "L16",
"L27"= "L17",
   "L11"="x1", 
  "x1"+(-.5,-.5)="x2",
  "x1"+(.5,-.5)="x3",
  "x1"+(0,-.25)="c",  
   "x1";"x2"**\crv{"x1" & "c" & "x2"+(0,.25)}, 
   "x1";"x3"**\crv{"x1" & "c" & "x3"+(0,.25)}, 
"x2"= "L21",
"x3"= "L22",
  "L12";"L12" +(0,-.5)="L23"**\dir{-},
  "L13";"L13" +(0,-.5)="L24"**\dir{-},
  "L14";"L14" +(0,-.5)="L25"**\dir{-},
  "L15";"L15" +(0,-.5)="L26"**\dir{-},
  "L16";"L16" +(0,-.5)="L27"**\dir{-},
  "L17";"L17" +(0,-.5)="L28"**\dir{-},
"L21"= "L11",
"L22"= "L12",
"L23"= "L13",
"L24"= "L14",
"L25"= "L15",
"L26"= "L16",
"L27"= "L17",
"L28"= "L18",
  "L11";"L11" +(0,-.5)="L11"**\dir{-},
  "L12";"L12" +(0,-.5)="L12"**\dir{-},
  "L13";"L13" +(0,-.5)="L13"**\dir{-},
  "L14";"L14" +(0,-.5)="L14"**\dir{-},
  "L15"="x1",       
  "L16"="x2",  , 
  "x1"+(0,-.5)="x3",
  "x2"+(0,-.5)="x4",
  "x1";"x2"**\dir{}?(.5)="m", 
  "m" + (0,-.25)="c", 
  "x1";"x4"**\crv{ "x1"+(0,-.25) & "c" & "x4"+(0,.25)},  
  "x2";"c"+<3pt,2pt>="c2"**\crv{"x2"+(0,-.25) & "c2"},
  "x3";"c"+<-3pt,-2pt>="c1"**\crv{"x3"+(0,.25) & "c1"}, 
  "c"+(0,.10)*!<0pt,-6pt>{\si},
"x3"="L15",
"x4"="L16",
  "L17";"L17" +(0,-.5)="L17"**\dir{-},
  "L18";"L18" +(0,-.5)="L18"**\dir{-},
  "L11";"L11" +(0,-.5)="L11"**\dir{-},
  "L12";"L12" +(0,-.5)="L12"**\dir{-},
  "L13";"L13" +(0,-.5)="L13"**\dir{-},
  "L14";"L14" +(0,-.5)="L14"**\dir{-},
  "L15";"L15" +(0,-.5)="L15"**\dir{-},
  "L16"="x1",       
  "L17"="x2",  , 
  "x1"+(0,-.5)="x3",
  "x2"+(0,-.5)="x4",
  "x1";"x2"**\dir{}?(.5)="m", 
  "m" + (0,-.25)="c", 
  "x1";"x4"**\crv{ "x1"+(0,-.25) & "c" & "x4"+(0,.25)},  
  "x2";"c"+<3pt,2pt>="c2"**\crv{"x2"+(0,-.25) & "c2"},
  "x3";"c"+<-3pt,-2pt>="c1"**\crv{"x3"+(0,.25) & "c1"}, 
  "c"+(0,.10)*!<0pt,-6pt>{\si},
"x3"="L16",
"x4"="L17",
  "L18";"L18" +(0,-.5)="L18"**\dir{-},
  "L11";"L11" +(0,-.5)="L11"**\dir{-},
  "L12";"L12" +(0,-.5)="L12"**\dir{-},
  "L13";"L13" +(0,-.5)="L13"**\dir{-},
  "L14";"L14" +(0,-.5)="L14"**\dir{-},
  "L15";"L15" +(0,-.5)="L15"**\dir{-},
  "L16";"L16" +(0,-.5)="L16"**\dir{-},
  "L17"="x1",       
  "L18"="x2",  , 
  "x1"+(0,-.5)="x3",
  "x2"+(0,-.5)="x4",
  "x1";"x2"**\dir{}?(.5)="m", 
  "m" + (0,-.25)="c", 
  "x1";"x4"**\crv{ "x1"+(0,-.25) & "c" & "x4"+(0,.25)},  
  "x2";"c"+<3pt,2pt>="c2"**\crv{"x2"+(0,-.25) & "c2"},
  "x3";"c"+<-3pt,-2pt>="c1"**\crv{"x3"+(0,.25) & "c1"}, 
  "c"+(0,.10)*!<0pt,-6pt>{\si},
"x3"="L17",
"x4"="L18",
  "L11";"L11" +(0,-.5)="L11"**\dir{-},
  "L12";"L12" +(0,-.5)="L12"**\dir{-},
  "L13"="x1",       
  "L14"="x2",  , 
  "x1"+(0,-.5)="x3",
  "x2"+(0,-.5)="x4",
  "x1";"x2"**\dir{}?(.5)="m", 
  "m" + (0,-.25)="c", 
  "x1";"x4"**\crv{ "x1"+(0,-.25) & "c" & "x4"+(0,.25)},  
  "x2";"c"+<3pt,2pt>="c2"**\crv{"x2"+(0,-.25) & "c2"},
  "x3";"c"+<-3pt,-2pt>="c1"**\crv{"x3"+(0,.25) & "c1"}, 
  "c"+(0,.10)*!<0pt,-6pt>{\si},
"x3"="L13",
"x4"="L14",
  "L11";"L11" +(0,-.5)="L11"**\dir{-},
  "L12";"L12" +(0,-.5)="L12"**\dir{-},
  "L13";"L13" +(0,-.5)="L13"**\dir{-},
  "L15";"L15" +(0,-1)="L15"**\dir{-},
  "L16";"L16" +(0,-1)="L16"**\dir{-},
  "L17";"L17" +(0,-1)="L17"**\dir{-},
  "L18";"L18" +(0,-1)="L18"**\dir{-},
  "L11"="x1", 
  "x1"+(0,-.25)="x2",
  "x1";"x2"**\dir{}?(.5)="c", 
  "x1";"x2"**\dir{-},   
"c"-(0,0)*!<-3pt,0pt>{\bullet \delta  },  
"x2"+(0,0)*!<0pt,2pt>{\circ},  
  "L12"-(0,.25)*!<3pt,0pt>{S \bullet }, 
  "L12";"L12" +(0,-.5)="L11"**\dir{-}, 
  "L13"-(0,.25)*!<-3pt,0pt>{\bullet S^2 }, 
  "L13";"L13" +(0,-.5)="L12"**\dir{-},
  "L14"="x1", 
  "x1"+(0,-.25)="x2",
  "x1";"x2"**\dir{}?(.5)="c", 
  "x1";"x2"**\dir{-},   
"c"-(0,0)*!<-3pt,0pt>{\bullet \varepsilon  },  
"x2"+(0,0)*!<0pt,2pt>{\circ},  
   "x2" + (0, -.5)="x1", 
  "x1"+(0,-.25)="x2",
  "x1";"x2"**\dir{}?(.5)="c", 
  "x1";"x2"**\dir{-},   
  "c"-(0,0)*!<-3pt,0pt>{\bullet \eta },  
  "x1"+(0,0)*!<0pt,-2pt>{\circ}, 
"x2"="L13",
  "L15"-(0,.25)*!<-3pt,0pt>{\bullet \widetilde{S} }, 
  "L15";"L15" +(0,-.5)="L14"**\dir{-},
  "L16"-(0,.25)*!<-3pt,0pt>{\bullet S }, 
  "L16";"L16" +(0,-.5)="L15"**\dir{-},
  "L17";"L17" +(0,-.5)="L16"**\dir{-},
   "L17" + (1, 0)="x1", 
  "x1"+(0,-.5)="x2",
  "x1";"x2"**\dir{}?(.5)="c", 
  "x1";"x2"**\dir{-},   
  "c"-(0,0)*!<-3pt,0pt>{\bullet \sigma },  
  "x1"+(0,0)*!<0pt,-2pt>{\circ}, 
"x2"="L17",
  "L18"="x1", 
  "x1"+(0,-.25)="x2",
  "x1";"x2"**\dir{}?(.5)="c", 
  "x1";"x2"**\dir{-},   
"c"-(0,0)*!<-6pt,0pt>{\bullet \delta S},  
"x2"+(0,0)*!<0pt,2pt>{\circ},  
  "L11"="x1",       
  "L12"="x2",  , 
  "x1"+(0,-.5)="x3",
  "x2"+(0,-.5)="x4",
  "x1";"x2"**\dir{}?(.5)="m", 
  "m" + (0,-.25)="c", 
  "x1";"x4"**\crv{ "x1"+(0,-.25) & "c" & "x4"+(0,.25)},  
  "x2";"c"+<3pt,2pt>="c2"**\crv{"x2"+(0,-.25) & "c2"},
  "x3";"c"+<-3pt,-2pt>="c1"**\crv{"x3"+(0,.25) & "c1"}, 
  "c"+(0,.10)*!<0pt,-6pt>{\si},
"x3"="L11",
"x4"="L12",
  "L13";"L13" +(0,-.5)="L13"**\dir{-},
  "L14";"L14" +(0,-.5)="L14"**\dir{-},
  "L15";"L15" +(0,-.5)="L15"**\dir{-},
  "L16";"L16" +(0,-.5)="L16"**\dir{-},
  "L17";"L17" +(0,-.5)="L17"**\dir{-},
  "L11";"L11" +(0,-.5)="L11"**\dir{-},
  "L12"="x1",       
  "L13"="x2",  , 
  "x1"+(0,-.5)="x3",
  "x2"+(0,-.5)="x4",
  "x1";"x2"**\dir{}?(.5)="m", 
  "m" + (0,-.25)="c", 
  "x1";"x4"**\crv{ "x1"+(0,-.25) & "c" & "x4"+(0,.25)},  
  "x2";"c"+<3pt,2pt>="c2"**\crv{"x2"+(0,-.25) & "c2"},
  "x3";"c"+<-3pt,-2pt>="c1"**\crv{"x3"+(0,.25) & "c1"}, 
  "c"+(0,.10)*!<0pt,-6pt>{\si},
"x3"="L12",
"x4"="L13",
  "L14";"L14" +(0,-.5)="L14"**\dir{-},
  "L15";"L15" +(0,-.5)="L15"**\dir{-},
  "L16";"L16" +(0,-.5)="L16"**\dir{-},
  "L17";"L17" +(0,-.5)="L17"**\dir{-},
  "L11";"L11" +(0,-.5)="L11"**\dir{-},
  "L12";"L12" +(0,-.5)="L12"**\dir{-},
  "L13"="x1",       
  "L14"="x2",  , 
  "x1"+(0,-.5)="x3",
  "x2"+(0,-.5)="x4",
  "x1";"x2"**\dir{}?(.5)="m", 
  "m" + (0,-.25)="c", 
  "x1";"x4"**\crv{ "x1"+(0,-.25) & "c" & "x4"+(0,.25)},  
  "x2";"c"+<3pt,2pt>="c2"**\crv{"x2"+(0,-.25) & "c2"},
  "x3";"c"+<-3pt,-2pt>="c1"**\crv{"x3"+(0,.25) & "c1"}, 
  "c"+(0,.10)*!<0pt,-6pt>{\si},
"x3"="L13",
"x4"="L14",
  "L15";"L15" +(0,-.5)="L15"**\dir{-},
  "L16";"L16" +(0,-.5)="L16"**\dir{-},
  "L17";"L17" +(0,-.5)="L17"**\dir{-},
  "L11";"L11" +(0,-.5)="L11"**\dir{-},
  "L12";"L12" +(0,-.5)="L12"**\dir{-},
  "L13";"L13" +(0,-.5)="L13"**\dir{-},
  "L14"="x1",       
  "L15"="x2",  , 
  "x1"+(0,-.5)="x3",
  "x2"+(0,-.5)="x4",
  "x1";"x2"**\dir{}?(.5)="m", 
  "m" + (0,-.25)="c", 
  "x1";"x4"**\crv{ "x1"+(0,-.25) & "c" & "x4"+(0,.25)},  
  "x2";"c"+<3pt,2pt>="c2"**\crv{"x2"+(0,-.25) & "c2"},
  "x3";"c"+<-3pt,-2pt>="c1"**\crv{"x3"+(0,.25) & "c1"}, 
  "c"+(0,.10)*!<0pt,-6pt>{\si},
"x3"="L14",
"x4"="L15",
  "L16";"L16" +(0,-.5)="L16"**\dir{-},
  "L17";"L17" +(0,-.5)="L17"**\dir{-},
  "L11";"L11" +(0,-.5)="L11"**\dir{-},
  "L12";"L12" +(0,-.5)="L12"**\dir{-},
  "L13";"L13" +(0,-.5)="L13"**\dir{-},
  "L14";"L14" +(0,-.5)="L14"**\dir{-},
  "L15"="x1",       
  "L16"="x2",  , 
  "x1"+(0,-.5)="x3",
  "x2"+(0,-.5)="x4",
  "x1";"x2"**\dir{}?(.5)="m", 
  "m" + (0,-.25)="c", 
  "x1";"x4"**\crv{ "x1"+(0,-.25) & "c" & "x4"+(0,.25)},  
  "x2";"c"+<3pt,2pt>="c2"**\crv{"x2"+(0,-.25) & "c2"},
  "x3";"c"+<-3pt,-2pt>="c1"**\crv{"x3"+(0,.25) & "c1"}, 
  "c"+(0,.10)*!<0pt,-6pt>{\si},
"x3"="L15",
"x4"="L16",
  "L17";"L17" +(0,-.5)="L17"**\dir{-},
  "L11";"L11" +(0,-.5)="L11"**\dir{-},
  "L12"="x1",       
  "L13"="x2",  , 
  "x1"+(0,-.5)="x3",
  "x2"+(0,-.5)="x4",
  "x1";"x2"**\dir{}?(.5)="m", 
  "m" + (0,-.25)="c", 
  "x1";"x4"**\crv{ "x1"+(0,-.25) & "c" & "x4"+(0,.25)},  
  "x2";"c"+<3pt,2pt>="c2"**\crv{"x2"+(0,-.25) & "c2"},
  "x3";"c"+<-3pt,-2pt>="c1"**\crv{"x3"+(0,.25) & "c1"}, 
  "c"+(0,.10)*!<0pt,-6pt>{\si},
"x3"="L12",
"x4"="L13",
  "L14";"L14" +(0,-.5)="L14"**\dir{-},
  "L15";"L15" +(0,-.5)="L15"**\dir{-},
  "L16";"L16" +(0,-.5)="L16"**\dir{-},
  "L17";"L17" +(0,-.5)="L17"**\dir{-},
  "L11";"L11" +(0,-.5)="L11"**\dir{-},
  "L12";"L12" +(0,-.5)="L12"**\dir{-},
  "L13"="x1",       
  "L14"="x2",  , 
  "x1"+(0,-.5)="x3",
  "x2"+(0,-.5)="x4",
  "x1";"x2"**\dir{}?(.5)="m", 
  "m" + (0,-.25)="c", 
  "x1";"x4"**\crv{ "x1"+(0,-.25) & "c" & "x4"+(0,.25)},  
  "x2";"c"+<3pt,2pt>="c2"**\crv{"x2"+(0,-.25) & "c2"},
  "x3";"c"+<-3pt,-2pt>="c1"**\crv{"x3"+(0,.25) & "c1"}, 
  "c"+(0,.10)*!<0pt,-6pt>{\si},
"x3"="L13",
"x4"="L14",
  "L15";"L15" +(0,-.5)="L15"**\dir{-},
  "L16";"L16" +(0,-.5)="L16"**\dir{-},
  "L17";"L17" +(0,-.5)="L17"**\dir{-},
  "L11"="x1",       
  "L12"="x2",  , 
  "x1"+(0,-.5)="x3",
  "x2"+(0,-.5)="x4",
  "x1";"x2"**\dir{}?(.5)="m", 
  "m" + (0,-.25)="c", 
  "x1";"x4"**\crv{ "x1"+(0,-.25) & "c" & "x4"+(0,.25)},  
  "x2";"c"+<3pt,2pt>="c2"**\crv{"x2"+(0,-.25) & "c2"},
  "x3";"c"+<-3pt,-2pt>="c1"**\crv{"x3"+(0,.25) & "c1"}, 
  "c"+(0,.10)*!<0pt,-6pt>{\si},
"x3"="L11",
"x4"="L12",
  "L13";"L13" +(0,-.5)="L13"**\dir{-},
  "L14";"L14" +(0,-.5)="L14"**\dir{-},
  "L15";"L15" +(0,-.5)="L15"**\dir{-},
  "L16";"L16" +(0,-.5)="L16"**\dir{-},
  "L17";"L17" +(0,-.5)="L17"**\dir{-},
  "L11";"L11" +(0,-.5)="L11"**\dir{-},
  "L12";"L12" +(0,-.5)="L12"**\dir{-},
  "L13";"L13" +(0,-.5)="L13"**\dir{-},
  "L14"="x1", 
  "L15"="x2",  
  "x1";"x2"**\dir{}?(.5)="m", 
  "m" + (0,-.25)="c", 
  "m" + (0,-.5)="x3", 
  "x1";"x3"**\crv{"x1"-(0,.25) & "c" & "x3"},
  "x2";"x3"**\crv{"x2"-(0,.25) & "c"  & "x3"}, 
"x3"="L14",
  "L16";"L16" +(0,-.5)="L15"**\dir{-},
  "L17";"L17" +(0,-.5)="L16"**\dir{-},
  "L11";"L11" +(0,-.5)="L11"**\dir{-},
  "L12";"L12" +(0,-.5)="L12"**\dir{-},
  "L13"="x1", 
  "L14"="x2",  
  "x1";"x2"**\dir{}?(.5)="m", 
  "m" + (0,-.25)="c", 
  "m" + (0,-.5)="x3", 
  "x1";"x3"**\crv{"x1"-(0,.25) & "c" & "x3"},
  "x2";"x3"**\crv{"x2"-(0,.25) & "c"  & "x3"}, 
"x3"="L13",
  "L15";"L15" +(0,-.5)="L14"**\dir{-},
  "L16";"L16" +(0,-.5)="L15"**\dir{-},
  "L11"="x1", 
  "L12"="x2",  
  "x1";"x2"**\dir{}?(.5)="m", 
  "m" + (0,-.25)="c", 
  "m" + (0,-.5)="x3", 
  "x1";"x3"**\crv{"x1"-(0,.25) & "c" & "x3"},
  "x2";"x3"**\crv{"x2"-(0,.25) & "c"  & "x3"}, 
"x3"="L11",
  "L13";"L13" +(0,-.5)="L12"**\dir{-},
  "L14"="x1", 
  "L15"="x2",  
  "x1";"x2"**\dir{}?(.5)="m", 
  "m" + (0,-.25)="c", 
  "m" + (0,-.5)="x3", 
  "x1";"x3"**\crv{"x1"-(0,.25) & "c" & "x3"},
  "x2";"x3"**\crv{"x2"-(0,.25) & "c"  & "x3"}, 
"x3"="L13",
  "L11"="x1",       
  "L12"="x2",  , 
  "x1"+(0,-.5)="x3",
  "x2"+(0,-.5)="x4",
  "x1";"x2"**\dir{}?(.5)="m", 
  "m" + (0,-.25)="c", 
  "x1";"x4"**\crv{ "x1"+(0,-.25) & "c" & "x4"+(0,.25)},  
  "x2";"c"+<3pt,2pt>="c2"**\crv{"x2"+(0,-.25) & "c2"},
  "x3";"c"+<-3pt,-2pt>="c1"**\crv{"x3"+(0,.25) & "c1"}, 
  "c"+(0,.10)*!<0pt,-6pt>{\si},
"x3"="L11",
"x4"="L12",
  "L13";"L13" +(0,-.5)="L13"**\dir{-},
  "L11";"L11" +(0,-.5)="L11"**\dir{-},
  "L12"="x1", 
  "L13"="x2",  
  "x1";"x2"**\dir{}?(.5)="m", 
  "m" + (0,-.25)="c", 
  "m" + (0,-.5)="x3", 
  "x1";"x3"**\crv{"x1"-(0,.25) & "c" & "x3"},
  "x2";"x3"**\crv{"x2"-(0,.25) & "c"  & "x3"}, 
"x3"="L12",
  "L11"-(0,.5)*!<0pt,-6pt>{H}, 
  "L12"-(0,.5)*!<0pt,-6pt>{H}, 
\endxy
\]

By the same method, using standard identities and properties, including:
\[(1,\si)(\si,1)(1,\si) =(\si,1)(1,\si)(\si,1) ,  \]
one reduces the above digram to the following one:

\[
\xy /r2.5pc/:,   
(0,0)="L11",
(2,0)="L12",
(3,0)="L13",
  "L11"+(0,.5)*!<0pt,6pt>{H},
  "L12"+(0,.5)*!<0pt,6pt>{H}, 
  "L13"+(0,.5)*!<0pt,6pt>{H},
"L11";"L11" + (0,-.5) = "L21" **\dir{-},  
   "L12"="x1", 
  "x1"+(-.5,-.5)="x2",
  "x1"+(.5,-.5)="x3",
  "x1"+(0,-.25)="c",  
   "x1";"x2"**\crv{"x1" & "c" & "x2"+(0,.25)}, 
   "x1";"x3"**\crv{"x1" & "c" & "x3"+(0,.25)}, 
"x2"= "L22",
"x3"= "L23",
"L13";"L13" + (0,-.5) = "L24" **\dir{-},  
   "L21"="x1", 
  "x1"+(-.5,-.5)="x2",
  "x1"+(.5,-.5)="x3",
  "x1"+(0,-.25)="c",  
   "x1";"x2"**\crv{"x1" & "c" & "x2"+(0,.25)}, 
   "x1";"x3"**\crv{"x1" & "c" & "x3"+(0,.25)}, 
"x2"= "L31",
"x3"= "L32",
 "L22"="x1", 
  "x1"+(0,-.5)="x2",
  "x1";"x2"**\dir{}?(.5)="c", 
  "x1";"x2"**\dir{-},  
  "c"-(0,0)*!<-4pt,0pt>{\bullet \widetilde{S} }, 
"x2"= "L33",
"L23"="x1", 
  "x1"+(0,-.5)="x2",
  "x1";"x2"**\dir{}?(.5)="c", 
  "x1";"x2"**\dir{-},  
  "c"-(0,0)*!<-4pt,0pt>{\bullet S }, 
"x2"= "L34",
"L24";"L24" + (0,-.5) = "L35" **\dir{-},  
"L24" +(1,-.25) = "x1", 
  "x1"+(0,-.25)="x2",
  "x1";"x2"**\dir{}?(.5)="c", 
  "x1";"x2"**\dir{-},   
  "c"-(0,0)*!<-3pt,0pt>{\bullet \sigma },
  "x1"+(0,0)*!<0pt,-2pt>{\circ},
"x2"= "L36",
"L31"="L21",
"L32"="L22",
"L33"="L23",
"L34"="L24",
"L35"="L25",
"L36"="L26",
"L21";"L21" + (0,-.5) = "L31" **\dir{-},  
"L22"="x1", 
  "x1"+(0,-.5)="x2",
  "x1";"x2"**\dir{}?(.5)="c", 
  "x1";"x2"**\dir{-},  
  "c"-(0,0)*!<-4pt,0pt>{\bullet S }, 
"x2"= "L32",
   "L31" ="x1", 
    "L32"="x2",   
  "x1";"x2"**\dir{}?(.5)="m", 
  "m" + (0,-.25)="c", 
  "m" + (0,-.5)="x3", 
  "x1";"x3"**\crv{"x1"-(0,.25) & "c" & "x3"},
  "x2";"x3"**\crv{"x2"-(0,.25) & "c"  & "x3"}, 
"x3"="L31",
"L31" ="x1", 
  "x1"+(0,-.25)="x2",
  "x1";"x2"**\dir{}?(.5)="c", 
  "x1";"x2"**\dir{-},    
  "c"-(0,0)*!<-3pt,0pt>{\bullet \delta },  
  "x2"+(0,0)*!<0pt,2pt>{\circ}, 
"L23";"L23" + (0,-1.5) = "L31" **\dir{-},  
"L24";"L24" + (0,-1.5) = "L32" **\dir{-},  
"L25";"L25" + (0,-1.5) = "L33" **\dir{-},  
"L26";"L26" + (0,-1.5) = "L34" **\dir{-},  
"L31"="L21",
"L32"="L22",
"L33"="L23",
"L34"="L24",
   "L21"="x1",       
   "L22"="x2", 
  "x1"+(0,-.5)="x3",
  "x2"+(0,-.5)="x4",
  "x1";"x2"**\dir{}?(.5)="m", 
  "m" + (0,-.25)="c", 
  "x1";"x4"**\crv{ "x1"+(0,-.25) & "c" & "x4"+(0,.25)},  
  "x2";"c"+<3pt,2pt>="c2"**\crv{"x2"+(0,-.25) & "c2"},
  "x3";"c"+<-3pt,-2pt>="c1"**\crv{"x3"+(0,.25) & "c1"}, 
  "c"+(0,.10)*!<0pt,-6pt>{\si},  
"x3"= "L31",
"x4"= "L32",
"L23";"L23" + (0,-.5) = "L33" **\dir{-},  
"L24";"L24" + (0,-.5) = "L34" **\dir{-},  
"L31"="L21",
"L32"="L22",
"L33"="L23",
"L34"="L24",
"L21";"L21" + (0,-.5) = "L31" **\dir{-},  
   "L22"="x1",      
   "L23"="x2",  
  "x1"+(0,-.5)="x3",
  "x2"+(0,-.5)="x4",
  "x1";"x2"**\dir{}?(.5)="m", 
  "m" + (0,-.25)="c", 
  "x1";"x4"**\crv{ "x1"+(0,-.25) & "c" & "x4"+(0,.25)},  
  "x2";"c"+<3pt,2pt>="c2"**\crv{"x2"+(0,-.25) & "c2"},
  "x3";"c"+<-3pt,-2pt>="c1"**\crv{"x3"+(0,.25) & "c1"}, 
  "c"+(0,.10)*!<0pt,-6pt>{\si},  
"x3"= "L32",
"x4"= "L33",
"L24";"L24" + (0,-.5) = "L34" **\dir{-},  
"L31"="L21",
"L32"="L22",
"L33"="L23",
"L34"="L24",
"L21";"L21" + (0,-.5) = "L31" **\dir{-},  
"L22";"L22" + (0,-.5) = "L32" **\dir{-},  
   "L23"="x1", 
   "L24"="x2",   
  "x1";"x2"**\dir{}?(.5)="m", 
  "m" + (0,-.25)="c", 
  "m" + (0,-.5)="x3", 
  "x1";"x3"**\crv{"x1"-(0,.25) & "c" & "x3"},
  "x2";"x3"**\crv{"x2"-(0,.25) & "c"  & "x3"}, 
"x3"="L33", 
"L31"="L21",
"L32"="L22",
"L33"="L23",
   "L21"="x1", 
   "L22"="x2",  
  "x1";"x2"**\dir{}?(.5)="m", 
  "m" + (0,-.25)="c", 
  "m" + (0,-.5)="x3", 
  "x1";"x3"**\crv{"x1"-(0,.25) & "c" & "x3"},
  "x2";"x3"**\crv{"x2"-(0,.25) & "c"  & "x3"}, 
"x3"="L31", 
"L23";"L23" + (0,-.5) = "L32" **\dir{-},  
"L31"="L21",
"L32"="L22",
 "L31"-(0,.5)*!<0pt,-6pt>{H},   
 "L32"-(0,.5)*!<0pt,-6pt>{H}, 
"L32"+(1.25,2.5)*!<0pt,-6pt>{\overset{~\eqref{antipode m(1,S)Del = eta epn}}{=}}, 
\endxy 
\]
\[
\xy /r2.5pc/:,   
(0,0)="L11",
(1,0)="L12",
(2.5,0)="L13",
  "L11"+(0,.5)*!<0pt,6pt>{H},
  "L12"+(0,.5)*!<0pt,6pt>{H}, 
  "L13"+(0,.5)*!<0pt,6pt>{H},
    "L11"="x1",
  "x1"+(0,-.25)="x2",
  "x1";"x2"**\dir{}?(.5)="c", 
  "x1";"x2"**\dir{-},   
  "c"-(0,0)*!<-3pt,0pt>{\bullet \varepsilon },  
  "x2"+(0,0)*!<0pt,2pt>{\circ}, 
"L12";"L12" + (0,-.5) = "L21" **\dir{-},  
"L13";"L13" + (0,-.5) = "L22" **\dir{-},  
   "L21"="x1", 
  "x1"+(-.5,-.5)="x2",
  "x1"+(.5,-.5)="x3",
  "x1"+(0,-.25)="c",  
   "x1";"x2"**\crv{"x1" & "c" & "x2"+(0,.25)}, 
   "x1";"x3"**\crv{"x1" & "c" & "x3"+(0,.25)}, 
"x2"= "L31",
"x3"= "L32",
"L22";"L22" + (0,-.5) = "L33" **\dir{-},  
     "L13" +(1, -.25) ="x1", 
  "x1"+(0,-.25)="x2",
  "x1";"x2"**\dir{}?(.5)="c", 
  "x1";"x2"**\dir{-},   
  "c"-(0,0)*!<-3pt,0pt>{\bullet \sigma },  
  "x1"+(0,0)*!<0pt,-2pt>{\circ},  
"x2" = "L34",
   "L31"="x1", 
  "x1"+(0,-.5)="x2",
  "x1";"x2"**\dir{}?(.5)="c", 
  "x1";"x2"**\dir{-},  
  "c"-(0,0)*!<4pt,0pt>{\widetilde{S} \bullet}, 
"x2" = "L31",
   "L32"="x1", 
  "x1"+(0,-.5)="x2",
  "x1";"x2"**\dir{}?(.5)="c", 
  "x1";"x2"**\dir{-},  
  "c"-(0,0)*!<-4pt,0pt>{\bullet S }, 
"x2" = "L32",
    "L31"="x1",   
   "L32"="x2", 
  "x1"+(0,-.5)="x3",
  "x2"+(0,-.5)="x4",
  "x1";"x2"**\dir{}?(.5)="m", 
  "m" + (0,-.25)="c", 
  "x1";"x4"**\crv{ "x1"+(0,-.25) & "c" & "x4"+(0,.25)},  
  "x2";"c"+<3pt,2pt>="c2"**\crv{"x2"+(0,-.25) & "c2"},
  "x3";"c"+<-3pt,-2pt>="c1"**\crv{"x3"+(0,.25) & "c1"}, 
  "c"+(0,.10)*!<0pt,-6pt>{\si}, 
"x3" = "L31",
"x4" = "L32",
"L33";"L33" + (0,-1) = "L33" **\dir{-},  
"L34";"L34" + (0,-1.5) = "L34" **\dir{-},  
"L31"="L21",
"L32"="L22",
"L33"="L23",
"L34"="L24",
"L21";"L21" + (0,-.5) = "L31" **\dir{-},  
    "L22"="x1",     
   "L23"="x2",  
  "x1"+(0,-.5)="x3",
  "x2"+(0,-.5)="x4",
  "x1";"x2"**\dir{}?(.5)="m", 
  "m" + (0,-.25)="c", 
  "x1";"x4"**\crv{ "x1"+(0,-.25) & "c" & "x4"+(0,.25)},  
  "x2";"c"+<3pt,2pt>="c2"**\crv{"x2"+(0,-.25) & "c2"},
  "x3";"c"+<-3pt,-2pt>="c1"**\crv{"x3"+(0,.25) & "c1"}, 
  "c"+(0,.10)*!<0pt,-6pt>{\si}, 
"x3" = "L32",
"x4" = "L33",
"L24";"L24" + (0,-.5) = "L34" **\dir{-},  
"L31"="L21",
"L32"="L22",
"L33"="L23",
"L34"="L24",
   "L21"="x1", 
   "L22"="x2",  
  "x1";"x2"**\dir{}?(.5)="m", 
  "m" + (0,-.25)="c", 
  "m" + (0,-.5)="x3", 
  "x1";"x3"**\crv{"x1"-(0,.25) & "c" & "x3"},
  "x2";"x3"**\crv{"x2"-(0,.25) & "c"  & "x3"}, 
"x3"="L31", 
   "L23"="x1", 
   "L24"="x2",  
  "x1";"x2"**\dir{}?(.5)="m", 
  "m" + (0,-.25)="c", 
  "m" + (0,-.5)="x3", 
  "x1";"x3"**\crv{"x1"-(0,.25) & "c" & "x3"},
  "x2";"x3"**\crv{"x2"-(0,.25) & "c"  & "x3"}, 
"x3"="L32", 
 "L31"-(0,.5)*!<0pt,-6pt>{H},   
 "L32"-(0,.5)*!<0pt,-6pt>{H}, 
  "L22"+(2.5,1)*!<0pt,0pt>{\overset{~\eqref{ Del ws}}{=} }, 
"L13" + (2,0)="L11",
"L13" + (3,0)="L12",
"L13" + (4,0)="L13",
  "L11"+(0,.5)*!<0pt,6pt>{H},
  "L12"+(0,.5)*!<0pt,6pt>{H}, 
  "L13"+(0,.5)*!<0pt,6pt>{H},
    "L11"="x1", 
  "x1"+(0,-.25)="x2",
  "x1";"x2"**\dir{}?(.5)="c", 
  "x1";"x2"**\dir{-},  
  "c"-(0,0)*!<-3pt,0pt>{\bullet \varepsilon }, 
  "x2"+(0,0)*!<0pt,2pt>{\circ},  
"L12";"L12" + (0,-.5) = "L21" **\dir{-},  
"L13";"L13" + (0,-1) = "L22" **\dir{-},  
   "L21"="x1", 
  "x1"+(0,-.5)="x2",
  "x1";"x2"**\dir{}?(.5)="c", 
  "x1";"x2"**\dir{-},  
  "c"-(0,0)*!<-4pt,0pt>{\bullet \widetilde{S} }, 
"x2" = "L31",
"L22";"L22" + (0,-.5) = "L33" **\dir{-},  
     "L22" +(1, -.25) ="x1", 
  "x1"+(0,-.25)="x2",
  "x1";"x2"**\dir{}?(.5)="c", 
  "x1";"x2"**\dir{-},  
  "c"-(0,0)*!<-3pt,0pt>{\bullet \sigma }, 
  "x1"+(0,0)*!<0pt,-2pt>{\circ},  
"x2" = "L34",
   "L31"="x1", 
  "x1"+(-.5,-.5)="x2",
  "x1"+(.5,-.5)="x3",
  "x1"+(0,-.25)="c",  
   "x1";"x2"**\crv{"x1" & "c" & "x2"+(0,.25)}, 
   "x1";"x3"**\crv{"x1" & "c" & "x3"+(0,.25)}, 
"x2"= "L31",
"x3"= "L32",
"L31"="L21",
"L32"="L22",
"L33"="L23",
"L34"="L24",
"L21";"L21" + (0,-.5) = "L31" **\dir{-},  
    "L22"="x1",     
   "L23"="x2",  
  "x1"+(0,-.5)="x3",
  "x2"+(0,-.5)="x4",
  "x1";"x2"**\dir{}?(.5)="m", 
  "m" + (0,-.25)="c", 
  "x1";"x4"**\crv{ "x1"+(0,-.25) & "c" & "x4"+(0,.25)},  
  "x2";"c"+<3pt,2pt>="c2"**\crv{"x2"+(0,-.25) & "c2"},
  "x3";"c"+<-3pt,-2pt>="c1"**\crv{"x3"+(0,.25) & "c1"}, 
  "c"+(0,.10)*!<0pt,-6pt>{\si}, 
"x3" = "L32",
"x4" = "L33",
"L24";"L24" + (0,-.5) = "L34" **\dir{-},  
"L31"="L21",
"L32"="L22",
"L33"="L23",
"L34"="L24",
   "L21"="x1",
   "L22"="x2",  
  "x1";"x2"**\dir{}?(.5)="m", 
  "m" + (0,-.25)="c", 
  "m" + (0,-.5)="x3", 
  "x1";"x3"**\crv{"x1"-(0,.25) & "c" & "x3"},
  "x2";"x3"**\crv{"x2"-(0,.25) & "c"  & "x3"}, 
"x3"="L31", 
   "L23"="x1", 
   "L24"="x2",  
  "x1";"x2"**\dir{}?(.5)="m", 
  "m" + (0,-.25)="c", 
  "m" + (0,-.5)="x3", 
  "x1";"x3"**\crv{"x1"-(0,.25) & "c" & "x3"},
  "x2";"x3"**\crv{"x2"-(0,.25) & "c"  & "x3"}, 
"x3"="L32", 
 "L31"-(0,.5)*!<0pt,-6pt>{H},   
 "L32"-(0,.5)*!<0pt,-6pt>{H}, 
\endxy
\]
$=(m,m)(1,\si,1)(\Delta \widetilde{S},1,\sigma)(\varepsilon,2)= \tau_2 \sigma_0 $\\
\end{proof}

\begin{lemma}
It is easy to verify that
\[m(S^2,\widetilde{S})\si \Delta = \eta \delta \]

\begin{equation}\label{m(S2,ws)si Del}
\xy /r2.5pc/:, 
(0,0) = "L11",
   "L11"="x1", 
  "x1"+(-.5,-.5)="x2",
  "x1"+(.5,-.5)="x3",
  "x1"+(0,-.25)="c",  
   "x1";"x2"**\crv{"x1" & "c" & "x2"+(0,.25)}, 
   "x1";"x3"**\crv{"x1" & "c" & "x3"+(0,.25)}, 
  "x1"+(0,.5)*!<0pt,6pt>{H},  
"x2"= "L11",
"x3"= "L12",
    "L11"="x1",      
   "L12"="x2",  
  "x1"+(0,-.5)="x3",
  "x2"+(0,-.5)="x4",
  "x1";"x2"**\dir{}?(.5)="m", 
  "m" + (0,-.25)="c", 
  "x1";"x4"**\crv{ "x1"+(0,-.25) & "c" & "x4"+(0,.25)},  
  "x2";"c"+<3pt,2pt>="c2"**\crv{"x2"+(0,-.25) & "c2"},
  "x3";"c"+<-3pt,-2pt>="c1"**\crv{"x3"+(0,.25) & "c1"}, 
  "c"+(0,.10)*!<0pt,-4pt>{\si}, 
"x3" = "L11",
"x4" = "L12",
  "L11"="x1", 
  "x1"+(0,-.5)="x2",
  "x1";"x2"**\dir{}?(.5)="c", 
  "x1";"x2"**\dir{-},   
  "c"-(0,0)*!<-4pt,0pt>{\bullet S^2 }, 
"x2"="L11",
  "L12"="x1", 
  "x1"+(0,-.5)="x2",
  "x1";"x2"**\dir{}?(.5)="c", 
  "x1";"x2"**\dir{-},   
  "c"-(0,0)*!<-4pt,0pt>{\bullet \widetilde{S} },  
"x2"="L12",
  "L11"="x1", 
  "L12"="x2",  
  "x1";"x2"**\dir{}?(.5)="m", 
  "m" + (0,-.25)="c", 
  "m" + (0,-.5)="x3", 
  "x1";"x3"**\crv{"x1"-(0,.25) & "c" & "x3"},
  "x2";"x3"**\crv{"x2"-(0,.25) & "c"  & "x3"}, 
"x3"="L11",
"x3"-(0,.5)*!<0pt,-6pt>{H},  
 "x3" + (1.25,1)*!<0pt,-6pt>{=}, 
    (2,0)="x1", 
   "x1"+(0,-.75)="x2",
  "x1";"x2"**\dir{}?(.5)="c", 
  "x1";"x2"**\dir{-},  
  "c"-(0,0)*!<-3pt,0pt>{\bullet \delta },  
  "x2"+(0,0)*!<0pt,2pt>{\circ}, 
  "x1"+(0,.5)*!<0pt,6pt>{H}, 
"x2"="L11",
   (2,-1)="x1",
  "x1"+(0,-1)="x2",
  "x1";"x2"**\dir{}?(.5)="c", 
  "x1";"x2"**\dir{-}, 
  "c"-(0,0)*!<-3pt,0pt>{\bullet \eta }, 
  "x1"+(0,0)*!<0pt,-2pt>{\circ}, 
   "x2"-(0,.5)*!<0pt,-6pt>{H}, 
"x2"="L11",
\endxy  
\end{equation} 
\end{lemma}

We use this lemma in the proof of next theorem.

\begin{theorem} \label{symm and ciclc cndtion}
Under the conditions of Theorem ~\eqref{brdd ver of CM thry in nonsymmetric}, $$\tau_2^3 = \si_{H,H}^2$$.
\end{theorem}

\begin{proof}   We have:

\begin{align*} 
\tau_2^3=& \tau_2 \tau_2^2\\
=&\tau_2(m_2)(\Delta \widetilde{S},1,\sigma)(m_2)(\Delta \widetilde{S},1,\sigma)\\
=&\tau_2(m,m)(1,\si ,1)(\si(\widetilde{S},S)\Delta,1,\sigma)(m,m)(1,\si ,1)(\si(\widetilde{S},S)\Delta,1,\sigma)\\
=&\tau_2(m,m)(1,\si ,1)(\si,2)(\widetilde{S},S,1,\sigma)(\Delta,1)(m,m)(1,\si ,1)(\si,2)\\
~~~~~&(\widetilde{S},S,1,\sigma)(\Delta,1) 
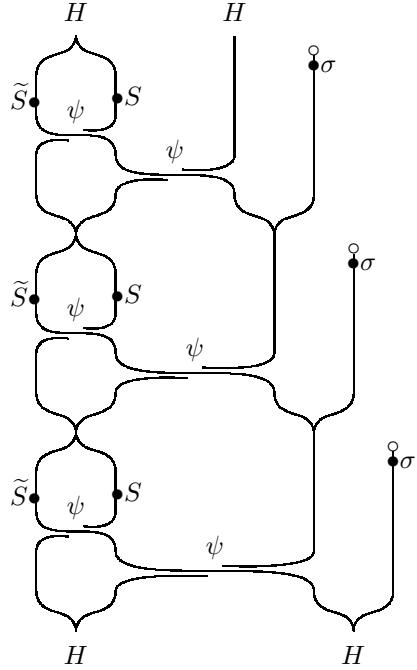
\overset{~\eqref{diagoftau2}}{=}
\end{align*}

\[
\xy /r2.5pc/:,   
(0,0)="L11",
(2,0)="L12",
  "L11"+(0,.5)*!<0pt,6pt>{H},
  "L12"+(0,.5)*!<0pt,6pt>{H}, 
   "L11"="x1", 
  "x1"+(-.5,-.5)="x2",
  "x1"+(.5,-.5)="x3",
  "x1"+(0,-.25)="c", 
   "x1";"x2"**\crv{"x1" & "c" & "x2"+(0,.25)}, 
   "x1";"x3"**\crv{"x1" & "c" & "x3"+(0,.25)}, 
"x2"="L21",
"x3"="L22",
 "L12"="x1", 
  "x1"+(0,-.5)="x2",
  "x1";"x2"**\dir{-},  
"x2"="L23",
  "L12"+ (1,-.25)="x1", 
  "x1"+(0,-.25)="x2",
  "x1";"x2"**\dir{}?(.5)="c", 
  "x1";"x2"**\dir{-},   
  "c"-(0,0)*!<-3pt,0pt>{\bullet \sigma },  
  "x1"+(0,0)*!<0pt,-2pt>{\circ}, 
"x2"="L24",
  "L21"="x1", 
  "x1"+(0,-.5)="x2",
  "x1";"x2"**\dir{}?(.5)="c", 
  "x1";"x2"**\dir{-},  
  "c"-(0,0)*!<4pt,0pt>{\widetilde{S} \bullet},  
"x2"="L31",
  "L22"="x1", 
  "x1"+(0,-.5)="x2",
  "x1";"x2"**\dir{}?(.5)="c", 
  "x1";"x2"**\dir{-},  
  "c"-(0,0)*!<-4pt,0pt>{\bullet S},  
"x2"="L32",
  "L23"="x1", 
  "x1"+(0,-.5)="x2",
  "x1";"x2"**\dir{-},   
"x2"="L33",
  "L24"="x1", 
  "x1"+(0,-.5)="x2",
  "x1";"x2"**\dir{-},   
"x2"="L34",
"L33";"L33" + (0,-.5) = "L33" **\dir{-},  
"L34";"L34" + (0,-.5) = "L34" **\dir{-},  
"L31"="L21",
"L32"="L22",
"L33"="L23",
"L34"="L24",
   "L21"="x1",      
     "L22"="x2",  
  "x1"+(0,-.5)="x3",
  "x2"+(0,-.5)="x4",
  "x1";"x2"**\dir{}?(.5)="m", 
  "m" + (0,-.25)="c", 
  "x1";"x4"**\crv{ "x1"+(0,-.25) & "c" & "x4"+(0,.25)},  
  "x2";"c"+<3pt,2pt>="c2"**\crv{"x2"+(0,-.25) & "c2"},
  "x3";"c"+<-3pt,-2pt>="c1"**\crv{"x3"+(0,.25) & "c1"}, 
  "c"+(0,.10)*!<0pt,-6pt>{\si},
"x3"="L31",
"x4"="L32",
"L23"+(0,-0)="L33",
"L24"+(0,-0)="L34",
"L31"="L21",
"L32"="L22",
"L33"="L23",
"L34"="L24",
"L21";"L21" + (0,-.5) = "L31" **\dir{-},  
  "L22"="x1",       
  "L23"="x2",  , 
  "x1"+(0,-.5)="x3",
  "x2"+(0,-.5)="x4",
  "x1";"x2"**\dir{}?(.5)="m", 
  "m" + (0,-.25)="c", 
  "x1";"x4"**\crv{ "x1"+(0,-.25) & "c" & "x4"+(0,.25)},  
  "x2";"c"+<3pt,2pt>="c2"**\crv{"x2"+(0,-.25) & "c2"},
  "x3";"c"+<-3pt,-2pt>="c1"**\crv{"x3"+(0,.25) & "c1"}, 
  "c"+(0,.10)*!<0pt,-6pt>{\si},
"x3"= "L32",
"x4"= "L33",
"L24";"L24" + (0,-.5) = "L34" **\dir{-},  
"L31"="L21",
"L32"="L22",
"L33"="L23",
"L34"="L24",
  "L21"="x1", 
  "L22"="x2",  
  "x1";"x2"**\dir{}?(.5)="m", 
  "m" + (0,-.25)="c", 
  "m" + (0,-.5)="x3", 
  "x1";"x3"**\crv{"x1"-(0,.25) & "c" & "x3"},
  "x2";"x3"**\crv{"x2"-(0,.25) & "c"  & "x3"}, 
"x3"="L31",
  "L23"="x1", 
  "L24"="x2",  
  "x1";"x2"**\dir{}?(.5)="m", 
  "m" + (0,-.25)="c", 
  "m" + (0,-.5)="x3", 
  "x1";"x3"**\crv{"x1"-(0,.25) & "c" & "x3"},
  "x2";"x3"**\crv{"x2"-(0,.25) & "c"  & "x3"}, 
"x3"="L32",
"L31"="L21",
"L32"="L22",
"L21"="L11",
"L22"="L12",
   "L11"="x1", 
  "x1"+(-.5,-.5)="x2",
  "x1"+(.5,-.5)="x3",
  "x1"+(0,-.25)="c",  
   "x1";"x2"**\crv{"x1" & "c" & "x2"+(0,.25)}, 
   "x1";"x3"**\crv{"x1" & "c" & "x3"+(0,.25)}, 
"x2"="L21",
"x3"="L22",
 "L12"="x1", 
  "x1"+(0,-.5)="x2",
  "x1";"x2"**\dir{-},   
"x2"="L23",
  "L12"+ (1,-.25)="x1", 
  "x1"+(0,-.25)="x2",
  "x1";"x2"**\dir{}?(.5)="c", 
  "x1";"x2"**\dir{-},   
  "c"-(0,0)*!<-3pt,0pt>{\bullet \sigma },  
  "x1"+(0,0)*!<0pt,-2pt>{\circ}, 
"x2"="L24",
  "L21"="x1",
  "x1"+(0,-.5)="x2",
  "x1";"x2"**\dir{}?(.5)="c", 
  "x1";"x2"**\dir{-},   
  "c"-(0,0)*!<4pt,0pt>{\widetilde{S} \bullet}, 
"x2"="L31",
  "L22"="x1", 
  "x1"+(0,-.5)="x2",
  "x1";"x2"**\dir{}?(.5)="c", 
  "x1";"x2"**\dir{-},  
  "c"-(0,0)*!<-4pt,0pt>{\bullet S},  
"x2"="L32",
  "L23"="x1", 
  "x1"+(0,-.5)="x2",
  "x1";"x2"**\dir{-}, 
"x2"="L33",
  "L24"="x1", 
  "x1"+(0,-.5)="x2",
  "x1";"x2"**\dir{-},  
"x2"="L34",
"L33";"L33" + (0,-.5) = "L33" **\dir{-},  
"L34";"L34" + (0,-.5) = "L34" **\dir{-},  
"L31"="L21",
"L32"="L22",
"L33"="L23",
"L34"="L24",
   "L21"="x1",      
     "L22"="x2",  
  "x1"+(0,-.5)="x3",
  "x2"+(0,-.5)="x4",
  "x1";"x2"**\dir{}?(.5)="m", 
  "m" + (0,-.25)="c", 
  "x1";"x4"**\crv{ "x1"+(0,-.25) & "c" & "x4"+(0,.25)},  
  "x2";"c"+<3pt,2pt>="c2"**\crv{"x2"+(0,-.25) & "c2"},
  "x3";"c"+<-3pt,-2pt>="c1"**\crv{"x3"+(0,.25) & "c1"}, 
  "c"+(0,.10)*!<0pt,-6pt>{\si},
"x3"="L31",
"x4"="L32",
"L23"+(0,-0)="L33",
"L24"+(0,-0)="L34",
"L31"="L21",
"L32"="L22",
"L33"="L23",
"L34"="L24",
"L21";"L21" + (0,-.5) = "L31" **\dir{-},  
  "L22"="x1",       
  "L23"="x2",  , 
  "x1"+(0,-.5)="x3",
  "x2"+(0,-.5)="x4",
  "x1";"x2"**\dir{}?(.5)="m", 
  "m" + (0,-.25)="c", 
  "x1";"x4"**\crv{ "x1"+(0,-.25) & "c" & "x4"+(0,.25)},  
  "x2";"c"+<3pt,2pt>="c2"**\crv{"x2"+(0,-.25) & "c2"},
  "x3";"c"+<-3pt,-2pt>="c1"**\crv{"x3"+(0,.25) & "c1"}, 
  "c"+(0,.10)*!<0pt,-6pt>{\si},
"x3"= "L32",
"x4"= "L33",
"L24";"L24" + (0,-.5) = "L34" **\dir{-},  
"L31"="L21",
"L32"="L22",
"L33"="L23",
"L34"="L24",
  "L21"="x1", 
  "L22"="x2",  
  "x1";"x2"**\dir{}?(.5)="m", 
  "m" + (0,-.25)="c", 
  "m" + (0,-.5)="x3", 
  "x1";"x3"**\crv{"x1"-(0,.25) & "c" & "x3"},
  "x2";"x3"**\crv{"x2"-(0,.25) & "c"  & "x3"}, 
"x3"="L31",
  "L23"="x1", 
  "L24"="x2", 
  "x1";"x2"**\dir{}?(.5)="m", 
  "m" + (0,-.25)="c", 
  "m" + (0,-.5)="x3", 
  "x1";"x3"**\crv{"x1"-(0,.25) & "c" & "x3"},
  "x2";"x3"**\crv{"x2"-(0,.25) & "c"  & "x3"}, 
"x3"="L32",
"L31"="L21",
"L32"="L22",
"L21"="L11",
"L22"="L12",
   "L11"="x1", 
  "x1"+(-.5,-.5)="x2",
  "x1"+(.5,-.5)="x3",
  "x1"+(0,-.25)="c",  
   "x1";"x2"**\crv{"x1" & "c" & "x2"+(0,.25)}, 
   "x1";"x3"**\crv{"x1" & "c" & "x3"+(0,.25)}, 
"x2"="L21",
"x3"="L22",
 "L12"="x1", 
  "x1"+(0,-.5)="x2",
  "x1";"x2"**\dir{-},   
"x2"="L23",
  "L12"+ (1,-.25)="x1", 
  "x1"+(0,-.25)="x2",
  "x1";"x2"**\dir{}?(.5)="c", 
  "x1";"x2"**\dir{-},  
  "c"-(0,0)*!<-3pt,0pt>{\bullet \sigma }, 
  "x1"+(0,0)*!<0pt,-2pt>{\circ}, 
"x2"="L24",
  "L21"="x1", 
  "x1"+(0,-.5)="x2",
  "x1";"x2"**\dir{}?(.5)="c", 
  "x1";"x2"**\dir{-},  
  "c"-(0,0)*!<4pt,0pt>{\widetilde{S} \bullet},  
"x2"="L31",
  "L22"="x1", 
  "x1"+(0,-.5)="x2",
  "x1";"x2"**\dir{}?(.5)="c", 
  "x1";"x2"**\dir{-},  
  "c"-(0,0)*!<-4pt,0pt>{\bullet S},  
"x2"="L32",
  "L23"="x1", 
  "x1"+(0,-.5)="x2",
  "x1";"x2"**\dir{-},   
"x2"="L33",
  "L24"="x1", 
  "x1"+(0,-.5)="x2",
  "x1";"x2"**\dir{-},   
"x2"="L34",
"L33";"L33" + (0,-.5) = "L33" **\dir{-},  
"L34";"L34" + (0,-.5) = "L34" **\dir{-},  
"L31"="L21",
"L32"="L22",
"L33"="L23",
"L34"="L24",
   "L21"="x1",       
     "L22"="x2",  
  "x1"+(0,-.5)="x3",
  "x2"+(0,-.5)="x4",
  "x1";"x2"**\dir{}?(.5)="m", 
  "m" + (0,-.25)="c", 
  "x1";"x4"**\crv{ "x1"+(0,-.25) & "c" & "x4"+(0,.25)},  
  "x2";"c"+<3pt,2pt>="c2"**\crv{"x2"+(0,-.25) & "c2"},
  "x3";"c"+<-3pt,-2pt>="c1"**\crv{"x3"+(0,.25) & "c1"}, 
  "c"+(0,.10)*!<0pt,-6pt>{\si},
"x3"="L31",
"x4"="L32",
"L23"+(0,-0)="L33",
"L24"+(0,-0)="L34",
"L31"="L21",
"L32"="L22",
"L33"="L23",
"L34"="L24",
"L21";"L21" + (0,-.5) = "L31" **\dir{-},  
  "L22"="x1",       
  "L23"="x2",  , 
  "x1"+(0,-.5)="x3",
  "x2"+(0,-.5)="x4",
  "x1";"x2"**\dir{}?(.5)="m", 
  "m" + (0,-.25)="c", 
  "x1";"x4"**\crv{ "x1"+(0,-.25) & "c" & "x4"+(0,.25)},  
  "x2";"c"+<3pt,2pt>="c2"**\crv{"x2"+(0,-.25) & "c2"},
  "x3";"c"+<-3pt,-2pt>="c1"**\crv{"x3"+(0,.25) & "c1"}, 
  "c"+(0,.10)*!<0pt,-6pt>{\si},
"x3"= "L32",
"x4"= "L33",
"L24";"L24" + (0,-.5) = "L34" **\dir{-},  
"L31"="L21",
"L32"="L22",
"L33"="L23",
"L34"="L24",
  "L21"="x1", 
  "L22"="x2", 
  "x1";"x2"**\dir{}?(.5)="m", 
  "m" + (0,-.25)="c", 
  "m" + (0,-.5)="x3", 
  "x1";"x3"**\crv{"x1"-(0,.25) & "c" & "x3"},
  "x2";"x3"**\crv{"x2"-(0,.25) & "c"  & "x3"}, 
  "x3"+(0,-.5)*!<0pt,-6pt>{H},  
"x3"="L31",
  "L23"="x1", 
  "L24"="x2",  
  "x1";"x2"**\dir{}?(.5)="m", 
  "m" + (0,-.25)="c", 
  "m" + (0,-.5)="x3", 
  "x1";"x3"**\crv{"x1"-(0,.25) & "c" & "x3"},
  "x2";"x3"**\crv{"x2"-(0,.25) & "c"  & "x3"}, 
  "x3"+(0,-.5)*!<0pt,-6pt>{H}, 
"x3"="L32",
"L31"="L21",
"L32"="L22",
\endxy
\]

\begin{align*} 
\overset{~\eqref{naturalityofsi},~\eqref{compatibilityDelm}}{=}&\tau_2(m,m)(1,\si ,1)(\si,2)(\widetilde{S},S,1,\sigma)(m,m,m)(1,\si,3)(\Delta,\Delta,2)\\
~~~~~&(S,1,\widetilde{S},\sigma)(1,\si)(\si,1)(\Delta,1)\\
=&\tau_2(m,m)(1,\si ,1)(\si,2)(\widetilde{S} m,Sm,m,\sigma)(1,\si,3)(\Delta S,\Delta,\widetilde{S},\sigma)\\
~~~~~&(1,\si)(\si,1)(\Delta,1) =
\end{align*}

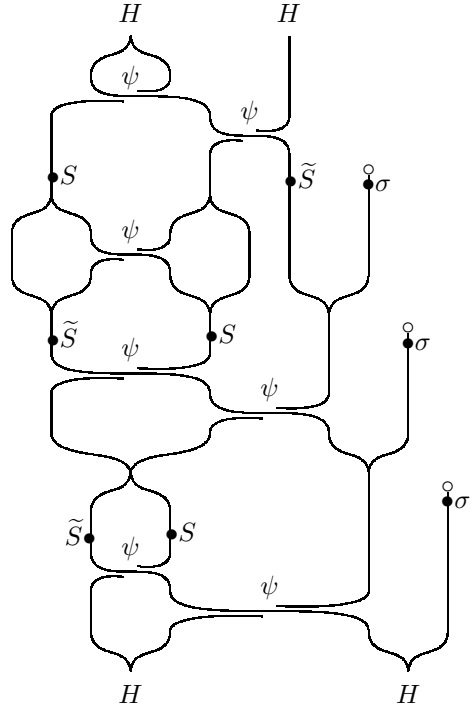
\[
\xy /r2.5pc/:,   
(7,0)="L11",
(9,0)="L12",
  "L11"+(0,.5)*!<0pt,6pt>{H},
  "L12"+(0,.5)*!<0pt,6pt>{H}, 
   "L11"="x1", 
  "x1"+(-.5,-.5)="x2",
  "x1"+(.5,-.5)="x3",
  "x1"+(0,-.25)="c", 
   "x1";"x2"**\crv{"x1" & "c" & "x2"+(0,.25)}, 
   "x1";"x3"**\crv{"x1" & "c" & "x3"+(0,.25)}, 
"x2"="L21",
"x3"="L22",
"L12";"L12" + (0,-.5) = "L23" **\dir{-},  
   "L21"="x1",    
     "L22"="x2", 
  "x1"+(-.5,-.5)="x3",
  "x2"+(.5,-.5)="x4",
  "x1";"x2"**\dir{}?(.5)="m", 
  "m" + (0,-.25)="c", 
  "x1";"x4"**\crv{ "x1"+(0,-.25) & "c" & "x4"+(0,.25)},  
  "x2";"c"+<3pt,2pt>="c2"**\crv{"x2"+(0,-.25) & "c2"},
  "x3";"c"+<-3pt,-2pt>="c1"**\crv{"x3"+(0,.25) & "c1"}, 
  "c"+(0,.10)*!<0pt,-4pt>{\si},
"x3"="L31",
"x4"="L32",
"L23";"L23" + (0,-.5) = "L33" **\dir{-},  
"L31"="L21",
"L32"="L22",
"L33"="L23",
"L21";"L21" + (0,-.5) = "L31" **\dir{-},  
   "L22"="x1",      
     "L23"="x2",  
  "x1"+(0,-.5)="x3",
  "x2"+(0,-.5)="x4",
  "x1";"x2"**\dir{}?(.5)="m", 
  "m" + (0,-.25)="c", 
  "x1";"x4"**\crv{ "x1"+(0,-.25) & "c" & "x4"+(0,.25)},  
  "x2";"c"+<3pt,2pt>="c2"**\crv{"x2"+(0,-.25) & "c2"},
  "x3";"c"+<-3pt,-2pt>="c1"**\crv{"x3"+(0,.25) & "c1"}, 
  "c"+(0,.10)*!<0pt,-6pt>{\si},
"x3"="L32",
"x4"="L33",
"L31"="L21",
"L32"="L22",
"L33"="L23",
  "L21"="x1", 
  "x1"+(0,-.5)="x2",
  "x1";"x2"**\dir{}?(.5)="c", 
  "x1";"x2"**\dir{-},  
  "c"-(0,0)*!<-4pt,0pt>{\bullet S },  
"x2"="L31",
"L22";"L22" + (0,-.5) = "L32" **\dir{-},  
  "L23"="x1",
  "x1"+(0,-.5)="x2",
  "x1";"x2"**\dir{}?(.5)="c", 
  "x1";"x2"**\dir{-},  
  "c"-(0,0)*!<-4pt,0pt>{\bullet \widetilde{S} }, 
"x2"="L33",
  "L23"+ (1,-.25)="x1", 
  "x1"+(0,-.25)="x2",
  "x1";"x2"**\dir{}?(.5)="c", 
  "x1";"x2"**\dir{-},  
  "c"-(0,0)*!<-3pt,0pt>{\bullet \sigma }, 
  "x1"+(0,0)*!<0pt,-2pt>{\circ}, 
"x2"="L34",
"L31"="L21",
"L32"="L22",
"L33"="L23",
"L34"="L24",
   "L21"="x1", 
  "x1"+(-.5,-.5)="x2",
  "x1"+(.5,-.5)="x3",
  "x1"+(0,-.25)="c", 
   "x1";"x2"**\crv{"x1" & "c" & "x2"+(0,.25)}, 
   "x1";"x3"**\crv{"x1" & "c" & "x3"+(0,.25)}, 
"x2"="L31",
"x3"="L32",
   "L22"="x1", 
  "x1"+(-.5,-.5)="x2",
  "x1"+(.5,-.5)="x3",
  "x1"+(0,-.25)="c", 
   "x1";"x2"**\crv{"x1" & "c" & "x2"+(0,.25)}, 
   "x1";"x3"**\crv{"x1" & "c" & "x3"+(0,.25)}, 
"x2"="L33",
"x3"="L34",
"L23";"L23" + (0,-.5) = "L35" **\dir{-},  
"L24";"L24" + (0,-.5) = "L36" **\dir{-},  
"L31"="L21",
"L32"="L22",
"L33"="L23",
"L34"="L24",
"L35"="L25",
"L36"="L26",
"L21";"L21" + (0,-.5) = "L31" **\dir{-},  
   "L22"="x1",    
     "L23"="x2",  
  "x1"+(0,-.5)="x3",
  "x2"+(0,-.5)="x4",
  "x1";"x2"**\dir{}?(.5)="m", 
  "m" + (0,-.25)="c", 
  "x1";"x4"**\crv{ "x1"+(0,-.25) & "c" & "x4"+(0,.25)},  
  "x2";"c"+<3pt,2pt>="c2"**\crv{"x2"+(0,-.25) & "c2"},
  "x3";"c"+<-3pt,-2pt>="c1"**\crv{"x3"+(0,.25) & "c1"}, 
  "c"+(0,.10)*!<0pt,-6pt>{\si},
"x3"="L32",
"x4"="L33",
"L24";"L24" + (0,-.5) = "L34" **\dir{-},  
"L25";"L25" + (0,-.5) = "L35" **\dir{-},  
"L26";"L26" + (0,-.5) = "L36" **\dir{-},  
"L31"="L21",
"L32"="L22",
"L33"="L23",
"L34"="L24",
"L35"="L25",
"L36"="L26",
  "L21"="x1", 
  "L22"="x2", 
  "x1";"x2"**\dir{}?(.5)="m", 
  "m" + (0,-.25)="c", 
  "m" + (0,-.5)="x3", 
  "x1";"x3"**\crv{"x1"-(0,.25) & "c" & "x3"},
  "x2";"x3"**\crv{"x2"-(0,.25) & "c"  & "x3"}, 
"x3"="L31",
  "L23"="x1", 
  "L24"="x2", 
  "x1";"x2"**\dir{}?(.5)="m", 
  "m" + (0,-.25)="c", 
  "m" + (0,-.5)="x3", 
  "x1";"x3"**\crv{"x1"-(0,.25) & "c" & "x3"},
  "x2";"x3"**\crv{"x2"-(0,.25) & "c"  & "x3"}, 
"x3"="L32",
  "L25"="x1", 
  "L26"="x2",   
  "x1";"x2"**\dir{}?(.5)="m", 
  "m" + (0,-.25)="c", 
  "m" + (0,-.5)="x3", 
  "x1";"x3"**\crv{"x1"-(0,.25) & "c" & "x3"},
  "x2";"x3"**\crv{"x2"-(0,.25) & "c"  & "x3"}, 
"x3"="L33",
"L31"="L21",
"L32"="L22",
"L33"="L23",
  "L21"="x1", 
  "x1"+(0,-.5)="x2",
  "x1";"x2"**\dir{}?(.5)="c", 
  "x1";"x2"**\dir{-},   
  "c"-(0,0)*!<-4pt,0pt>{\bullet \widetilde{S} },  
"x2"="L31",
  "L22"="x1", 
  "x1"+(0,-.5)="x2",
  "x1";"x2"**\dir{}?(.5)="c", 
  "x1";"x2"**\dir{-},  
  "c"-(0,0)*!<-4pt,0pt>{\bullet S}, 
"x2"="L32",
"L23";"L23" + (0,-.5) = "L33" **\dir{-},  
  "L23"+ (1,-.25)="x1", 
  "x1"+(0,-.25)="x2",
  "x1";"x2"**\dir{}?(.5)="c", 
  "x1";"x2"**\dir{-}, 
  "c"-(0,0)*!<-3pt,0pt>{\bullet \sigma },  
  "x1"+(0,0)*!<0pt,-2pt>{\circ}, 
"x2"="L34",
"L31"="L21",
"L32"="L22",
"L33"="L23",
"L34"="L24",
   "L21"="x1",    
     "L22"="x2",  
  "x1"+(0,-.5)="x3",
  "x2"+(0,-.5)="x4",
  "x1";"x2"**\dir{}?(.5)="m", 
  "m" + (0,-.25)="c", 
  "x1";"x4"**\crv{ "x1"+(0,-.25) & "c" & "x4"+(0,.25)},  
  "x2";"c"+<3pt,2pt>="c2"**\crv{"x2"+(0,-.25) & "c2"},
  "x3";"c"+<-3pt,-2pt>="c1"**\crv{"x3"+(0,.25) & "c1"}, 
  "c"+(0,.10)*!<0pt,-6pt>{\si},
"x3"="L31",
"x4"="L32",
"L23";"L23" + (0,-.5) = "L33" **\dir{-},  
"L24";"L24" + (0,-.5) = "L34" **\dir{-},  
"L31"="L21",
"L32"="L22",
"L33"="L23",
"L34"="L24",
"L21";"L21" + (0,-.5) = "L31" **\dir{-},  
  "L22"="x1",       
  "L23"="x2",  , 
  "x1"+(0,-.5)="x3",
  "x2"+(0,-.5)="x4",
  "x1";"x2"**\dir{}?(.5)="m", 
  "m" + (0,-.25)="c", 
  "x1";"x4"**\crv{ "x1"+(0,-.25) & "c" & "x4"+(0,.25)},  
  "x2";"c"+<3pt,2pt>="c2"**\crv{"x2"+(0,-.25) & "c2"},
  "x3";"c"+<-3pt,-2pt>="c1"**\crv{"x3"+(0,.25) & "c1"}, 
  "c"+(0,.10)*!<0pt,-6pt>{\si},
"x3"= "L32",
"x4"= "L33",
"L24";"L24" + (0,-.5) = "L34" **\dir{-},  
"L31"="L21",
"L32"="L22",
"L33"="L23",
"L34"="L24",
  "L21"="x1", 
  "L22"="x2",  
  "x1";"x2"**\dir{}?(.5)="m", 
  "m" + (0,-.25)="c", 
  "m" + (0,-.5)="x3", 
  "x1";"x3"**\crv{"x1"-(0,.25) & "c" & "x3"},
  "x2";"x3"**\crv{"x2"-(0,.25) & "c"  & "x3"}, 
"x3"="L31",
  "L23"="x1", 
  "L24"="x2",   
  "x1";"x2"**\dir{}?(.5)="m", 
  "m" + (0,-.25)="c", 
  "m" + (0,-.5)="x3", 
  "x1";"x3"**\crv{"x1"-(0,.25) & "c" & "x3"},
  "x2";"x3"**\crv{"x2"-(0,.25) & "c"  & "x3"}, 
"x3"="L32",
"L31"="L21",
"L32"="L22",
"L21"="L11",
"L22"="L12",
   "L11"="x1", 
  "x1"+(-.5,-.5)="x2",
  "x1"+(.5,-.5)="x3",
  "x1"+(0,-.25)="c",  
   "x1";"x2"**\crv{"x1" & "c" & "x2"+(0,.25)}, 
   "x1";"x3"**\crv{"x1" & "c" & "x3"+(0,.25)}, 
"x2"="L21",
"x3"="L22",
 "L12"="x1", 
  "x1"+(0,-.5)="x2",
  "x1";"x2"**\dir{-},  
"x2"="L23",
  "L12"+ (1,-.25)="x1", 
  "x1"+(0,-.25)="x2",
  "x1";"x2"**\dir{}?(.5)="c", 
  "x1";"x2"**\dir{-},  
  "c"-(0,0)*!<-3pt,0pt>{\bullet \sigma }, 
  "x1"+(0,0)*!<0pt,-2pt>{\circ},  
"x2"="L24",
  "L21"="x1", 
  "x1"+(0,-.5)="x2",
  "x1";"x2"**\dir{}?(.5)="c", 
  "x1";"x2"**\dir{-},  
  "c"-(0,0)*!<4pt,0pt>{\widetilde{S} \bullet}, 
"x2"="L31",
  "L22"="x1", 
  "x1"+(0,-.5)="x2",
  "x1";"x2"**\dir{}?(.5)="c", 
  "x1";"x2"**\dir{-},  
  "c"-(0,0)*!<-4pt,0pt>{\bullet S}, 
"x2"="L32",
  "L23"="x1", 
  "x1"+(0,-.5)="x2",
  "x1";"x2"**\dir{-},  
"x2"="L33",
  "L24"="x1", 
  "x1"+(0,-.5)="x2",
  "x1";"x2"**\dir{-},  
"x2"="L34",
"L33";"L33" + (0,-.5) = "L33" **\dir{-},  
"L34";"L34" + (0,-.5) = "L34" **\dir{-},  
"L31"="L21",
"L32"="L22",
"L33"="L23",
"L34"="L24",
   "L21"="x1",     
     "L22"="x2",  
  "x1"+(0,-.5)="x3",
  "x2"+(0,-.5)="x4",
  "x1";"x2"**\dir{}?(.5)="m", 
  "m" + (0,-.25)="c", 
  "x1";"x4"**\crv{ "x1"+(0,-.25) & "c" & "x4"+(0,.25)},  
  "x2";"c"+<3pt,2pt>="c2"**\crv{"x2"+(0,-.25) & "c2"},
  "x3";"c"+<-3pt,-2pt>="c1"**\crv{"x3"+(0,.25) & "c1"}, 
  "c"+(0,.10)*!<0pt,-6pt>{\si},
"x3"="L31",
"x4"="L32",
"L23"+(0,-0)="L33",
"L24"+(0,-0)="L34",
"L31"="L21",
"L32"="L22",
"L33"="L23",
"L34"="L24",
"L21";"L21" + (0,-.5) = "L31" **\dir{-},  
  "L22"="x1",       
  "L23"="x2",  , 
  "x1"+(0,-.5)="x3",
  "x2"+(0,-.5)="x4",
  "x1";"x2"**\dir{}?(.5)="m", 
  "m" + (0,-.25)="c", 
  "x1";"x4"**\crv{ "x1"+(0,-.25) & "c" & "x4"+(0,.25)},  
  "x2";"c"+<3pt,2pt>="c2"**\crv{"x2"+(0,-.25) & "c2"},
  "x3";"c"+<-3pt,-2pt>="c1"**\crv{"x3"+(0,.25) & "c1"}, 
  "c"+(0,.10)*!<0pt,-6pt>{\si},
"x3"= "L32",
"x4"= "L33",
"L24";"L24" + (0,-.5) = "L34" **\dir{-},  
"L31"="L21",
"L32"="L22",
"L33"="L23",
"L34"="L24",
  "L21"="x1", 
  "L22"="x2",  
  "x1";"x2"**\dir{}?(.5)="m", 
  "m" + (0,-.25)="c", 
  "m" + (0,-.5)="x3", 
  "x1";"x3"**\crv{"x1"-(0,.25) & "c" & "x3"},
  "x2";"x3"**\crv{"x2"-(0,.25) & "c"  & "x3"}, 
  "x3"+(0,-.5)*!<0pt,-6pt>{H},  
"x3"="L31",
  "L23"="x1", 
  "L24"="x2",  
  "x1";"x2"**\dir{}?(.5)="m", 
  "m" + (0,-.25)="c", 
  "m" + (0,-.5)="x3", 
  "x1";"x3"**\crv{"x1"-(0,.25) & "c" & "x3"},
  "x2";"x3"**\crv{"x2"-(0,.25) & "c"  & "x3"}, 
  "x3"+(0,-.5)*!<0pt,-6pt>{H}, 
"x3"="L32",
"L31"="L21",
"L32"="L22",

\endxy 
\]

\begin{align*} 
\overset{~\eqref{Del S}, ~\eqref{S m}, ~\eqref{ws m}}{=}&\tau_2(m,m)(1,\si ,1)(\si,2)(m(\widetilde{S},\widetilde{S})\si,m(S,S)\si,m,\sigma)\\
~~~~~&(1,\si,3)((S,S)\si \Delta ,\Delta,\widetilde{S},\sigma)(1,\si)(\si,1)(\Delta,1)=
\end{align*}

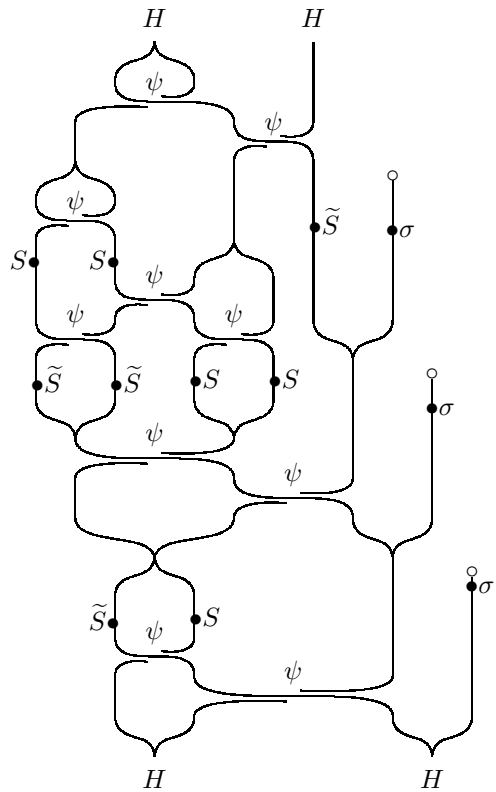
\[
\xy /r2.5pc/:,   
(7,0)="L11",
(9,0)="L12",
  "L11"+(0,.5)*!<0pt,6pt>{H},
  "L12"+(0,.5)*!<0pt,6pt>{H}, 
   "L11"="x1", 
  "x1"+(-.5,-.5)="x2",
  "x1"+(.5,-.5)="x3",
  "x1"+(0,-.25)="c", 
   "x1";"x2"**\crv{"x1" & "c" & "x2"+(0,.25)}, 
   "x1";"x3"**\crv{"x1" & "c" & "x3"+(0,.25)}, 
"x2"="L21",
"x3"="L22",
"L12";"L12" + (0,-.5) = "L23" **\dir{-},  
   "L21"="x1",   
     "L22"="x2", 
  "x1"+(-.5,-.5)="x3",
  "x2"+(.5,-.5)="x4",
  "x1";"x2"**\dir{}?(.5)="m", 
  "m" + (0,-.25)="c", 
  "x1";"x4"**\crv{ "x1"+(0,-.25) & "c" & "x4"+(0,.25)},  
  "x2";"c"+<3pt,2pt>="c2"**\crv{"x2"+(0,-.25) & "c2"},
  "x3";"c"+<-3pt,-2pt>="c1"**\crv{"x3"+(0,.25) & "c1"}, 
  "c"+(0,.10)*!<0pt,-4pt>{\si},
"x3"="L31",
"x4"="L32",
"L23";"L23" + (0,-.5) = "L33" **\dir{-},  
"L31"="L21",
"L32"="L22",
"L33"="L23",
"L21";"L21" + (0,-.5) = "L31" **\dir{-},  
   "L22"="x1",    
     "L23"="x2",  
  "x1"+(0,-.5)="x3",
  "x2"+(0,-.5)="x4",
  "x1";"x2"**\dir{}?(.5)="m", 
  "m" + (0,-.25)="c", 
  "x1";"x4"**\crv{ "x1"+(0,-.25) & "c" & "x4"+(0,.25)},  
  "x2";"c"+<3pt,2pt>="c2"**\crv{"x2"+(0,-.25) & "c2"},
  "x3";"c"+<-3pt,-2pt>="c1"**\crv{"x3"+(0,.25) & "c1"}, 
  "c"+(0,.10)*!<0pt,-4pt>{\si},
"x3"="L32",
"x4"="L33",
"L31"="L11",
"L32"="L12",
"L33"="L13",
   "L11"="x1", 
  "x1"+(-.5,-.5)="x2",
  "x1"+(.5,-.5)="x3",
  "x1"+(0,-.25)="c", 
   "x1";"x2"**\crv{"x1" & "c" & "x2"+(0,.25)}, 
   "x1";"x3"**\crv{"x1" & "c" & "x3"+(0,.25)}, 
"x2"= "L21",
"x3"= "L22",
    "L21"="x1",   
   "L22"="x2",  
  "x1"+(0,-.5)="x3",
  "x2"+(0,-.5)="x4",
  "x1";"x2"**\dir{}?(.5)="m", 
  "m" + (0,-.25)="c", 
  "x1";"x4"**\crv{ "x1"+(0,-.25) & "c" & "x4"+(0,.25)},  
  "x2";"c"+<3pt,2pt>="c2"**\crv{"x2"+(0,-.25) & "c2"},
  "x3";"c"+<-3pt,-2pt>="c1"**\crv{"x3"+(0,.25) & "c1"}, 
  "c"+(0,.10)*!<0pt,-4pt>{\si}, 
"x3" = "L21",
"x4" = "L22",
   "L21"="x1", 
  "x1"+(0,-.5)="x2",
  "x1";"x2"**\dir{}?(.5)="c", 
  "x1";"x2"**\dir{-},   
  "c"-(0,0)*!<4pt,0pt>{S \bullet}, 
"x2" = "L31",
   "L22"="x1", 
  "x1"+(0,-.5)="x2",
  "x1";"x2"**\dir{}?(.5)="c", 
  "x1";"x2"**\dir{-},  
  "c"-(0,0)*!<4pt,0pt>{S \bullet}, 
"x2" = "L32",
"L12";"L12" + (0,-1) = "L33" **\dir{-},  
   "L33"="x1", 
  "x1"+(-.5,-.5)="x2",
  "x1"+(.5,-.5)="x3",
  "x1"+(0,-.25)="c", 
   "x1";"x2"**\crv{"x1" & "c" & "x2"+(0,.25)}, 
   "x1";"x3"**\crv{"x1" & "c" & "x3"+(0,.25)}, 
"x2"="L33",
"x3"="L34",
  "L13"="x1", 
  "x1"+(0,-1.5)="x2",
  "x1";"x2"**\dir{}?(.5)="c", 
  "x1";"x2"**\dir{-},  
  "c"-(0,0)*!<-4pt,0pt>{\bullet \widetilde{S} }, 
"x2"="L35",
  "L13"+ (1,-.25)="x1", 
  "x1"+(0,-1.25)="x2",
  "x1";"x2"**\dir{}?(.5)="c", 
  "x1";"x2"**\dir{-},  
  "c"-(0,0)*!<-3pt,0pt>{\bullet \sigma }, 
  "x1"+(0,0)*!<0pt,-2pt>{\circ}, 
"x2"="L36",
"L31"="L21",
"L32"="L22",
"L33"="L23",
"L34"="L24",
"L35"="L25",
"L36"="L26",
"L21";"L21" + (0,-.5) = "L31" **\dir{-},  
   "L22"="x1",    
     "L23"="x2", 
  "x1"+(0,-.5)="x3",
  "x2"+(0,-.5)="x4",
  "x1";"x2"**\dir{}?(.5)="m", 
  "m" + (0,-.25)="c", 
  "x1";"x4"**\crv{ "x1"+(0,-.25) & "c" & "x4"+(0,.25)},  
  "x2";"c"+<3pt,2pt>="c2"**\crv{"x2"+(0,-.25) & "c2"},
  "x3";"c"+<-3pt,-2pt>="c1"**\crv{"x3"+(0,.25) & "c1"}, 
  "c"+(0,.10)*!<0pt,-6pt>{\si},
"x3"="L32",
"x4"="L33",
"L24";"L24" + (0,-.5) = "L34" **\dir{-},  
"L25";"L25" + (0,-.5) = "L35" **\dir{-},  
"L26";"L26" + (0,-.5) = "L36" **\dir{-},  
"L31"="L11",
"L32"="L12",
"L33"="L13",
"L34"="L14",
"L35"="L15",
"L36"="L16",
  "L11"="x1",       
  "L12"="x2",  , 
  "x1"+(0,-.5)="x3",
  "x2"+(0,-.5)="x4",
  "x1";"x2"**\dir{}?(.5)="m", 
  "m" + (0,-.25)="c", 
  "x1";"x4"**\crv{ "x1"+(0,-.25) & "c" & "x4"+(0,.25)},  
  "x2";"c"+<3pt,2pt>="c2"**\crv{"x2"+(0,-.25) & "c2"},
  "x3";"c"+<-3pt,-2pt>="c1"**\crv{"x3"+(0,.25) & "c1"}, 
  "c"+(0,.10)*!<0pt,-6pt>{\si}, 
"x3"="L21",
"x4"="L22",
  "L21"="x1", 
  "x1"+(0,-.5)="x2",
  "x1";"x2"**\dir{}?(.5)="c", 
  "x1";"x2"**\dir{-},  
  "c"-(0,0)*!<-4pt,0pt>{\bullet \widetilde{S}  }, 
"x2"="L21",
  "L22"="x1", 
  "x1"+(0,-.5)="x2",
  "x1";"x2"**\dir{}?(.5)="c", 
  "x1";"x2"**\dir{-},  
  "c"-(0,0)*!<-4pt,0pt>{\bullet \widetilde{S}  }, 
"x2"="L22",
  "L21"="x1", 
  "L22"="x2", 
  "x1";"x2"**\dir{}?(.5)="m", 
  "m" + (0,-.25)="c", 
  "m" + (0,-.5)="x3", 
  "x1";"x3"**\crv{"x1"-(0,.25) & "c" & "x3"},
  "x2";"x3"**\crv{"x2"-(0,.25) & "c"  & "x3"}, 
"x3"="L31",
  "L13"="x1",       
  "L14"="x2",  , 
  "x1"+(0,-.5)="x3",
  "x2"+(0,-.5)="x4",
  "x1";"x2"**\dir{}?(.5)="m", 
  "m" + (0,-.25)="c", 
  "x1";"x4"**\crv{ "x1"+(0,-.25) & "c" & "x4"+(0,.25)},  
  "x2";"c"+<3pt,2pt>="c2"**\crv{"x2"+(0,-.25) & "c2"},
  "x3";"c"+<-3pt,-2pt>="c1"**\crv{"x3"+(0,.25) & "c1"}, 
  "c"+(0,.10)*!<0pt,-6pt>{\si}, 
"x3"="L21",
"x4"="L22",
  "L21"="x1", 
  "x1"+(0,-.5)="x2",
  "x1";"x2"**\dir{}?(.5)="c", 
  "x1";"x2"**\dir{-},  
  "c"-(0,0)*!<-4pt,0pt>{\bullet S  }, 
"x2"="L21",
  "L22"="x1", 
  "x1"+(0,-.5)="x2",
  "x1";"x2"**\dir{}?(.5)="c", 
  "x1";"x2"**\dir{-},  
  "c"-(0,0)*!<-4pt,0pt>{\bullet S  }, 
"x2"="L22",
  "L21"="x1", 
  "L22"="x2",  
  "x1";"x2"**\dir{}?(.5)="m", 
  "m" + (0,-.25)="c", 
  "m" + (0,-.5)="x3", 
  "x1";"x3"**\crv{"x1"-(0,.25) & "c" & "x3"},
  "x2";"x3"**\crv{"x2"-(0,.25) & "c"  & "x3"}, 
"x3"="L32",
  "L15"="x1", 
  "L16"="x2",  
  "x1";"x2"**\dir{}?(.5)="m", 
  "m" + (0,-.25)="c", 
  "m" + (0,-.5)="x3", 
  "x1";"x3"**\crv{"x1"-(0,.25) & "c" & "x3"},
  "x2";"x3"**\crv{"x2"-(0,.25) & "c"  & "x3"}, 
"x3"="L33",
"L33";"L33" + (0,-1) = "L33" **\dir{-},  
  "x3"+ (1,-.25)="x1", 
  "x1"+(0,-.75)="x2",
  "x1";"x2"**\dir{}?(.5)="c", 
  "x1";"x2"**\dir{-},  
  "c"-(0,0)*!<-3pt,0pt>{\bullet \sigma },  
  "x1"+(0,0)*!<0pt,-2pt>{\circ}, 
"x2"="L34",
"L31"="L21",
"L32"="L22",
"L33"="L23",
"L34"="L24",
   "L21"="x1",    
     "L22"="x2",  
  "x1"+(0,-.5)="x3",
  "x2"+(0,-.5)="x4",
  "x1";"x2"**\dir{}?(.5)="m", 
  "m" + (0,-.25)="c", 
  "x1";"x4"**\crv{ "x1"+(0,-.25) & "c" & "x4"+(0,.25)},  
  "x2";"c"+<3pt,2pt>="c2"**\crv{"x2"+(0,-.25) & "c2"},
  "x3";"c"+<-3pt,-2pt>="c1"**\crv{"x3"+(0,.25) & "c1"}, 
  "c"+(0,.10)*!<0pt,-6pt>{\si},
"x3"="L31",
"x4"="L32",
"L23";"L23" + (0,-.5) = "L33" **\dir{-},  
"L24";"L24" + (0,-.5) = "L34" **\dir{-},  
"L31"="L21",
"L32"="L22",
"L33"="L23",
"L34"="L24",
"L21";"L21" + (0,-.5) = "L31" **\dir{-},  
  "L22"="x1",       
  "L23"="x2",  , 
  "x1"+(0,-.5)="x3",
  "x2"+(0,-.5)="x4",
  "x1";"x2"**\dir{}?(.5)="m", 
  "m" + (0,-.25)="c", 
  "x1";"x4"**\crv{ "x1"+(0,-.25) & "c" & "x4"+(0,.25)},  
  "x2";"c"+<3pt,2pt>="c2"**\crv{"x2"+(0,-.25) & "c2"},
  "x3";"c"+<-3pt,-2pt>="c1"**\crv{"x3"+(0,.25) & "c1"}, 
  "c"+(0,.10)*!<0pt,-6pt>{\si},
"x3"= "L32",
"x4"= "L33",
"L24";"L24" + (0,-.5) = "L34" **\dir{-},  
"L31"="L21",
"L32"="L22",
"L33"="L23",
"L34"="L24",
  "L21"="x1", 
  "L22"="x2",  
  "x1";"x2"**\dir{}?(.5)="m", 
  "m" + (0,-.25)="c", 
  "m" + (0,-.5)="x3", 
  "x1";"x3"**\crv{"x1"-(0,.25) & "c" & "x3"},
  "x2";"x3"**\crv{"x2"-(0,.25) & "c"  & "x3"}, 
"x3"="L31",
  "L23"="x1", 
  "L24"="x2", 
  "x1";"x2"**\dir{}?(.5)="m", 
  "m" + (0,-.25)="c", 
  "m" + (0,-.5)="x3", 
  "x1";"x3"**\crv{"x1"-(0,.25) & "c" & "x3"},
  "x2";"x3"**\crv{"x2"-(0,.25) & "c"  & "x3"}, 
"x3"="L32",
"L31"="L21",
"L32"="L22",
"L21"="L11",
"L22"="L12",
   "L11"="x1", 
  "x1"+(-.5,-.5)="x2",
  "x1"+(.5,-.5)="x3",
  "x1"+(0,-.25)="c", 
   "x1";"x2"**\crv{"x1" & "c" & "x2"+(0,.25)}, 
   "x1";"x3"**\crv{"x1" & "c" & "x3"+(0,.25)}, 
"x2"="L21",
"x3"="L22",
 "L12"="x1", 
  "x1"+(0,-.5)="x2",
  "x1";"x2"**\dir{-},   
"x2"="L23",
  "L12"+ (1,-.25)="x1", 
  "x1"+(0,-.25)="x2",
  "x1";"x2"**\dir{}?(.5)="c", 
  "x1";"x2"**\dir{-},  
  "c"-(0,0)*!<-3pt,0pt>{\bullet \sigma }, 
  "x1"+(0,0)*!<0pt,-2pt>{\circ}, 
"x2"="L24",
  "L21"="x1", 
  "x1"+(0,-.5)="x2",
  "x1";"x2"**\dir{}?(.5)="c", 
  "x1";"x2"**\dir{-},  
  "c"-(0,0)*!<4pt,0pt>{\widetilde{S} \bullet},  
"x2"="L31",
  "L22"="x1", 
  "x1"+(0,-.5)="x2",
  "x1";"x2"**\dir{}?(.5)="c", 
  "x1";"x2"**\dir{-}, 
  "c"-(0,0)*!<-4pt,0pt>{\bullet S},  
"x2"="L32",
  "L23"="x1", 
  "x1"+(0,-.5)="x2",
  "x1";"x2"**\dir{-},  
"x2"="L33",
  "L24"="x1", 
  "x1"+(0,-.5)="x2",
  "x1";"x2"**\dir{-}, 
"x2"="L34",
"L33";"L33" + (0,-.5) = "L33" **\dir{-},  
"L34";"L34" + (0,-.5) = "L34" **\dir{-},  
"L31"="L21",
"L32"="L22",
"L33"="L23",
"L34"="L24",
   "L21"="x1",  
     "L22"="x2", 
  "x1"+(0,-.5)="x3",
  "x2"+(0,-.5)="x4",
  "x1";"x2"**\dir{}?(.5)="m", 
  "m" + (0,-.25)="c", 
  "x1";"x4"**\crv{ "x1"+(0,-.25) & "c" & "x4"+(0,.25)},  
  "x2";"c"+<3pt,2pt>="c2"**\crv{"x2"+(0,-.25) & "c2"},
  "x3";"c"+<-3pt,-2pt>="c1"**\crv{"x3"+(0,.25) & "c1"}, 
  "c"+(0,.10)*!<0pt,-6pt>{\si},
"x3"="L31",
"x4"="L32",
"L23"+(0,-0)="L33",
"L24"+(0,-0)="L34",
"L31"="L21",
"L32"="L22",
"L33"="L23",
"L34"="L24",
"L21";"L21" + (0,-.5) = "L31" **\dir{-},  
  "L22"="x1",       
  "L23"="x2",  , 
  "x1"+(0,-.5)="x3",
  "x2"+(0,-.5)="x4",
  "x1";"x2"**\dir{}?(.5)="m", 
  "m" + (0,-.25)="c", 
  "x1";"x4"**\crv{ "x1"+(0,-.25) & "c" & "x4"+(0,.25)},  
  "x2";"c"+<3pt,2pt>="c2"**\crv{"x2"+(0,-.25) & "c2"},
  "x3";"c"+<-3pt,-2pt>="c1"**\crv{"x3"+(0,.25) & "c1"}, 
  "c"+(0,.10)*!<0pt,-6pt>{\si},
"x3"= "L32",
"x4"= "L33",
"L24";"L24" + (0,-.5) = "L34" **\dir{-},  
"L31"="L21",
"L32"="L22",
"L33"="L23",
"L34"="L24",
  "L21"="x1", 
  "L22"="x2",  
  "x1";"x2"**\dir{}?(.5)="m", 
  "m" + (0,-.25)="c", 
  "m" + (0,-.5)="x3", 
  "x1";"x3"**\crv{"x1"-(0,.25) & "c" & "x3"},
  "x2";"x3"**\crv{"x2"-(0,.25) & "c"  & "x3"}, 
  "x3"+(0,-.5)*!<0pt,-6pt>{H}, 
"x3"="L31",
  "L23"="x1", 
  "L24"="x2",  
  "x1";"x2"**\dir{}?(.5)="m", 
  "m" + (0,-.25)="c", 
  "m" + (0,-.5)="x3", 
  "x1";"x3"**\crv{"x1"-(0,.25) & "c" & "x3"},
  "x2";"x3"**\crv{"x2"-(0,.25) & "c"  & "x3"}, 
  "x3"+(0,-.5)*!<0pt,-6pt>{H}, 
"x3"="L32",
"L31"="L21",
"L32"="L22",
\endxy
\]

Again, by the same method, using standard identities and properties, including:
\[(1,\si)(\si,1)(1,\si) =(\si,1)(1,\si)(\si,1), \]
one reduces the above digram to the following one:

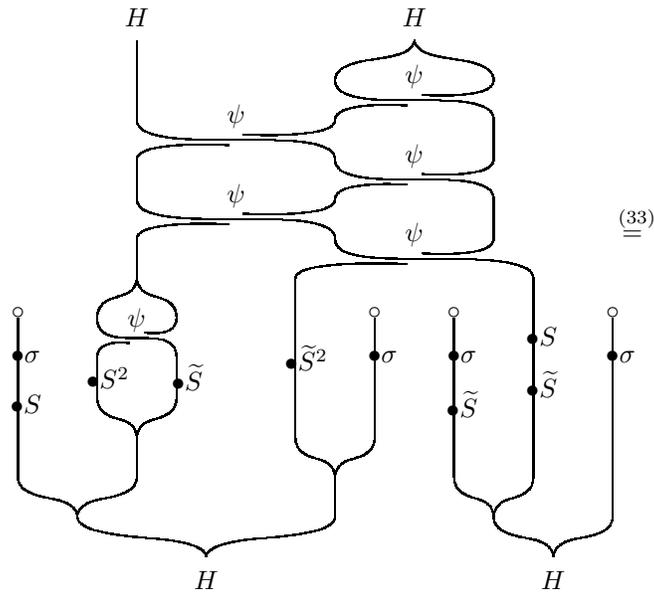
\[
\xy /r2.5pc/:,  
(0,0)="L11",
(3.5,0)="L12",
  "L11"+(0,.5)*!<0pt,6pt>{H},
  "L12"+(0,.5)*!<0pt,6pt>{H}, 
"L11";"L11" + (0,-.5) = "L11" **\dir{-},  
   "L12"="x1", 
  "x1"+(-1,-.5)="x2",
  "x1"+(1,-.5)="x3",
  "x1"+(0,-.25)="c",  
   "x1";"x2"**\crv{"x1" & "c" & "x2"+(0,.25)}, 
   "x1";"x3"**\crv{"x1" & "c" & "x3"+(0,.25)}, 
"x2"= "L12",
"x3"= "L13",
"L11";"L11" + (0,-.5) = "L11" **\dir{-},  
    "L12"="x1",    
   "L13"="x2", 
  "x1"+(0,-.5)="x3",
  "x2"+(0,-.5)="x4",
  "x1";"x2"**\dir{}?(.5)="m", 
  "m" + (0,-.25)="c", 
  "x1";"x4"**\crv{ "x1"+(0,-.25) & "c" & "x4"+(0,.25)},  
  "x2";"c"+<3pt,2pt>="c2"**\crv{"x2"+(0,-.25) & "c2"},
  "x3";"c"+<-3pt,-2pt>="c1"**\crv{"x3"+(0,.25) & "c1"}, 
  "c"+(0,.10)*!<0pt,-6pt>{\si}, 
"x3" = "L12",
"x4" = "L13",
    "L11"="x1",   
   "L12"="x2", 
  "x1"+(0,-.5)="x3",
  "x2"+(0,-.5)="x4",
  "x1";"x2"**\dir{}?(.5)="m", 
  "m" + (0,-.25)="c", 
  "x1";"x4"**\crv{ "x1"+(0,-.25) & "c" & "x4"+(0,.25)},  
  "x2";"c"+<3pt,2pt>="c2"**\crv{"x2"+(0,-.25) & "c2"},
  "x3";"c"+<-3pt,-2pt>="c1"**\crv{"x3"+(0,.25) & "c1"}, 
  "c"+(0,.10)*!<0pt,-6pt>{\si}, 
"x3" = "L11",
"x4" = "L12",
"L13";"L13" + (0,-.5) = "L13" **\dir{-},  
"L11";"L11" + (0,-.5) = "L11" **\dir{-},  
    "L12"="x1",     
   "L13"="x2",  
  "x1"+(0,-.5)="x3",
  "x2"+(0,-.5)="x4",
  "x1";"x2"**\dir{}?(.5)="m", 
  "m" + (0,-.25)="c", 
  "x1";"x4"**\crv{ "x1"+(0,-.25) & "c" & "x4"+(0,.25)},  
  "x2";"c"+<3pt,2pt>="c2"**\crv{"x2"+(0,-.25) & "c2"},
  "x3";"c"+<-3pt,-2pt>="c1"**\crv{"x3"+(0,.25) & "c1"}, 
  "c"+(0,.10)*!<0pt,-6pt>{\si}, 
"x3" = "L12",
"x4" = "L13",
    "L11"="x1",    
   "L12"="x2", 
  "x1"+(0,-.5)="x3",
  "x2"+(0,-.5)="x4",
  "x1";"x2"**\dir{}?(.5)="m", 
  "m" + (0,-.25)="c", 
  "x1";"x4"**\crv{ "x1"+(0,-.25) & "c" & "x4"+(0,.25)},  
  "x2";"c"+<3pt,2pt>="c2"**\crv{"x2"+(0,-.25) & "c2"},
  "x3";"c"+<-3pt,-2pt>="c1"**\crv{"x3"+(0,.25) & "c1"}, 
  "c"+(0,.10)*!<0pt,-6pt>{\si}, 
"x3" = "L11",
"x4" = "L12",
"L13";"L13" + (0,-.5) = "L13" **\dir{-},  
"L11";"L11" + (0,-.5) = "L11" **\dir{-},  
    "L12"="x1",   
   "L13"="x2", 
  "x1"+(-.5,-.5)="x3",
  "x2"+(.5,-.5)="x4",
  "x1";"x2"**\dir{}?(.5)="m", 
  "m" + (0,-.25)="c", 
  "x1";"x4"**\crv{ "x1"+(0,-.25) & "c" & "x4"+(0,.25)},  
  "x2";"c"+<3pt,2pt>="c2"**\crv{"x2"+(0,-.25) & "c2"},
  "x3";"c"+<-3pt,-2pt>="c1"**\crv{"x3"+(0,.25) & "c1"}, 
  "c"+(0,.10)*!<0pt,-6pt>{\si}, 
"x3" = "L12",
"x4" = "L13",
"L11" + (-1.5,0)="L21",
"L11" ="L22",
"L12" ="L23",
"L12" + (.5,0)="L24",
"L12" + (1,0)="L25",
"L13" ="L26",
"L13" + (1,0)="L27",
   "L21"  + (0,-.5) ="x1", 
  "x1"+(0,-1.5)="x2",
  "x1";"x2"**\dir{}?(.5)="c", 
  "x1";"x2"**\dir{-},  
  "c"-(0,0)*!<-3pt,-8pt>{\bullet \sigma }, 
  "c"+(0,0)*!<-3pt,10pt>{\bullet S }, 
  "x1"+(0,0)*!<0pt,-2pt>{\circ}, 
"x2" ="L11",
   "L22"="x1", 
  "x1"+(-.5,-.5)="x2",
  "x1"+(.5,-.5)="x3",
  "x1"+(0,-.25)="c", 
   "x1";"x2"**\crv{"x1" & "c" & "x2"+(0,.25)}, 
   "x1";"x3"**\crv{"x1" & "c" & "x3"+(0,.25)}, 
"x2"= "L12",
"x3"= "L13",
    "L12"="x1",    
   "L13"="x2", 
  "x1"+(0,-.5)="x3",
  "x2"+(0,-.5)="x4",
  "x1";"x2"**\dir{}?(.5)="m", 
  "m" + (0,-.25)="c", 
  "x1";"x4"**\crv{ "x1"+(0,-.25) & "c" & "x4"+(0,.25)},  
  "x2";"c"+<3pt,2pt>="c2"**\crv{"x2"+(0,-.25) & "c2"},
  "x3";"c"+<-3pt,-2pt>="c1"**\crv{"x3"+(0,.25) & "c1"}, 
  "c"+(0,.10)*!<0pt,-4pt>{\si}, 
"x3" = "L12",
"x4" = "L13",
  "L12"="x1", 
  "x1"+(0,-.5)="x2",
  "x1";"x2"**\dir{}?(.5)="c", 
  "x1";"x2"**\dir{-}, 
  "c"-(0,0)*!<-4pt,0pt>{\bullet S^2 },  
"x2"="L12",
  "L13"="x1", 
  "x1"+(0,-.5)="x2",
  "x1";"x2"**\dir{}?(.5)="c", 
  "x1";"x2"**\dir{-},   
  "c"-(0,0)*!<-4pt,0pt>{\bullet \widetilde{S} },  
"x2"="L13",
  "L12"="x1", 
  "L13"="x2",  
  "x1";"x2"**\dir{}?(.5)="m", 
  "m" + (0,-.25)="c", 
  "m" + (0,-.5)="x3", 
  "x1";"x3"**\crv{"x1"-(0,.25) & "c" & "x3"},
  "x2";"x3"**\crv{"x2"-(0,.25) & "c"  & "x3"}, 
"x3"="L12",
  "L23"="x1", 
  "x1"+(0,-2)="x2",
  "x1";"x2"**\dir{}?(.5)="c", 
  "x1";"x2"**\dir{-},   
  "c"-(0,0)*!<-4pt,0pt>{\bullet \widetilde{S}^2 },  
"x2"="L13",
   "L24" + (.5,-.5)  ="x1", 
  "x1"+(0,-1.5)="x2",
  "x1";"x2"**\dir{}?(.5)="c", 
  "x1";"x2"**\dir{-},  
  "c"-(0,0)*!<-3pt,-8pt>{\bullet \sigma }, 
  "x1"+(0,0)*!<0pt,-2pt>{\circ}, 
"x2" ="L14",
   "L25" + (1,-.5)  ="x1", 
  "x1"+(0,-1.5)="x2",
  "x1";"x2"**\dir{}?(.5)="c", 
  "x1";"x2"**\dir{-}, 
  "c"-(0,0)*!<-3pt,-8pt>{\bullet \sigma }, 
  "c"+(0,0)*!<-3pt,10pt>{\bullet \widetilde{S} }, 
  "x1"+(0,0)*!<0pt,-2pt>{\circ}, 
"x2" ="L15",
   "L26" ="x1", 
  "x1"+(0,-2)="x2",
  "x1";"x2"**\dir{}?(.5)="c", 
  "x1";"x2"**\dir{-},  
  "c"-(0,0)*!<-3pt,-8pt>{\bullet S }, 
  "c"+(0,0)*!<-3pt,10pt>{\bullet \widetilde{S} }, 
"x2" ="L16",
   "L27" + (0,-.5) ="x1", 
  "x1"+(0,-1.5)="x2",
  "x1";"x2"**\dir{}?(.5)="c", 
  "x1";"x2"**\dir{-},  
  "c"-(0,0)*!<-3pt,-8pt>{\bullet \sigma }, 
  "x1"+(0,0)*!<0pt,-2pt>{\circ}, 
"x2" ="L17",
"L11";"L11" + (0,-.5) = "L11" **\dir{-},  
"L12";"L12" + (0,-.5) = "L12" **\dir{-},  
  "L13"="x1",
  "L14"="x2",   
  "x1";"x2"**\dir{}?(.5)="m", 
  "m" + (0,-.25)="c", 
  "m" + (0,-.5)="x3", 
  "x1";"x3"**\crv{"x1"-(0,.25) & "c" & "x3"},
  "x2";"x3"**\crv{"x2"-(0,.25) & "c"  & "x3"}, 
"x3"="L13",
"L15";"L15" + (0,-.5) = "L14" **\dir{-},  
"L16";"L16" + (0,-.5) = "L15" **\dir{-},  
"L17";"L17" + (0,-.5) = "L16" **\dir{-},  
  "L11"="x1", 
  "L12"="x2", 
  "x1";"x2"**\dir{}?(.5)="m", 
  "m" + (0,-.25)="c", 
  "m" + (0,-.5)="x3", 
  "x1";"x3"**\crv{"x1"-(0,.25) & "c" & "x3"},
  "x2";"x3"**\crv{"x2"-(0,.25) & "c"  & "x3"}, 
"x3"="L21",
"L13";"L13" + (0,-.5) = "L22" **\dir{-},  
  "L21"="x1",
  "L22"="x2", 
  "x1";"x2"**\dir{}?(.5)="m", 
  "m" + (0,-.25)="c", 
  "m" + (0,-.5)="x3", 
  "x1";"x3"**\crv{"x1"-(0,.25) & "c" & "x3"},
  "x2";"x3"**\crv{"x2"-(0,.25) & "c"  & "x3"}, 
"x3"="L11",
  "L14"="x1", 
  "L15"="x2",   
  "x1";"x2"**\dir{}?(.5)="m", 
  "m" + (0,-.25)="c", 
  "m" + (0,-.5)="x3", 
  "x1";"x3"**\crv{"x1"-(0,.25) & "c" & "x3"},
  "x2";"x3"**\crv{"x2"-(0,.25) & "c"  & "x3"}, 
"x3"="L21",
"L16";"L16" + (0,-.5) = "L22" **\dir{-},  
  "L21"="x1", 
  "L22"="x2",  
  "x1";"x2"**\dir{}?(.5)="m", 
  "m" + (0,-.25)="c", 
  "m" + (0,-.5)="x3", 
  "x1";"x3"**\crv{"x1"-(0,.25) & "c" & "x3"},
  "x2";"x3"**\crv{"x2"-(0,.25) & "c"  & "x3"}, 
"x3"="L12",
"L11"+(0,-.5)*!<0pt,-6pt>{H}, 
"L12"+(0,-.5)*!<0pt,-6pt>{H}, 
 "x3" + (1,4)*!<0pt,-6pt>{\overset{~\eqref{m(S2,ws)si Del}}{=}}, 
\endxy
\]

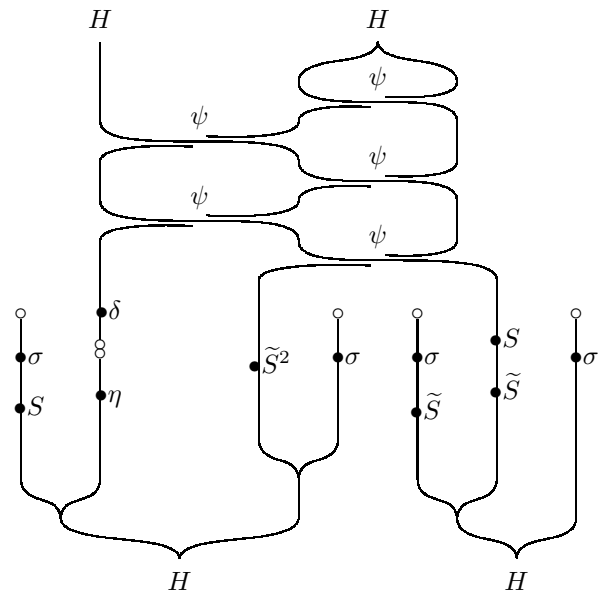
\[
\xy /r2.5pc/:,  
(0,0)="L11",
(3.5,0)="L12",
  "L11"+(0,.5)*!<0pt,6pt>{H},
  "L12"+(0,.5)*!<0pt,6pt>{H}, 
"L11";"L11" + (0,-.5) = "L11" **\dir{-},  
   "L12"="x1",
  "x1"+(-1,-.5)="x2",
  "x1"+(1,-.5)="x3",
  "x1"+(0,-.25)="c", 
   "x1";"x2"**\crv{"x1" & "c" & "x2"+(0,.25)}, 
   "x1";"x3"**\crv{"x1" & "c" & "x3"+(0,.25)}, 
"x2"= "L12",
"x3"= "L13",
"L11";"L11" + (0,-.5) = "L11" **\dir{-},  
    "L12"="x1",    
   "L13"="x2",  
  "x1"+(0,-.5)="x3",
  "x2"+(0,-.5)="x4",
  "x1";"x2"**\dir{}?(.5)="m", 
  "m" + (0,-.25)="c", 
  "x1";"x4"**\crv{ "x1"+(0,-.25) & "c" & "x4"+(0,.25)},  
  "x2";"c"+<3pt,2pt>="c2"**\crv{"x2"+(0,-.25) & "c2"},
  "x3";"c"+<-3pt,-2pt>="c1"**\crv{"x3"+(0,.25) & "c1"}, 
  "c"+(0,.10)*!<0pt,-6pt>{\si}, 
"x3" = "L12",
"x4" = "L13",
    "L11"="x1",   
   "L12"="x2", 
  "x1"+(0,-.5)="x3",
  "x2"+(0,-.5)="x4",
  "x1";"x2"**\dir{}?(.5)="m", 
  "m" + (0,-.25)="c", 
  "x1";"x4"**\crv{ "x1"+(0,-.25) & "c" & "x4"+(0,.25)},  
  "x2";"c"+<3pt,2pt>="c2"**\crv{"x2"+(0,-.25) & "c2"},
  "x3";"c"+<-3pt,-2pt>="c1"**\crv{"x3"+(0,.25) & "c1"}, 
  "c"+(0,.10)*!<0pt,-6pt>{\si}, 
"x3" = "L11",
"x4" = "L12",
"L13";"L13" + (0,-.5) = "L13" **\dir{-},  
"L11";"L11" + (0,-.5) = "L11" **\dir{-},  
    "L12"="x1",    
   "L13"="x2",  
  "x1"+(0,-.5)="x3",
  "x2"+(0,-.5)="x4",
  "x1";"x2"**\dir{}?(.5)="m", 
  "m" + (0,-.25)="c", 
  "x1";"x4"**\crv{ "x1"+(0,-.25) & "c" & "x4"+(0,.25)},  
  "x2";"c"+<3pt,2pt>="c2"**\crv{"x2"+(0,-.25) & "c2"},
  "x3";"c"+<-3pt,-2pt>="c1"**\crv{"x3"+(0,.25) & "c1"}, 
  "c"+(0,.10)*!<0pt,-6pt>{\si}, 
"x3" = "L12",
"x4" = "L13",
    "L11"="x1",    
   "L12"="x2",  
  "x1"+(0,-.5)="x3",
  "x2"+(0,-.5)="x4",
  "x1";"x2"**\dir{}?(.5)="m", 
  "m" + (0,-.25)="c", 
  "x1";"x4"**\crv{ "x1"+(0,-.25) & "c" & "x4"+(0,.25)},  
  "x2";"c"+<3pt,2pt>="c2"**\crv{"x2"+(0,-.25) & "c2"},
  "x3";"c"+<-3pt,-2pt>="c1"**\crv{"x3"+(0,.25) & "c1"}, 
  "c"+(0,.10)*!<0pt,-6pt>{\si}, 
"x3" = "L11",
"x4" = "L12",
"L13";"L13" + (0,-.5) = "L13" **\dir{-},  
"L11";"L11" + (0,-.5) = "L11" **\dir{-},  
    "L12"="x1",    
   "L13"="x2",  
  "x1"+(-.5,-.5)="x3",
  "x2"+(.5,-.5)="x4",
  "x1";"x2"**\dir{}?(.5)="m", 
  "m" + (0,-.25)="c", 
  "x1";"x4"**\crv{ "x1"+(0,-.25) & "c" & "x4"+(0,.25)},  
  "x2";"c"+<3pt,2pt>="c2"**\crv{"x2"+(0,-.25) & "c2"},
  "x3";"c"+<-3pt,-2pt>="c1"**\crv{"x3"+(0,.25) & "c1"}, 
  "c"+(0,.10)*!<0pt,-6pt>{\si}, 
"x3" = "L12",
"x4" = "L13",
"L11" + (-1,0)="L21",
"L11" ="L22",
"L12" ="L23",
"L12" + (.5,0)="L24",
"L12" + (1,0)="L25",
"L13" ="L26",
"L13" + (1,0)="L27",
   "L21"  + (0,-.5) ="x1", 
  "x1"+(0,-1.5)="x2",
  "x1";"x2"**\dir{}?(.5)="c", 
  "x1";"x2"**\dir{-},  
  "c"-(0,0)*!<-3pt,-8pt>{\bullet \sigma }, 
  "c"+(0,0)*!<-3pt,10pt>{\bullet S}, 
  "x1"+(0,0)*!<0pt,-2pt>{\circ}, 
"x2" ="L11",
    "L22"="x1", 
   "x1"+(0,-.75)="x2",
  "x1";"x2"**\dir{}?(.5)="c", 
  "x1";"x2"**\dir{-},   
  "c"-(0,0)*!<-3pt,0pt>{\bullet \delta }, 
  "x2"+(0,0)*!<0pt,2pt>{\circ}, 
"x2"="L31",
   "L31" + (0,-.25)="x1", 
  "x1"+(0,-1)="x2",
  "x1";"x2"**\dir{}?(.5)="c", 
  "x1";"x2"**\dir{-},  
  "c"-(0,0)*!<-3pt,0pt>{\bullet \eta },  
  "x1"+(0,0)*!<0pt,-2pt>{\circ},  
"x2"="L12",
  "L23"="x1", 
  "x1"+(0,-2)="x2",
  "x1";"x2"**\dir{}?(.5)="c", 
  "x1";"x2"**\dir{-},  
  "c"-(0,0)*!<-4pt,0pt>{\bullet \widetilde{S}^2 },  
"x2"="L13",
   "L24" + (.5,-.5)  ="x1", 
  "x1"+(0,-1.5)="x2",
  "x1";"x2"**\dir{}?(.5)="c", 
  "x1";"x2"**\dir{-},   
  "c"-(0,0)*!<-3pt,-8pt>{\bullet \sigma }, 
  "x1"+(0,0)*!<0pt,-2pt>{\circ}, 
"x2" ="L14",
   "L25" + (1,-.5)  ="x1", 
  "x1"+(0,-1.5)="x2",
  "x1";"x2"**\dir{}?(.5)="c", 
  "x1";"x2"**\dir{-},  
  "c"-(0,0)*!<-3pt,-8pt>{\bullet \sigma }, 
  "c"+(0,0)*!<-3pt,10pt>{\bullet \widetilde{S} }, 
  "x1"+(0,0)*!<0pt,-2pt>{\circ}, 
"x2" ="L15",
   "L26" ="x1", 
  "x1"+(0,-2)="x2",
  "x1";"x2"**\dir{}?(.5)="c", 
  "x1";"x2"**\dir{-},   
  "c"-(0,0)*!<-3pt,-8pt>{\bullet S }, 
  "c"+(0,0)*!<-3pt,10pt>{\bullet \widetilde{S} }, 
"x2" ="L16",
   "L27" + (0,-.5) ="x1", 
  "x1"+(0,-1.5)="x2",
  "x1";"x2"**\dir{}?(.5)="c", 
  "x1";"x2"**\dir{-},   
  "c"-(0,0)*!<-3pt,-8pt>{\bullet \sigma }, 
  "x1"+(0,0)*!<0pt,-2pt>{\circ}, 
"x2" ="L17",
"L11";"L11" + (0,-.5) = "L11" **\dir{-},  
"L12";"L12" + (0,-.5) = "L12" **\dir{-},  
  "L13"="x1", 
  "L14"="x2",  
  "x1";"x2"**\dir{}?(.5)="m", 
  "m" + (0,-.25)="c", 
  "m" + (0,-.5)="x3", 
  "x1";"x3"**\crv{"x1"-(0,.25) & "c" & "x3"},
  "x2";"x3"**\crv{"x2"-(0,.25) & "c"  & "x3"}, 
"x3"="L13",
"L15";"L15" + (0,-.5) = "L14" **\dir{-},  
"L16";"L16" + (0,-.5) = "L15" **\dir{-},  
"L17";"L17" + (0,-.5) = "L16" **\dir{-},  
  "L11"="x1",
  "L12"="x2",  
  "x1";"x2"**\dir{}?(.5)="m", 
  "m" + (0,-.25)="c", 
  "m" + (0,-.5)="x3", 
  "x1";"x3"**\crv{"x1"-(0,.25) & "c" & "x3"},
  "x2";"x3"**\crv{"x2"-(0,.25) & "c"  & "x3"}, 
"x3"="L21",
"L13";"L13" + (0,-.5) = "L22" **\dir{-},  
  "L21"="x1",
  "L22"="x2",  
  "x1";"x2"**\dir{}?(.5)="m", 
  "m" + (0,-.25)="c", 
  "m" + (0,-.5)="x3", 
  "x1";"x3"**\crv{"x1"-(0,.25) & "c" & "x3"},
  "x2";"x3"**\crv{"x2"-(0,.25) & "c"  & "x3"}, 
"x3"="L11",
  "L14"="x1", 
  "L15"="x2",  
  "x1";"x2"**\dir{}?(.5)="m", 
  "m" + (0,-.25)="c", 
  "m" + (0,-.5)="x3", 
  "x1";"x3"**\crv{"x1"-(0,.25) & "c" & "x3"},
  "x2";"x3"**\crv{"x2"-(0,.25) & "c"  & "x3"}, 
"x3"="L21",
"L16";"L16" + (0,-.5) = "L22" **\dir{-},  
  "L21"="x1", 
  "L22"="x2", 
  "x1";"x2"**\dir{}?(.5)="m", 
  "m" + (0,-.25)="c", 
  "m" + (0,-.5)="x3", 
  "x1";"x3"**\crv{"x1"-(0,.25) & "c" & "x3"},
  "x2";"x3"**\crv{"x2"-(0,.25) & "c"  & "x3"}, 
"x3"="L12",
"L11"+(0,-.5)*!<0pt,-6pt>{H}, 
"L12"+(0,-.5)*!<0pt,-6pt>{H}, 
\endxy
\]

\begin{align*} 
=&(m(m,1),m(m,1))(2,m,3)(S \sigma,\eta\delta,\widetilde{S}^2,\sigma,\widetilde{S}\sigma,\widetilde{S} S,\sigma)\\
~~~~~&(1,\si)(\si,1)(1,\si)(\si,1)(1,\si)(1,\Delta)\\
=&(m(m,1),m(m,1))(1,\eta,4)(1,m,3)(S \sigma,\widetilde{S}^2,\sigma,\widetilde{S}\sigma,\widetilde{S} S,\sigma)\\
~~~~~&(\si)(\delta,2)(\si,1)(1,\si)(\si,1)(1,\si)(1,\Delta)\\
=&(m(m(1,\eta),1),m(m,1))(1,m,3)(S \sigma,\widetilde{S}^2,\sigma,\widetilde{S}\sigma,\widetilde{S} S,\sigma)(\si)\\
~~~~~&(2,\delta)(\si,1)(1,\si)(1,\Delta)\\
=&(m,m(m,1))(1,m,3)(S \sigma,\widetilde{S}^2,\sigma,\widetilde{S}\sigma,\widetilde{S} S,\sigma)(\si)(\si)(1,\delta,1)(1,\Delta)\\
=&(m,m(m,1))(1,m,3)(S \sigma,\widetilde{S}^2,\sigma,\widetilde{S}\sigma,\widetilde{S} ,\sigma)(1,S)\si^2(1,\delta,1)(1,\Delta)\\
=&(m,m(m,1))(1,m,3)(S \sigma,\widetilde{S}^2,\sigma,\widetilde{S}\sigma,\widetilde{S} ,\sigma)\si^2(1,S)(1,\delta,1)(1,\Delta)\\
=&(m,m(m,1))(1,m,3)(S \sigma,\widetilde{S}^2,\sigma,\widetilde{S}\sigma,\widetilde{S} ,\sigma)\si^2(1,(\delta,S)\Delta)\\
=&(m,m(m,1))(1,m,3)(S \sigma,\widetilde{S}^2,\sigma,\widetilde{S}\sigma,\widetilde{S} ,\sigma)\si^2(1,\widetilde{S})\\
=&(m(1,m),m(m,1))(S \sigma,\widetilde{S}^2,\sigma,\widetilde{S}\sigma,\widetilde{S} ,\sigma)(1,\widetilde{S})\si^2\\
=&(m(1,m)(S \sigma,\widetilde{S}^2,\sigma),m(m,1)(\widetilde{S}\sigma,\widetilde{S}^2 ,\sigma))\si^2 =
\end{align*}

\[
\xy /r2.5pc/:, 
(0,0)="L11",
(2.5,0)="L12",
  "L11"+(0,.5)*!<0pt,6pt>{H},
  "L12"+(0,.5)*!<0pt,6pt>{H}, 
    "L11"="x1",     
   "L12"="x2",  
  "x1"+(0,-.5)="x3",
  "x2"+(0,-.5)="x4",
  "x1";"x2"**\dir{}?(.5)="m", 
  "m" + (0,-.25)="c", 
  "x1";"x4"**\crv{ "x1"+(0,-.25) & "c" & "x4"+(0,.25)},  
  "x2";"c"+<3pt,2pt>="c2"**\crv{"x2"+(0,-.25) & "c2"},
  "x3";"c"+<-3pt,-2pt>="c1"**\crv{"x3"+(0,.25) & "c1"}, 
  "c"+(0,.10)*!<0pt,-6pt>{\si}, 
"x3" = "L11",
"x4" = "L12",
    "L11"="x1",    
   "L12"="x2", 
  "x1"+(0,-.5)="x3",
  "x2"+(0,-.5)="x4",
  "x1";"x2"**\dir{}?(.5)="m", 
  "m" + (0,-.25)="c", 
  "x1";"x4"**\crv{ "x1"+(0,-.25) & "c" & "x4"+(0,.25)},  
  "x2";"c"+<3pt,2pt>="c2"**\crv{"x2"+(0,-.25) & "c2"},
  "x3";"c"+<-3pt,-2pt>="c1"**\crv{"x3"+(0,.25) & "c1"}, 
  "c"+(0,.10)*!<0pt,-3pt>{\si}, 
"x3" = "L11",
"x4" = "L12",
"L11" + (-1,0)= "L21",
"L11" = "L22",
"L11" + (1,0)= "L23",
   "L21"  + (0,-.5) ="x1", 
  "x1"+(0,-1.5)="x2",
  "x1";"x2"**\dir{}?(.5)="c", 
  "x1";"x2"**\dir{-},  
  "c"-(0,0)*!<-3pt,-8pt>{\bullet \sigma }, 
  "c"+(0,0)*!<-3pt,10pt>{\bullet S}, 
  "x1"+(0,0)*!<0pt,-2pt>{\circ}, 
"x2" ="L21",
   "L22"="x1", 
  "x1"+(0,-2)="x2",
  "x1";"x2"**\dir{}?(.5)="c", 
  "x1";"x2"**\dir{-},   
  "c"-(0,0)*!<-4pt,0pt>{\bullet \widetilde{S}^2 }, 
"x2" ="L22",
   "L23"  + (-.25,-.5) ="x1", 
  "x1"+(0,-1.5)="x2",
  "x1";"x2"**\dir{}?(.5)="c", 
  "x1";"x2"**\dir{-},  
  "c"-(0,0)*!<-3pt,-8pt>{\bullet \sigma }, 
  "x1"+(0,0)*!<0pt,-2pt>{\circ}, 
"x2" ="L23",
  "L21"="x1", 
  "L22"="x2",   
  "x1";"x2"**\dir{}?(.5)="m", 
  "m" + (0,-.25)="c", 
  "m" + (0,-.5)="x3", 
  "x1";"x3"**\crv{"x1"-(0,.25) & "c" & "x3"},
  "x2";"x3"**\crv{"x2"-(0,.25) & "c"  & "x3"}, 
"x3"="L21",
"L23";"L23" + (0,-.5) = "L22" **\dir{-},  
  "L21"="x1", 
  "L22"="x2",  
  "x1";"x2"**\dir{}?(.5)="m", 
  "m" + (0,-.25)="c", 
  "m" + (0,-.5)="x3", 
  "x1";"x3"**\crv{"x1"-(0,.25) & "c" & "x3"},
  "x2";"x3"**\crv{"x2"-(0,.25) & "c"  & "x3"}, 
"x3"="L21",
"L21"+(0,-.5)*!<0pt,-6pt>{H}, 
"L12" + (-1,0)= "L21",
"L12" = "L22",
"L12" + (1,0)= "L23",
   "L21"  + (0,-.5) ="x1", 
  "x1"+(0,-1.5)="x2",
  "x1";"x2"**\dir{}?(.5)="c", 
  "x1";"x2"**\dir{-},  
  "c"-(0,0)*!<-3pt,-8pt>{\bullet \sigma }, 
  "c"+(0,0)*!<-3pt,10pt>{\bullet \widetilde{S} }, 
  "x1"+(0,0)*!<0pt,-2pt>{\circ}, 
"x2" ="L21",
   "L22"="x1", 
  "x1"+(0,-2)="x2",
  "x1";"x2"**\dir{}?(.5)="c", 
  "x1";"x2"**\dir{-},  
  "c"-(0,0)*!<-4pt,0pt>{\bullet \widetilde{S}^2 }, 
"x2" ="L22",
   "L23"  + (0,-.5) ="x1", 
  "x1"+(0,-1.5)="x2",
  "x1";"x2"**\dir{}?(.5)="c", 
  "x1";"x2"**\dir{-},   
  "c"-(0,0)*!<-3pt,-8pt>{\bullet \sigma }, 
  "x1"+(0,0)*!<0pt,-2pt>{\circ}, 
"x2" ="L23",
  "L21"="x1", 
  "L22"="x2",  
  "x1";"x2"**\dir{}?(.5)="m", 
  "m" + (0,-.25)="c", 
  "m" + (0,-.5)="x3", 
  "x1";"x3"**\crv{"x1"-(0,.25) & "c" & "x3"},
  "x2";"x3"**\crv{"x2"-(0,.25) & "c"  & "x3"}, 
"x3"="L21",
"L23";"L23" + (0,-.5) = "L22" **\dir{-},  
  "L21"="x1", 
  "L22"="x2", 
  "x1";"x2"**\dir{}?(.5)="m", 
  "m" + (0,-.25)="c", 
  "m" + (0,-.5)="x3", 
  "x1";"x3"**\crv{"x1"-(0,.25) & "c" & "x3"},
  "x2";"x3"**\crv{"x2"-(0,.25) & "c"  & "x3"}, 
"x3"="L21",
"L21"+(0,-.5)*!<0pt,-6pt>{H}, 
 "x3" + (1.5,3)*!<0pt,-6pt>{\overset{~\eqref{BMPI}}{=}}, 
(5.5,0)="L11",
(8,0)="L12",
  "L11"+(0,.5)*!<0pt,6pt>{H},
  "L12"+(0,.5)*!<0pt,6pt>{H}, 
    "L11"="x1",      
   "L12"="x2",  
  "x1"+(0,-.5)="x3",
  "x2"+(0,-.5)="x4",
  "x1";"x2"**\dir{}?(.5)="m", 
  "m" + (0,-.25)="c", 
  "x1";"x4"**\crv{ "x1"+(0,-.25) & "c" & "x4"+(0,.25)},  
  "x2";"c"+<3pt,2pt>="c2"**\crv{"x2"+(0,-.25) & "c2"},
  "x3";"c"+<-3pt,-2pt>="c1"**\crv{"x3"+(0,.25) & "c1"}, 
  "c"+(0,.10)*!<0pt,-6pt>{\si}, 
"x3" = "L11",
"x4" = "L12",
    "L11"="x1",     
   "L12"="x2",  
  "x1"+(0,-.5)="x3",
  "x2"+(0,-.5)="x4",
  "x1";"x2"**\dir{}?(.5)="m", 
  "m" + (0,-.25)="c", 
  "x1";"x4"**\crv{ "x1"+(0,-.25) & "c" & "x4"+(0,.25)},  
  "x2";"c"+<3pt,2pt>="c2"**\crv{"x2"+(0,-.25) & "c2"},
  "x3";"c"+<-3pt,-2pt>="c1"**\crv{"x3"+(0,.25) & "c1"}, 
  "c"+(0,.10)*!<0pt,-3pt>{\si}, 
"x3" = "L11",
"x4" = "L12",
"L11"+(0,-.5)*!<0pt,-6pt>{H},
"L12"+(0,-.5)*!<0pt,-6pt>{H}, 
\endxy
\]

\[ =(1,1)\si^2=\si^2 \]
\end{proof}

\begin{remark} \label{gnral taun} In general we have for all $n \neq 0$ :
$$\tau_n^{n+1}= (\si_{H^{(n-1)},H})^n $$
which is equal to $id$, for all $n \neq 0$, iff $\si^2 =id$.
\end{remark}

For $n=3$ this statement means:
\[\tau_3^{4}= (\si_{H^{2},H})^3, \]
which is visualized in the following picture:

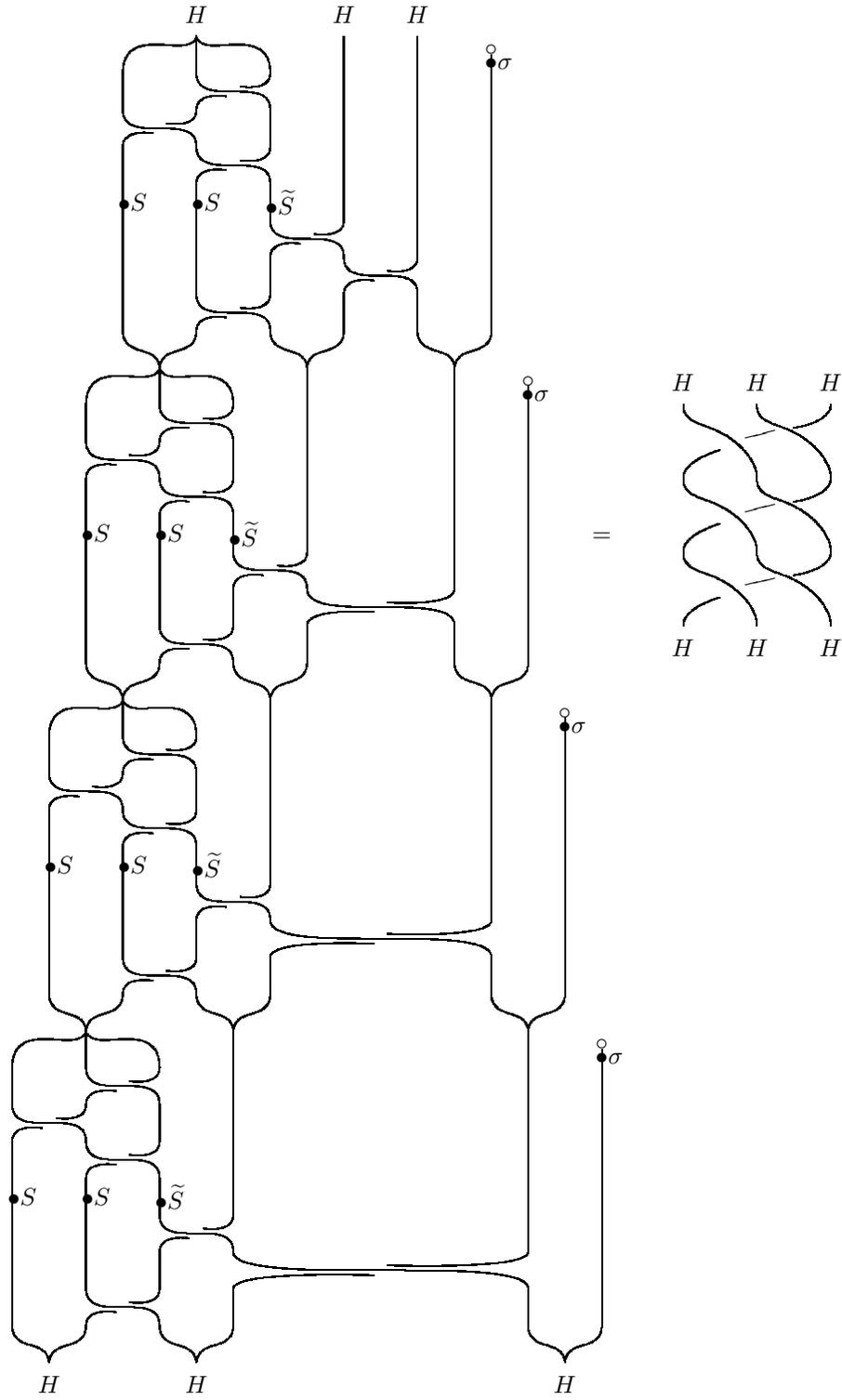
\[
\xy /r2.5pc/:,   
(0,0)="L11",
(2,0)="L12",
(3,0)="L13",
  "L11"+(0,.5)*!<0pt,6pt>{H},
  "L12"+(0,.5)*!<0pt,6pt>{H}, 
  "L13"+(0,.5)*!<0pt,6pt>{H},
   "L11"="x1", 
  "x1"+(-1,-.5)="x2",
  "x1"+(0,-.5)="x3",
  "x1"+(1,-.5)="x4",
     "x1";"x3"**\dir{}?(.5)="c", 
   "x1";"x2"**\crv{"x1" & "c" & "x2"+(0,.5)}, 
   "x1";"x3"**\crv{"x1" & "c" & "x3"+(0,.5)}, 
   "x1";"x4"**\crv{"x1" & "c" & "x4"+(0,.5)}, 
"x2"="L21",
"x3"="L22",
"x4"="L23",
 "L12"="x1", 
  "x1"+(0,-.5)="x2",
  "x1";"x2"**\dir{-},  
"x2"="L24",
 "L13"="x1", 
  "x1"+(0,-.5)="x2",
  "x1";"x2"**\dir{-},   
"x2"="L25",
  "L13"+ (1,-.25)="x1", 
  "x1"+(0,-.25)="x2",
  "x1";"x2"**\dir{}?(.5)="c", 
  "x1";"x2"**\dir{-},  
  "c"-(0,0)*!<-3pt,0pt>{\bullet \sigma },  
  "x1"+(0,0)*!<0pt,-2pt>{\circ},  
"x2"="L26",
  "L21"="x1", 
  "x1"+(0,-.5)="x2",
  "x1";"x2"**\dir{-},   
"x2"="L31",
   "L22"="x1",    
     "L23"="x2", 
  "x1"+(0,-.5)="x3",
  "x2"+(0,-.5)="x4",
  "x1";"x2"**\dir{}?(.5)="m", 
  "m" + (0,-.25)="c", 
  "x1";"x4"**\crv{ "x1"+(0,-.25) & "c" & "x4"+(0,.25)},  
  "x2";"c"+<3pt,2pt>="c2"**\crv{"x2"+(0,-.25) & "c2"},
  "x3";"c"+<-3pt,-2pt>="c1"**\crv{"x3"+(0,.25) & "c1"}, 
  "c"+(0,.10)*!<0pt,-6pt>{},
"x3"="L32",
"x4"="L33",
"L24"+(0,-0)="L34",
"L25"+(0,-0)="L35",
"L26"+(0,-0)="L36",
   "L31"="x1",     
     "L32"="x2",  
  "x1"+(0,-.5)="x3",
  "x2"+(0,-.5)="x4",
  "x1";"x2"**\dir{}?(.5)="m", 
  "m" + (0,-.25)="c", 
  "x1";"x4"**\crv{ "x1"+(0,-.25) & "c" & "x4"+(0,.25)},  
  "x2";"c"+<3pt,2pt>="c2"**\crv{"x2"+(0,-.25) & "c2"},
  "x3";"c"+<-3pt,-2pt>="c1"**\crv{"x3"+(0,.25) & "c1"}, 
  "c"+(0,.10)*!<0pt,-6pt>{},
"x3"="L41",
"x4"="L42",
  "L33"="x1", 
  "x1"+(0,-.5)="x2",
  "x1";"x2"**\dir{-},   
"x2"="L43",
  "L24"="x1", 
  "x1"+(0,-.5)="x2",
  "x1";"x2"**\dir{-},   
"x2"="L44",
  "L25"="x1", 
  "x1"+(0,-.5)="x2",
  "x1";"x2"**\dir{-},   
"x2"="L45",
  "L26"="x1", 
  "x1"+(0,-.5)="x2",
  "x1";"x2"**\dir{-},   
"x2"="L46",
  "L41"="x1", 
  "x1"+(0,-.5)="x2",
  "x1";"x2"**\dir{-}, 
"x2"="L51",
   "L42"="x1",    
     "L43"="x2",  
  "x1"+(0,-.5)="x3",
  "x2"+(0,-.5)="x4",
  "x1";"x2"**\dir{}?(.5)="m", 
  "m" + (0,-.25)="c", 
  "x1";"x4"**\crv{ "x1"+(0,-.25) & "c" & "x4"+(0,.25)},  
  "x2";"c"+<3pt,2pt>="c2"**\crv{"x2"+(0,-.25) & "c2"},
  "x3";"c"+<-3pt,-2pt>="c1"**\crv{"x3"+(0,.25) & "c1"}, 
  "c"+(0,.10)*!<0pt,-6pt>{},
"x3"="L52",
"x4"="L53",
  "L44"="x1", 
  "x1"+(0,-.5)="x2",
  "x1";"x2"**\dir{-},  
"x2"="L54",
  "L45"="x1", 
  "x1"+(0,-.5)="x2",
  "x1";"x2"**\dir{-},   
"x2"="L55",
  "L46"="x1", 
  "x1"+(0,-.5)="x2",
  "x1";"x2"**\dir{-},   
"x2"="L56",
  "L51"="x1", 
  "x1"+(0,-.5)="x2",
  "x1";"x2"**\dir{}?(.5)="c", 
  "x1";"x2"**\dir{-},  
  "c"-(0,0)*!<-4pt,0pt>{\bullet S},  
"x2"="L61",
  "L52"="x1", 
  "x1"+(0,-.5)="x2",
  "x1";"x2"**\dir{}?(.5)="c", 
  "x1";"x2"**\dir{-},   
  "c"-(0,0)*!<-4pt,0pt>{\bullet S}, 
"x2"="L62",
  "L53"="x1", 
  "x1"+(0,-.5)="x2",
  "x1";"x2"**\dir{}?(.5)="c", 
  "x1";"x2"**\dir{-},  
  "c"-(0,0)*!<-4pt,0pt>{\bullet \widetilde{S } }, 
"x2"="L63",
  "L54"="x1", 
  "x1"+(0,-.5)="x2",
  "x1";"x2"**\dir{-},   
"x2"="L64",
  "L55"="x1", 
  "x1"+(0,-.5)="x2",
  "x1";"x2"**\dir{-},  
"x2"="L65",
  "L56"="x1", 
  "x1"+(0,-.5)="x2",
  "x1";"x2"**\dir{-},  
"x2"="L66",
"L61"="L21",
"L62"="L22",
"L63"="L23",
  "L64";"L64" + (0,-.5)="L24"**\dir{-}, 
  "L65";"L65" + (0,-.5)="L25"**\dir{-}, 
  "L66";"L66" + (0,-.5)="L26"**\dir{-}, 
  "L21"="x1", 
  "x1"+(0,-.5)="x2",
  "x1";"x2"**\dir{-},   
"x2"="L31",
  "L22"="x1", 
  "x1"+(0,-.5)="x2",
  "x1";"x2"**\dir{-},   
"x2"="L32",
   "L23"="x1",    
   "L24"="x2",  
  "x1"+(0,-.5)="x3",
  "x2"+(0,-.5)="x4",
  "x1";"x2"**\dir{}?(.5)="m", 
  "m" + (0,-.25)="c", 
  "x1";"x4"**\crv{ "x1"+(0,-.25) & "c" & "x4"+(0,.25)},  
  "x2";"c"+<3pt,2pt>="c2"**\crv{"x2"+(0,-.25) & "c2"},
  "x3";"c"+<-3pt,-2pt>="c1"**\crv{"x3"+(0,.25) & "c1"},   
  "c"+(0,.10)*!<0pt,-6pt>{},
"x3"="L33",
"x4"="L34",
  "L25"="x1", 
  "x1"+(0,-.5)="x2",
  "x1";"x2"**\dir{-},  
"x2"="L35",
  "L26"="x1", 
  "x1"+(0,-.5)="x2",
  "x1";"x2"**\dir{-},   
"x2"="L36",
"L31"="L21",
"L32"="L22",
"L33"="L23",
"L34"="L24",
"L35"="L25",
"L36"="L26",
  "L21"="x1", 
  "x1"+(0,-.5)="x2",
  "x1";"x2"**\dir{-},   
"x2"="L31",
  "L22"="x1", 
  "x1"+(0,-.5)="x2",
  "x1";"x2"**\dir{-},  
"x2"="L32",
  "L23"="x1", 
  "x1"+(0,-.5)="x2",
  "x1";"x2"**\dir{-},  
"x2"="L33",
   "L24"="x1",      
   "L25"="x2",  
  "x1"+(0,-.5)="x3",
  "x2"+(0,-.5)="x4",
  "x1";"x2"**\dir{}?(.5)="m", 
  "m" + (0,-.25)="c", 
  "x1";"x4"**\crv{ "x1"+(0,-.25) & "c" & "x4"+(0,.25)},  
  "x2";"c"+<3pt,2pt>="c2"**\crv{"x2"+(0,-.25) & "c2"},
  "x3";"c"+<-3pt,-2pt>="c1"**\crv{"x3"+(0,.25) & "c1"}, 
  "c"+(0,.10)*!<0pt,-6pt>{},
"x3"="L34",
"x4"="L35",
  "L26"="x1", 
  "x1"+(0,-.5)="x2",
  "x1";"x2"**\dir{-},   
"x2"="L36",
"L31"="L21",
"L32"="L22",
"L33"="L23",
"L34"="L24",
"L35"="L25",
"L36"="L26",
  "L21"="x1", 
  "x1"+(0,-.5)="x2",
  "x1";"x2"**\dir{-},  
"x2"="L31",
   "L22"="x1",   
   "L23"="x2", 
  "x1"+(0,-.5)="x3",
  "x2"+(0,-.5)="x4",
  "x1";"x2"**\dir{}?(.5)="m", 
  "m" + (0,-.25)="c", 
  "x1";"x4"**\crv{ "x1"+(0,-.25) & "c" & "x4"+(0,.25)},  
  "x2";"c"+<3pt,2pt>="c2"**\crv{"x2"+(0,-.25) & "c2"},
  "x3";"c"+<-3pt,-2pt>="c1"**\crv{"x3"+(0,.25) & "c1"}, 
  "c"+(0,.10)*!<0pt,-6pt>{},
"x3"="L32",
"x4"="L33",
  "L24"="x1", 
  "x1"+(0,-.5)="x2",
  "x1";"x2"**\dir{-},   
"x2"="L34",
  "L25"="x1", 
  "x1"+(0,-.5)="x2",
  "x1";"x2"**\dir{-},   
"x2"="L35",
  "L26"="x1", 
  "x1"+(0,-.5)="x2",
  "x1";"x2"**\dir{-},   
"x2"="L36",
"L31"="L21",
"L32"="L22",
"L33"="L23",
    "L34"="L24",
"L35"="L25",
"L36"="L26",
  "L21"="x1", 
  "L22"="x2",  
  "x1";"x2"**\dir{}?(.5)="m", 
  "m" + (0,-.25)="c", 
  "m" + (0,-.5)="x3", 
  "x1";"x3"**\crv{"x1"-(0,.25) & "c" & "x3"},
  "x2";"x3"**\crv{"x2"-(0,.25) & "c"  & "x3"}, 
"x3"= "L31",
  "L23"="x1", 
  "L24"="x2",   
  "x1";"x2"**\dir{}?(.5)="m", 
  "m" + (0,-.25)="c", 
  "m" + (0,-.5)="x3", 
  "x1";"x3"**\crv{"x1"-(0,.25) & "c" & "x3"},
  "x2";"x3"**\crv{"x2"-(0,.25) & "c"  & "x3"}, 
"x3"= "L32",
  "L25"="x1", 
  "L26"="x2",  
  "x1";"x2"**\dir{}?(.5)="m", 
  "m" + (0,-.25)="c", 
  "m" + (0,-.5)="x3", 
  "x1";"x3"**\crv{"x1"-(0,.25) & "c" & "x3"},
  "x2";"x3"**\crv{"x2"-(0,.25) & "c"  & "x3"}, 
"x3"= "L33",
"L31"="L21",
"L32"="L22",
"L33"="L23",
"L21"="L11",
"L22"="L12",
"L23"="L13",
   "L11"="x1", 
  "x1"+(-1,-.5)="x2",
  "x1"+(0,-.5)="x3",
  "x1"+(1,-.5)="x4",
     "x1";"x3"**\dir{}?(.5)="c", 
   "x1";"x2"**\crv{"x1" & "c" & "x2"+(0,.5)}, 
   "x1";"x3"**\crv{"x1" & "c" & "x3"+(0,.5)}, 
   "x1";"x4"**\crv{"x1" & "c" & "x4"+(0,.5)}, 
"x2"="L21",
"x3"="L22",
"x4"="L23",
 "L12"="x1", 
  "x1"+(0,-.5)="x2",
  "x1";"x2"**\dir{-},  
"x2"="L24",
 "L13"="x1", 
  "x1"+(0,-.5)="x2",
  "x1";"x2"**\dir{-},  
"x2"="L25",
  "L13"+ (1,-.25)="x1", 
  "x1"+(0,-.25)="x2",
  "x1";"x2"**\dir{}?(.5)="c", 
  "x1";"x2"**\dir{-},  
  "c"-(0,0)*!<-3pt,0pt>{\bullet \sigma },  
  "x1"+(0,0)*!<0pt,-2pt>{\circ},  
"x2"="L26",
  "L21"="x1", 
  "x1"+(0,-.5)="x2",
  "x1";"x2"**\dir{-},  
"x2"="L31",
   "L22"="x1",    
     "L23"="x2",  
  "x1"+(0,-.5)="x3",
  "x2"+(0,-.5)="x4",
  "x1";"x2"**\dir{}?(.5)="m", 
  "m" + (0,-.25)="c", 
  "x1";"x4"**\crv{ "x1"+(0,-.25) & "c" & "x4"+(0,.25)},  
  "x2";"c"+<3pt,2pt>="c2"**\crv{"x2"+(0,-.25) & "c2"},
  "x3";"c"+<-3pt,-2pt>="c1"**\crv{"x3"+(0,.25) & "c1"}, 
  "c"+(0,.10)*!<0pt,-6pt>{},
"x3"="L32",
"x4"="L33",
"L24"+(0,-0)="L34",
"L25"+(0,-0)="L35",
"L26"+(0,-0)="L36",
   "L31"="x1",  
     "L32"="x2", 
  "x1"+(0,-.5)="x3",
  "x2"+(0,-.5)="x4",
  "x1";"x2"**\dir{}?(.5)="m", 
  "m" + (0,-.25)="c", 
  "x1";"x4"**\crv{ "x1"+(0,-.25) & "c" & "x4"+(0,.25)},  
  "x2";"c"+<3pt,2pt>="c2"**\crv{"x2"+(0,-.25) & "c2"},
  "x3";"c"+<-3pt,-2pt>="c1"**\crv{"x3"+(0,.25) & "c1"}, 
  "c"+(0,.10)*!<0pt,-6pt>{},
"x3"="L41",
"x4"="L42",
  "L33"="x1", 
  "x1"+(0,-.5)="x2",
  "x1";"x2"**\dir{-},  
"x2"="L43",
  "L24"="x1", 
  "x1"+(0,-.5)="x2",
  "x1";"x2"**\dir{-},  
"x2"="L44",
  "L25"="x1", 
  "x1"+(0,-.5)="x2",
  "x1";"x2"**\dir{-},   
"x2"="L45",
  "L26"="x1", 
  "x1"+(0,-.5)="x2",
  "x1";"x2"**\dir{-},   
"x2"="L46",
  "L41"="x1", 
  "x1"+(0,-.5)="x2",
  "x1";"x2"**\dir{-},  
"x2"="L51",
   "L42"="x1",    
     "L43"="x2",  
  "x1"+(0,-.5)="x3",
  "x2"+(0,-.5)="x4",
  "x1";"x2"**\dir{}?(.5)="m", 
  "m" + (0,-.25)="c", 
  "x1";"x4"**\crv{ "x1"+(0,-.25) & "c" & "x4"+(0,.25)},  
  "x2";"c"+<3pt,2pt>="c2"**\crv{"x2"+(0,-.25) & "c2"},
  "x3";"c"+<-3pt,-2pt>="c1"**\crv{"x3"+(0,.25) & "c1"}, 
  "c"+(0,.10)*!<0pt,-6pt>{},
"x3"="L52",
"x4"="L53",
  "L44"="x1", 
  "x1"+(0,-.5)="x2",
  "x1";"x2"**\dir{-},   
"x2"="L54",
  "L45"="x1", 
  "x1"+(0,-.5)="x2",
  "x1";"x2"**\dir{-},   
"x2"="L55",
  "L46"="x1", 
  "x1"+(0,-.5)="x2",
  "x1";"x2"**\dir{-},  
"x2"="L56",
  "L51"="x1", 
  "x1"+(0,-.5)="x2",
  "x1";"x2"**\dir{}?(.5)="c", 
  "x1";"x2"**\dir{-},  
  "c"-(0,0)*!<-4pt,0pt>{\bullet S},  
"x2"="L61",
  "L52"="x1",
  "x1"+(0,-.5)="x2",
  "x1";"x2"**\dir{}?(.5)="c", 
  "x1";"x2"**\dir{-}, 
  "c"-(0,0)*!<-4pt,0pt>{\bullet S}, 
"x2"="L62",
  "L53"="x1", 
  "x1"+(0,-.5)="x2",
  "x1";"x2"**\dir{}?(.5)="c", 
  "x1";"x2"**\dir{-},   
  "c"-(0,0)*!<-4pt,0pt>{\bullet \widetilde{S } }, 
"x2"="L63",
  "L54"="x1", 
  "x1"+(0,-.5)="x2",
  "x1";"x2"**\dir{-},  
"x2"="L64",
  "L55"="x1", 
  "x1"+(0,-.5)="x2",
  "x1";"x2"**\dir{-},  
"x2"="L65",
  "L56"="x1", 
  "x1"+(0,-.5)="x2",
  "x1";"x2"**\dir{-},  
"x2"="L66",
"L61"="L21",
"L62"="L22",
"L63"="L23",
  "L64";"L64" + (0,-.5)="L24"**\dir{-},  
  "L65";"L65" + (0,-.5)="L25"**\dir{-}, 
  "L66";"L66" + (0,-.5)="L26"**\dir{-},  
  "L21"="x1", 
  "x1"+(0,-.5)="x2",
  "x1";"x2"**\dir{-},   
"x2"="L31",
  "L22"="x1", 
  "x1"+(0,-.5)="x2",
  "x1";"x2"**\dir{-},   
"x2"="L32",
   "L23"="x1",    
   "L24"="x2",  
  "x1"+(0,-.5)="x3",
  "x2"+(0,-.5)="x4",
  "x1";"x2"**\dir{}?(.5)="m", 
  "m" + (0,-.25)="c", 
  "x1";"x4"**\crv{ "x1"+(0,-.25) & "c" & "x4"+(0,.25)},  
  "x2";"c"+<3pt,2pt>="c2"**\crv{"x2"+(0,-.25) & "c2"},
  "x3";"c"+<-3pt,-2pt>="c1"**\crv{"x3"+(0,.25) & "c1"},   
  "c"+(0,.10)*!<0pt,-6pt>{},
"x3"="L33",
"x4"="L34",
  "L25"="x1", 
  "x1"+(0,-.5)="x2",
  "x1";"x2"**\dir{-},   
"x2"="L35",
  "L26"="x1", 
  "x1"+(0,-.5)="x2",
  "x1";"x2"**\dir{-},   
"x2"="L36",
"L31"="L21",
"L32"="L22",
"L33"="L23",
"L34"="L24",
"L35"="L25",
"L36"="L26",
  "L21"="x1", 
  "x1"+(0,-.5)="x2",
  "x1";"x2"**\dir{-},   
"x2"="L31",
  "L22"="x1", 
  "x1"+(0,-.5)="x2",
  "x1";"x2"**\dir{-},   
"x2"="L32",
  "L23"="x1", 
  "x1"+(0,-.5)="x2",
  "x1";"x2"**\dir{-},   
"x2"="L33",
   "L24"="x1",     
   "L25"="x2",  
  "x1"+(0,-.5)="x3",
  "x2"+(0,-.5)="x4",
  "x1";"x2"**\dir{}?(.5)="m", 
  "m" + (0,-.25)="c", 
  "x1";"x4"**\crv{ "x1"+(0,-.25) & "c" & "x4"+(0,.25)},  
  "x2";"c"+<3pt,2pt>="c2"**\crv{"x2"+(0,-.25) & "c2"},
  "x3";"c"+<-3pt,-2pt>="c1"**\crv{"x3"+(0,.25) & "c1"}, 
  "c"+(0,.10)*!<0pt,-6pt>{},
"x3"="L34",
"x4"="L35",
  "L26"="x1", 
  "x1"+(0,-.5)="x2",
  "x1";"x2"**\dir{-},   
"x2"="L36",
"L31"="L21",
"L32"="L22",
"L33"="L23",
"L34"="L24",
"L35"="L25",
"L36"="L26",
  "L21"="x1", 
  "x1"+(0,-.5)="x2",
  "x1";"x2"**\dir{-},  
"x2"="L31",
   "L22"="x1",   
   "L23"="x2",  
  "x1"+(0,-.5)="x3",
  "x2"+(0,-.5)="x4",
  "x1";"x2"**\dir{}?(.5)="m", 
  "m" + (0,-.25)="c", 
  "x1";"x4"**\crv{ "x1"+(0,-.25) & "c" & "x4"+(0,.25)},  
  "x2";"c"+<3pt,2pt>="c2"**\crv{"x2"+(0,-.25) & "c2"},
  "x3";"c"+<-3pt,-2pt>="c1"**\crv{"x3"+(0,.25) & "c1"}, 
  "c"+(0,.10)*!<0pt,-6pt>{},
"x3"="L32",
"x4"="L33",
  "L24"="x1", 
  "x1"+(0,-.5)="x2",
  "x1";"x2"**\dir{-}, 
"x2"="L34",
  "L25"="x1", 
  "x1"+(0,-.5)="x2",
  "x1";"x2"**\dir{-},   
"x2"="L35",
  "L26"="x1", 
  "x1"+(0,-.5)="x2",
  "x1";"x2"**\dir{-},  
"x2"="L36",
"L31"="L21",
"L32"="L22",
"L33"="L23",
    "L34"="L24",
"L35"="L25",
"L36"="L26",
  "L21"="x1", 
  "L22"="x2",  
  "x1";"x2"**\dir{}?(.5)="m", 
  "m" + (0,-.25)="c", 
  "m" + (0,-.5)="x3", 
  "x1";"x3"**\crv{"x1"-(0,.25) & "c" & "x3"},
  "x2";"x3"**\crv{"x2"-(0,.25) & "c"  & "x3"}, 
"x3"= "L31",
  "L23"="x1", 
  "L24"="x2",   
  "x1";"x2"**\dir{}?(.5)="m", 
  "m" + (0,-.25)="c", 
  "m" + (0,-.5)="x3", 
  "x1";"x3"**\crv{"x1"-(0,.25) & "c" & "x3"},
  "x2";"x3"**\crv{"x2"-(0,.25) & "c"  & "x3"}, 
"x3"= "L32",
  "L25"="x1", 
  "L26"="x2",  
  "x1";"x2"**\dir{}?(.5)="m", 
  "m" + (0,-.25)="c", 
  "m" + (0,-.5)="x3", 
  "x1";"x3"**\crv{"x1"-(0,.25) & "c" & "x3"},
  "x2";"x3"**\crv{"x2"-(0,.25) & "c"  & "x3"}, 
"x3"= "L33",
"L31"="L21",
"L32"="L22",
"L33"="L23",
  "L23"+(1.5,2)*!<0pt,-6pt>{=}, 
"L21"="L11",
"L22"="L12",
"L23"="L13",
"L11"="x1", 
  "x1"+(-1,-.5)="x2",
  "x1"+(0,-.5)="x3",
  "x1"+(1,-.5)="x4",
     "x1";"x3"**\dir{}?(.5)="c", 
   "x1";"x2"**\crv{"x1" & "c" & "x2"+(0,.5)}, 
   "x1";"x3"**\crv{"x1" & "c" & "x3"+(0,.5)}, 
   "x1";"x4"**\crv{"x1" & "c" & "x4"+(0,.5)}, 
"x2"="L21",
"x3"="L22",
"x4"="L23",
 "L12"="x1", 
  "x1"+(0,-.5)="x2",
  "x1";"x2"**\dir{-},  
"x2"="L24",
 "L13"="x1", 
  "x1"+(0,-.5)="x2",
  "x1";"x2"**\dir{-},   
"x2"="L25",
  "L13"+ (1,-.25)="x1", 
  "x1"+(0,-.25)="x2",
  "x1";"x2"**\dir{}?(.5)="c", 
  "x1";"x2"**\dir{-},  
  "c"-(0,0)*!<-3pt,0pt>{\bullet \sigma },  
  "x1"+(0,0)*!<0pt,-2pt>{\circ},  
"x2"="L26",
  "L21"="x1", 
  "x1"+(0,-.5)="x2",
  "x1";"x2"**\dir{-},   
"x2"="L31",
   "L22"="x1",    
     "L23"="x2", 
  "x1"+(0,-.5)="x3",
  "x2"+(0,-.5)="x4",
  "x1";"x2"**\dir{}?(.5)="m", 
  "m" + (0,-.25)="c", 
  "x1";"x4"**\crv{ "x1"+(0,-.25) & "c" & "x4"+(0,.25)},  
  "x2";"c"+<3pt,2pt>="c2"**\crv{"x2"+(0,-.25) & "c2"},
  "x3";"c"+<-3pt,-2pt>="c1"**\crv{"x3"+(0,.25) & "c1"}, 
  "c"+(0,.10)*!<0pt,-6pt>{},
"x3"="L32",
"x4"="L33",
"L24"+(0,-0)="L34",
"L25"+(0,-0)="L35",
"L26"+(0,-0)="L36",
   "L31"="x1",     
     "L32"="x2",  
  "x1"+(0,-.5)="x3",
  "x2"+(0,-.5)="x4",
  "x1";"x2"**\dir{}?(.5)="m", 
  "m" + (0,-.25)="c", 
  "x1";"x4"**\crv{ "x1"+(0,-.25) & "c" & "x4"+(0,.25)},  
  "x2";"c"+<3pt,2pt>="c2"**\crv{"x2"+(0,-.25) & "c2"},
  "x3";"c"+<-3pt,-2pt>="c1"**\crv{"x3"+(0,.25) & "c1"}, 
  "c"+(0,.10)*!<0pt,-6pt>{},
"x3"="L41",
"x4"="L42",
  "L33"="x1", 
  "x1"+(0,-.5)="x2",
  "x1";"x2"**\dir{-},   
"x2"="L43",
  "L24"="x1", 
  "x1"+(0,-.5)="x2",
  "x1";"x2"**\dir{-},   
"x2"="L44",
  "L25"="x1", 
  "x1"+(0,-.5)="x2",
  "x1";"x2"**\dir{-},   
"x2"="L45",
  "L26"="x1", 
  "x1"+(0,-.5)="x2",
  "x1";"x2"**\dir{-},   
"x2"="L46",
  "L41"="x1", 
  "x1"+(0,-.5)="x2",
  "x1";"x2"**\dir{-}, 
"x2"="L51",
   "L42"="x1",    
     "L43"="x2",  
  "x1"+(0,-.5)="x3",
  "x2"+(0,-.5)="x4",
  "x1";"x2"**\dir{}?(.5)="m", 
  "m" + (0,-.25)="c", 
  "x1";"x4"**\crv{ "x1"+(0,-.25) & "c" & "x4"+(0,.25)},  
  "x2";"c"+<3pt,2pt>="c2"**\crv{"x2"+(0,-.25) & "c2"},
  "x3";"c"+<-3pt,-2pt>="c1"**\crv{"x3"+(0,.25) & "c1"}, 
  "c"+(0,.10)*!<0pt,-6pt>{},
"x3"="L52",
"x4"="L53",
  "L44"="x1", 
  "x1"+(0,-.5)="x2",
  "x1";"x2"**\dir{-},  
"x2"="L54",
  "L45"="x1", 
  "x1"+(0,-.5)="x2",
  "x1";"x2"**\dir{-},   
"x2"="L55",
  "L46"="x1", 
  "x1"+(0,-.5)="x2",
  "x1";"x2"**\dir{-},   
"x2"="L56",
  "L51"="x1", 
  "x1"+(0,-.5)="x2",
  "x1";"x2"**\dir{}?(.5)="c", 
  "x1";"x2"**\dir{-},  
  "c"-(0,0)*!<-4pt,0pt>{\bullet S },  
"x2"="L61",
  "L52"="x1", 
  "x1"+(0,-.5)="x2",
  "x1";"x2"**\dir{}?(.5)="c", 
  "x1";"x2"**\dir{-},   
  "c"-(0,0)*!<-4pt,0pt>{\bullet S}, 
"x2"="L62",
  "L53"="x1", 
  "x1"+(0,-.5)="x2",
  "x1";"x2"**\dir{}?(.5)="c", 
  "x1";"x2"**\dir{-},  
  "c"-(0,0)*!<-4pt,0pt>{\bullet \widetilde {S} }, 
"x2"="L63",
  "L54"="x1", 
  "x1"+(0,-.5)="x2",
  "x1";"x2"**\dir{-},   
"x2"="L64",
  "L55"="x1", 
  "x1"+(0,-.5)="x2",
  "x1";"x2"**\dir{-},  
"x2"="L65",
  "L56"="x1", 
  "x1"+(0,-.5)="x2",
  "x1";"x2"**\dir{-},  
"x2"="L66",
"L61"="L21",
"L62"="L22",
"L63"="L23",
  "L64";"L64" + (0,-.5)="L24"**\dir{-}, 
  "L65";"L65" + (0,-.5)="L25"**\dir{-}, 
  "L66";"L66" + (0,-.5)="L26"**\dir{-}, 
  "L21"="x1", 
  "x1"+(0,-.5)="x2",
  "x1";"x2"**\dir{-},   
"x2"="L31",
  "L22"="x1", 
  "x1"+(0,-.5)="x2",
  "x1";"x2"**\dir{-},   
"x2"="L32",
   "L23"="x1",    
   "L24"="x2",  
  "x1"+(0,-.5)="x3",
  "x2"+(0,-.5)="x4",
  "x1";"x2"**\dir{}?(.5)="m", 
  "m" + (0,-.25)="c", 
  "x1";"x4"**\crv{ "x1"+(0,-.25) & "c" & "x4"+(0,.25)},  
  "x2";"c"+<3pt,2pt>="c2"**\crv{"x2"+(0,-.25) & "c2"},
  "x3";"c"+<-3pt,-2pt>="c1"**\crv{"x3"+(0,.25) & "c1"},   
  "c"+(0,.10)*!<0pt,-6pt>{},
"x3"="L33",
"x4"="L34",
  "L25"="x1", 
  "x1"+(0,-.5)="x2",
  "x1";"x2"**\dir{-},  
"x2"="L35",
  "L26"="x1", 
  "x1"+(0,-.5)="x2",
  "x1";"x2"**\dir{-},   
"x2"="L36",
"L31"="L21",
"L32"="L22",
"L33"="L23",
"L34"="L24",
"L35"="L25",
"L36"="L26",
  "L21"="x1", 
  "x1"+(0,-.5)="x2",
  "x1";"x2"**\dir{-},   
"x2"="L31",
  "L22"="x1", 
  "x1"+(0,-.5)="x2",
  "x1";"x2"**\dir{-},  
"x2"="L32",
  "L23"="x1", 
  "x1"+(0,-.5)="x2",
  "x1";"x2"**\dir{-},  
"x2"="L33",
   "L24"="x1",      
   "L25"="x2",  
  "x1"+(0,-.5)="x3",
  "x2"+(0,-.5)="x4",
  "x1";"x2"**\dir{}?(.5)="m", 
  "m" + (0,-.25)="c", 
  "x1";"x4"**\crv{ "x1"+(0,-.25) & "c" & "x4"+(0,.25)},  
  "x2";"c"+<3pt,2pt>="c2"**\crv{"x2"+(0,-.25) & "c2"},
  "x3";"c"+<-3pt,-2pt>="c1"**\crv{"x3"+(0,.25) & "c1"}, 
  "c"+(0,.10)*!<0pt,-6pt>{},
"x3"="L34",
"x4"="L35",
  "L26"="x1", 
  "x1"+(0,-.5)="x2",
  "x1";"x2"**\dir{-},   
"x2"="L36",
"L31"="L21",
"L32"="L22",
"L33"="L23",
"L34"="L24",
"L35"="L25",
"L36"="L26",
  "L21"="x1", 
  "x1"+(0,-.5)="x2",
  "x1";"x2"**\dir{-},  
"x2"="L31",
   "L22"="x1",   
   "L23"="x2", 
  "x1"+(0,-.5)="x3",
  "x2"+(0,-.5)="x4",
  "x1";"x2"**\dir{}?(.5)="m", 
  "m" + (0,-.25)="c", 
  "x1";"x4"**\crv{ "x1"+(0,-.25) & "c" & "x4"+(0,.25)},  
  "x2";"c"+<3pt,2pt>="c2"**\crv{"x2"+(0,-.25) & "c2"},
  "x3";"c"+<-3pt,-2pt>="c1"**\crv{"x3"+(0,.25) & "c1"}, 
  "c"+(0,.10)*!<0pt,-6pt>{},
"x3"="L32",
"x4"="L33",
  "L24"="x1", 
  "x1"+(0,-.5)="x2",
  "x1";"x2"**\dir{-},   
"x2"="L34",
  "L25"="x1", 
  "x1"+(0,-.5)="x2",
  "x1";"x2"**\dir{-},   
"x2"="L35",
  "L26"="x1", 
  "x1"+(0,-.5)="x2",
  "x1";"x2"**\dir{-},   
"x2"="L36",
"L31"="L21",
"L32"="L22",
"L33"="L23",
    "L34"="L24",
"L35"="L25",
"L36"="L26",
  "L21"="x1", 
  "L22"="x2",  
  "x1";"x2"**\dir{}?(.5)="m", 
  "m" + (0,-.25)="c", 
  "m" + (0,-.5)="x3", 
  "x1";"x3"**\crv{"x1"-(0,.25) & "c" & "x3"},
  "x2";"x3"**\crv{"x2"-(0,.25) & "c"  & "x3"}, 
"x3"= "L31",
  "L23"="x1", 
  "L24"="x2",   
  "x1";"x2"**\dir{}?(.5)="m", 
  "m" + (0,-.25)="c", 
  "m" + (0,-.5)="x3", 
  "x1";"x3"**\crv{"x1"-(0,.25) & "c" & "x3"},
  "x2";"x3"**\crv{"x2"-(0,.25) & "c"  & "x3"}, 
"x3"= "L32",
  "L25"="x1", 
  "L26"="x2",  
  "x1";"x2"**\dir{}?(.5)="m", 
  "m" + (0,-.25)="c", 
  "m" + (0,-.5)="x3", 
  "x1";"x3"**\crv{"x1"-(0,.25) & "c" & "x3"},
  "x2";"x3"**\crv{"x2"-(0,.25) & "c"  & "x3"}, 
"x3"= "L33",
"L31"="L21",
"L32"="L22",
"L33"="L23",
"L21"="L11",
"L22"="L12",
"L23"="L13",
   "L11"="x1", 
  "x1"+(-1,-.5)="x2",
  "x1"+(0,-.5)="x3",
  "x1"+(1,-.5)="x4",
     "x1";"x3"**\dir{}?(.5)="c", 
   "x1";"x2"**\crv{"x1" & "c" & "x2"+(0,.5)}, 
   "x1";"x3"**\crv{"x1" & "c" & "x3"+(0,.5)}, 
   "x1";"x4"**\crv{"x1" & "c" & "x4"+(0,.5)}, 
"x2"="L21",
"x3"="L22",
"x4"="L23",
 "L12"="x1", 
  "x1"+(0,-.5)="x2",
  "x1";"x2"**\dir{-},  
"x2"="L24",
 "L13"="x1", 
  "x1"+(0,-.5)="x2",
  "x1";"x2"**\dir{-},  
"x2"="L25",
  "L13"+ (1,-.25)="x1", 
  "x1"+(0,-.25)="x2",
  "x1";"x2"**\dir{}?(.5)="c", 
  "x1";"x2"**\dir{-},  
  "c"-(0,0)*!<-3pt,0pt>{\bullet \sigma },  
  "x1"+(0,0)*!<0pt,-2pt>{\circ},  
"x2"="L26",
  "L21"="x1", 
  "x1"+(0,-.5)="x2",
  "x1";"x2"**\dir{-},  
"x2"="L31",
   "L22"="x1",    
     "L23"="x2",  
  "x1"+(0,-.5)="x3",
  "x2"+(0,-.5)="x4",
  "x1";"x2"**\dir{}?(.5)="m", 
  "m" + (0,-.25)="c", 
  "x1";"x4"**\crv{ "x1"+(0,-.25) & "c" & "x4"+(0,.25)},  
  "x2";"c"+<3pt,2pt>="c2"**\crv{"x2"+(0,-.25) & "c2"},
  "x3";"c"+<-3pt,-2pt>="c1"**\crv{"x3"+(0,.25) & "c1"}, 
  "c"+(0,.10)*!<0pt,-6pt>{},
"x3"="L32",
"x4"="L33",
"L24"+(0,-0)="L34",
"L25"+(0,-0)="L35",
"L26"+(0,-0)="L36",
   "L31"="x1",  
     "L32"="x2", 
  "x1"+(0,-.5)="x3",
  "x2"+(0,-.5)="x4",
  "x1";"x2"**\dir{}?(.5)="m", 
  "m" + (0,-.25)="c", 
  "x1";"x4"**\crv{ "x1"+(0,-.25) & "c" & "x4"+(0,.25)},  
  "x2";"c"+<3pt,2pt>="c2"**\crv{"x2"+(0,-.25) & "c2"},
  "x3";"c"+<-3pt,-2pt>="c1"**\crv{"x3"+(0,.25) & "c1"}, 
  "c"+(0,.10)*!<0pt,-6pt>{},
"x3"="L41",
"x4"="L42",
  "L33"="x1", 
  "x1"+(0,-.5)="x2",
  "x1";"x2"**\dir{-},  
"x2"="L43",
  "L24"="x1", 
  "x1"+(0,-.5)="x2",
  "x1";"x2"**\dir{-},  
"x2"="L44",
  "L25"="x1", 
  "x1"+(0,-.5)="x2",
  "x1";"x2"**\dir{-},   
"x2"="L45",
  "L26"="x1", 
  "x1"+(0,-.5)="x2",
  "x1";"x2"**\dir{-},   
"x2"="L46",
  "L41"="x1", 
  "x1"+(0,-.5)="x2",
  "x1";"x2"**\dir{-},  
"x2"="L51",
   "L42"="x1",    
     "L43"="x2",  
  "x1"+(0,-.5)="x3",
  "x2"+(0,-.5)="x4",
  "x1";"x2"**\dir{}?(.5)="m", 
  "m" + (0,-.25)="c", 
  "x1";"x4"**\crv{ "x1"+(0,-.25) & "c" & "x4"+(0,.25)},  
  "x2";"c"+<3pt,2pt>="c2"**\crv{"x2"+(0,-.25) & "c2"},
  "x3";"c"+<-3pt,-2pt>="c1"**\crv{"x3"+(0,.25) & "c1"}, 
  "c"+(0,.10)*!<0pt,-6pt>{},
"x3"="L52",
"x4"="L53",
  "L44"="x1", 
  "x1"+(0,-.5)="x2",
  "x1";"x2"**\dir{-},   
"x2"="L54",
  "L45"="x1", 
  "x1"+(0,-.5)="x2",
  "x1";"x2"**\dir{-},   
"x2"="L55",
  "L46"="x1", 
  "x1"+(0,-.5)="x2",
  "x1";"x2"**\dir{-},  
"x2"="L56",
  "L51"="x1", 
  "x1"+(0,-.5)="x2",
  "x1";"x2"**\dir{}?(.5)="c", 
  "x1";"x2"**\dir{-},  
  "c"-(0,0)*!<-4pt,0pt>{\bullet S},  
"x2"="L61",
  "L52"="x1",
  "x1"+(0,-.5)="x2",
  "x1";"x2"**\dir{}?(.5)="c", 
  "x1";"x2"**\dir{-}, 
  "c"-(0,0)*!<-4pt,0pt>{\bullet S}, 
"x2"="L62",
  "L53"="x1", 
  "x1"+(0,-.5)="x2",
  "x1";"x2"**\dir{}?(.5)="c", 
  "x1";"x2"**\dir{-},   
  "c"-(0,0)*!<-4pt,0pt>{\bullet \widetilde{S } }, 
"x2"="L63",
  "L54"="x1", 
  "x1"+(0,-.5)="x2",
  "x1";"x2"**\dir{-},  
"x2"="L64",
  "L55"="x1", 
  "x1"+(0,-.5)="x2",
  "x1";"x2"**\dir{-},  
"x2"="L65",
  "L56"="x1", 
  "x1"+(0,-.5)="x2",
  "x1";"x2"**\dir{-},  
"x2"="L66",
"L61"="L21",
"L62"="L22",
"L63"="L23",
  "L64";"L64" + (0,-.5)="L24"**\dir{-},  
  "L65";"L65" + (0,-.5)="L25"**\dir{-}, 
  "L66";"L66" + (0,-.5)="L26"**\dir{-},  
  "L21"="x1", 
  "x1"+(0,-.5)="x2",
  "x1";"x2"**\dir{-},   
"x2"="L31",
  "L22"="x1", 
  "x1"+(0,-.5)="x2",
  "x1";"x2"**\dir{-},   
"x2"="L32",
   "L23"="x1",    
   "L24"="x2",  
  "x1"+(0,-.5)="x3",
  "x2"+(0,-.5)="x4",
  "x1";"x2"**\dir{}?(.5)="m", 
  "m" + (0,-.25)="c", 
  "x1";"x4"**\crv{ "x1"+(0,-.25) & "c" & "x4"+(0,.25)},  
  "x2";"c"+<3pt,2pt>="c2"**\crv{"x2"+(0,-.25) & "c2"},
  "x3";"c"+<-3pt,-2pt>="c1"**\crv{"x3"+(0,.25) & "c1"},   
  "c"+(0,.10)*!<0pt,-6pt>{},
"x3"="L33",
"x4"="L34",
  "L25"="x1", 
  "x1"+(0,-.5)="x2",
  "x1";"x2"**\dir{-},   
"x2"="L35",
  "L26"="x1", 
  "x1"+(0,-.5)="x2",
  "x1";"x2"**\dir{-},   
"x2"="L36",
"L31"="L21",
"L32"="L22",
"L33"="L23",
"L34"="L24",
"L35"="L25",
"L36"="L26",
  "L21"="x1", 
  "x1"+(0,-.5)="x2",
  "x1";"x2"**\dir{-},   
"x2"="L31",
  "L22"="x1", 
  "x1"+(0,-.5)="x2",
  "x1";"x2"**\dir{-},   
"x2"="L32",
  "L23"="x1", 
  "x1"+(0,-.5)="x2",
  "x1";"x2"**\dir{-},   
"x2"="L33",
   "L24"="x1",     
   "L25"="x2",  
  "x1"+(0,-.5)="x3",
  "x2"+(0,-.5)="x4",
  "x1";"x2"**\dir{}?(.5)="m", 
  "m" + (0,-.25)="c", 
  "x1";"x4"**\crv{ "x1"+(0,-.25) & "c" & "x4"+(0,.25)},  
  "x2";"c"+<3pt,2pt>="c2"**\crv{"x2"+(0,-.25) & "c2"},
  "x3";"c"+<-3pt,-2pt>="c1"**\crv{"x3"+(0,.25) & "c1"}, 
  "c"+(0,.10)*!<0pt,-6pt>{},
"x3"="L34",
"x4"="L35",
  "L26"="x1", 
  "x1"+(0,-.5)="x2",
  "x1";"x2"**\dir{-},   
"x2"="L36",
"L31"="L21",
"L32"="L22",
"L33"="L23",
"L34"="L24",
"L35"="L25",
"L36"="L26",
  "L21"="x1", 
  "x1"+(0,-.5)="x2",
  "x1";"x2"**\dir{-},  
"x2"="L31",
   "L22"="x1",   
   "L23"="x2",  
  "x1"+(0,-.5)="x3",
  "x2"+(0,-.5)="x4",
  "x1";"x2"**\dir{}?(.5)="m", 
  "m" + (0,-.25)="c", 
  "x1";"x4"**\crv{ "x1"+(0,-.25) & "c" & "x4"+(0,.25)},  
  "x2";"c"+<3pt,2pt>="c2"**\crv{"x2"+(0,-.25) & "c2"},
  "x3";"c"+<-3pt,-2pt>="c1"**\crv{"x3"+(0,.25) & "c1"}, 
  "c"+(0,.10)*!<0pt,-6pt>{},
"x3"="L32",
"x4"="L33",
  "L24"="x1", 
  "x1"+(0,-.5)="x2",
  "x1";"x2"**\dir{-}, 
"x2"="L34",
  "L25"="x1", 
  "x1"+(0,-.5)="x2",
  "x1";"x2"**\dir{-},   
"x2"="L35",
  "L26"="x1", 
  "x1"+(0,-.5)="x2",
  "x1";"x2"**\dir{-},  
"x2"="L36",
"L31"="L21",
"L32"="L22",
"L33"="L23",
    "L34"="L24",
"L35"="L25",
"L36"="L26",
  "L21"="x1", 
  "L22"="x2",  
  "x1";"x2"**\dir{}?(.5)="m", 
  "m" + (0,-.25)="c", 
  "m" + (0,-.5)="x3", 
  "x1";"x3"**\crv{"x1"-(0,.25) & "c" & "x3"},
  "x2";"x3"**\crv{"x2"-(0,.25) & "c"  & "x3"}, 
"x3"= "L31",
  "L23"="x1", 
  "L24"="x2",   
  "x1";"x2"**\dir{}?(.5)="m", 
  "m" + (0,-.25)="c", 
  "m" + (0,-.5)="x3", 
  "x1";"x3"**\crv{"x1"-(0,.25) & "c" & "x3"},
  "x2";"x3"**\crv{"x2"-(0,.25) & "c"  & "x3"}, 
"x3"= "L32",
  "L25"="x1", 
  "L26"="x2",  
  "x1";"x2"**\dir{}?(.5)="m", 
  "m" + (0,-.25)="c", 
  "m" + (0,-.5)="x3", 
  "x1";"x3"**\crv{"x1"-(0,.25) & "c" & "x3"},
  "x2";"x3"**\crv{"x2"-(0,.25) & "c"  & "x3"}, 
"x3"= "L33",
"L31"="L21",
"L32"="L22",
"L33"="L23",
  "L21"-(0,.5)*!<0pt,-6pt>{H},  
  "L22"-(0,.5)*!<0pt,-6pt>{H}, 
  "L23"-(0,.5)*!<0pt,-6pt>{H}, 
\endxy
\quad \quad  %%%%%%%%%%%%%%%%%%%%%%%%%%%%%%%%%%%%%%%
\xy /r2.5pc/:, 
(1,-5)="L11",
(2,-5)="L12",
(3,-5)="L13",
"L11"="a1",
"L12"="a2",
"L13"="a3",
"a1"+(0,.5)*!<0pt,6pt>{H}, 
"a2"+(0,.5)*!<0pt,6pt>{H}, 
"a3"+(0,.5)*!<0pt,6pt>{H},
"a1" + (0,-1)="b1",
"a2" + (0,-1)="b2",
"a3" + (0,-1)="b3",
"a1";"a2"**\dir{}?(.5)="m12", 
  "m12" + (0,-.25)="c12", 
"a2";"a3"**\dir{}?(.5)="m23", 
  "m23" + (0,-.25)="c23", 
  "a1";"b2"**\crv{ "a1"+(0,-.25) & "c12" & "b2"+(0,.25)},  
  "a2";"b3"**\crv{ "a2"+(0,-.25) & "c23" & "b3"+(0,.25)}, 
                "b1";"c12"+<-0pt,-11pt>="c21"**\crv{"b1"+(0,.25) & "c21"},         
               "a3";"c23"+<-0pt,-2pt>="c32"**\crv{"a3"+(0,-.25) & "c32"},
"c12" +(.35,-.20);"c23" +(-.25,-.10)**\dir{-}, 
"b1"="L11",
"b2"="L12",
"b3"="L13",
"L11"="a1",
"L12"="a2",
"L13"="a3",
"a1" + (0,-1)="b1",
"a2" + (0,-1)="b2",
"a3" + (0,-1)="b3",
"a1";"a2"**\dir{}?(.5)="m12", 
  "m12" + (0,-.25)="c12", 
"a2";"a3"**\dir{}?(.5)="m23", 
  "m23" + (0,-.25)="c23", 
  "a1";"b2"**\crv{ "a1"+(0,-.25) & "c12" & "b2"+(0,.25)},  
  "a2";"b3"**\crv{ "a2"+(0,-.25) & "c23" & "b3"+(0,.25)}, 
                "b1";"c12"+<-0pt,-11pt>="c21"**\crv{"b1"+(0,.25) & "c21"},         
               "a3";"c23"+<-0pt,-2pt>="c32"**\crv{"a3"+(0,-.25) & "c32"},
"c12" +(.35,-.20);"c23" +(-.25,-.10)**\dir{-}, 
"b1"="L11",
"b2"="L12",
"b3"="L13",
"L11"="a1",
"L12"="a2",
"L13"="a3",
"a1" + (0,-1)="b1",
"a2" + (0,-1)="b2",
"a3" + (0,-1)="b3",
"a1";"a2"**\dir{}?(.5)="m12", 
  "m12" + (0,-.25)="c12", 
"a2";"a3"**\dir{}?(.5)="m23", 
  "m23" + (0,-.25)="c23", 
  "a1";"b2"**\crv{ "a1"+(0,-.25) & "c12" & "b2"+(0,.25)},  
  "a2";"b3"**\crv{ "a2"+(0,-.25) & "c23" & "b3"+(0,.25)}, 
                "b1";"c12"+<-0pt,-11pt>="c21"**\crv{"b1"+(0,.25) & "c21"},         
               "a3";"c23"+<-0pt,-2pt>="c32"**\crv{"a3"+(0,-.25) & "c32"},
"c12" +(.35,-.20);"c23" +(-.25,-.10)**\dir{-},
"b1"-(0,.5)*!<0pt,-6pt>{H},
"b2"-(0,.5)*!<0pt,-6pt>{H}, 
"b3"-(0,.5)*!<0pt,-6pt>{H},   
\endxy 
\]%\label{tau4intermsos si}

Now we proceed to our last example, the Hopf cyclic cohomology for 
quasitriangular quasi-Hopf algebras, which was one of the main motivations 
for this work. One knows that a quasitriangular quasi-Hopf algebra $(H,\,R,\,\Phi,
\, \alpha ,\, \beta)$ is a Hopf algebra in the braided monoidal
category of (left) $H$-modules ~\cite{bn,maj} . This braided Hopf
algebra $\underline{H}$ has the following structure. As a vector
space $\underline{H} = H$, with $H$-module structure given by conjugation
\[ a \rhd h =a^{(1)}h S(a^{(2)}). \]
The Hopf algebra structure on $\underline{H}$ is given by ~\cite{bn}:

\[ \underline{m}(a,b)=a \underline{.} b = X_1 a S(x_1 X_2)\alpha x_2 X_3^{(1)} b S(x_3 X_3^{(2)}), \] 
with unit $\beta$.

\begin{multline*}  
\underline{\Delta} (h)= h^{(\underline{1})} \ot h^{(\underline{2})}=\\ 
x_1 X_1 h^{(1)}g_1 S(x_2 R_2y_3 X_3^{(2)}) \ot x_3 R_1 \rhd y_1X_2 h^{(2)}g_2S(y_2X_3^{(1)}) ,
\end{multline*}
with counit $\underline{\varepsilon} = \varepsilon$, and antipode

\[ \underline{S} (h)=x_1R_2p_2S(q_1(X_2R_1p_1 \rhd h)S(q_2)X_3) .  \]

Here we have used the following notations:

\[\Phi=X_1\ot X_2\ot X_3  ,\]
\[ \Phi^{-1}=x_1\ot x_2\ot x_3 = y_1\ot y_2\ot y_3  ,\]
\[p_1\ot p_2=x_1\ot x_2 \beta S(x_3)  ,\]
\[q_1\ot q_2=X_1\ot S^{-1}(\alpha X_3)X_2  ,\]
and
\[g_1\ot g_2=\Delta(S(x_1)\alpha x_2)\xi(S\ot S)(\Delta^{\text{op}}(x_3)),  \]
 where,
 \[\xi=B_1\beta S(B_4)\ot B_2 \beta S(B_3),   \]
\[ B_1 \ot B_2 \ot B_3\ot B_4 = (\Delta \ot id \ot id)(\Phi)(\Phi^{-1}\ot 1)  .\]

Let $(\delta, \sigma)$ be a braided modular pair in involution for
$\underline{H}$. Using Theorem ~\eqref{brdd ver of CM thry in nonsymmetric} we
can associate a para-cocyclic module to $\underline{H}$. By passing to
 subspaces $\text{ker}(1- \tau_n^{n+1})$ we obtain a cocyclic module.
Therefore we have the following corollary.

\begin{corollary}
For any quasitriangular quasi-Hopf algebra $(H,\,R,\,\Phi,\,\alpha ,\, \beta)$ 
endowed with a braided modular pair in involution $(\delta, \sigma)$, 
the complex obtained in Theorem ~\eqref{brdd ver of CM thry in nonsymmetric} 
defines a Hopf cyclic cohomology for $H$, denoted by 
 $HC_{(\delta,\sigma)}^* (H,\,R,\,\Phi,\,\alpha ,\, \beta)$.
\end{corollary}

\end{document}